\documentclass[titlepage,11pt]{article}
\usepackage[USenglish]{babel}
\usepackage{multicol}
\usepackage{a4wide}
\usepackage[all]{xy}
\usepackage{amssymb}
\usepackage{amsmath}
\usepackage{latexsym}
\usepackage{pifont}
\usepackage{pdfsync}
\usepackage{color}
\usepackage{mathrsfs}
\usepackage{enumerate}
\usepackage{float,fancyvrb,verbatim,multicol}
\usepackage{ifthen}
\usepackage[parfill]{parskip}    % Activate to begin paragraphs with an empty line rather than an indent
\usepackage{xr}
\usepackage{makeidx}
\usepackage[utf8]{inputenc} 
\title{Measures and all that --- A Tutorial}
\author{Ernst-Erich Doberkat\\Math ++ Software, Bochum\\\texttt{eed@doberkat.de}}
\date{\today}
\newcommand{\labelImpl}[2]{\ensuremath{\ref{#1}~\Rightarrow~\ref{#2}}}

\newcommand{\Inff}[1]{\ensuremath{#1_\infty}}

\newcommand{\RelOp}{\ensuremath{\ltimes}}

\newcommand{\Paths}{\ensuremath{\textsc{Paths}}}

\newcommand{\Steady}[2][\RelOp p]{\ensuremath{\mathcal{S}_{#1}({#2})}}

\newcommand{\PathQuant}[2][\RelOp p]{\ensuremath{\mathcal{P}_{#1}({#2})}}

\newcommand{\PathQuantSenza}{\ensuremath{\mathcal{P}}}
\newcommand{\Next}[2][I]{\ensuremath{\mathcal{X}^{#1}\:{#2}}}

\newcommand{\Until}[3][I]{\ensuremath{{#2}\:\mathcal{U}^{#1}\:{#3}}}
\newcommand{\Pfeil}[1][t]{\ensuremath{\overset{#1}{\longrightarrow}}}
\newcommand{\At}[2][\sigma]{\ensuremath{{#1}@{#2}}}

\newcommand{\Klasse}[2]{\left[#1\right]_{#2}}
\newcommand{\Faktor}[2]{{#1}/{#2}}
\newcommand{\fMap}[1]{\eta_{#1}}
\newcommand{\Bild}[2]{{#1}\left[#2\right]}
\newcommand{\InvBild}[2]{\Bild{#1^{-1}}{#2}}
\newcommand{\Kern}[1]{\mathsf{ker}\left(#1\right)}
\newcommand{\Folge}[1]{(#1_n)_{n \in \Nat}}

\newcommand{\supp}{\mathsf{supp}}

\newcommand{\spaceFont}[1]{\mathfrak{#1}}

\newcommand{\Prob}[1]{\spaceFont{P}\left(#1\right)}
\newcommand{\SubProb}[1]{\spaceFont{S}\left(#1\right)}

\newcommand{\Category}[1]{\ensuremath{\mathbf{#1}}}
\newcommand{\SubProbSenza}{\spaceFont{S}}
\newcommand{\ProbSenza}{\spaceFont{P}}
\newcommand{\PowerSet}[1]{\ensuremath{\mathcal{P}\left(#1\right)}}

\newcommand{\Borel}[1]{\ensuremath{{\mathcal B}(#1)}}

\edef\LinkeKlammer{\lbrack\!\lbrack}
\edef\RechteKlammer{\rbrack\!\rbrack}
\newcommand{\Gilt}[1][\phi]{\ensuremath{\LinkeKlammer#1\RechteKlammer}}

\newcommand{\Trans}{\rightsquigarrow}

\newcommand{\compositionSymbol}[1]{\mathbf{#1}}

\newcommand{\klComp}{\compositionSymbol{\ast}}

\newcommand{\Closure}[1]{\ensuremath{\mathsf{cl}(#1)}}

\newcommand{\Dach}[1]{{}^{\Box}{#1}}

\newtheorem{definition}{Definition}[section]
\newcommand{\BeginDefinition}[1]{%
  \begin{definition}\label{#1}
}
\newcommand{\EndDefinition}{\end{definition}}

\newtheorem{example}[definition]{Example}%[section]
\newcommand{\BeginExample}[1]{%
  \begin{example}\label{#1}\rm
}
\newcommand{\EndExample}{--- \end{example}}

\newtheorem{observation}[definition]{Observation}%[section]
\newcommand{\BeginObservation}[1]{
  \begin{observation}\label{#1}\rm
}
\newcommand{\EndObservation}{--- \end{observation}}

\newtheorem{theorem}[definition]{Theorem}%[section]
\newcommand{\BeginTheorem}[1]{%
  \begin{theorem}\label{#1}
}
\newcommand{\EndTheorem}{\end{theorem}}

\newtheorem{corollary}[definition]{Corollary}%[section]
\newcommand{\BeginCorollary}[1]{
  \begin{corollary}\label{#1}
}

\newtheorem{proposition}[definition]{Proposition}%[section]
\newcommand{\BeginProposition}[1]{%
  \begin{proposition}\label{#1}
}
\newcommand{\EndProposition}{\end{proposition}}

\newcommand{\EndCorollary}{\end{corollary}}
\newtheorem{lemma}[definition]{Lemma}%[section]
\newcommand{\BeginLemma}[1]{%
  \begin{lemma}\label{#1}
}
\newcommand{\EndLemma}{\end{lemma}}

\newtheorem{claim}{Claim}
\newcommand{\BeginClaim}[1]{%
  \begin{claim}\label{#1}
}
\newcommand{\EndClaim}{\end{claim}}

\newenvironment{proof}{\textbf{Proof\ }}{\ensuremath{\QED}}
\newcommand{\BeginProof}{\begin{proof}}
\newcommand{\EndProof}{\end{proof}}

\newenvironment{remark}{\textbf{Remark:\ }}{}
\newcommand{\BeginRemark}{\begin{remark}}
\newcommand{\EndRemark}{\QED\end{remark}}
\newcommand{\QED}{%
\ensuremath{\dashv}
}

\newcommand{\Real}{\mathbb{R}}
\newcommand{\pReal}{\mathbb{R}_{+}}
\newcommand{\Nat}{\mathbb{N}}
\newcommand{\Rational}{\mathbb{Q}}

\newcommand{\unit}{\ensuremath{\eta}}

\def\theta{\vartheta}
\def\dots{\ldots}
\def\Dach#1{\vec{#1}}

\makeatletter
\def\@axx#1{\ensuremath{\mathbb{(#1)}}}
\def\AC{\@axx{AC}}
\def\WO{\@axx{WO}}
\def\ZL{\@axx{ZL}}
\def\MP{\@axx{MP}}
\def\MI{\@axx{MI}}
\def\AD{\@axx{AD}}
\makeatother

\def\Nats{\Nat^+}
\def\NatInf{\Nat^{\infty}}
\def\sA{\ensuremath{\mathscr{A}}}

\newcommand{\isEquiv}[3]{\ensuremath{{#1}\ {#3}\ {#2}}}
\newenvironment{theExercises}
{\begin{exercise}\rm}
{\end{exercise}}
\newtheorem{exercise}{Exercise}%[chapter]
\newcommand{\BeginExercise}[1]{%
  \begin{theExercises}\label{#1}
}
\newcommand{\EndExercise}{\end{theExercises}}
\newcommand{\Exercise}[2]{\BeginExercise{#1}{#2}\EndExercise}

\def\endEx{{\Large\ding{44}}}
\renewcommand{\EndExample}{\endEx \end{example}}

\def\CatFont{\mathbf}
\def\Category#1{\ensuremath{\CatFont{#1}}}
\renewcommand{\spaceFont}[1]{\CatFont{#1}}
\def\Closed{\ensuremath{\mathbb{F}}}
\def\catK{\Category{K}}

\makeatletter
\def\@theProbSenza#1{\ensuremath{\mathbb{#1}}}
\def\SubProbSenza{\@theProbSenza{S}}
\def\ProbSenza{\@theProbSenza{P}}
\def\FinSenza{\@theProbSenza{M}}
\def\SigmaSenza{\@theProbSenza{M}_{\ensuremath{\sigma}}}
\makeatother
\renewcommand{\Prob}[1]{\ProbSenza\left(#1\right)}
\def\SubProb#1{\SubProbSenza(#1)}
\def\FinM#1{\FinSenza(#1)}
\def\SigmaM#1{\SigmaSenza(#1)}
\def\schwach#1{\ensuremath{\pmb{\wp}(#1)}}
\makeatletter
\def\@beta{\pmb{\beta}}
\def\betaSenza#1{\@beta_{#1}}
\def\theBeta#1#2#3{\ensuremath{\betaSenza{#3}(#1, #2)}}
\def\supp{\ensuremath{\mathrm{supp}}}
\makeatother

\renewcommand{\Borel}[1]{\ensuremath{{\mathcal B}(#1)}}
\newcommand{\Baire}[1]{\ensuremath{{\mathcal Ba}(#1)}}
\def\o{\varnothing}

\def\RealT{\widetilde{\Real}}
\def\MeasbFnct#1{\ensuremath{{\cal F}(#1)}}
\def\MeasbFnctT#1{\ensuremath{\widetilde{{\cal F}}(#1)}}
\def\MeasbFnctP#1{\ensuremath{{\cal F}_+(#1)}}
\def\MeasbFnctTP#1{\ensuremath{\widetilde{{\cal F}}_+(#1)}}
\def\StepFnct#1{\ensuremath{{\cal T}(#1)}}
\def\StepFnctP#1{\ensuremath{{\cal T}_{+}(#1)}}
\def\Graph#1{\ensuremath{\mathsf{graph}(#1)}}
\def\DefSect{
\ifthenelse{\boolean{isBook}}{
\def\Section{\chapter}
\def\Subsection{\section}
\def\Subsubsection{\subsection}
}{
\def\Section{\section}
\def\Subsection{\subsection}
\def\Subsubsection{\subsubsection}
}
}

\def\Rand#1{\begin{color}{red}\fbox{#1}\end{color}}
\def\Closure#1{\ensuremath{{#1}^{a}}}

\makeatletter
\def\@Tut#1#2{\cite[#1]{#2}}
\def\CategCite#1{\@Tut{#1}{EED-Categs}}
\def\SetCite#1{\@Tut{#1}{EED-Tut_sets}}
\makeatother

\makeatletter
\def\@norm#1#2#3{\ensuremath{||{#1}||_{#2}^{#3}}}
\def\infNorm#1#2{\@norm{#1}{\infty}{#2}}
\def\aNorm#1{\@norm{#1}{}{}}
\makeatother

\makeatletter
\def\@conv#1{\stackrel{#1}{\longrightarrow}} 
\def\aeC{\@conv{a.e.}} 
\def\nmC{\@conv{i.m.}}
\makeatother

\makeatletter
\newcommand{\@theL}[3]{\ensuremath{{#3}_{#1}{#2}}}
\newcommand{\cLp}[2][p]{\@theL{#1}{(#2)}{{\cal L}}}
\newcommand{\rLp}[2][p]{\@theL{#1}{(#2)}{L}}
\newcommand\cLpS[1][p]{\@theL{#1}{}{{\cal L}}}
\newcommand\rLpS[1][p]{\@theL{#1}{}{L}}
\makeatother

\newcommand{\Int}[3]{\ensuremath{\int_{#1}{#2}\ d{#3}}}
\newcommand{\dInt}[3][X]{\Int{#1}{#2}{#3}}
\def\absCont{\ensuremath{<\negthickspace<}}
\def\Paths{\ensuremath{(S \times \pReal)^\infty}}
\def\Inff#1{\ensuremath{#1_{\infty}}}
\newcommand{\AllFormulas}[1][AP]{\ensuremath{\mathfrak{L}_{#1}}}
\def\RelOp{\bowtie}
\def\comp#1{\ensuremath{\overline{#1}}}
\renewcommand{\Bild}[2]{{#1}\bigl[#2\bigr]}
\renewcommand{\Folge}[2][n]{\ensuremath{({#2}_{#1})_{#1\in\Nat}}}

\newboolean{isBook}
\setboolean{isBook}{false}
\DefSect
\def\Rand#1{\marginpar{\begin{color}{red}#1\end{color}}}
\def\phi{\varphi}

\usepackage{fancyhdr}
\makeindex
\begin{document}
\maketitle
\begin{abstract}
  This tutorial gives an overview of some of the basic techniques of
  measure theory. It includes a study of Borel sets and their
  generators for Polish and for analytic spaces, the weak topology on
  the space of all finite positive measures including its metrics, as
  well as measurable selections. Integration is covered, and product
  measures are introduced, both for finite and for arbitrary factors,
  with an application to projective systems. Finally, the duals of the
  Lp-spaces are discussed, together with the Radon-Nikodym Theorem and
  the Riesz Representation Theorem. Case studies include applications
  to stochastic Kripke models, to bisimulations, and to quotients for
  transition kernels.
\end{abstract}
%%% Local Variables: 
%%% mode: latex
%%% TeX-master: "../Mskr3"
%%% End: 

%\input{Disclaimer}
%\newpage
\tableofcontents\newpage
%\input{Categs/Hilf}
%\input{Probabs/ProbabsMain}
%!TEX root =  ../Mskr3.tex
%\def\Folder{Probabs}
%\def\Input{\input}
%\Input{\Folder/Overview}
\Section{Overview}
\label{sec:meas-overview}

Markov transition systems are based on transition probabilities on a
measurable space. This is a generalization of discrete spaces, where
certain sets are declared to be measurable. So, in contrast to assuming
that we know the probability for the transition between two states, we
have to model the probability of a transition going from one state to
a set of states: point-to-point probabilities are no longer available
due to working in a comparatively large space. Measurable spaces
are the domains of the probabilities involved. This approach has the
advantage of being more general than finite or countable spaces, but
now one deals with a fairly involved mathematical structure;
all of a sudden the dictionary has to be extended with words like
``universally measurable'' or ``sub-$\sigma$-algebra''. Measure
theory becomes an area where one has to find answers to questions
which did not appear to be particularly involved before, in the much simpler world
of discrete measures (the impression should not arise that this author
thinks that discrete measures are kiddie stuff, they are sometimes
difficult enough to handle. The continuous case, as it is called
sometimes, offers questions which simply do not arise in the discrete
context). Many arguments in this area are of a measure theoretic nature; this
tutorial makes an attempt to introduce the necessary tools and
techniques. 

It starts off with a discussion of $\sigma$-algebras,
which have already been met in~\SetCite{Section
  1.6}\footnote{\cite{EED-Tut_sets} and~\cite{EED-Categs} are other
  tutorials in this series. This is not an installment of mysteries,
  so one can be read quite independently from the others.}. We look at
the structure of $\sigma$-algebras, in particular at its generators;
it turns out that the underlying space has something to say about
it. In particular we will deal with Polish spaces and their
brethren. Two aspects deserve to be singled out. The $\sigma$-algebra
on the base space determines a $\sigma$-algebra on the space of all
finite measures, and, if this space has a topology, it determines also
a topology, the Alexandrov topology. These constructions are studied,
since they also affect the applications in logic, and in transition
systems, in which measures are vital. Second, we show that we can
construct measurable selections, which then enable constructions which
are interesting from a categorical point of view~\CategCite{2.4.2,
  2.6.1, 2.7.3}. 

After having laid the groundwork with a discussion of
$\sigma$-algebras as the domains of measures, we show that the
integral of a measurable function can be constructed through an
approximation process, very much in the tradition of the Riemann
integral, but with a larger scope. We also go the other way: given an
integral, we construct a measure from it. This is the elegant way
P. J. Daniell did propose for constructing measures, and it can be
brought to fruit in this context for a fairly simple proof of the
Riesz Representation Theorem on compact metric spaces.

Having all these tools at our disposal, we look at product measures, which
can be introduced now through a kind of line sweeping --- if you want
to measure an area in the plane, you measure the line length as you
sweep over the area; this produces a function of the abscissa, which
then yields the area. One of the main tools here is Fubini's
Theorem. The product measure is not confined to two factors, we
discuss the general case. This includes a discussion of projective
systems, which may be considered as a generalization of sequences of
products. A case study shows that projective systems arise easily in
the study of continuous time stochastic logics.

Now that integrals are available, we turn back and have a look at
topologies on spaces of measures; one suggests itself --- the weak
topology which is induced by the continuous functions. This is related
to the Alexandrov topology. It is shown that there is a particularly
handy metric for the weak topology, and that the space of all finite
measures is complete with this metric, so that we now have a Polish
space. This is capitalized on when discussing selections for
set-valued maps into it, which are helpful in showing that Polish
spaces are closed under bisimulations. We use measurable selections
for an investigation into the structure of quotients in the Kleisli
monad, providing another example for the interplay of arguments from
measure theory and categories.

Finally, we take up a true classic: $\rLpS$-spaces. We start from
Hilbert spaces, apply the representation of linear functionals on
$\rLpS[2]$ to obtain the Radon-Nikodym Theorem through von Neumann's
ingenious proof, and derive from it the representation of the dual
spaces. This is applied to disintegration, where we show that a
measure on a product can be decomposed into a projection and a
transition kernel (on the surface this does not look like an
application area for $\rLpS$-spaces; the relationship derives from the
Radon-Nikodym Theorem).

Because we are driven by applications to Markov transition systems and
similar objects, we did not strive for the most general approach to
measure and integral. In particular, we usually formulate the results
for finite or $\sigma$-finite measures, leaving the more general cases
outside of our focus. This means also that we did not deal with
complex measures (and the associated linear spaces over the complex
numbers), but things are discussed in the realm of real numbers;
we show, however, in which way one could start to deal with complex
measures when the occasion arises. Of course, a lot of things had
to be left out, among them a more careful study of the Borel hierarchy
and applications to Descriptive Set Theory, as well as martingales.

%%% Local Variables: 
%%% mode: latex
%%% TeX-master: "../Mskr3"
%%% End: 

%\Input{\Folder/MeasSets}
%spell checked - 24Aug14
\Section{Measurable Sets and Functions}
\label{sec:measurable-sets-and-functions}

This section contains a systematic study of measurable spaces and
measurable functions with a view towards later developments. A brief overview is in order.

The measurable structure is
lifted to the space of finite measures, which form a measurable set
under the weak-*-$\sigma$-algebra. This is studied in
Section~\ref{sec:sigma-algebra-on-measures}. If the underlying space
carries a topology, the topological structure is handed down to finite
measures through the Alexandrov topology. We will have a look at it in
Section~\ref{sec:topology-on-measures}. The measurable functions from
a measurable space to the reals form a vector space, which is also a
lattice, and we will show that the step functions,i.e., those
functions which take only finite number of values, are dense with
respect to pointwise convergence. This mode of convergence is relaxed
in the presence of a measure in various ways to almost uniform
convergence, convergence almost everywhere, and to convergence in measure
(Sections~\ref{sec:ess-bounded-fncts} and~\ref{sec:conv-a-e}), from
which also various (pseudo-) metrics and norms may be derived.

If the underlying measurable spaces are the Borel sets of a
metric space, and if the metric has a countable dense set, then the
Borel sets are countably generated as well. But the irritating
observation is that being countably generated is not hereditary~---~a
sub-$\sigma$-algebra of a countable $\sigma$-algebra needs not be
countably generated. So countably generated $\sigma$-algebras deserve
a separate look, which is what we will do in
Section~\ref{sec:count-gen}. The very important class of Polish spaces
will be studied in this context as well, and we will show to
manipulate a Polish topology into making certain measurable functions
continuous. Polish spaces generalize to analytic spaces in a most
natural manner, for example when taking the factor of a countably
generated equivalence relation in a Polish space; we will study the
relationship in Section~\ref{sec:borel-and-analytic}.The most
important tool here is Souslin's Separation Theorem. This discussion
leads quickly to a discussion of the abstract Souslin operation in
Section~\ref{sec:souslin-op}, through which analytic sets may be
generated in a Polish space. From there it is but a small step to introducing
universally measurable sets in Section~\ref{sec:uni-meas}, which turn
out to be closed under Souslin's operation in general measurable
spaces. Two applications of these techniques are given: Lubin's
Theorem extends a measure from a countably generated
sub-$\sigma$-algebra of the Borel sets of an analytic space to the
Borel sets proper, and we show that a transition kernel can be extended to
the universal completion (Sections~\ref{sec:lubin-extension} and
\ref{sec:ext-stoch-rel-to-compl}). Lubin's Theorem is established
through von Neumann's Selection Theorem, which provides a universally
measurable right inverse to a surjective measurable map from an
analytic space to a separable measurable space. The topic of
selections is taken up in Section~\ref{sec:meas-selections}, where the
selection theorem due to Kuratowski and Ryll-Nardzewski is in the
center of attention. It gives conditions under which a map which takes
values in the closed non-empty subsets of a Polish space has a
measurable selector. This is of interest, e.g., when it comes to
establish the existence of bisimulations for Markov transition
systems, or for identifying the quotient structure of transition
kernels.

\Subsection{Measurable Sets}
\label{sec:measurable-sets}

Recall that a measurable space\index{space!measurable} $(X, \mathcal{A})$
consists of a set $X$ with a $\sigma$-algebra $\mathcal{A}$, which
is an Boolean algebra of subsets of $X$ that is closed under countable
unions (hence countable intersections or countable disjoint unions).
If $\mathcal{A}_0$ is a family of subsets of $X$, then
\begin{equation*}
\sigma\left(\mathcal{A}_0\right) = \bigcap\{\mathcal{B} \mid
\mathcal{B} \text{ is a  $\sigma$-algebra on $M$ with }
\mathcal{A}_0 \subseteq \mathcal{A}\}
\end{equation*}
is the smallest $\sigma$-algebra on $M$ which contains
$\mathcal{A}_0$. This construction works since the power set
$\mathcal{P}(X)$ is a $\sigma$-algebra on $X$. Take for example as a
generator $\mathcal{I}$ all open intervals in the real numbers
$\Real$, then $\sigma(\mathcal{I})$ is the $\sigma$-algebra of
real \emph{Borel sets}\index{space!Polish!Borel sets}. These Borel sets
are denoted by $\Borel{\Real}$, and, since each open subset of $\Real$
can be represented as a countable union of open intervals,
$\Borel{\Real}$ is the smallest $\sigma$-algebra which contains the
open sets of $\Real$. Unless otherwise stated, the real numbers are
equipped with the $\sigma$-algebra $\Borel{\Real}$.

In general, if $(X, \tau)$ is a topological space, the
$\sigma$-algebra $\Borel{\tau} := \sigma(\tau)$ are called its
\emph{Borel sets}\index{Borel sets}\index{Borel
  sets}\index{space!topological!Borel sets}\index{topology!Borel
  sets}. They will be discussed extensively in the context of Polish
spaces. This is, however, not the only $\sigma$-algebra of interest on
a topological space.

\BeginExample{ex-baire-sets}
Call $F\subseteq
X$ \emph{functionally closed} iff $F = \InvBild{f}{\{0\}}$ for some continuous
function $f: X\to \Real$, $G\subseteq X$ is called functionally open
iff $G = X\setminus F$ with $F$ functionally closed. The \emph{Baire
  sets} $\Baire{X}$ of
$(X, \tau)$ are the $\sigma$-algebra generated by the functionally
closed sets of the space. \index{Baire
  sets}\index{space!topological!Baire sets}\index{topology!Baire
  sets} 

If $(X, d)$ is a metric space, let $F\subseteq X$ be closed, then 
\begin{equation*}
d(x, F) := \inf\{d(x, y) \mid  y\in F\}
\end{equation*}
is the distance of $x$ to $F$ with $x\in F$ iff $d(x, F) = 0$.
Moreover, $d(\cdot, F)$ is continuous, thus $F = \InvBild{d(\cdot,
  F)}{\{0\}}$ is functionally closed, hence the Baire and the Borel
sets coincide for metric spaces.

Note that 
$
|d(x, F) - d(y, F)| \geq d(x, y),
$
so that $d(\cdot, F)$ is even uniformly continuous.
\EndExample

The next example constructs a $\sigma$-algebra which comes up quite naturally in the study of stochastic nondeterminism.

\BeginExample{hit-sigma-algebra}
Let ${\cal A}\subseteq \PowerSet{X}$ for some set $X$, the family of
hit sets, and ${\cal G}$ a distinguished subsets of
$\PowerSet{X}$. Define the \index{hit!$\sigma$-algebra}\emph{hit-$\sigma$-algebra} ${\cal
  H}_{{\cal A}}({\cal G})$ as the smallest
$\sigma$-algebra on ${\cal G}$ which contains all the sets $H_{A}$ with $A\in{\cal
  A}$, where $H_{A}$ is the hit set associated with $A$, i.e., $H_{A}
:= \{B\in{\cal G}\mid B\cap A \not=\emptyset\}$. 
\EndExample

Rather than working with a closure operation $\sigma(\cdot)$, one
sometimes can adjoin additional elements to obtain a $\sigma$-algebra
from a given one, see also
Exercise~\ref{ex-extend-top-and-sigma}. This is demonstrated for a
$\sigma$-ideal through the following construction, which will be
helpful when completing a measure space. Recall that ${\cal
  N}\subseteq\PowerSet{X}$ is a $\sigma$-ideal iff it is an order
ideal which is closed under countable unions~\SetCite{Definition
  2.91}.

\BeginLemma{adjoin-sigma-ideal}
Let ${\cal A}$ be a $\sigma$-algebra on a set $X$, ${\cal
  N}\subseteq\PowerSet{X}$ a $\sigma$-ideal. Then 
\begin{equation*}
  {\cal A}_{{\cal N}} := \{A\Delta N \mid A\in {\cal A}, N\in{\cal N}\}
\end{equation*}
is the smallest $\sigma$-algebra containing both ${\cal A}$ and ${\cal N}$.
\EndLemma

\BeginProof
Is is sufficient to demonstrate that ${\cal A}_{{\cal N}}$ is a
$\sigma$-algebra. Because 
\begin{equation*}
X\setminus (A\Delta N) = X\Delta(A\Delta N)
= (X\Delta A)\Delta N = (X\setminus A)\Delta N,
\end{equation*}
we see that ${\cal A}_{{\cal N}}$ is closed under complementation. Now
let $\bigl(A_{n}\Delta N_{n}\bigr)_{n\in\Nat}$ be a sequence of sets
with $\Folge{A}$ in ${\cal A}$ and $\Folge{N}$ in ${\cal N}$, we have 
\begin{equation*}
  \bigcup_{n\in\Nat}(A_{n}\Delta N_{n}) =
  \bigl(\bigcup_{n\in\Nat}A_{n}\bigr)\Delta N
\end{equation*}
with 
\begin{equation*}
  N = \bigcup_{n\in\Nat}(A_{n}\Delta N_{n})\Delta
  \bigl(\bigcup_{n\in\Nat}A_{n}\bigr)
\stackrel{(\ddag)}{\subseteq}\bigcup_{n\in\Nat}\bigl((A_{n}\Delta N_{n})\Delta A_{n}\bigr)
= \bigcup_{n\in\Nat}N_{n},
\end{equation*}
using Exercise~\ref{char-fnct} in $(\ddag)$.
Because ${\cal N}$ is a $\sigma$-ideal, we conclude that $N\in{\cal
  N}$. Thus ${\cal A}_{{\cal N}}$ is also closed under countable
unions. Since $\emptyset, X\in{\cal A}_{{\cal N}}$, we conclude
that this set is a $\sigma$-algebra indeed. 
\EndProof
%%%

{%\setboolean{isBook}{true}
\ifthenelse{\boolean{isBook}}{
}{
It turns out to be
most convenient to have a closer look at the construction of
$\sigma$-algebras when the family of sets we start from has already
some structure. This gives the occasion to introduce Dynkin's
\index{$\pi$-$\lambda$-Theorem.}$\pi$-$\lambda$-Theorem. This is an important tool, which eases
sometimes the task of identifying the $\sigma$-algebra generated from
some family of sets.

\BeginTheorem{Pi-Lambda}{\textbf{($\pi$-$\lambda$-Theorem)}}
Let $\mathcal{P}$ be a family of
subsets of $S$ that is closed under finite intersections (this is
called a
\emph{$\pi$-class}). Then $\sigma(\mathcal{P})$ is the smallest
$\lambda$-class containing $\mathcal{P}$, where a family
$\mathcal{L}$ of subsets of $S$ is called a \emph{$\lambda$-class}
iff it is closed under complements and countable disjoint unions.
\EndTheorem

\BeginProof
1.
Let $\mathcal{L}$ be the smallest $\lambda$-class containing $P$, then
we show that $\mathcal{L}$ is a $\sigma$-algebra.

2.
We show first that
it is an algebra. Being a $\lambda$-class, $\mathcal{L}$ is closed under
complementation. Let
$A \subseteq S$, then
$
\mathcal{L}_A := \{B \subseteq S \mid A \cap B \in \mathcal{L}\}
$
is a $\lambda$-class again: if $A \cap B \in \mathcal{L}$, then
$$
A \cap \left(S\setminus B\right) = A\setminus B = S\setminus ((A \cap B) \cup (S\setminus A)),
$$
which is in $\mathcal{L}$, since $(A \cap B) \cap S\setminus A = \emptyset$, and
since $\mathcal{L}$ is closed under disjoint unions.

If $A \in \mathcal{P}$, then $\mathcal{P} \subseteq \mathcal{L}_A$, because
$\mathcal{P}$ is closed under intersections. Because $\mathcal{L}_A$ is a
$\lambda$-system, this implies
$
\mathcal{L} \subseteq\mathcal{L}_A
$
for all $A \in \mathcal{P}$. Now take $B \in \mathcal{L}$, then the
preceding argument shows that $\mathcal{P} \subseteq \mathcal{L}_B$,
and again we may conclude that $\mathcal{L} \subseteq \mathcal{L}_B$. Thus
we have shown that $A \cap B \in \mathcal{L}$, provided $A, B \in \mathcal{L}$, so that
$\mathcal{L}$ is closed under finite intersections. Thus $\mathcal{L}$ is a Boolean algebra.

3.
$\mathcal{L}$ is a $\sigma$-algebra as well. It is enough to show that
$\mathcal{L}$ is closed under countable unions. But since
$$
\bigcup_{n \in \Nat} A_n = \bigcup_{n \in \Nat} \left(A_n \setminus \bigcup_{i=1}^{n-1}A_i\right),
$$
this follows immediately.
\EndProof

Consider an immediate and fairly typical application. It states that
two finite measures are equal on a $\sigma$-algebra, provided they are
equal on a generator which is closed under finite intersections. The
proof technique is worth noting: We collect all sets for which the
assertion holds into one family of sets and investigate its
properties, starting from an originally given set. If we find that the
family has the desired property, then we look at the corresponding
closure. To be specific, have a look at the  proof of the following statement.

\BeginLemma{are-equal}
Let $\mu, \nu$ be finite measures on a $\sigma$-algebra $\sigma({\cal B})$,
where ${\cal B}$ is a family of sets which
is closed under finite intersections. Then
$\mu(A) = \nu(A)$ for all $A\in\sigma({\cal B})$, provided $\mu(B) = \nu(B)$ for
all $B\in{\cal B}$. 
\EndLemma

\BeginProof
We have a look at all sets for which the assertion is true, and
investigate this set. Put
\begin{equation*}
  {\cal G} := \{A\in\sigma({\cal B}) \mid \mu(A) = \nu(A)\},
\end{equation*}
then ${\cal G}$ has these properties:
\begin{itemize}
\item ${\cal B}\subseteq{\cal G}$ by assumption.
\item Since ${\cal B}$ is closed under finite intersections,
  $S\in{\cal B}\subseteq{\cal G}.$
\item ${\cal G}$ is closed under complements.
\item ${\cal G}$ is closed under countable disjoint unions; in fact,
  let $\Folge{A}$ be a sequence of mutually disjoint sets in ${\cal
    G}$ and $A := \bigcup_{n\in\Nat} A_n$, then
  \begin{equation*}
    \mu(A) = \sum_{n\in\Nat} \mu(A_n) = \sum_{n\in\Nat} \nu(A_n) = \nu(A),
  \end{equation*}
hence $A\in{\cal G}$. 
\end{itemize}
But this means that ${\cal G}$ is a $\lambda$-class containing ${\cal
  B}$. But the smallest $\lambda$-class containing ${\cal G}$ is
$\sigma({\cal B})$ by Theorem~\ref{Pi-Lambda}, so that we have now
\begin{equation*}
  \sigma({\cal B}) \subseteq {\cal G} \subseteq \sigma({\cal B}),
\end{equation*}
the last inclusion coming from the definition of ${\cal G}$. Thus we
may conclude that ${\cal G} = \sigma({\cal B})$, hence all sets in
$\sigma({\cal B})$ have the desired property.
\EndProof
}
}

If $(Y, \mathcal{B})$ is another measurable space, then a map $ f: X
\rightarrow Y $ is called
\emph{$\mathcal{A}$-$\mathcal{B}$-measurable}\index{map!measurable}
iff the inverse image under $f$ of each set in $\mathcal{B}$ is a
member of $\mathcal{A}$, hence iff $ \InvBild{f}{G} \in \mathcal{A}
$ holds for all $G \in \mathcal{B}$.

Checking measurability is made easier by the observation that
it suffices for the inverse images of a generator to be measurable sets.
\BeginLemma{ElemGen}
Let $(X, \mathcal{A})$ and $(Y, \mathcal{B})$ be measurable spaces,
and assume that $\mathcal{B} = \sigma(\mathcal{B}_0)$ is generated by
a family $\mathcal{B}_0$ of subsets of $Y$. Then $f: X \rightarrow Y$ is
$\mathcal{A}$-$\mathcal{B}$-measurable iff $ \InvBild{f}{G} \in
\mathcal{A} $ holds for all $G \in \mathcal{B}_0.$
\EndLemma

\BeginProof
Clearly, if $f$ is $\mathcal{A}$-$\mathcal{B}$-measurable, then $ \InvBild{f}{G} \in
\mathcal{A} $ holds for all $G \in \mathcal{B}_0.$

Conversely, suppose $ \InvBild{f}{G} \in
\mathcal{A} $ holds for all $G \in \mathcal{B}_0,$ then we need to show
that $\InvBild{f}{G} \in \mathcal{A}$ for all $G \in \mathcal{B}$.
In fact, consider the set $\mathcal{G}$ for which the assertion is true,
\begin{equation*}
\mathcal{G} := \{G \in \mathcal{B} \mid \InvBild{f}{G} \in \mathcal{A}\}.
\end{equation*}
An elementary calculation shows that the empty set and $Y$ are both members of $\mathcal{G}$,
and since $\InvBild{f}{Y\setminus G} = X\setminus \InvBild{f}{G}$,
$\mathcal{G}$ is closed under complementation. Because
\begin{equation*}\textstyle
\InvBild{f}{\bigcup_{i \in I} G_i} = \bigcup_{i \in I} \InvBild{f}{G_i}
\end{equation*}
holds for any index set $I$, $\mathcal{G}$ is closed under finite and countable
unions. Thus $\mathcal{G}$ is a $\sigma$-algebra, so that
$\sigma(\mathcal{G}) = \mathcal{G}$ holds. By assumption, $\mathcal{B}_0 \subseteq \mathcal{G}$,
so that
\begin{equation*}
\mathcal{A} = \sigma(\mathcal{B}_0) \subseteq \sigma(\mathcal{G}) = \mathcal{G}\subseteq \mathcal{A}
\end{equation*}
is inferred. Thus all elements of $\mathcal{B}$ have their inverse image in $\mathcal{A}$.
\EndProof

An example is furnished by a real valued function $f: X \rightarrow
\Real$ on $X$ which is $\mathcal{A}$-$\Borel{\Real}$-measurable iff $
\{x \in X \mid f(x) \RelOp t\} \in \mathcal{A} $ holds for each $t \in
\Real$; the relation $\RelOp$ may be taken from $<, \leq, \geq, >.$ We
infer in particular that a function $f$ from an topological space $(X,
\tau)$ which is upper or lower \index{map!semicontinuous}semicontinuous
(i.e., for which in the \emph{upper semicontinuous} case the set
$\{x\in X \mid f(x) < c\}$ is open, and in the \emph{lower
  semicontinuous} case the set $\{x\in X\mid f(x) > c\}$ is open,
$c\in\Real$ being arbitrary), is Borel measurable. Hence a continuous
function is Borel measurable. A continuous function $f: X\to Y$ into a
metric space $Y$ is Baire measurable (Exercise~\ref{ex-cont-baire-meas}).

These observations will be used frequently.

The proof's strategy is to have a look at all objects that have the
desired property, and to show that this \emph{set of good guys} is a
$\sigma$-algebra (this is why this approach is sometimes called the
\emph{\index{principle of good sets}principle of good
  sets}~\cite{Elstrodt}). It is similar to showing in a proof by
induction that the set of all natural numbers having a certain
property is closed under constructing the successor.  Then we show
that the generator of the $\sigma$-algebra is contained in the good
guys, which is rather similar to begin the induction.  Taking both
steps together then yields the desired properties for both cases. We
will encounter this pattern of proof over and over again.

An example is furnished by the equivalence relation induced by a family of sets.

\BeginExample{prob-equiv-gen}
Given a subset ${\cal C}\subseteq\PowerSet{X}$ for a set $X$, define the equivalence relation $\equiv_{{\cal C}}$ on $X$ upon setting 
\begin{equation*}
  \isEquiv{x}{x'}{\equiv_{{\cal C}}} \text{ iff } \forall C\in {\cal C}: x\in C \Leftrightarrow x'\in C.
\end{equation*}
Thus $\isEquiv{x}{x'}{\equiv_{\cal C}}$ iff ${\cal C}$ cannot separate the elements $x$ and $x'$; call $\equiv_{{\cal C}}$ the equivalence relation generated by ${\cal C}$. 

Now let ${\cal A}$ be a $\sigma$-algebra on $X$ with ${\cal A} = \sigma({\cal A}_{0})$. Then ${\cal A}$ and ${\cal A}_{0}$ generate the same equivalence relation, i.e., $\equiv_{{\cal A}}\ =\ \equiv_{{\cal A}_{0}}$. In fact, define for $x, x'\in X$ with $\isEquiv{x}{x'}{\equiv_{{\cal A}_{0}}}$
\begin{equation*}
  {\cal B} := \{A\in {\cal A}\mid x\in A \Leftrightarrow x'\in A\}
\end{equation*}
Then ${\cal B}$ is a $\sigma$-algebra with ${\cal A}_{0}\subseteq{\cal B}$, hence $\sigma({\cal A}_{0})\subseteq{\cal B}\subseteq{\cal A}$, so that ${\cal A} = {\cal B}$. Thus $\isEquiv{x}{x'}{\equiv_{{\cal A}_{0}}}$ implies $\isEquiv{x}{x'}{\equiv_{\cal A}}$; since the reverse implication is obvious, the claim is established. 
\EndExample

If $(X, \mathcal{A})$ is a measurable space and
$
f: X \rightarrow Y
$
is a map, then
\begin{equation*}
\mathcal{B} := \{D \subseteq Y \mid  \InvBild{f}{D} \in
\mathcal{A}\}
\end{equation*}
is the largest $\sigma$-algebra $\mathcal{B}_0$ on $N$ that renders$f$ $\mathcal{A}$-$\mathcal{B}_0$-measurable; then $\mathcal{B}$ is called the
\emph{final} $\sigma$-algebra\index{$\sigma$-algebra!final} with respect to 
$f$. In fact, because the inverse set operator $f^{-1}$ is compatible
with the Boolean operations, it is immediate that $\mathcal{B}$ is closed under
the operations for a $\sigma$-algebra, and a little moment's reflection
shows that this is also the largest $\sigma$-algebra with this property.

Symmetrically, let $ g: P \rightarrow X $ be a map, then
\begin{equation*}
\InvBild{g}{\mathcal{A}} := \{\InvBild{g}{E}\mid  E \in
\mathcal{A}\}
\end{equation*}
is the smallest $\sigma$-algebra $\mathcal{P}_0$ on $P$ that renders
$g: \mathcal{P}_0 \rightarrow \mathcal{A}$ measurable; accordingly,
$\InvBild{g}{\mathcal{M}}$ is called
\emph{initial}\index{$\sigma$-algebra!initial} with respect to $f$. Similarly,
$\InvBild{g}{\mathcal{A}}$ is a $\sigma$-algebra, and it is fairly
clear that this is the smallest one with the desired property. In particular,
the inclusion $i_Q: Q \rightarrow X$ becomes measurable for a subset
$Q \subseteq X$ when $Q$ is endowed with the $\sigma$-algebra
$
\{Q \cap B \mid B \in \mathcal{A}\}.
$
It is called the \emph{trace of $\mathcal{A}$ on $Q$}\index{$\sigma$-algebra!trace}
and is denoted ---~in a slight abuse of notation~--- by $\mathcal{A} \cap Q$.

Initial and final $\sigma$-algebras generalize in an obvious way to families of maps.
For example,
$
\sigma\left(\bigcup_{i \in I} \InvBild{g_i}{\mathcal{A}_i}\right)
$
is the smallest $\sigma$-algebra $\mathcal{P}_0$ on $P$ which makes
all the maps $g_i: P \rightarrow X_i$ $\mathcal{P}_0$-$\mathcal{A}_i$-measurable
for a family $\left((X_i, \mathcal{A}_i)\right)_{i \in I}$ of measurable spaces.

This is an intrinsic, universal characterization of the initial $\sigma$-algebra
for a single map.

\BeginLemma{final-char}
Let $(X, \mathcal{A})$ be a measurable space and $f: X \rightarrow Y$ be a map.
The following conditions are equivalent:
\begin{enumerate}
\item\label{L-final-char-1}The  $\sigma$-algebra $\mathcal{B}$ on $Y$ is final with respect to $f$.
\item\label{L-final-char-2}If $(P, \mathcal{P})$ is a measurable space, and
$g: Y \rightarrow P$ is a map, then the $\mathcal{A}$-$\mathcal{P}$-measurability of $g \circ f$
implies the $\mathcal{B}$-$\mathcal{P}$-measurability of $g$.
\end{enumerate}
\EndLemma

\BeginProof
1.
Taking care of $\labelImpl{L-final-char-1}{L-final-char-2}$, we note that
\begin{equation*}
\InvBild{(g \circ f)}{\mathcal{P}} = \InvBild{f}{\InvBild{g}{\mathcal{P}}} \subseteq \mathcal{A}.
\end{equation*}
Consequently, $\InvBild{g}{\mathcal{P}}$ is one of the $\sigma$-algebras $\mathcal{B}_0$ with
$\InvBild{f}{\mathcal{B}_0} \subseteq \mathcal{A}$. Since $\mathcal{B}$ is the largest of them, we have
$\InvBild{g}{\mathcal{P}} \subseteq \mathcal{B}$. Hence $g$ is $\mathcal{B}$-$\mathcal{P}$-measurable.

2.
In order to establish $\labelImpl{L-final-char-2}{L-final-char-1}$, we have to show that
$\mathcal{B}_0 \subseteq \mathcal{B}$ whenever $\mathcal{B}_0$ is a
$\sigma$-algebra on $\mathcal{Y}$ with $\InvBild{f}{\mathcal{B}_0} \subseteq \mathcal{A}$. Put
$(P, \mathcal{P}) := (Y, \mathcal{B}_0)$, and let $g$ be the identity $id_Y$. Because
$\InvBild{f}{\mathcal{B}_0} \subseteq \mathcal{A}$, we see that $id_Y \circ f$ is
$\mathcal{B}_0$-$\mathcal{A}$-measurable. Thus $id_Y$ is
$\mathcal{B}$-$\mathcal{B}_0$-measurable.
But this means $\mathcal{B}_0 \subseteq \mathcal{B}$.
\EndProof

We will use the final $\sigma$-algebra mainly for factoring through
an equivalence relation. In fact, let $\alpha$ be an equivalence relation on a set
$X$, where $(X, \mathcal{A})$ is a measurable space. Then the factor map
\begin{equation*}
\fMap{\alpha}:
\begin{cases}
  X & \rightarrow \Faktor{X}{\alpha}\\
  x & \mapsto \Klasse{x}{\alpha}
\end{cases}
\end{equation*}
that maps each element to its class can be made a measurable map by
taking the final $\sigma$-algebra $\Faktor{\mathcal{A}}{\alpha}$
with respect to $\fMap{\alpha}$ and $\mathcal{A}$
as the $\sigma$-algebra on $\Faktor{X}{\alpha}$.

Dual to Lemma~\ref{final-char}, the initial $\sigma$-algebra is characterized.
\BeginLemma{initial-char}
Let $(Y, \mathcal{B})$ be a measurable space and $f: X \rightarrow Y$ be a map.
The following conditions are equivalent:
\begin{enumerate}
\item\label{L-initial-char-1}The  $\sigma$-algebra $\mathcal{A}$ on $X$ is initial with respect to $f$.
\item\label{L-initial-char-2}If $(P, \mathcal{P})$ is a measurable space, and
$g: P \rightarrow X$ is a map, then the $\mathcal{P}$-$\mathcal{B}$-measurability of $f \circ g$
implies the $\mathcal{P}$-$\mathcal{A}$-measurability of $g$.
\end{enumerate}
\QED
\EndLemma

Let $\left((A_i, \mathcal{A}_i)\right)_{i \in I}$ be a family of measurable
spaces, then the product-$\sigma$-algebra $ \bigotimes_{i \in
I}\mathcal{A}_i $ denotes that initial $\sigma$-algebra
on $ \prod_{i \in I} X_i $ for the
projections
\begin{equation*}
\pi_j: \langle m_i \mid  i \in I\rangle \mapsto m_j.
\end{equation*}
It is not difficult to see that\index{$\sigma$-algebra!product}
$
\bigotimes_{i \in I} \mathcal{A}_i = \sigma(\mathcal{Z})
$
with\label{cylinder-sets}
\begin{equation*}
\mathcal{Z} := \{\prod_{i \in I} E_i
\mid  \forall i \in I: E_i \in \mathcal{M}_i, E_i = M_i \text{ for
almost all indices}\}
\end{equation*}
as the collection of \emph{cylinder sets}\index{set!cylinder}
(use Theorem~\ref{Pi-Lambda} and the observation that $\mathcal{Z}$ is closed under intersection).
% we will make frequent use of cylinders when dealing with infinite products
% for interpreting continuous time stochastic logics in Chapter~\ref{sec:interpr-modal-temp}.

For $I = \{1, 2\}$, the $\sigma$-algebra $
\mathcal{A}_1 \otimes \mathcal{A}_2 $ is generated from the set of
\emph{measurable rectangles}\index{measurable!rectangle}
$$
\{E_1 \times E_2\mid  E_1 \in \mathcal{A}_1, E_2 \in
\mathcal{A}_2\}.
$$

Dually, the sum\index{$\sigma$-algebra!sum} $(X_1 + X_2,
\mathcal{A}_1 + \mathcal{A}_2)$ of the measurable spaces $(X_1,
\mathcal{A}_1)$ and $(X_2, \mathcal{A}_2)$ is defined through the
final $\sigma$-algebra on the sum $X_1 + X_2$ for the injections
$
X_i \rightarrow X_1 + X_2.
$
This is the special case of the coproduct
$
\bigoplus_{i \in I} (X_i, \mathcal{A}_i),
$
where the $\sigma$-algebra $\bigoplus_{i \in I}\mathcal{A}_i$
is initial with respect to the injections.

\Subsubsection{A $\sigma$-Algebra On Spaces Of Measures}
\label{sec:sigma-algebra-on-measures}

We will now introduce a $\sigma$-algebra on the space of all $\sigma$-finite measures. It is induced by evaluating measures at fixed events. Note the inversion: instead of observing a measure assigning a real number to a set, we take a set and have it act on measures. This approach is fairly natural for many applications.

In addition to $\SubProbSenza$ resp. $\ProbSenza$, the functors which
assign to each measurable space its subprobabilities and its
probabilities (see~\CategCite{Section 1.4.2}), we introduce the space of finite resp. $\sigma$-finite
measures. Denote by \index{$\FinM{X, {\cal A}}$}$\FinM{X, {\cal A}}$
the set of all finite measures on $(X, {\cal A})$, the set of all
$\sigma$-finite measures is denoted by \index{$\SigmaM{X, {\cal
      A}}$}$\SigmaM{X, {\cal A}}$. Each set $A\in{\cal A}$ gives rise
to the evaluation map $ev_{A}: \mu\mapsto\mu(A)$; the
\emph{weak-*-$\sigma$-algebra} \index{$\schwach{X, {\cal A}}$}$\schwach{X,
  {\cal A}}$ on $\FinM{X, {\cal A}}$ is the initial $\sigma$-algebra
with respect to the family $\{ev_{A}\mid A\in{\cal A}\}$ (actually, it
suffices to consider a generator ${\cal A}_{0}$ of ${\cal A}$, see
Exercise~\ref{ex-weak-generator}). It is clear that we have
\begin{equation*}
  \schwach{X, {\cal A}} = \sigma(\{\theBeta{A}{\bowtie q}{{\cal A}}\mid A\in{\cal A}, q\in\pReal\})
\end{equation*}
when we define
\begin{equation*}
\theBeta{A}{\bowtie q}{{\cal A}} := \{\mu \in\FinM{X, {\cal A}}\mid \mu(A) \bowtie q\}.
\end{equation*}
Here $\bowtie$ is one of the relational operators $\leq, <, \geq, >$,
and it apparent that $q$ may be taken from the rationals. We will use
the same symbol\index{ $\betaSenza{{\cal A}}$} $\betaSenza{{\cal A}}$
when we refer to probabilities or subprobabilities, if no confusion
arises. Thus the base space from which the weak-*-$\sigma$-algebra
will be constructed should be clear from the context.

Let $(Y, {\cal B})$ be another measurable space, and let $f: X\to Y$ be ${\cal A}$-${\cal B}$-measurable. Define 
\begin{equation*}
  \FinM{f}(\mu)(B) := \mu(\InvBild{f}{B})
\end{equation*}
for $\mu\in\FinM{X, {\cal A}}$ and for $B\in{\cal B}$, then
$\FinM{f}(\mu)\in\FinM{Y, {\cal B}}$, hence $\FinM{f}: \FinM{X, {\cal
    A}}\to \FinM{Y, {\cal B}}$ is a map, and since 
\begin{equation*}
\InvBild{(\FinM{f})}{\theBeta{B}{\bowtie q}{{\cal B}}} =
\theBeta{\InvBild{f}{B}}{\bowtie q}{{\cal A}},
\end{equation*}
this map is
$\schwach{{\cal A}}$-$\schwach{{\cal B}}$-measurable. Thus $\FinSenza$ is
an endofunctor on the category of measurable spaces.

Measurable maps into $\SigmaM{\cdot}$ deserve special attention.

\BeginDefinition{transition-kernels}
Given measurable spaces $(X, {\cal A})$ and $(Y, {\cal B})$, an ${\cal A}$-$\schwach{{\cal B}}$ measurable map $K: X \to \SigmaM{Y, {\cal B}}$ is called a \emph{\index{kernel!transition}transition kernel} and denoted by $K: (X, {\cal A})\Trans (Y, {\cal B})$. 
\EndDefinition

A transition kernel $K: (X, {\cal A})\Trans (Y, {\cal B})$ models a situation in which each $x\in X$ is associated with a $\sigma$-finite measure $K(x)$ on $(Y, {\cal B})$. In a probabilistic setting, this may be interpreted as the probability that a system reacts on input $x$ with $K(x)$ as the probability distribution of its responses. For example, if $(X, {\cal A}) = (Y, {\cal B})$ is the state space of a probabilistic transition system, then $K(x)(B)$ is often interpreted as the probability that the next state is a member of measurable set $B$ after a transition from $x$. 

This is an immediate characterization of transition kernels.

\BeginLemma{char-trans-kernel}
$K: (X, {\cal A})\Trans (Y, {\cal B})$ is a transition kernel iff these conditions are satisfied
\begin{enumerate}
\item $K(x)$ is a $\sigma$-finite measure on $(Y, {\cal B})$ for each $x\in X$.
\item $x\mapsto K(x)(B)$ is a measurable function for each $B\in{\cal B}$. 
\end{enumerate}
\EndLemma

\BeginProof
If $K: (X, {\cal A})\Trans (Y, {\cal B})$, then $K(x)$ is a $\sigma$-finite measure on $(Y, {\cal B})$, and 
\begin{equation*}
  \{x\in X| K(x)(B) > q\} = \InvBild{K}{\theBeta{B}{> q}{{\cal B}}} \in {\cal A}.
\end{equation*}
Thus $x\mapsto K(x)(B)$ is measurable for all $B\in{\cal B}$. Conversely, if $x\mapsto K(x)(B)$ is measurable for $b\in{\cal B}$, then the above equation shows that $\InvBild{K}{\theBeta{B}{> q}{{\cal B}}} \in {\cal A}$, so $K: (X, {\cal A})\to \SigmaM{Y, {\cal B}}$ is  ${\cal A}$-$\schwach{{\cal B}}$ measurable by Lemma~\ref{ElemGen}. 
\EndProof

A special case of transition kernels are
\index{kernel!Markov}\emph{Markov kernels}, sometimes also called
\emph{\index{stochastic relations}stochastic relations}. These are
kernels the image of which is in $\SubProbSenza$ or in $\ProbSenza$, whatever
the case may be.

\BeginExample{modal-logic}
Transition kernels may be used for interpreting \index{logic!modal}modal logics. Consider this grammar for formulas
\begin{equation*}
  \phi ::= \top \mid \phi_{1}\wedge\phi_{2}\mid \Diamond_{q}\phi
\end{equation*}
with $q\in\Rational, q \geq 0$. The informal interpretation in a probabilistic transition system is that $\top$ always holds, and that $\Diamond_{q}\phi$ holds with probability not smaller that $q$ after a transition in a state in which formula $\phi$ holds. Now let $M: (X, {\cal A})\Trans (X, {\cal A})$ be a transition kernel, and define inductively
\begin{align*}
  \Gilt[\top]_{M} :=&\ X\\
\Gilt[\phi_{1}\wedge\phi_{2}]_{M} :=&\ \Gilt[\phi_{1}]_{M}\cap\Gilt[\phi_{2}]_{M}\\
\Gilt[\Diamond_{q}\phi]_{M} :=&\ \{x\in X \mid M(x)(\Gilt_{M})\geq q\}\\
=&\ \InvBild{M}{{\theBeta{B}{\geq q}{{\cal A}}}}
\end{align*}
It is easy to show by induction on the structure of the formula that the sets $\Gilt_{M}$ are measurable, since $M$ is a transition kernel, for a generalization, see Example~\ref{modal-logic-rev}. 
\EndExample

\Subsubsection{The Alexandrov Topology On Spaces of Measures}
\label{sec:topology-on-measures}

Given a topological space $(X, \tau)$, the Borel sets $\Borel{\tau} =
\sigma(\tau)$ and the Baire sets $\Baire{X}$ come for free as
measurable structures: $\Borel{\tau}$ the smallest $\sigma$-algebra on
$X$ that contains the open sets; measurability of maps with respect to
the Borel sets is referred to as \emph{Borel
  measurability}\index{Borel measurability}.  $\Baire{X}$ is the
smallest $\sigma$-algebra on $X$ which contains the functionally
closed sets; they provide yet another measurable structure on $(X,
\tau)$, this time involving the continuous real valued
functions. Since $\Borel{X} = \Baire{X}$ for a metric space by
Example~\ref{ex-baire-sets}, the distinction between these
$\sigma$-algebras vanishes, and the Borel sets as the $\sigma$-algebra
generated by the open sets dominate the scene.

We will now define a topology of spaces of measures on a topological
space in a similar way, and relate this topology to the weak-*-$\sigma$-algebra, for the
time being in a special case. Fix a Hausdorff space $(X, \tau)$; the
space will be specialized as the discussion proceeds. Define for the
functionally open sets $G_{1}, \dots, G_{n}$, the functionally close
sets $F_{1}, \dots, F_{n}$ and $\epsilon>0$ for $\mu_{0}\in\FinM{X,
  \Baire{X}}$ the sets
\begin{align*}
  W_{G_{1}, \dots, G_{n}, \epsilon}(\mu_{0}) 
& := 
\{\mu\in\FinM{X, \Baire{X}}\mid \mu(G_{i}) > \mu_{0}(G_{i}) -
\epsilon\text{ for } 1\leq i \leq n, |\mu(X)-\mu_{0}(X)|< \epsilon\},\\
  W_{F_{1}, \dots, F_{n}, \epsilon}(\mu_{0}) 
& := 
\{\mu\in\FinM{X, \Baire{X}}\mid \mu(F_{i}) < \mu_{0}(F_{i}) +
\epsilon\text{ for } 1\leq i \leq n, |\mu(X)-\mu_{0}(X)|< \epsilon\}
\end{align*}

The topology which has the sets $W_{G_{1},\dots, G_{n},
  \epsilon}(\mu_{0})$ as a basis is called the \emph{Alexandrov
  topology} or
\emph{\index{topology!Alexandrov}A-topology}~\cite[8.10~(iv)]{Bogachev}. The
A-topology is defined in terms of Baire sets rather than Borel sets of
$(X, \tau)$. This is so because the Baire sets provide a scenario
which take the continuous functions on $(X, \tau)$ directly into
account. This is in general not the case with the Borel sets, which
are defined purely in terms of set theoretic operations. But the
distinction vanishes when we turn to metric spaces, see
Example~\ref{ex-baire-sets}. Note also that we deal with finite
measures here.

\BeginLemma{a-top-is-hausdorff}
The A-topology on $\FinM{X, \Baire{X}}$ is Hausdorff. 
\EndLemma

\BeginProof
The family of functionally closed sets of $X$ is closed under finite
intersections, hence if two measure coincide on the functionally
closed sets, they must coincide on the Baire sets $\Baire{X}$ of $X$ by
the $\pi$-$\lambda$-Theorem~\ref{Pi-Lambda}. 
\EndProof

Convergence in the A-topology is easily characterized in terms of
functionally open or closed sets. Recall that for a sequence
$\Folge{c}$ of real numbers the statements $\limsup_{n\to \infty}
c\leq c$ is equivalent to $\inf_{n\in\Nat}\sup_{k\geq n}c_{k}\leq c$
which in turn is equivalent to $\forall \epsilon>0\exists
n\in\Nat\forall k\geq n: c_{k} < c + \epsilon$. Similarly for
$\liminf_{n\to \infty}c_{n}$. This proves:

\BeginProposition{convergence-a-topology}
Let $\Folge{\mu}$ be a sequence of measures in $\FinM{X,
  \Baire{X}}$, then the following statements are equivalent.
\begin{enumerate}
\item\label{conv-A-1} $\mu_{n}\to \mu$ in the A-topology.
\item\label{conv-A-2} $\limsup_{n\to \infty}\mu_{n}(F) \leq \mu(F)$ for each
  functionally closed set $F$, and $\mu_{n}(X)\to \mu(X)$.
\item\label{conv-A-3} $\liminf_{n\to \infty}\mu_{n}(G) \geq \mu(G)$ for each
  functionally open set $G$, and $\mu_{n}(X)\to \mu(X)$.
\end{enumerate}
\QED
\EndProposition

This criterion is sometimes a little impractical, since it deals with
inequalities. We could have equality in the limit for all those sets for which the
boundary has $\mu$-measure zero, but, alas, the boundary may not be
Baire measurable. So we try with an approximation --- we approximate a
Baire set from within by a functionally open set (corresponding to the
interior) and from the outside by a closed set (corresponding to the
closure). This is discussed in some detail now.

Given $\mu\in\FinM{X, \Baire{X}}$, define by ${\cal R}_{\mu}$ all
those Baire sets which have a functional boundary of vanishing
$\mu$-measure, formally
\begin{equation*}
  {\cal R}_{\mu} := \{E\in\Baire{X}\mid G\subseteq E \subseteq F,
  \mu(F\setminus G) = 0, G\text{ functionally open}, F\text{
    functionally closed}\}.
\end{equation*}
Hence if $X$ is a metric space, $E\in{\cal R}_{\mu}$ iff $\mu(\partial
E) = 0$ for the boundary
$\partial E$ of $E$.

This is another criterion for convergence in the A-topology.

\BeginCorollary{conv-a-top-boundary}
Let $\Folge{\mu}$ be a sequence of Baire measures. Then $\mu_{n}\to
\mu$ in the A-topology iff $\mu_{n}(E) \to \mu(E)$ for all $E\in{\cal R}_{\mu}$. 
\EndCorollary

\BeginProof
The condition is necessary by
Proposition~\ref{convergence-a-topology}. Assume, on the other hand,
that $\mu_{n}(E) \to \mu(E)$ for all $E\in{\cal R}_{\mu}$, and take a
functionally open set $G$. We find $f: X\to  \Real$ continuous such
that $G = \{x\in X\mid f(x) > 0\}$. Fix $\epsilon>0$, then we can find
$c>0$ such that 
\begin{align*}
\mu(G) & < \mu(\{x\in X \mid f(x) > c\}) + \epsilon,\\
\mu(\{x\in X\mid f(x) > c\}) & = \mu(\{x\in X\mid f(x)\geq
c\}).
\end{align*}
Hence $E := \{x\in X\mid f(x) > c\}\in{\cal R}_{\mu}$, since
$E$ is open and $F := \{x\in X\mid f(x)\geq c\}$ is closed with
$\mu(F\setminus E) = 0$. So $\mu_{n}(E) \to \mu(E)$, by assumption,
and
\begin{equation*}
  \liminf_{n\to \infty}\mu_{n}(G) \geq \lim_{n\to \infty}\mu_{n}(E) =
  \mu(E) > \mu(G) - \epsilon.
\end{equation*}
Since $\epsilon>0$ was arbitrary, we infer
$
\liminf_{n\to \infty}\mu_{n}(G) \geq \mu(G).
$
Because $G$ was an arbitrary functionally open set, we infer from
Proposition~\ref{convergence-a-topology} that $\Folge{\mu}$ converges
in the A-topology to $\mu$. 
\EndProof

The family ${\cal R}_{\mu}$ has some interesting properties, which will
be of use later on, because, as we will show in a moment, it contains
the basis for the topology. This holds whenever there are enough continuous
functions to separate points from closed sets not containing
them. Before we can state this property, which will be helpful in the
analysis of the A-topology below, we introduce $\mu$-atoms, which are
of interest for themselves (we will define later, in
Definition~\ref{probs-atom}, atoms on a strictly order theoretic
basis, without reference to measures).

\BeginDefinition{measure-atom}
A set $A\in{\cal A}$ is called an $\mu$-atom iff $\mu(A) > 0$, and if
$\mu(B)\in\{0, \mu(A)\}$ for every $B\in{\cal A}$ with $B\subseteq A$.
\EndDefinition

Thus a $\mu$-atom does not permit values other than $0$ and $\mu(A)$
for its measurable subsets, so two different $\mu$-atoms $A$ and $A'$
are essentially disjoint, since $\mu(A\cap A') = 0$.

\BeginLemma{only-countably-many-atoms}
For the finite measure space $(X, {\cal A}, \mu)$ there exists an at
most countable set $\{A_{i}\mid i\in I\}$ of atoms such that
$X\setminus\bigcup_{i\in I}A_{i}$ is free of $\mu$-atoms.
\EndLemma

\BeginProof
If we do not have any atoms, we are done. Otherwise, let $A_{1}$ be an
arbitrary atom. This is the beginning; proceeding inductively, assume
that the atoms $A_{1},\dots, A_{n}$ are already selected, and let 
$
{\cal A}_{n} := \{A\in {\cal A} \mid A\subseteq
X\setminus\bigcup_{i=1}^{n}A_{i}\text{ is an atom}\}.  
$
If ${\cal A}_{n}=\emptyset$, we are done. Otherwise select the atom $A_{n+1}\in{\cal A}_{n}$ with 
$\mu(A_{n+1}) \geq \frac{1}{2}\cdot\sup_{A\in{\cal A}_{n}}\mu(A).$ Observe that $A_{1}, \dots, A_{n+1}$ are mutually disjoint. 

Let $\{A_{i}\mid i\in I\}$ be the set of atoms selected in this way, after the selection has terminated. Assume that $A\subseteq X\setminus\bigcup_{i\in I}A_{i}$ is an atom, then the index set $I$ must be infinite, and $\mu(A_{i})\geq \mu(A)$ for all $i\in I$. But since 
$
\sum_{i\in I}\mu(A_{i}) \leq \mu(X) < \infty,
$
we conclude that $\mu(A_{i})\to 0$, consequently, $\mu(A) = 0$, hence $A$ cannot be a $\mu$-atom.
\EndProof

This is a useful consequence.

\BeginCorollary{countably-many-point-values}
Let $f: X\to \Real$ be a continuous function. Then there are at most countably many $r\in\Real$ such that $\mu(\{x\in X \mid f(x) = r\})>0$.
\EndCorollary

\BeginProof
Consider the image measure $\FinM{f}(\mu): B \mapsto
\mu(\InvBild{f}{B})$ on $\Borel{\Real}$. If $\mu(\{x\in X\mid f(x) =
r\}) > 0$, then $\{r\}$ is a $\FinM{f}(\mu)$-atom. By
Lemma~\ref{only-countably-many-atoms} there are only countably many $\FinM{f}(\mu)$-atoms.
\EndProof

Returning to ${\cal R}_{\mu}$ we are now in a position to have a
closer look at its structure.

\BeginProposition{without-boundary-has-basis}
${\cal R}_{\mu}$ is a Boolean algebra. If $(X, \tau)$ is completely
regular, then ${\cal R}_{\mu}$ contains a basis for the topology
$\tau$.
\EndProposition

\BeginProof
It is immediate that ${\cal R}_{\mu}$ is closed under complementation,
and it is easy to see that it is closed under finite unions.

Let $f: X\to \Real$ be continuous, and define $U(f, r) := \{x\in X
\mid f(x) > r\}$, then $U(f, r)$ is open, and $\partial U(f, r)
\subseteq \{x\in X \mid f(x) = r\}$, thus $M_{f} := \{r\in \Real \mid
\mu(\partial U(f, r)) > 0\}$ is at most countable, such that the sets
$U(f, r)\in{\cal R}_{\mu}$, whenever $r\not\in M_{f}$. 

Now let $x\in X$ and $G$ be an open neighborhood of $x$. Because $X$ is
completely regular, we can find $f: x\to [0, 1]$ continuous such that
$f(y) = 1$ for all $y\not\in G$, and $f(x) = 0$. Hence we can find
$r\not\in M_{f}$ such that $x\in U(f, r)\subseteq G$. So ${\cal
  R}_{\mu}$ is in fact a basis for the topology.
\EndProof

Under the conditions above, ${\cal R}_{\mu}$ contains a base for $\tau$, we lift this base to $\FinM{X, \Baire{X}}$ in the hope of obtaining a base for the A-topology. This works, as we will show now.  

\BeginCorollary{a-top-basis-without-bound}
Let $X$ be a completely regular topological space, then the A-topology has a basis consisting of sets of the form
\begin{equation*}
  Q_{A_{1}, \dots, A_{n}, \epsilon}(\mu) := \{\nu\in\FinM{X, \Baire{X}} \mid  |\mu(A_{i}) - \nu(A_{i})| < \epsilon\text{ for } i = 1, \dots, n\}
\end{equation*}
with $\epsilon>0$, $n\in\Nat$ and $A_{1}, \dots, A_{n}\in{\cal R}_{\mu}$. 
\EndCorollary

\BeginProof
Let $W_{G_{1}, \dots, G_{n}, \epsilon}(\mu)$ with functionally open sets $G_{1}, \dots, G_{n}$ and $\epsilon>0$ be given. Select $A_{i}\in{\cal R}_{\mu}$ functionally open with $A_{i}\subseteq G_{i}$ and $\mu(A_{i}) > \mu(G_{i}) - \epsilon/2$, then it is easy to see that 
$
Q_{A_{1}, \dots, A_{n}, \epsilon/2}(\mu) \subseteq W_{G_{1}, \dots, G_{n}, \epsilon}(\mu).
$
\EndProof

We will specialize the discussion now to metric spaces. So fix a
metric space $(X, d)$, which we may assume to be bounded (otherwise we
switch to the equivalent metric $\langle x, y\rangle \mapsto d(x,
y)/(1-d(x, y))$). Recall that the $\epsilon$-neighborhood
$B^{\epsilon}$ of a set $B\subseteq X$ is defined as $ B^{\epsilon} :=
\{x\in X \mid d(x, B) < \epsilon\}.  $ Thus $B^{\epsilon}$ is always
an open set. Since the Baire and the Borel sets coincide in a metric
space (see Example~\ref{ex-baire-sets}), the A-topology is defined on
$\FinM{X, \Borel{X}}$, and we will relate it to a metric now.

Define the
\emph{\index{distance!Lévy-Prohorov}Lévy-Prohorov distance}
$d_{P}(\mu, \nu)$ of the
measures $\mu, \nu\in\FinM{X, \Borel{X}}$ through
\begin{equation*}
  d_{P}(\mu, \nu) := \inf\bigl\{\epsilon>0\mid \nu(B) \leq
  \mu(B^{\epsilon}) + \epsilon, \mu(B) \leq \nu(B^{\epsilon}) +
  \epsilon\text{ for all } B\in\Borel{X}\bigr\}
\end{equation*}

We note first that $d_{P}$ defines a metric, and that we can find a
metrically exact copy of the base space $X$ in the space $\FinM{X,
  \Borel{X}}$.

\BeginLemma{levy-prohorov-is-metric}
$d_{P}$ is a metric on $\FinM{X, \Borel{X}}$. $X$ is isometrically
isomorphic to the set $\{\delta_{x}\mid  x\in X\}$ of Dirac measures. 
\EndLemma

\BeginProof
It is clear that $d_{P}(\mu, \nu) = d_{P}(\nu, \mu)$. Let $d_{P}(\mu,
\nu) = 0$, then $\mu(F) \leq \nu(F^{1/n}) + 1/n$ and $\nu(F) \leq
\mu(F^{1/n})+1/n$ for each closed set $F\subseteq X$, hence $\nu(F) =
\mu(F)$ (note that $F^{1}\supseteq F^{1/2}\supseteq F^{1/3}\supseteq \dots$ and $F = \bigcap_{n\in\Nat}F^{1/n}$). Thus $\mu =
\nu$. If we have for all $B\in\Borel{X}$ that $\nu(B) \leq \mu(B^{\epsilon}) + \epsilon$, $\mu(B) \leq
\nu(B^{\epsilon})+\epsilon$ and $\mu(B) \leq \rho(B^{\delta}) +
\delta$, $\rho(B) \leq m(B^{\delta}) +\delta$, then $\nu(B) \leq
\rho(B^{\epsilon+\delta})+\epsilon+\delta$ and
$\rho(B) \leq \nu(B^{\epsilon+\delta}) + \epsilon+\delta$, thus
$d_{P}(\mu, \nu) \leq d_{P}(\mu, \rho) + d_{P}(\rho, \nu)$. We also
have $d_{P}(\delta_{x}, \delta_{y}) = d(x, y)$, from which the
isometry derives. 
\EndProof

We will relate the metric topology to the A-topology now. Without additional assumptions this relationship can be stated:

\BeginProposition{A-coarser-than-metric}
Each open set in the A-topology is also metrically open, hence A-topology is coarser than the metric topology.
\EndProposition

\BeginProof
Let $W_{F_{1}, \dots, F_{n}, \epsilon}(\mu)$ be an open basic
neighborhood of $\mu$ in the A-topology with $F_{1}, \dots, F_{n}$
closed. We want to find an open metric neighborhood with center $\mu$ which is contained in this A-neighborhood.

Because $(F^{1/n})_{n\in \Nat}$ is a decreasing sequence with
$\inf_{n\in\Nat}\mu(F_{n}) = \mu(F)$, whenever $F$ is closed, we can
find $\delta>0$ such that $\mu(F_{i}^{\delta}) < \mu(F_{i})+\epsilon/2$ for $1
\leq i \leq n$ and $0 < \delta < \epsilon/2$. Thus, if $d_{P}(\mu,
\nu) < \delta$, we have for $i = 1, \dots, n$ that
$
\nu(F_{i}) < \mu(F_{i}^{\delta}) + \delta < \mu(F_{i}) +\epsilon.
$
But this means that $\nu\in W_{F_{1}, \dots, F_{n},
  \epsilon}(\mu)$. 

Thus each neighborhood in the A-topology contains in fact an
open ball for the $d_{P}$-metric. 
\EndProof

The converse of Proposition~\ref{A-coarser-than-metric} can only be established under additional conditions,
which, however, are met for separable metric spaces. It is a
generalization of $\sigma$-continuity: while the latter deals
with sequences of sets, the concept of $\tau$-regularity deals with
the more general notion of directed families of open sets (recall that a
family ${\cal M}$ of sets is called \emph{directed} iff given $M_{1},
M_{2}\in{\cal M}$ there exists $M'\in{\cal M}$ with $M_{1}\cup
M_{2}\subseteq M'$).  

\BeginDefinition{tau-continuous}
A measure $\mu\in \FinM{X, \Borel{X}}$ is called
\emph{\index{$\tau$-regular}$\tau$-regular} iff 
\begin{equation*}
  \mu(\bigcup {\cal G}) = \sup_{G\in{\cal G}}\mu(G)
\end{equation*}
for each directed family ${\cal G}$ of open sets.
\EndDefinition

It is clear that we restrict our attention to open sets, because the
union of a directed family of arbitrary measurable sets is not necessarily
measurable. It is also clear that the condition above is satisfied for
countable increasing sequences of open sets, so that $\tau$-regularity
generalizes $\sigma$-continuity.

It turns out that finite measures on separable metric spaces are $\tau$-regular. Roughly speaking, this is due to the fact that countably many open sets determine the family of open sets, so that the space cannot be too large when looked at as a measure space. 

\BeginLemma{regular-tau-in-separable}
Let $(X, d)$ be a separable metric space, then each $\mu\in\FinM{X,
  \Borel{X}}$ is $\tau$-regular. 
\EndLemma

\BeginProof
Let ${\cal G}_{0}$ be a countable basis for the metric topology. If
${\cal G}$ is a directed family of open sets, we can find for each
$G\in {\cal G}$ a countable cover $(G_{i})_{i\in I_{G}}$ from ${\cal
  G}_{0}$ with $G = \bigcup_{i\in I_{G}}G_{i}$ and $\mu(G) =
\sup_{i\in I_{G}}\mu(G_{i})$. Thus
\begin{equation*}\textstyle
  \mu(\bigcup{\cal G}) = \sup\mu\bigl(\{\mu(G) \mid G\in{\cal G}_{0},
  G\subseteq \bigcup{\cal G}\}\bigr) = \sup_{G\in{\cal G}}\mu(G).
\end{equation*}
\EndProof

As a trivial consequence it is observed that $\mu(\bigcup{\cal G}) = 0$, where ${\cal G}$ is the family of all open sets $G$ with $\mu(G) = 0$.  

The important observation for our purposes is that a $\tau$-regular
measure is supported by a closed set which in terms of $\mu$ can be
chosen as being as tightly fitting as possible.

\BeginLemma{support-for-regular}
Let $(X, d)$ be a separable metric space. Given
$\mu\in\FinM{X, \Borel{X}}$ with $\mu(X) > 0$, there exists a smallest closed set $C_{\mu}$ such that
$\mu(C_{\mu}) = \mu(X)$. $C_{\mu}$
is called the \emph{\index{support}support of $\mu$} and is denoted by \index{$\supp(\mu)$}$\supp(\mu)$. 
\EndLemma

\BeginProof
Let ${\cal F}$ be the family of all closed sets $F$ with $\mu(F) =
\mu(X)$, then $\{X\setminus F\mid F\in{\cal F}\}$ is a directed family
of open sets, hence 
$
\mu(\bigcap{\cal F}) = \inf_{F\in{\cal F}}\mu(F) = \mu(X).
$
Define $\supp(\mu) := \bigcap{\cal F}$, then $\supp(\mu)$ is closed
with $\mu(\supp(\mu)) = \mu(X)$;
if $F\subseteq X$ is a closed set with $\mu(F) = \mu(X)$, then
$F\in{\cal F}$, hence $\supp(\mu)\subseteq F$. 
\EndProof

We can characterize the support of $\mu$ also in terms of open sets;
this is but a simple consequence of Lemma~\ref{support-for-regular}.

\BeginCorollary{member-of-support}
Under the assumptions of Lemma~\ref{support-for-regular} we have $x\in
\supp(\mu)$ iff $\mu(U)>0$ for each open neighborhood $U$ of $x$.
\QED
\EndCorollary

After all these preparations (with some interesting vistas to the
landscape of measures), we are in a position to show that the metric
topology on $\FinM{X, \Borel{X}}$ coincides with the A-topology for
$X$ separable metric. The following lemma will be the central
statement; it is formulated and proved separately, because its proof
is somewhat technical. Recall that the \emph{diameter} $\mathsf{diam}(Q)$\index{diameter} of $Q \subseteq X$
as
\begin{equation*}
\mathsf{diam}(Q) := \sup \{d(x_1, x_2) \mid x_1, x_2 \in Q\}.
\end{equation*}

\BeginLemma{ball-contains-neighborhood}
Every $d_{P}$-ball with center $\mu\in\FinM{X, \Borel{X}}$ contains a neighborhood of $\mu$ of the A-topology, if $(X, d)$ is separable metric.
\EndLemma

\BeginProof
Fix $\mu\in\FinM{X, \Borel{X}}$ and $\epsilon>0$, pick $\delta>0$ with
$4\cdot\delta<\epsilon$; it is no loss of generality to assume that
$\mu(X) = 1$. Because $X$ is separable metric, the support $S :=
\supp(\mu)$ is defined by Lemma~\ref{support-for-regular}. Because $S$
is closed, we can cover $S$ with a countable number $\Folge{V}$ of
open sets the diameter of which is less that $\delta$ and
$\mu(\partial V_{n}) = 0$ by
Proposition~\ref{without-boundary-has-basis}. Define 
\begin{align*}
  A_{1} & := V_{1}, \\
A_{n} & := \bigcup_{i=1}^{n}V_{i}\setminus\bigcup_{j=1}^{n-1}V_{j},
\end{align*}
then $\Folge{A}$ is a mutually disjoint family of sets which cover $S$, and for which $\mu(\partial A_{n}) = 0$ holds for all $n\in\Nat$. We can find an index $k$ such that $\mu(\bigcup_{i=1}^{k}) > 1-\delta$. Let $T_{1}, \dots, T_{\ell}$ be all sets which are a union of some of the sets $A_{1}, \dots, A_{k}$, then 
\begin{equation*}
  W := W_{T_{1}, \dots, T_{\ell}, \epsilon}(\mu)
\end{equation*}
is a neighborhood of $\mu$ in the A-topology by Corollary~\ref{a-top-basis-without-bound}. We claim that $d_{P}(\mu, \nu) < \epsilon$ for all $\nu\in W$. In fact, let $B\in\Borel{X}$ be arbitrary, and put 
\begin{equation*}
  A := \bigcup\{A_{i}\mid 1\leq i \leq k, A\cap B\not=\emptyset\},
\end{equation*}
then $A$ is among the $T$s just constructed, and $B\cap S \subseteq A\cup \bigcup_{i=k+1}^{\infty}A_{i}$.Moreover, we know that $A\subseteq B^{\delta}$, because each $A_{i}$ has a diameter less that $\delta$. This yields
\begin{equation*}
  \mu(B) = \mu(B\cap S) \leq \mu(A) + \delta < \nu(A) + 2\cdot\delta \leq \nu(B^{\delta}) + 2\cdot\delta.
\end{equation*}
On the other hand, we have%\footnote{Noch einmal nachrechnen.}
\begin{align*}
  \nu(B) & = \nu(B\cap S) + \nu(B\cap (X\setminus S))\\
& \leq \nu(A\cap \bigcup_{i=k+1}^{\infty}A_{i}) + 3\cdot\delta\\
& \leq \nu(A) + 3\cdot\delta\\
& \leq \mu(A) + 3\cdot\delta\\
& \leq \mu(B^{\delta}) + 4\cdot\delta.
\end{align*}
Hence $d_{P}(\mu, \nu) < 4\cdot\delta < \epsilon$. Thus $W$ is contained in the open ball around $\mu$ with radius smaller $\epsilon$. 
\EndProof

We have established

\BeginTheorem{a-top-is-metrizable}
The A-topology on $\FinM{X, \Borel{X}}$ is metrizable by the Lévy-Skohorod metric $d_{P}$, provided $(X, d)$ is a separable metric space. \QED 
\EndTheorem

We will see later that $d_{P}$ is not the only metric for this
topology, and that the corresponding metric space has interesting and
helpful properties. Some of these properties are best derived through an integral representation, for which a careful study of real-valued functions is required. This is what we are going to investigate in Section~\ref{sec:real-valu-funct}. But before doing this, we have a brief and tentative look at the relation between the Borel sets for A-topology and weak-*-$\sigma$-algebra.

\BeginLemma{weak-star-vs-a-borel}
Let $X$ be a metric space, then the weak-*-$\sigma$-algebra is contained in the Borel sets of the A-topology. If the A-topology has a countable basis, both $\sigma$-algebras are equal.
\EndLemma

\BeginProof
Denote by ${\cal C}$ the Borel sets of the A-topology on $\FinM{X, \Borel{X}}$. 

Since $X$ is metric, the Baire sets and the Borel sets coincide. For each closed set $F$, the evaluation map $ev_{F}: \mu\mapsto\mu(F)$ is upper semi-continuous by Proposition~\ref{convergence-a-topology}, so that the set 
\begin{equation*}
  {\cal G} := \{A\in\Borel{X} \mid ev_{A}\text{ is } {\cal C}-\text{measurable}\}
\end{equation*}
contains all closed sets. Because ${\cal G}$ is closed under complementation and countable disjoint unions, we conclude that ${\cal G}$ contains $\Borel{X}$. Hence $\schwach{\Borel{X}}\subseteq {\cal C}$ by minimality of $\schwach{\Borel{X}}$. 

2.
Assume that the A-topology has a countable basis, then each open set can be represented as a countable union of sets of the form $  W_{G_{1}, \dots, G_{n}, \epsilon}(\mu_{0}) $ with $G_1, \dots, G_n$ open. But $W_{G_{1}, \dots, G_{n}, \epsilon}(\mu_{0})\in\schwach{X, \Borel{X}}$, so that each open set is a member of $\schwach{X, \Borel{X}}$.  This implies the other inclusion.
\EndProof

We will investigate the A-topology further in Section~\ref{sec:weak-top} and turn to real-valued functions now. 

\Subsection{Real-Valued Functions} 
\label{sec:real-valu-funct}

We discuss the set of all measurable and bounded functions into the
real line now. We show first that the set of all these functions is
closed under the usual algebraic operations, so that it is a vector
space, and that it is also closed under finite infima and suprema,
rendering it a distributive lattice; in fact, algebraic operations and
order are compatible. Then we show that the measurable step functions
are dense with respect to pointwise convergence. This is an important
observation, which will help us later on to transfer linear properties
from indicator functions (a.k.a. measurable sets) to general
measurable functions. This prepares the stage for discussing
convergence of functions in the presence of a measure; we will deal
with convergence almost everywhere, which neglects a set of measure
zero for the purposes of convergence, and convergence in measure,
which is defined in terms of a pseudo metric, but surprisingly turns
out to be related to convergence almost everywhere through
subsequences of subsequences (this sounds a bit mysterious, so carry
on).

\BeginLemma{prob-lincomb-fnct}
Let $f, g: X\to \Real$ be ${\cal A}$-$\Borel{\Real}$-measurable
functions for the measurable space $(X, {\cal A})$. Then $f\wedge g$,
$f\vee g$ and $\alpha\cdot f + \beta\cdot g$ are ${\cal
  A}$-$\Borel{\Real}$-measurable for $\alpha, \beta\in\Real$. 
\EndLemma

\BeginProof
If $f$ is measurable, $\alpha\cdot f$ is. This follows immediately from Lemma~\ref{ElemGen}. From
\begin{equation*}
  \{x\in X\mid  f(x)+g(x) < q\} = 
\bigcup_{r_{1}, r_{2}\in\Rational, r_{1}+r_{2}\leq q}\bigl(\{x\mid f(x)< r_{1}\}\cap\{x\mid g(x)<r_{2}\}\bigr)
\end{equation*}
we obtain that the sum of measurable functions is measurable again. Since 
\begin{align*}
  \{x\in X \mid (f\wedge g)(x) < q\} & = \{x \mid f(x) < q\}\cup\{x \mid
  g(x) < q\}\\
 \{x\in X \mid (f\vee g)(x) < q\} & = \{x \mid f(x) < q\}\cap\{x \mid
  g(x) < q\},
\end{align*}
we see that both $f\wedge g$ and $f\vee g$ are measurable. 
\EndProof

\BeginCorollary{cor-prob-lincomb-fnct}
If $f: X\to \Real$ is ${\cal A}$-$\Borel{\Real}$-measurable, so is $|f|$. 
\EndCorollary

\BeginProof
Write $|f| = f^{+}-f^{-}$ with $f^{+} := f\vee 0$ and $f^{-} := (-f)\vee 0$. 
\EndProof

The consequence is that for a measurable space $(X, {\cal A})$ the set
\begin{equation*}
\MeasbFnct{X, \mathcal{A}} := \{f: N \rightarrow \Real \mid f \text{
is } \mathcal{A}-\Borel{\Real}\text{ measurable and bounded}\}
\end{equation*}
is both a vector space and a distributive lattice; in fact, it is a
vector lattice, see Definition~\ref{vector-lattice} on
page~\pageref{vector-lattice}. Assume that $\Folge{f}\subseteq
\MeasbFnct{X, \mathcal{A}} $ is a sequence of bounded measurable
functions such that $f: x \mapsto \liminf_{n\to \infty} f_{n}(x)$ is a
bounded function, then $f\in\MeasbFnct{X, \mathcal{A}}$. This is so
because
\begin{align*}
  \{x\in X \mid \liminf_{n\to \infty}f_{n}(x)\leq q\}
& = \{x \mid  \sup_{n\in \Nat}\inf_{k\geq n}f_{k}(x) \leq q\}\\
& = \bigcap_{n\in\Nat}\{x \mid \inf_{k\geq n}f_{k}(x) \leq q\}\\
& = \bigcap_{n\in\Nat}\bigcap_{\ell\in\Nat}\{x \mid \inf_{k\geq n}f_{k}(x) < q + 1/\ell\}\\
& = \bigcap_{n\in\Nat}\bigcap_{\ell\in \Nat}\bigcup_{k\geq n}\{x \mid f_{k}(x) < q + 1/\ell\}
\end{align*}
Similarly, if $x\mapsto \limsup_{n\to \infty} f_{n}(x)$ defines a bounded function, then it is measurable as well. Consequently, if the sequence $(f_{n}(x))_{n\in\Nat}$ converges to a bounded function $f$, then $f\in\MeasbFnct{X, \mathcal{A}}$. 

Hence we have shown
\BeginProposition{closed-under-limits}
Let $\Folge{f} \subseteq \MeasbFnct{X, \mathcal{A}}$ be a sequence of bounded measurable functions. Then
\begin{itemize}
  \item If $f_*(x) := \liminf_{n\to \infty} f_{n}(x)$ defines a bounded function, then $f_*\in \MeasbFnct{X, \mathcal{A}}$,
  \item if $f^*(x) := \limsup_{n\to \infty} f_{n}(x)$ defines a bounded function, then $f^*\in \MeasbFnct{X, \mathcal{A}}$.
\end{itemize}
\QED
\EndProposition

We use occasionally the representation of sets through
indicator functions. Recall for $A \subseteq X$ its
\emph{indicator function}\index{function!indicator}\index{indicator function}
\begin{equation*}
\chi_A(x) :=
\begin{cases}
  1, & \text{if } x \in A\\
  0, & \text{if } x \notin A.
\end{cases}
\end{equation*}
Clearly, if $\mathcal{A}$ is a $\sigma$-algebra on $X$, then $A \in
\mathcal{A}$ iff $\chi_A$ is a
$\mathcal{A}$-$\Borel{\Real}$-measurable function. This is so since
we have for the inverse image of an interval under $\chi_A$
\begin{equation*}
\InvBild{\chi_A}{[0, q]} =
\begin{cases}
  \emptyset, & \text{if } q < 0,\\
  X\setminus A, & \text{if } 0 \leq q < 1,\\
  X, & \text{if } q \geq 1.
\end{cases}
\end{equation*}

A measurable
\emph{step function}\index{function!step}
\begin{equation*}
f = \sum_{i = 1}^n \alpha_i \cdot \chi_{A_i}
\end{equation*}
is a linear combination of indicator functions with $A_i \in
\mathcal{N}$. Since $\chi_{A}\in\MeasbFnct{X, \mathcal{A}} $ for $A\in{\cal A}$, measurable step functions are indeed measurable functions. 

\BeginProposition{ApproxStepFncts}
Let
$
(X, \mathcal{A})
$
be a measurable space. Then
\begin{enumerate}
 \item For $f \in \MeasbFnct{X, \mathcal{A}}$ with $f \geq 0$ there
  exists an increasing sequence  $\Folge{f}$ of
  step functions $f_n \in \MeasbFnct{X, \mathcal{A}}$ with
  \begin{equation*}
  f(x) = \sup_{n \in \Nat} f_n(x)
\end{equation*}
  for all $x \in X.$
 \item For $f \in \MeasbFnct{X, \mathcal{A}}$  there
  exists a sequence  $\Folge{f}$ of
  step functions $f_n \in \MeasbFnct{N, \mathcal{A}}$ with
  \begin{equation*}
  f(x) = \lim_{n \rightarrow \infty} f_n(x)
\end{equation*}
  for all $x \in X.$
\end{enumerate}
\EndProposition

\BeginProof
1.
Take $f\geq 0$, and assume without loss of generality that $f\leq 1$ (otherwise, if $0\leq f \leq m$, consider $f/m$). Put 
\begin{equation*}
  A_{i, n} := \{x\in X \mid i/n \leq f(x) < (i+1)/n\},
\end{equation*}
for $n\in\Nat$, $0\leq i \leq n$, then $A_{i, n}\in {\cal A}$, since $f$ is measurable. Define 
\begin{equation*}
  f_{n}(x) := \sum_{0\leq i < 2^{n}}i\cdot 2^{-n}\chi_{A_{i, 2^{n}}}.
\end{equation*}
Then $f_{n}$ is a measurable step function, and $f_{n}\leq f$, moreover $\Folge{f}$ is increasing. This is so because given $n\in \Nat, x\in X$, we can find $i$ such that $x\in A_{i, 2^{n}} = A_{2i, 2^{n+1}}\cup A_{2i+1, 2^{n+1}}$. If $f(x) < (2i+1)/2^{n+1}$, we have $x\in A_{2i, 2^{n+1}}$ with $f_{n}(x) = f_{n+1}(x)$, if, however, $(2i+1)/2^{n+1}\leq f(x)$, we have $f_{n}(x) < f_{n+1}(x)$. 

Given $\epsilon > 0$, choose $n_{0}\in \Nat$ with $2^{-n}< \epsilon$ for $n\geq n_{0}$. Let $x\in X, n\geq n_{0}$, then $x\in A_{i, 2^{n}}$ for some $i$, hence $|f_{n}(x) - f(x)| = f(x) - i2^{-n} < 2^{-n} < \epsilon$.  Thus $f = \sup_{n\in\Nat} f_{n}$. 

2.
Given $f\in\MeasbFnct{X, \mathcal{A}}$, write $f_{1} := f\wedge 0$ and $f_{2} := f\vee 0$, then $f = f_{1}+ f_{2}$ with $f_{1}\leq 0$ and $f_{2}\geq 0$ as measurable and bounded functions. Hence $f_{2} = \sup_{n\in\Nat} g_{n} = \lim_{n\to \infty}g_{n}$ and $-f_{1} = -\sup_{n\in\Nat}h_{n} = -\lim_{n\to \infty}h_{n}$ for increasing sequences of step functions $\Folge{g}$ and $\Folge{h}$. Thus $f = \lim_{n\to \infty}(g_{n}+h_{n})$, and $g_{n}+h_{n}$ is a step function for each $n\in\Nat$.  
\EndProof

Given $f: X\to \Real$ with $f\geq 0$, the set $\{\langle x, q\rangle
\in X\times\Real \mid 0\leq f(x) \leq q\}$ can be visualized as the
area between the $X$-axis and the graph of the function. We obtain as
a consequence that this set is measurable, provided $f$ is measurable.
This gives an example of a product measurable set. To be
specific

\BeginCorollary{pre-choquet}
Let $f: X\to \Real$ with $f\geq 0$ be a bounded measurable function for a measurable space $(X, {\cal A})$, and define
\begin{equation*}
  C_{\Join}(f) := \{\langle x, q\rangle \mid 0\leq q \Join f(x)\}\subseteq X\times\Real
\end{equation*}
for the relational operator $\Join$ taken from $\{\geq, <, =, \not=, >, \geq\}$. Then $C_{\Join}(f)\in{\cal A}\otimes\Borel{\Real}$. 
\EndCorollary

\BeginProof
We prove the assertion for $C(f) := C_{<}(f)$, from which the other cases may easily be derived, e.g., 
\begin{equation*}
  C_{\leq}(f) = \bigcap_{k\in \Nat}\{\langle x, q\rangle \mid  f(x) < q+1/k\} 
= \bigcap_{k\in\Nat}C_{<}(f-1/k).
\end{equation*}
Consider these cases.
\begin{enumerate}
  \item If $f = \chi_A$ with $A \in \mathcal{A}$, then
$
C(f) = X\setminus A \times \{0\} \cup A \times [0, 1[ \in \mathcal{A}\otimes \Borel{\Real}.
$
  \item If $f$ is represented as a step function with a finite
  number of mutually disjoint steps, say, $f = \sum_{i=1}^k r_i\cdot \chi_{A_i}$
  with $r_i \geq 0$ and all $A_i \in \mathcal{A}$, then
  \begin{equation*}
  C(f) = \left(X\setminus\bigcup_{i=1}^k A_i\right) \times \{0\} \cup \bigcup_{i=1}^k A_i \times [0, r_i[
  \in \mathcal{A}\otimes \Borel{\Real}.
\end{equation*}
  \item If $f$ is represented as a monotone limit of step function
  $\Folge{f}$ with $f_n \geq 0$ according to Proposition~\ref{ApproxStepFncts}, then
  $
  C(f) = \bigcup_{n \in \Nat} C(f_n),
  $
  thus $C(f) \in \mathcal{A}\otimes \Borel{\Real}$.
\end{enumerate}
\EndProof

\BeginExample{modal-logic-rev}
Consider the simple modal logic in Example~\ref{modal-logic},
interpreted through a transition kernel $M: (X, {\cal A})\Trans (X,
{\cal A})$. Given a formula $\phi$, the set $ \{\langle x, r\rangle
\mid M(x)(\Gilt_{M}) \geq r\} $ is a member of ${\cal
  A}\otimes\Borel{\Real}$. Note that $\Gilt[\Diamond_{q}\phi]_{M}$ is
the cut of this set at $q$. Hence this observation generalizes
measurability of $\Gilt[\cdot]_{M}$, one of the cornerstones for
interpreting modal logics probabilistically.
\EndExample

We will turn now to the interplay of measurable functions and measures
and have a look at different modes of convergence for sequences of
measurable functions in the presence of a (finite) measure.

\Subsubsection{Essentially Bounded Functions}
\label{sec:ess-bounded-fncts}
Fix for this section a finite measure space $(X, {\cal A}, \mu)$. We say that a measurable property holds \emph{$\mu$-almost everywhere} (abbreviated as \emph{\index{$\mu$-a.e.}$\mu$-a.e.}) iff the set on which the property does not hold has $\mu$-measure zero. 

The measurable function $f\in\MeasbFnct{X, \mathcal{A}}$ is called \emph{\index{$\mu$-essentially bounded}$\mu$-essentially bounded} iff 
\begin{equation*}
  \infNorm{f}{\mu} := \inf\{a\in\Real \mid |f|\leq_{\mu} a\} < \infty,
\end{equation*}
where $f\leq_{\mu}a$ indicates that $f\leq a$ holds $\mu$-a.e. Thus a $\mu$-essentially bounded function may occasionally take arbitrary large values, but the set of these values must be negligible in terms of $\mu$.

The set 
\begin{equation*}
\index{${\cal L}_{\infty}(\mu)$}{\cal L}_{\infty}(\mu) := {\cal L}_{\infty}(X, {\cal A},\mu) := \{f\in\MeasbFnct{X, \mathcal{A}} \mid \infNorm{f}{\mu} < \infty\}
\end{equation*}
of all $\mu$-essentially bounded functions is a real vector space, and we have for $\infNorm{\cdot}{\mu}$ these properties.

\BeginLemma{prelim-prop-inf-norm}
Let $f, g\in \MeasbFnct{X, \mathcal{A}}$ essentially bounded, $\alpha, \beta\in\Real$, then $\infNorm{\cdot }{\mu}$ is a \emph{\index{pseudo-norm}pseudo-norm} on $\MeasbFnct{X, \mathcal{A}}$, i.e., 
\begin{enumerate}
\item If $\infNorm{f}{\mu} = 0$, then $f=_{\mu} 0$. 
\item $\infNorm{\alpha\cdot f}{\mu} = |\alpha|\cdot\infNorm{f}{\mu}$,
\item $\infNorm{ f + g}{\mu} \leq \infNorm{f}{\mu} + \infNorm{g}{\mu}$.
\end{enumerate}
\EndLemma

\BeginProof
If $\infNorm{f}{\mu} = 0$, we have $|f|\leq_{\mu}1/n$ for all
$n\in\Nat$, so that 
\begin{equation*}
\{x\in X \mid |f(x)| \not= 0\}\subseteq
\bigcup_{n\in\Nat} \{x\in X\mid |f(x)| \leq 1/n\},
\end{equation*}
consequently,
$f=_\mu 0$. The converse is trivial.  The second property follows from
$|f|\leq_{\mu} a$ iff $|\alpha\cdot f| \leq_{\mu} |\alpha|\cdot a$, the
third one from the observation that $|f|\leq_{\mu} a$ and
$|g|\leq_{\mu}b$ implies $|f + g| \leq |f| + |g| \leq_{\mu} a + b$.
\EndProof

So $\infNorm{\cdot}{\mu}$ \emph{nearly} a norm, but the
crucial property that the norm for a vector is zero only if the vector is
zero is missing. We factor ${\cal L}_{\infty}(X, {\cal A},\mu)$ with
respect to the equivalence relation $=_{\mu}$, then the set
\begin{equation*}
  L_{\infty}(\mu) := L_{\infty}(X, {\cal A},\mu) := \{[f] \mid f\in {\cal L}_{\infty}(X, {\cal A},\mu)\}
\end{equation*}
of all equivalence classes $[f]$ of $\mu$-essentially bounded measurable functions is a vector space again. This is so because $\isEquiv{f}{g}{=_{\mu}}$ and  $\isEquiv{f'}{g'}{=_{\mu}}$together imply $\isEquiv{f+f'}{g+g'}{=_{\mu}}$, and $\isEquiv{f}{g}{=_{\mu}}$ implies $\isEquiv{\alpha\cdot f}{\alpha\cdot g}{=_{\mu}}$ for all $\alpha\in\Real$. Moreover,
\begin{equation*}
\infNorm{[f]}{\mu} := \infNorm{f}{\mu}
\end{equation*}
defines a norm on this space. For easier reading we will identify in the sequel $f$ with its class $[f]$.

We obtain in this way a normed vector space, which is complete with respect to this norm.

\BeginProposition{linf-is-banach-space}
$(L_{\infty}(\mu), \infNorm{\cdot}{\mu})$ is a Banach space.
\EndProposition

\BeginProof
Let $\Folge{f}$ be a Cauchy sequence in $ L_{\infty}(X, {\cal A},\mu)$, and define 
\begin{equation*}
  N := \bigcup_{n_{1}, n_{2}\in\Nat}\{x\in X\mid |f_{n_{1}}(x) - f_{n_{2}}(x)| > \infNorm{f_{n_{1}} - f_{n_{2}}}\mu\},
\end{equation*}
then $\mu(N) = 0$. Put $g_{n} := \chi_{X\setminus N}\cdot f_{n}$, then $\Folge{g}$ converges uniformly with respect to the supremum norm $\infNorm{\cdot}{}$ to some element $g\in \MeasbFnct{X, \mathcal{A}}$, hence also $\infNorm{f_{n}- g}{\mu}\to 0$. Clearly, $g$ is bounded. 
\EndProof

This is the first instance of a vector space intimately connected with a measure space. We will discuss several of these spaces later on, when integration is at our disposal. 

The convergence of a sequence of measurable functions into $\Real$ in the presence of a finite measure is discussed now. Without a measure, we may use pointwise or uniform convergence for modelling approximations. Recall that \emph{\index{convergence!pointwise}pointwise convergence} of a sequence $\Folge{f}$ of functions to a function $f$ is given by 
\begin{equation}
\label{conv-everywhere}
  \forall x\in X: \lim_{n\to \infty}f_{n}(x) = f(x),
\end{equation}
and the stronger form of \emph{\index{convergence!uniform}uniform convergence} through 
\begin{equation*}
  \lim_{n\to \infty}\infNorm{f_{n}-f}{} = 0,
\end{equation*}
with $\infNorm{\cdot}{}$ as the supremum norm, given by
\begin{equation*}
  \infNorm{f}{} := \sup_{x\in X}|f(x)|. 
\end{equation*}
We will weaken the first condition~(\ref{conv-everywhere}), so that it holds not everywhere
but almost everywhere, thus the set on which it does not hold will be
a set of measure zero. This leads to the notion of convergence almost
everywhere, which will turn out to be quite close to uniform
convergence, as we will see when discussing Egorov's Theorem. Convergence almost everywhere will be weakened to convergence in measure, for which we will define a pseudo metric. This in turn gives rise to another Banach space upon factoring. 

\Subsubsection{Convergence \emph{almost everywhere} and \emph{in measure}}
\label{sec:conv-a-e} 

Recall that we work in a finite measure space $(X, {\cal A}, \mu)$. The sequence
$\Folge{f}$ of measurable functions $f_{n}\in\MeasbFnct{X,
  \mathcal{A}}$ is said to \emph{converge \index{convergence!almost
    everywhere}almost everywhere} to a function $f\in\MeasbFnct{X,
  \mathcal{A}}$ (written as $f_{n}\index{$\aeC$}\aeC f$) iff the sequence
$(f_{n}(x))_{n\in \Nat}$ converges pointwise to $f(x)$ for every $x$
outside a set of measure zero. Thus we have $\mu(X\setminus K) = 0$,
where $K := \{x\in X\mid f_{n}(x) \to f(x)\}$. Because
\begin{equation*}
  K = \bigcap_{n\in\Nat}\bigcup_{m\in\Nat}\bigcap_{\ell\geq m}\{x\in X\mid |f_{\ell}(x)- f(x)| < 1/n\},
\end{equation*}
$K$ is a measurable set. It is clear that $f_{n}\aeC f$ and $f_{n}\aeC f'$ implies that $f=_{\mu} f'$ holds. 

The next lemma shows that convergence everywhere is compatible with
the common algebraic operations on $\MeasbFnct{X, \mathcal{A}}$ like
addition, scalar multiplication and the lattice operations. Since
these functions can be represented as continuous function of several
variables, we formulate this closure property abstractly in terms of
compositions with continuous functions.

\BeginLemma{ae-compatible-with-algebraic}
Let $f_{i, n}\aeC f_{i}$ for $1 \leq i \leq k$, and assume that $g: \Real^{k}\to \Real$ is continuous. Then $g\circ (f_{1, n}, \dots, f_{k, n})\aeC g\circ (f_{1}, \dots, f_{k})$. 
\EndLemma

\BeginProof
Put $h_{n} := g\circ (f_{1, n}, \dots, f_{k, n})$. 
Since $g$ is continuous, we have 
\begin{equation*}
  \{x\in X\mid  \bigl(h_{n}(x)\bigr)_{n\in\Nat} \text{ does not converge}\}
\subseteq\bigcup_{j=1}^{k}{}\{x\in X\mid  \bigl(f_{j, n}(x)\bigr)_{n\in\Nat} \text{ does not converge}\},
\end{equation*}
hence the set on the left hand side has measure zero.
\EndProof

Intuitively, convergence almost everywhere means that the measure of the set
\begin{equation*}
  \bigcup_{n\geq k}\{x\in X \mid |f_{n}(x) - f(x)| > \epsilon\}
\end{equation*}
tends to zero, as $k\to \infty$, so we are coming closer and closer to
the limit function, albeit on a set the measure of which becomes
smaller and smaller. We show that this intuitive understanding yields
an adequate model for this kind of convergence.

\BeginLemma{con-ae-as-limit}
Let $\Folge{f}$ be a sequence of functions in $\MeasbFnct{X, \mathcal{A}}$ and $f\in\MeasbFnct{X, \mathcal{A}}$. Then the following conditions are equivalent
\begin{enumerate}
\item\label{con-ae-as-limit:1} $f_{n}\aeC f$.
\item\label{con-ae-as-limit:2} $\lim_{k\to \infty}\mu\bigl(\bigcup_{n\geq k}\{x\in X \mid |f_{n}(x) - f(x)| > \epsilon\}\bigr) = 0$ for every $\epsilon>0$.
\end{enumerate}
\EndLemma

\BeginProof
Let $\epsilon>0$ be given, then there exists $k\in\Nat$ with $1/k<\epsilon$, so that 
\begin{align*}
\lim_{k\to \infty}\mu\bigl(\bigcup_{n\geq k} \{x\in X\mid |f_{n}(x) - f(x)| > \epsilon\}\bigr)
& \stackrel{(*)}{=} 
 \mu\bigl(\bigcap_{k\in\Nat}\bigcup_{n\geq k} \{x\in X\mid |f_{n}(x) - f(x)| > \epsilon\}\bigr)\\
& \leq \mu(\{x\in X\mid (f_{n}(x))_{n\in\Nat}\text{ does not converge}\}).
\end{align*}
Now assume that $f_{n}\aeC f$, then the implication \labelImpl{con-ae-as-limit:1}{con-ae-as-limit:2} is immediate. If, however, $f_{n}\aeC f$ is false, then we find for each $\epsilon>0$ so that for all $k\in\Nat$ there exists $n\geq k$ with 
$
\mu(\{x\in X\mid |f_{n}(x)-f(x)|\geq \epsilon\}) > 0$. Thus \ref{con-ae-as-limit:2} cannot hold.
\EndProof

Note that the statement above requires a finite measure space, because
the measure of a decreasing sequence of sets is the infimum of the
individual measures, used in the equation marked $(*)$. This is not
necessarily valid for non-finite measure space.

The characterization implies that a.e.-Cauchy sequences converge.

\BeginCorollary{ae-cauchy-converges}
Let $\Folge{f}$ be an a.e.-Cauchy sequence in $\MeasbFnct{X, \mathcal{A}}$. Then $\Folge{f}$ converges.
\EndCorollary

\BeginProof
Because $\Folge{f}$ is an a.e.-Cauchy sequence, we have that $\mu(X\setminus K_{\epsilon}) = 0$ for every $\epsilon>0$, where
\begin{equation*}
  K_{\epsilon} := \bigcap_{k\in\Nat}\bigcup_{n, m\geq k}\{x\in X\mid |f_{n}(x)-f_{m}(x)|>\epsilon\}.
\end{equation*}
Put 
\begin{align*}
  N & := \bigcup_{k\in\Nat} K_{1/k},\\
g_{n} & := f_{n}\cdot\chi_{X\setminus N},
\end{align*}
then $\Folge{g}$ is a Cauchy sequence in $\MeasbFnct{X, \mathcal{A}}$ which converges pointwise to some $f\in\MeasbFnct{X, \mathcal{A}}$. Since $\mu(X\setminus N) = 0$, $f_{n}\aeC f$ follows. 
\EndProof

Convergence a.e. is very nearly uniform convergence, where \emph{very
  nearly} serves to indicate that
the set on which uniform convergence does not happen is arbitrarily
small. To be specific, we can find for each threshold a set the complement
of which has a measure smaller than this bound, on which convergence
is uniform. This is what
\index{theorem!Egorov}\emph{Egorov's Theorem} says.

\BeginProposition{egorov}
Let $f_{n}\aeC f$ for $f_{n}, f\in\MeasbFnct{X, \mathcal{A}}$. Given $\epsilon>0$, there exists $A\in{\cal A}$ such that 
\begin{enumerate}
\item $\sup_{x\in A}|f_{n}(x)-f(x)|\to 0$,
\item $\mu(X\setminus A) < \epsilon$.
\end{enumerate}
\EndProposition

The idea of the proof is that we look at each $x$ for which uniform convergence is spoiled by $1/k$. This set can be made arbitrary small in terms of $\mu$, so the union of all these sets can be made as small as we want. Outside this set we have uniform convergence. Let's look at a more formal treatment now. 

\BeginProof
Fix $\epsilon>0$, then there exists for each $k\in\Nat$ an index $n_{k}\in\Nat$ such that $\mu(B_{k}) < \epsilon/2^{k+1}$ with 
\begin{equation*}
  B_{k} := \bigcup_{m\geq n_{k}}\{x\in X\mid |f_{m}(x)-f(x)|>1/k\}.
\end{equation*}
Now put $A := \bigcap_{k\in\Nat}(X\setminus B_{k})$, then 
\begin{equation*}
  \mu(X\setminus A) \leq \sum_{k\in\Nat}\mu(B_{k}) \leq \epsilon,
\end{equation*}
and we have for all $k\in\Nat$
\begin{equation*}
  \sup_{x\in A}|f_{n}(x) - f(x)| \leq \sup_{x\not\in B_{k}}|f_{n}(x) -f(x)| \leq 1/k
\end{equation*}
for $n\geq n_{k}$. Thus 
\begin{equation*}
  \lim_{n\to \infty}\sup_{x\in A} |f_{n}(x) - f(x)| = 0,
\end{equation*}
as claimed. 
\EndProof

Convergence almost everywhere makes sure that the set on which a
sequence of functions does not converge has measure zero, and Egorov's
Theorem shows that this is \emph{almost} uniform convergence.

Convergence in measure for a finite measure space $(X,
{\cal A}, \mu)$ takes another approach: fix $\epsilon>0$, and consider
the set $\{x\in X \mid |f_{n}(x) - f(x)| > \epsilon\}$. If the measure
of this set (for a fixed, but arbitrary $\epsilon$) tends to zero, as $n\to \infty$, 
then we say that $\Folge{f}$ \emph{\index{convergence!in
    measure}converges in measure} to $f$, and write $f_{n}\index{$\nmC$}\nmC f$. In
order to have a closer look at this notion of convergence, we note
that it is invariant against equality almost everywhere: if
$f_{n}=_{\mu} g_{n}$ and $f=_{\mu}g$, then $f_{n}\nmC f$ implies
$g_{n}\nmC g$, and vice versa. 

We will introduce a pseudo metric $\delta$ on $\MeasbFnct{X, \mathcal{A}}$ first:
\begin{equation*}
  \delta(f, g) := \inf\bigl\{\epsilon >0 \mid \mu(\{x\in X \mid |f(x)-g(x)|>\epsilon\}\leq\epsilon\bigr\}.
\end{equation*}

These are some elementary properties of $\delta$:

\BeginLemma{properties-pseudo-metric}
Let $f, g, h\in\MeasbFnct{X, \mathcal{A}}$, then we have
\begin{enumerate}
\item\label{item-pseudo-metric:1} $\delta(f, g) = 0$ iff $f=_\mu g$,
\item\label{item-pseudo-metric:2} $\delta(f, g) = \delta(g, f)$,
\item\label{item-pseudo-metric:3} $\delta(f, g) \leq \delta(f, h) + \delta(h, g)$.
\end{enumerate}
\EndLemma

\BeginProof
If $\delta(f, g) = 0$, but $f\not=_{\mu}g$, there exists $k$ with
$\mu(\{x\in X\mid |f(x)-g(x)|>1/k\})>1/k$. This is a
contradiction. The other direction is trivial. Symmetry of $\delta$ is
also trivial, so the triangle inequality remains to be shown. If
$|f(x)-g(x)| > \epsilon_{1} + \epsilon_{2}$, then $|f(x)-h(x)| >
\epsilon_{1}$ or $|h(x) - g(x)|> \epsilon_{2}$, thus 
\begin{equation*}
\mu(\{x\in X
\mid |f(x)-g(x)| > \epsilon_{1} + \epsilon_{2}\}) \leq \mu(\{x\in
X\mid |f(x)-h(x)| > \epsilon_{1}\}) + \mu(\{x\in X\mid |h(x)-g(x)| >
\epsilon_{2}\}).
\end{equation*}
This implies the third property.
\EndProof

This, then, is the formal definition of convergence in measure:

\BeginDefinition{def-conv-in-measure}
 The sequence $\Folge{f}$ in $\MeasbFnct{X, \mathcal{A}}$ is said to  \emph{converge in measure }to $f\in\MeasbFnct{X, \mathcal{A}}$  (written as $f_{n}\nmC f$) iff $\delta(f_{n}, f) \to 0$, as $n\to \infty$. 
\EndDefinition

We can express convergence in measure in terms of convergence almost everywhere.

\BeginProposition{conv-in-measure-subseq}
$\Folge{f}$ converges in measure to $f$ iff each subsequence of $\Folge{f}$ contains a subsequence $\Folge{h}$ with $h_{n}\aeC f$. 
\EndProposition

 \BeginProof
 ``$\Rightarrow$'': Assume $f_{n}\nmC f$, and let $\epsilon>0$ be arbitrary but fixed. Let $\Folge{g}$ be a subsequence of $\Folge{f}$. We find a sequence of indices $n_{1} < n_{2} < \dots$ such that $\mu(\{x\in X \mid  |g_{n_{k}}(x) - f(x)|>\epsilon\}) < 1/k^{2}$. Let $h_{k} := g_{n_{k}}$, then we obtain 
 \begin{equation*}
   \mu(\bigcup_{k\geq \ell}\{x\in X \mid |h_{k} - f| > \epsilon\}) \leq \sum_{k\geq\ell}\frac{1}{k^{2}} \to 0,
 \end{equation*}
 as $\ell\to \infty$. Hence $h_{k}\aeC f$. 

 ``$\Leftarrow$'': If $\delta(f_{n}, f)\not\to 0$, we can find a subsequence $(f_{n_{k}})_{k\in\Nat}$ and $r>0$ such that $\mu(\{x\in X\mid |f_{n_{k}}(x) - f(x)|>r\})>r$ for all $k\in\Nat$. Let $\Folge{g}$ be a subsequence of this subsequence, then 
 \begin{equation*}
   \lim_{n\to \infty}\mu(\{x\in X\mid |g_{n}-f|>r\})
 \leq \lim_{n\to \infty}\mu(\bigcup_{m\geq n}\{x\in X\mid |g_{m}-f|>r\}) = 0
 \end{equation*}
 by Lemma~\ref{con-ae-as-limit}. This is a contradiction.
 \EndProof

Hence convergence almost everywhere implies convergence in measure. Just for the record:

\BeginCorollary{conv-ae-impl-conv-im}
If $\Folge{f}$ converges almost everywhere to $f$, then the sequence converges also in measure to $f$.
\QED
\EndCorollary

The converse relationship is a bit more involved. Intuitively, a sequence which converges in measure need not converge almost everywhere.

\BeginExample{im-but-not-ae}
Let $A_{i, n} := [(i-1)/n, i/n]$ for $n\in\Nat$ and $1 \leq i \leq n$, and consider the sequence 
$\Folge{f} := \langle \chi_{A_{1, 1}}, \chi_{A_{1, 2}}, \chi_{A_{2, 2}}, \chi_{A_{3, 1}}, \chi_{A_{3,2}}, \chi_{A_{3, 3}}, \dots\rangle,$ so that in general $\chi_{A_{1, n}}, \dots, \chi_{A_{n, n}}$ is followed by $\chi_{A_{1, n+1}}, \dots, \chi_{A_{n+1, n+1}}$. Let $\mu$ be Lebesgue measure $\lambda$ on $\Borel{[0, 1]}$. Given $\epsilon>0$, $\lambda(\{x\in[0, 1]\mid f_{n}(x)>\epsilon\})$ can be made arbitrarily small for any given $\epsilon>0$, hence $f_{n}\nmC 0$. On the other hand, $\bigl(f_{n}(x)\bigr)_{n\in\Nat}$ fails to converge or any $x\in[0, 1]$, so $f_{n}\aeC 0$ is false. 
\EndExample
 
We have, however, this observation, which draws atom into our game.

\BeginProposition{conv-im-impl-conv-ae-atom}
Let $(A_{i})_{i\in I}$ be the at most countable collection of $\mu$-atoms according to Lemma~\ref{only-countably-many-atoms} such that $B := X\setminus\bigcup_{i\in I}A_{i}$ does not contain any atoms. Then these conditions are equivalent:
\begin{enumerate}
\item\label{aes:1} Convergence in measure implies convergence almost everywhere. 
\item\label{aes:2} $\mu(B) = 0$.
\end{enumerate}
\EndProposition

\BeginProof
\labelImpl{aes:1}{aes:2}: Assume that $\mu(B)\leq 0$, then we know
that for each $k\in \Nat$ there exist mutually disjoint measurable
subsets $B_{1, k}, \dots, B_{k, k}$ of $B$ such that $\mu(B_{i, k}) =
1/k\cdot \mu(B)$ and $B = \bigcup_{1\leq i \leq k}B_{i, k}$. This is
so because $B$ does not contain any atoms. Put as above $\Folge{f} :=
\langle \chi_{B_{1, 1}}, \chi_{B_{1, 2}}, \chi_{B_{2, 2}}, \chi_{B_{3,
    1}}, \chi_{B_{3,2}}, \chi_{B_{3, 3}}, \dots\rangle,$ so that in
general $\chi_{B_{1, n}}, \dots, \chi_{B_{n, n}}$ is followed by
$\chi_{B_{1, n+1}}, \dots, \chi_{B_{n+1, n+1}}$. Because $\mu(\{x\in X
\mid f_{n}(x) > \epsilon\}$ can be made arbitrarily small for any
positive $\epsilon$, we find $f_{n}\nmC 0$. If we assume that
convergence in measure implies convergence almost everywhere, we have
$f_{n}\aeC 0$, but this is false, because $\liminf_{n\to \infty}f_{n}
= 0$ and $\limsup_{n\to \infty}f_{n} = \chi_{B}$. This is a
contradiction.

\labelImpl{aes:2}{aes:1}:
Let $\Folge{f}$ be a sequence with $f_{n}\nmC f$. Fix an atom $A_{i}$, then 
$
\mu(\{x\in A_{i}\mid |f_{n}(x) - f(x)|> 1/k\}) = 0
$
for all $n\geq n_{k}$ with $n_{k}$ suitably chosen; this is so because $A_{i}$ is an atom, hence measurable subsets of $A_{i}$ take only the values $0$ and $\mu(A_{i})$. Put 
\begin{equation*}
g := \inf_{n\in\Nat}\sup_{n_{1}, n_{2}\geq n} |f_{n_{1}}-f_{n_{2}}|,
\end{equation*}
then $g(x) \not= 0$ iff $(f_{n}(x))_{n\in\Nat}$ does not converge to $f(x)$. We infer
$\mu(\{x\in A_{i}\mid g(x)\geq 2/k\}) = 0$. Because the family $(A_{i})_{i\in I}$ is mutually disjoint, we conclude that $\mu(\{x\in X \mid  g(x)\geq 2/k\}) = 0$ for all $k\in \Nat$. But now look at this
\begin{equation*}
\mu(\{x\in X \mid  \liminf_{n\to \infty}f_{n}(x) < \limsup_{n\to \infty}f_{n}(x)\} = \mu(\{x\in X \mid g(x) > 0\}) = 0.
\end{equation*}
Consequently, $f_{n}\aeC f$. 
\EndProof

Again we want to be sure that convergence in measure is preserved by the usual algebraic operations like addition or taking the infimum, so we state as a counterpart to Lemma~\ref{ae-compatible-with-algebraic} now as an easy consequence.

 \BeginLemma{nm-compatible-with-algebraic}
 Let $f_{i, n}\nmC f_{i}$ for $1 \leq i \leq k$, and assume that $g: \Real^{k}\to \Real$ is continuous. Then $g\circ (f_{1, n}, \dots, f_{k, n})\nmC g\circ (f_{1}, \dots, f_{k})$. 
 \EndLemma

 \BeginProof
 By iteratively selecting subsequences, we can find subsequences $(h_{i, n})_{n\in\Nat}$ such that $h_{i, n}\aeC f_{i}$, as $n\to \infty$ for $1 \leq i \leq k$. Then apply Lemma~\ref{ae-compatible-with-algebraic} and Proposition~\ref{conv-in-measure-subseq}.
 \EndProof

 Let $F(X, \mathcal{A})$ be the factor space $\Faktor{\MeasbFnct{X,
     \mathcal{A}}}{=_{\mu}}$ of the space of all measurable functions
 with respect to $=_{\mu}$. Then this is a real vector space again,
 because the algebraic operations on the equivalence classes are well
 defined. Note that we have $\delta(f, g) = \delta(f', g')$, provided
 $f=_\mu g$ and $f'=_{\mu}g'$. We identify again the class
 $\Klasse{f}{=_\mu}$ with $f$. Define
 \begin{equation*}
   \aNorm{f}{} := \delta(f, 0)
 \end{equation*}
 for $f\in F(X, \mathcal{A})$.

 \BeginProposition{f-is-banach-space}
 $(F(X, \mathcal{A}), \aNorm{\cdot}{})$ is a Banach space.
 \EndProposition

\BeginProof
 1.
 It follows from Lemma~\ref{properties-pseudo-metric} and the observation $\delta(f, 0) = 0$ iff $f =_{\mu} 0$ that $\aNorm{\cdot}{}$ is a norm, so we have to show that $F(X, \mathcal{A})$ is complete with this norm. 

 2.
 Let $\Folge{f}$ be a Cauchy sequence in $F(X, \mathcal{A})$, then we can find a strictly increasing sequence $\Folge{\ell}$ of integers such that $\delta(f_{\ell_{n}}, f_{\ell_{n+1}})\leq 1/n^{2}$, hence 
 \begin{equation*}
 \mu(\{x\in X \mid |f_{\ell_{n}}(x) - f_{\ell_{n+1}}(x)| > 1/n^{2}\})\leq 1/n^{2}.
 \end{equation*}
 Let $\epsilon>0$ be given, then there exists $r\in\Nat$ with $\sum_{n\geq r}1/n^{2}< \epsilon$, hence we have 
 \begin{equation*}
   \bigcap_{n\in\Nat}\bigcup_{m, k \geq n}\{x\in X \mid |f_{\ell_{m}}(x) - f_{\ell_{k}}(x)|>\epsilon\}
 \subseteq 
 \bigcup_{n\geq k}\{x\in X\mid |f_{\ell_{n}}(x) - f_{\ell_{n+1}}(x)|<1/n^{2}\},
 \end{equation*}
 if $k\geq r$. Thus
 \begin{equation*}
   \mu(\bigcap_{n\in\Nat}\bigcup_{m, k \geq n}\{x\in X \mid |f_{\ell_{m}}(x) - f_{\ell_{k}}(x)|>\epsilon\})
 \leq \sum_{n\geq k}1/n^{2}\to 0,
 \end{equation*}
 as $k\to \infty$. Hence $(f_{\ell_{n}})_{n\in\Nat}$ is an a.e. Cauchy sequence which
 converges a.e. to some $f\in F(X, \mathcal{A})$, which by
 Proposition~\ref{conv-in-measure-subseq} implies that $f_{n}\nmC f$.
\EndProof

A consequence of $(F(X, \mathcal{A}), \aNorm{\cdot}{})$ being a Banach space is that $\MeasbFnct{X, \mathcal{A}}$ is complete with respect to convergence in measure for any finite measure $\mu$ on ${\cal A}$. Thus for any sequence $\Folge{f}$ of functions such that given $\epsilon>0$ there exists $n_{0}$ such that $\mu(\{x\in X\mid |f_{n}(x) - f_{m}(x)|>\epsilon\})<\epsilon$ for all $n, m\geq n_{0}$ we can find $f\in \MeasbFnct{X, \mathcal{A}}$ such that $f_{n}\nmC f$ with respect to $\mu$.

We will deal with measurable real valued functions again and in greater detail in Section~\ref{lp-spaces}; then we will have integration as a powerful tool at our disposal, and we will know more about Hilbert spaces. 

Now we turn to the study of $\sigma$-algebras and focus on those which have a countable set as their generator.

\Subsection{Countably Generated $\sigma$-Algebras}
\label{sec:count-gen}
Fix a measurable space $(X, {\cal A})$. The $\sigma$-algebra ${\cal A}$ is said to be \emph{\index{$\sigma$-algebra!countably generated}countably generated} iff there exists countable ${\cal A}_{0}$ such that ${\cal A} = \sigma({\cal A}_{0})$.

\BeginExample{probs-countably-gen}
Let $(X, \tau)$ be a topological space with a countable basis. Then $\Borel{X}$ is countably generated. In fact, if $\tau_{0}$ is the countable basis for $\tau$, then each open set $G$ can be written as $G = \bigcup_{n\in\Nat}G_{n}$ with $\Folge{G}\subseteq\tau_{0}$, thus each open set is an element of $\sigma(\tau_{0})$, consequently, $\Borel{X} = \sigma(\tau_{0})$. 
\EndExample

The observation in Example~\ref{probs-countably-gen} implies that the Borel sets for a separable metric space, in particular for a Polish space, is countably generated. 

Having  a countable dense subset for a metric space, we can use the
corresponding base for a fairly helpful characterization of the Borel sets. The
next Lemma says that the Borel sets are in this case countably generated.

\BeginLemma{DenseWithBorel}
Let $Y$ be a separable metric space with metric $d$. Denote by $B_{r}(y) := \{y'\in Y\mid d(y, y') < r\}$ the open ball with radius $r$ and center $y$. Then
\begin{equation*}
\Borel{Y} = \sigma(\{B_r(d) \mid r > 0\text{ rational}, d \in D\}),
\end{equation*}
where $D$ is countable and dense.
\EndLemma

\BeginProof
Because an open ball is an open set, we infer that
\begin{equation*}
\sigma(\{B_r(d) \mid r > 0\text{ rational}, d \in D\}) \subseteq \Borel{Y}.
\end{equation*}
Conversely, let $G$ be open. Then there exists a sequence $\Folge{B}$ of
open balls with rational radii such that
$
\bigcup_{n \in \Nat} B_n = G,
$
accounting for the other inclusion.
\EndProof

% This representation implies that the Borel sets  $\Borel{X}$ of our Polish
% space $X$ are countably generated.

Also the characterization of Borel sets in a metric space as the closure of the open (closed) sets under countable unions and countable intersections will be occasionally helpful.

\BeginLemma{simple-Borel}
The Borel sets in a metric space $Y$ are the smallest collection of
sets that contains the open (closed) sets and that are closed under countable
unions and countable intersections.
\EndLemma

\BeginProof
The smallest collection $\mathcal{G}$ of
sets that contains the open sets and that is closed under countable
unions and countable intersections is closed under complementation. This is
so since each closed set is a $G_\delta$ by Theorem~\ref{Polish-is-G-delta}.
Thus $\Borel{Y} \subseteq \mathcal{G}$; on the other hand
$\mathcal{G}\subseteq \Borel{Y}$ by construction.
\EndProof

The property of being countably generated is, however, not hereditary for a $\sigma$-algebra~---~a sub-$\sigma$-algebra of a countably generated $\sigma$-algebra is not necessarily countably generated. This is demonstrated by the following example. Incidentally, we will see in Example~\ref{inters-non-cg} that the intersection of two countably generated $\sigma$-algebras need not be countably generated again. This indicates that having a countable generator is a fickle property which has to be observed closely.

\BeginExample{probs-count-cocount}
Let 
\begin{equation*}
  {\cal C} := \{A\subseteq\Real \mid A\text{ or }\Real\setminus A\text{ is countable}\}
\end{equation*}
This $\sigma$-algebra is usually referred to the \emph{\index{$\sigma$-algebra!countable-cocountable}countable-cocountable $\sigma$-algebra}. Clearly, ${\cal C}\subseteq\Borel{\Real}$, and $\Borel{\Real}$ is countably generated by Example~\ref{probs-countably-gen}. But ${\cal C}$ is not countably generated. Assume that it is, so let ${\cal C}_{0}$ be a countable generator for ${\cal C}$; we may assume that every element of ${\cal C}_{0}$ is countable. Put $A := \bigcup{\cal C}_{0}$, then $A\in{\cal C}$, since $A$ is countable. But 
\begin{equation*}
{\cal D}
:= \{B\subseteq\Real \mid  B\subseteq A\text{ or }B\subseteq\Real\setminus A\} 
\end{equation*}
is a $\sigma$-algebra, and ${\cal D}= \sigma({\cal C}_{0})$.
On the other hand there exists $a\in\Real$ with $a\not\in A$, thus $A\cup\{a\}\in{\cal C}$ but $A\cup\{a\}\not\in{\cal D}$, a contradiction. 
\EndExample

Although the entire $\sigma$-algebra may not be countably generated, we may find for each element of a $\sigma$-algebra a countable generator:

\BeginLemma{exists-countable-gen}
Let ${\cal A}$ be a $\sigma$-algebra on a set $X$ which is generated by family ${\cal G}$ of subsets. Then we can find for each $A\in{\cal A}$ a countable subset ${\cal G}_{0}\subseteq{\cal G}$ such that $A\in\sigma({\cal G}_{0})$. 
\EndLemma

\BeginProof
Let ${\cal D}$ be the set of all $A\in{\cal A}$ for which the assertion is true, then ${\cal D}$ is closed under complements, and ${\cal G}\subseteq{\cal A}$. Moreover, ${\cal D}$ is closed under countable unions, since the union of a countable family of countable sets is countable again. Hence ${\cal D}$ is a $\sigma$-algebra which contains ${\cal G}$, hence it contains ${\cal A} = \sigma({\cal G})$. 
\EndProof

This has a fairly interesting and somewhat unexpected consequence, which will be of use later on. Recall that ${\cal A}\otimes{\cal B}$ is the smallest $\sigma$-algebra on $X\times Y$ which contains for measurable spaces $(X, {\cal A})$ and $(Y, {\cal B})$ all measurable rectangles $A\times B$ with $A\in{\cal A}$ and $B\in{\cal B}$. In particular, $\PowerSet{X}\otimes\PowerSet{X}$ is generated by $\{A\times B \mid A, B\subseteq X\}$. One may be tempted to assume that this $\sigma$-algebra is the same as $\PowerSet{X\times X}$, but this is not always the case, because we have

\BeginProposition{diagonal-in-powerset}
Denote by $\Delta_X$ the diagonal $\{\langle x, x\rangle \mid x\in X\}$ for a set $X$. Then $\Delta_X\in \PowerSet{X}\otimes\PowerSet{X}$ implies that the cardinality of $X$ does not exceed that of $\PowerSet{\Nat}$.
\EndProposition

\BeginProof
Assume $\Delta_X\in\PowerSet{X}\otimes\PowerSet{X}$, then there exists a
countable family ${\cal C}\subseteq\PowerSet{X}$ such that
$\Delta_X\in\sigma(\{A\times B \mid A, B\in{\cal C}\})$. The map $q: x
\mapsto \{C\in{\cal C}\mid x\in C\}$ from $X$ to $\PowerSet{{\cal C}}$
is injective. In fact, suppose it is not, then there exists $x\not=x'$
with $x\in C \Leftrightarrow x'\in C$ for all $C\in{\cal C}$, so we
have for all $C\in{\cal C}$ that either $\{x, x'\}\subseteq C$ or
$\{x, x'\}\cap C = \emptyset$, so that the pairs $\langle x, x\rangle$
and $\langle x', x'\rangle$ never occur alone in any $A\times B$ with
$A, B\in{\cal C}$. Hence $\Delta_X$ cannot be a member of
$\sigma(\{A\times B \mid A, B\in{\cal C}\})$, a contradiction. As a
consequence, $X$ cannot have more elements that $\PowerSet{\Nat}$.
\EndProof

Among the countably generated measurable spaces those are of interest
which permit to \emph{separate points}, so that if $x\not=x'$, we can
find $A\in{\cal C}$ with $x\in A$ and $x'\not\in A$; they are called
separable. Formally

\BeginDefinition{probs-sep-meas}
The $\sigma$-algebra ${\cal A}$ is called \emph{\index{separable}\index{$\sigma$-algebra!separable}separable} iff it is countably generated, and if for any two different elements of $X$ there exists a measurable set $A\in{\cal A}$ which contains one, but not the other. The measurable space $(X, {\cal A})$ is called separable iff its $\sigma$-algebra ${\cal A}$ is separable. 
\EndDefinition

The argumentation from Proposition~\ref{diagonal-in-powerset} yields

\BeginCorollary{probs-sep-meas-cor}
Let ${\cal A}$ be a separable $\sigma$-algebra over the set $X$ with ${\cal A} = \sigma({\cal A}_{0})$ for ${\cal A}_{0}$ countable. Then ${\cal A}_{0}$ separates points, and $\Delta_X\in{\cal A}\otimes{\cal A}$. 
\EndCorollary

\BeginProof
Because ${\cal A}$ separates points, we obtain from Example~\ref{prob-equiv-gen} that $\equiv_{{\cal A}_{0}}\ =\ \Delta_X$, where $\equiv_{{\cal A}_{0}}$ is the equivalence relation defined by ${\cal A}_{0}$. So ${\cal A}_{0}$ separates points. The representation 
\begin{equation*}
  X\times X\setminus \Delta_X = 
\bigcup_{A\in{\cal A}_{0}}A\times(X\setminus A)\cup(X\setminus A)\times A.
\end{equation*}
now yields $\Delta_X\in {\cal A}\otimes{\cal A}$. 
\EndProof

In fact, we can say even more.

\BeginProposition{sep-is-sep-metr}
A separable measurable space $(X, {\cal A})$ is isomorphic to $(X, \Borel{X})$ with the Borel sets coming from a metric $d$ on $X$ such that $(X, d)$ has is a separable metric space. 
\EndProposition

\BeginProof
1.
Let ${\cal A}_{0} = \{A_{n}\mid n\in\Nat\}$ be the countable generator for ${\cal A}$ which separates points. Define 
\begin{equation*}
  (M, {\cal M}) := (\{0, 1\}^{\Nat},
\bigotimes_{n\in\Nat}\PowerSet{\{0, 1\}})
\end{equation*}
as the product of countable many copies of the discrete space $(\{0, 1\}, \PowerSet{\{0, 1\}})$. Then $\bigotimes_{n\in\Nat}\PowerSet{\{0, 1\}}$ has as a basis the cylinder sets 
$
\{Z_{v}\mid v\in\{0, 1\}^{k}\text{ for some }k\in\Nat\}
$
with
$
Z_{v} := \{\Folge{t}\in M\mid \langle m_{1},\dots, m_{k}\rangle = v\}
$
for $v\in\{0, 1\}^{k}$, see~\pageref{cylinder-sets}. Define $f: X\to M$ through $f(x) := (\chi_{A_{n}}(x))_{n\in\Nat}$, then $f$ is injective, because ${\cal A}_{0}$ separates points. Put $Q := \Bild{f}{X}$, and ${\cal Q} := {\cal M}\cap Q$, the trace of ${\cal M}$ on $Q$. 

Now let $Y_{v} := Z_{v}\cap Q$ be an element of the generator for ${\cal Q}$ with $v = \langle m_{1}, \dots, m_{k}\rangle$ , then $\InvBild{f}{Y_{v}} = \bigcap_{j=1}^{k} C_{j}$ with $C_{j} := A_{j}$, if $m_{j} = 1$, and $C_{j} := X\setminus A_{j}$ otherwise. Consequently, $f: X\to Q$ is ${\cal A}$-${\cal Q}$-measurable. 

2.
Put for $x, y\in X$
\begin{equation*}
  d(x, y) := \sum_{n\in\Nat}2^{-n}\cdot \bigl|\chi_{A_{n}}(x) - \chi_{A_{n}}(y)\bigr|,
\end{equation*}
then $d$ is a metric on $X$ which has 
\begin{equation*}
{\cal G} := \bigl\{\bigcap_{j\in F}B_{j}\mid B_{j} \in{\cal A}_{0}\text{ or }X\setminus B_{j}\in{\cal A}_{0}, F\subseteq\Nat\text{ is finite}\bigr\}
\end{equation*}
as a countable basis. In fact, let $G\subseteq X$ be open; given $x\in
G$, there exists $\epsilon>0$ such that the open ball $
B_{\epsilon}(x) := \{x'\in X\mid d(x, x') < \epsilon\} $ with center
$x$ and radius $\epsilon$ is contained in $G$. Now choose $k$ with
$2^{-k}< \epsilon$, and put $v := \langle x_{1},\dots, x_{k}\rangle$,
then $x\in \bigcap_{j=1}^{k}B_{j}\subseteq B_{\epsilon}(x)$. This
argument shows also that ${\cal A} = \Borel{X}$. 

3.
Because $(X, d)$ has
a countable basis, it is a separable metric space. The map $f: X\to Q$
is a bijection which is measurable, and $f^{-1}$ is measurable as
well. This is so because $ \{A\in{\cal A}\mid \Bild{f}{A}\in{\cal Q}\}
$ is a $\sigma$-algebra which contains the basis ${\cal G}$.
\EndProof

This representation, which is due to \index{Mackey}Mackey, gives the representation
of separable measurable spaces as subspaces of the countable product
of the discrete space $(\{0, 1\}, \PowerSet{\{0, 1\}}$. This space is
also a compact metric space, so we may say that a separable measurable
space is isomorphic to a subspace of a compact metric space. We will
make use of this observation later on.

By the way, this innocently looking statement has some remarkable
consequences for our context. Just as an appetizer:

\BeginCorollary{Separable} Let $(X, \mathcal{A})$ be a separable
measurable space. Then
\begin{enumerate}
 \item \label{Separable-Eins} The diagonal $\Delta_{X}$ is measurable in the product, i.e.,
 \item \label{Separable-Zwei}
  If $f_i: X_i \rightarrow X$ is $\mathcal{A}_i - \mathcal{A}$-measurable,
  where $(X_i, \mathcal{A}_i)$ is a measurable space ($i =
  1, 2$), then
  $
  f_1^{-1}\left[\mathcal{A}\right] \otimes
  f_2^{-1}\left[\mathcal{A}\right] = (f_1 \times
  f_2)^{-1}\left[\mathcal{A} \otimes \mathcal{A}\right].
  $
\end{enumerate}
\EndCorollary

\BeginProof
1. Let $\Folge{A}$ be a generator for $X$ that separates point,
then
\begin{equation*}
(X \times X)\setminus\Delta_{X}
=
\bigcup_{n \in \Nat} \left(A_n \times X\setminus A_n \cup X\setminus A_n \times A_n\right),
\end{equation*}
which is a member of $\mathcal{A} \otimes \mathcal{A}$.

2.
The product $\sigma$-algebra
$
\mathcal{A} \otimes \mathcal{A}
$
is generated by the rectangles
$
B_1 \times B_2
$
with
$
B_i
$
taken from some generator
$
\mathcal{B}_0
$
for
$
\mathcal{B}\ (i = 1, 2).
$
Since
$
\left(f_1 \times f_2\right)^{-1}\left[B_1 \times B_2\right] =
f_1^{-1}\left[B_1\right] \times f_2^{-1}\left[B_2\right],
$
we see that
$
  (f_1 \times  f_2)^{-1}\left[\mathcal{B} \otimes \mathcal{B}\right]
  \subseteq f_1^{-1}\left[\mathcal{B}\right] \otimes
  f_2^{-1}\left[\mathcal{B}\right].
$
This is true without the assumption of separability. Now let
$\tau$ be a second countable metric topology on $Y$ with $
\mathcal{B} = \Borel{\tau} $ and let $\tau_0$ be a
countable base for the topology. Then
\begin{equation*}
\tau_p := \{T_1 \times T_2 \mid  T_1, T_2 \in \tau_0\}
\end{equation*}
is a countable base for the product topology
$
\tau \otimes \tau,
$
and (this is the crucial property)
\begin{equation*}
\mathcal{B} \otimes \mathcal{B} = \Borel{Y \times Y, \tau
\otimes \tau}
\end{equation*}
holds: because the projections from $X \times Y$ to $X$ and to $Y$
are measurable, we observe
$\mathcal{B} \otimes \mathcal{B} \subseteq \Borel{Y \times Y, \tau
\otimes \tau}$;
because $\tau_p$ is a countable base for the product
topology $\tau \otimes \tau$, we infer the other inclusion.

3.
Since for $T_1, T_2 \in
\tau_0$ clearly $$ f_1^{-1}\left[T_1\right] \times
f_2^{-1}\left[T_2\right] \in (f_1 \times
f_2)^{-1}\left[\tau_p\right] \subseteq (f_1 \times
f_2)^{-1}\left[\mathcal{B} \otimes \mathcal{B}\right] $$ holds, the
nontrivial inclusion is inferred from the fact that the smallest
$\sigma$-algebra containing $ \{f_1^{-1}\left[T_1\right] \times
f_2^{-1}\left[T_2\right] \mid  T_1, T_2 \in \tau_0\} $ equals
$ f_1^{-1}\left[\mathcal{B}\right] \otimes
  f_2^{-1}\left[\mathcal{B}\right].
$
\EndProof

Given a measurable function into a separable measurable space, we find that its kernel yields a measurable subset in the product of its domain. We will use the kernel for many a construction, so this little observation is quite helpful. 

\BeginCorollary{Kernel-is-Borel}
Let $f: X \rightarrow Y$ be a $\mathcal{A}$-$\mathcal{B}$-measurable map, where
$(X, \mathcal{A})$ and $(Y, \mathcal{B})$ are measurable spaces, the latter being separable. Then
the \emph{kernel} of $f$\index{map!kernel}
\begin{equation*}
\Kern{f} := \{\langle x_1, x_2\rangle \mid f(x_1) = f(x_2)\}
\end{equation*}
is a member of $\mathcal{A}\otimes\mathcal{A}$.
\EndCorollary

\BeginProof
Exercise~\ref{ex-prod-meas}.
\EndProof

The observation, made in the proof of Proposition~\ref{diagonal-in-powerset}, that it may not always be possible to separate two different elements in a measurable space through a measurable set led there to a contradiction. Nevertheless it leads to an interesting notion.

\BeginDefinition{probs-atom}
The set $A\in{\cal A}$ is called an \emph{\index{atom}atom} of ${\cal A}$ iff $B\subseteq A$ implies $B = \emptyset$ or $B = A$ for all $B\in{\cal A}$. 
\EndDefinition

For example, each singleton set $\{x\}$ is an atom for the $\sigma$-algebra $\PowerSet{X}$. Clearly, being an atom depends also on the $\sigma$-algebra. If $A$ is an atom, we have alternatively $B\subseteq A$ or $B\cap A = \emptyset$ for all $B\in{\cal A}$; this is more radical than being a $\mu$-atom, which merely restricts the values of $\mu(B)$ for measurable $B\subseteq A$ to $0$ or $\mu(A)$. Certainly, if $A$ is an atom, and $\mu(A)>0$, then $A$ is a $\mu$-atom. 

For a countably generated $\sigma$-algebra, atoms are easily identified.

\BeginProposition{probs-char-atom}
Let ${\cal A}_{0} = \{A_{n}\mid n\in\Nat\}$ be a countable generator of ${\cal A}$, and define 
\begin{equation*}
  A_{\alpha} := \bigcap_{n\in\Nat}A_{n}^{\alpha_{n}},
\end{equation*}
for $\alpha\in\{0, 1\}^{\Nat}$, where $A^{0} := A, A^{1} := X\setminus A$. Then 
$
\{A_{\alpha}\mid \alpha\in\{0, 1\}^{\Nat}, A_{\alpha}\not=\emptyset\}
$
is the set of all atoms of ${\cal A}$. 
\EndProposition

\BeginProof
Assume that there exist in ${\cal A}$ two different non-empty subsets $B_{1}, B_{2}$ of $A_{\alpha}$, and take $y_{1}\in B_{1}, y_{2}\in B_{2}$. Then $\isEquiv{y_{1}}{y_{2}}{\equiv_{{\cal A}_{0}}}$, but $\isEquiv{y_{1}}{y_{2}}{\not\equiv_{{\cal A}}}$, contradicting the observation in Example~\ref{prob-equiv-gen}. Hence $A_{\alpha}$ is an atom. Let $x\in A_{\alpha}$, then $A_{\alpha}$ is the equivalence class of $x$ with respect to the equivalence relation $\equiv_{{\cal A}_{0}}$, hence with respect to ${\cal A}$. Thus each atom is given by some $A_{\alpha}$. 
\EndProof

Incidentally, this gives another proof that the countable-cocountable $\sigma$-algebra over $\Real$ is not countably generated. Assume it is generated by $\{A_{n}\mid n\in\Nat\}$, then 
\begin{equation*}
  H := \bigcap\{A_{n}\mid A_{n}\text{ is cocountable}\}\cap\bigcap\{\Real\setminus A_{n}\mid A_{n}\text{ is countable}\}
\end{equation*}
is an atom, but $H$ is also cocountable. This is a contradiction to
$H$ being an atom. 

We relate atoms to measurable maps:

\BeginLemma{meas-const-on-atoms}
Let $f: X\to \Real$ be ${\cal A}$-$\Borel{\Real}$-measurable. If $A\in{\cal A}$ is an atom of ${\cal A}$, then $f$ is constant on $A$.
\EndLemma

\BeginProof
Assume that we can find $x_{1}, x_{2}\in A$ with $f(x_{1}) \not= f(x_{2})$, say, $f(x_{1}) < c < f(x_{2})$. Then $\{x\in A\mid f(x) < c\}$ and $\{x\in A\mid f(x) > c\}$ are two non-empty disjoint measurable subsets of $A$. This contradicts $A$ being an atom. 
\EndProof

We will specialize now our view of measurable spaces to the Borel sets of Polish spaces and their more general cousins, analytic sets.

%%% Local Variables: 
%%% mode: latex
%%% TeX-master: "../Mskr3"
%%% End: 

%\Input{\Folder/BorelSets}
%spell checked - 24Aug14
\Subsubsection{Borel Sets in Polish and Analytic Spaces}
\label{sec:borel-and-analytic}

General measurable spaces and even separable metric spaces are sometimes too general for supporting
specific structures. We deal with Polish  and analytic spaces which
are general enough to support interesting applications, but have
specific properties which help establishing vital properties.
We remind the reader first of some basic facts and
provide then some helpful tools for working with Polish spaces, and
their more general cousins, analytic spaces.

Fix for the time being $(X, \tau)$ as a topological space.  Recall
that a family $\mathcal{B}\subseteq\tau$ of open subsets of $X$ is called a
\emph{base} for topology $\tau$ iff each element of
$\tau$ can be represented as the union of elements of
$\mathcal{B}$.  This is equivalent to saying
that~\index{space!topological!base}\index{topology!base} $ \bigcup \{B
\mid B \in \mathcal{B}\} = X, $ and that we can find for each $x \in
B_1 \cap B_2$ with $B_1, B_2 \in \mathcal{B}$ an element $B_3 \in
\mathcal{B}$ with $x \in B_3 \subseteq B_1 \cap B_2$. A \emph{subbase}
$\mathcal{S}$ for $\tau$ has the property that the
set~\index{space!topological!subbase}\index{topology!subbase} $
\{\bigcap \mathcal{F} \mid \mathcal{F} \subseteq \mathcal{S} \text{
  finite}\} $ of finite intersections of elements of $\mathcal{S}$
forms a base for $\tau$.

Given another topological space $(Y, \theta)$, a map
$f: X \rightarrow Y$ is called \emph{$\tau$-$\theta$-continuous} iff
\index{space!topological!continuity}\index{continuity}\index{map!continuous}
the inverse image of an open set from $Y$ is open in $X$ again, i.e.,
iff $\InvBild{f}{\theta} \subseteq \tau$. The topological spaces
$(X,\tau)$ and $(Y,\theta)$  are called \emph{homeomorphic}
\index{space!topological!homeomophism}\index{homeomorphism}
iff there
exists a $\tau$-$\theta$-continuous bijection $f: X \rightarrow Y$ the inverse of
which is $\theta$-$\tau$-continuous.

Proceeding in analogy to measurable spaces, a topology $\tau$ on a
set $X$ is called \emph{initial}\index{topology!initial}
for a map $f: X \rightarrow Y$ with a topological space $(Y,\theta)$ iff
$\tau$ is the smallest topology $\tau_0$ on $X$
rendering $f$ a $\tau_0$-$\theta$-continuous
map. For example, if $Y \subseteq X$ is a subset, then the topological
subspace $(Y, \{Y \cap G \mid G \in \tau\})$ is just the initial
topology with respect to the inclusion map $i_Y: Y \rightarrow
X$\index{topology!subspace}.

Dually, if $(X,\tau)$ is a topological space and $f: X \rightarrow Y$
is a map, then the \emph{final topology}\index{topology!initial} $\mathcal{S}$
on $Y$ is the largest topology
$\mathcal{S}_0$ on $Y$ making $f$  $\mathcal{T}$-$\mathcal{S}_0$-continuous.
Both initial and final topologies generalize to families of spaces and maps.

The \emph{topological product}
$\prod_{i \in I} (X_i, \tau_i)$\index{topology!product}
of the topological spaces
$\left((X_i, \tau_i)\right)_{i \in I}$ is the Cartesian product $\prod_{i \in I}X_i$ endowed with
the initial topology with respect to the projections, and the \emph{topological sum}
$\coprod_{i \in I} (X_i, \tau_i)$\index{topology!sum} of the topological spaces
$\left((X_i, \tau_i)\right)_{i \in I}$ is the direct $\coprod_{i \in I}X_i$ endowed with
the final topology with respect to the injections.

An immediate consequence of Lemma~\ref{ElemGen} is that continuity implies Borel measurability.

\BeginLemma{ContImplBorel}
Let $(X_1, \tau_1)$ and $(X_2, \tau_2)$ be topological
spaces. Then $f: X_1 \rightarrow X_2$ is
$\Borel{\tau_1}$-$\Borel{\tau_2}$ measurable, provided $f$ is
$\tau_1$-$\tau_2$-continuous. \QED
\EndLemma

We note for later use that the limit of a sequence of measurable functions into a metric space is measurable again, see Exercise~\ref{ex-limit-is-measurable}.

\BeginProposition{limit-is-measurable}
Let $(X, \mathcal{A})$ be a measurable, $(Y, d)$ a metric space, and $\Folge{f}$ be a sequence of $\mathcal{A}$-$\Borel{Y}$-measurable functions $f_n: X \to Y$.  Then
\begin{itemize}
  \item the set $C := \{x \in X \mid (f_n(x))_{n\in \Nat}\text{ exists}\}$ is measurable,
  \item $f(x) := \lim_{n\to\infty} f_n(x)$ defines a $\mathcal{A}\cap C$-$\Borel{Y}$-measurable map $f: C\to Y$
\end{itemize}
\EndProposition

Neither general topological spaces nor metric spaces offer a structure
rich enough for the study of the transition systems that we will enter
into. We need to restrict the class of topological spaces to a
particularly interesting class of spaces that are traditionally called
\emph{Polish}.

As far as notation goes, we will write down a topological or a metric
space without its adornment through a topology or a metric, unless
this becomes really necessary.

Remember that a metric space $(X, d)$ is called \emph{complete} iff
each $d$-Cauchy sequence has a
limit\index{space!metric!complete}. Recall also that completeness is
really a property of the metric rather than the underlying topological
space, so a metrizable space may be complete with one metric and
incomplete with another one. In contrast, having a countable base is a
topological property which is invariant under the different metrics
the topology may admit.

\BeginDefinition{PolishSpace}
A \emph{Polish space}\index{space!Polish} $X$ is a topological space
the topology of which is metrizable through a complete metric, and
which has a countable base, or, equivalently, a countable dense
subset.
\EndDefinition

Familiar spaces are Polish, as these examples show.
\BeginExample{AllPolish-1}
The real $\Real$ with their usual topology, which is
induced by the open intervals, are a Polish space.
\EndExample

\BeginExample{AllPolish-2}
The open unit interval $]0, 1[$ with the usual topology
induced by the open intervals form a Polish space.

This comes probably as a surprise, because $]0, 1[$ is known
not to be complete with the usual metric. But all we need is
a dense subset (take the rationals
$
\Rational \cap ]0, 1[
$),
and a metric that generates the topology, and that is complete.
Define
\begin{equation*}
d(x, y) := \left|\ln\frac{x}{1-x} - \ln\frac{y}{1-y}\right|,
\end{equation*}
then this is a complete metric for $]0, 1[$. This is so since
$x \mapsto \ln(x/(1-x))$ is a continuous bijection from $]0, 1[$
to $\Real$, and the inverse $y \mapsto e^y/(1 + e^y)$
is also a continuous bijection.
\EndExample

\BeginLemma{ClosedIsPolish}
Let $X$ be a Polish space, and assume that $F \subseteq X$ is closed, then the
subspace $F$ is Polish as well.
\EndLemma

\BeginProof
Because $F$ is closed, each Cauchy sequence in $F$ has its limit
in $F$, so $F$ is complete. The topology that $F$ inherits
from $X$ has a countable base and is metrizable, so $F$ has a
countable dense subset, too.
\EndProof

\BeginLemma{Products+Coproducts}
Let $\Folge{X}$ be a sequence of Polish spaces, then
the product $\prod_{n \in \Nat} X_n$ and the coproduct $\coprod_{n \in \Nat} X_n$
are Polish spaces.
\EndLemma

\BeginProof
Assume that the topology $\tau_n$ on $X_n$ is metrized through metric
$d_n$, where it may be assumed that $d_n \leq 1$ holds (otherwise
use for $\tau_n$ the complete metric $d_n(x, y)/(1 + d_n(x, y))$). Then
\begin{equation*}
d(\Folge{x}, \Folge{y}) := \sum_{n \in \Nat} 2^{-n} d_n(x_n, y_n)
\end{equation*}
is a complete metric for the product topology $\prod_{n \in \Nat} \tau_n$.
For the coproduct, define the complete metric
\begin{equation*}
d(x, y) :=
\begin{cases}
  2, & \text{if } x \in X_n, y \in X_m, n\not= m\\
  d_n(x, y), & \text{if } x, y \in X_n.
\end{cases}
\end{equation*}

All this is established through standard arguments.
\EndProof

\BeginExample{Inf-Nat}
The set $\Nat$ of natural numbers with the discrete topology is a Polish space on
account of being the topological sum of its elements. Thus the set
$\Nat^\infty$ of all infinite sequences is a Polish space. The sets
\begin{equation*}
\Sigma_\alpha := \{\tau \in \Nat^\infty \mid \alpha\text{ is an initial piece of }\tau\}
\end{equation*}
for $\alpha \in \Nat^*$, the free monoid generated by $\Nat$, constitute
a base for the product topology.
\EndExample

This last example will be discussed in much greater detail later on.
It permits sometimes reducing the discussion of properties for general
Polish spaces to an investigation of the corresponding properties
of $\Nat^\infty$, the structure of the latter space being
more easily accessible than that of a general space.
We apply Example~\ref{Inf-Nat} directly to show that
all open subsets of a metric space $X$ with a countable
base can be represented through a single closed set in $\Nat^\infty \times X$.

\Rand{Bekannt.}
Recall that for $ D \subseteq X \times Y$ the \emph{vertical
cut}\index{cut!vertical} $D_{x}$ is defined through 
$
D_x := \{y \in Y \mid  \langle x, y\rangle \in D\}
$
and the \emph{horizontal cut}\index{cut!horizontal} $D^{y}$ is
$
D^y := \{x \in X \mid  \langle x, y\rangle \in D\}.
$
Note that $\bigl((X \times Y)\setminus D\bigl)_x = Y\setminus D_x$.

\BeginProposition{universal-open}
Let $X$ be a separable metric space.
Then there exists an open set $U \subseteq \Nat^\infty \times X$ and a closed set
$F \subseteq \Nat^\infty \times X$ with these properties:
\begin{enumerate}[a.]
  \item \label{universal-open-a}For each open set $G \subseteq X$ there exists $t \in \Nat^\infty$
such that $G = U_t$.
  \item \label{universal-open-b}For each closed set $C \subseteq X$ there exists $t \in \Nat^\infty$
such that $C = F_t$.
\end{enumerate}
\EndProposition

\BeginProof
0.
It is enough to establish the property for open sets; taking complements
will prove it for closed ones.

1.
Let $\Folge{V}$ be a basis for the open sets in $X$ with $V_n \not= \emptyset$
for all $n \in \Nat$. Define
\begin{equation*}
U := \{\langle t, x \rangle \mid t\in \NatInf, x \in \bigcup_{n \in \Nat} V_{t_n}\},
\end{equation*}
then $U \subseteq \Nat^\infty \times X$ is open. In fact, let $\langle t, x\rangle \in U$,
then there exists $n \in \Nat$ with $x \in V_n$, thus
$
\langle t, x\rangle \in \Sigma_n \times V_n \subseteq U,
$
and $\Sigma_n \times V_n$ is open in the product.

2.
Let $G \subseteq X$ be open. Because $\Folge{V}$ is a basis for the
topology,  there exists a sequence $t \in \Nat^\infty$
with $G = \bigcup_{n \in \Nat} V_{t_n} = U_t$.
\EndProof

The set $U$ is usually called a \emph{universal open set},
similar\index{universal set} for $F$. These universal sets will be
used rather heavily when we discuss analytic sets.

We have seen that a closed subset of a Polish space is a Polish
space in its own right; a similar argument shows that an open
subset of a Polish space is Polish as well. Both observations turn out to
be special cases of the characterization of Polish subspaces through $G_\delta$-sets.

We need for this characterization an auxiliary statement due to Kuratowski which permits the extension
of a continuous map from a subspace to a $G_\delta$-set containing it~---~just far
enough to be interesting to us.
Denote by $\Closure{A}$ the topological
closure of a set $A$\index{closure of a set}.

\BeginLemma{Kuratowski}
Let $Y$ be a complete metrizable space, $W$ a metric space, then a continuous
map $f: A \rightarrow Y$ can be extended to a  continuous map
$
f_*: G \rightarrow Y
$
with $G$ a $G_\delta$-set such that $A \subseteq G \subseteq \Closure{A}$.
\EndLemma

\BeginProof
1.
We may and do assume that the complete metric $d$ for $Y$ is bounded
by $1$, otherwise we move to the equivalent and complete metric $\langle x, y\rangle \mapsto d(x, y)/(1 + d(x, y))$, see Exercise~\ref{ex-equiv-compl-metric}. 
The \emph{oscillation} ${\o}_f(x)$\index{oscillation} of $f$ at $x \in \Closure{A}$ is defined as the smallest
diameter of the image of an open neighborhood of $x$, formally,
\begin{equation*}
{\o}_f(x) := \inf\{\mathsf{diam}(\Bild{f}{A \cap V}) \mid x \in V, V \text{ open}\}.
\end{equation*}
Because $f$ is continuous on $A$, we have ${\o}_f(x) = 0$ for each element $x$ of $A$. In fact, let $\epsilon>0$  be given, then there exists $\delta>0$ such that $\mathsf{diam}(\Bild{f}{A \cap V}) < \epsilon$, whenever $V$ is a neighborhood of $x$ of diameter less than $\delta$. Thus ${\o}_{f}(x) < \epsilon$; since $\epsilon>0$ was chosen to be arbitrary, the claim follows. 

2.
Put
$
G := \{x \in \Closure{A} \mid {\o}_f(x) = 0\},
$
then $A \subseteq G \subseteq \Closure{A}$, and $G$ is a $G_\delta$ in $W$. In fact,
represent $G$ as
\begin{equation*}
G = \bigcap_{n \in \Nat} \{x \in \Closure{A} \mid {\o}_f(x) < \frac{1}{n}\},
\end{equation*}
so we have to show that $\{x \in \Closure{A} \mid {\o}_f(x) < q\}$ is open in $\Closure{A}$ fr any $q>0$. But we have 
\begin{equation*}
\{x \in \Closure{A} \mid {\o}_f(x) < q\}
=
\bigcup \{V \cap \Closure{A} \mid \mathsf{diam}(\Bild{f}{V \cap A}) < q\}.
\end{equation*}
This is the union of sets open in $\Closure{A}$, hence is an open set itself. Note that $\Closure{A}$ is ---~as a closed set ~---~a $G_\delta$ in $W$.

3.
Now take an element $x \in G \subseteq \Closure{A}$. Then there exists a sequence
$\Folge{x}$ of elements $x_n \in A$ with $x_n \rightarrow x$. Given $\epsilon > 0$,
we find a neighborhood $V$ of $x$ with $\mathsf{diam}(\Bild{f}{A \cap V}) < \epsilon$, since the oscillation of $f$ at $x$ is $0$. Because
$x_n \rightarrow x$, we know that we can find an index $n_{\epsilon}\in\Nat$ such that $x_m \in V \cap A$ for all $m > n_\epsilon$. This implies that
the sequence $(f(x_n))_{n \in \Nat}$ is a Cauchy sequence in $Y$. It converges
because $Y$ is complete. Put
\begin{equation*}
f_*(x) := \lim_{n \rightarrow \infty} f(x_n).
\end{equation*}

4.
We have to show now that 
\begin{itemize}
\item $f_{*}$ is well-defined.
\item $f_{*}$ extends $f$.
\item $f_{*}$ is continuous.
\end{itemize}
Assume that we can find $x\in G$ such that $\Folge{x}$ and $\Folge{x'}$ are sequences in $A$ with $x_{n}\to x$ and $x'_{n}\to x$, but $\lim_{n\to \infty}f(x_{n}) \not= \lim_{n\to \infty}f(x'_{n})$. Thus we find some $\eta>0$ such that $d(f(x_{n}), f(x'_{n})) \geq \eta$ infinitely often. Then the oscillation of $f$ at $x$ is at least $\eta>0$, a contradiction. This implies that $f_{*}$ is well-defined, and it implies also that $f_{*}$ extends $f$. Now let $x\in G$. If $\epsilon>0$ is given, we find a neighborhood $V$ of $x$  with $\mathsf{diam}(\Bild{f}{A\cap V})<\epsilon$. Thus, if $x'\in G\cap V$, then $d(f_{*}(x), f_{*}(x')) < \epsilon$. Hence $f_{*}$ is continuous. 
\EndProof

This technical Lemma is an important step in establishing a far reaching characterization
of subspaces of Polish spaces that are Polish in their own right. We will show now that a subset $X$ of a Polish space is a Polish space in its own right iff it is a $G_{\delta}$-set. We will present Kuratowski's proof for it. It is not difficult to show that $X$ must be a $G_{\delta}$-set, using Lemma~\ref{Kuratowski}. The tricky part is the converse, and at its very center is the following idea: assume that we have represented $X = \bigcap_{k\in\Nat}G_{k}$ with each $G_{k}$ open, and assume that we have a Cauchy sequence $\Folge{x}\subseteq X$ with $x_{n}\to x$. How do we prevent $x$ from being outside $X$? Well, what we will do is to set up \index{Kuratowski's trap}Kuratowski's trap, preventing the sequence to wander off. The trap is a new complete and equivalent metric $D$, which is makes it impossible for the sequence to behave bad. So if $x$ is trapped to be an element of $X$, we may conclude that $X$ is complete, and the assertion may be established. 

Before we begin with the easier half, we fix a Polish space $Y$ and a complete metric $d$ on $Y$.

\BeginLemma{polish-is-g-delta-easy}
If $X\subseteq Y$ is a Polish space, then $X$ is a $G_{\delta}$-set.
\EndLemma

\BeginProof
$X$ is complete, hence closed in $Y$. The identity $id_{X}: X\to Y$ can be extended continuously by Lemma~\ref{Kuratowski} to a $G_{\delta}$-set $G$ with $X \subseteq G \subseteq \Closure{X}$, thus $G = X$, so $X$ is a $G_{\delta}$-set. 
\EndProof

Now let $X = \bigcap_{k\in\Nat}G_{k}$ with $G_{k}$ open for all $k\in\Nat$. In order to prepare for Kuratowski's trap, we define 
\begin{equation*}
  f_{k}(x, x') := \bigl|\frac{1}{d(x, Y\setminus G_{k})} - \frac{1}{d(x', Y\setminus G_{k})}\bigr|
\end{equation*}
for $x, x'\in X$. Because $G_{k}$ is open, we have $x\in G_{k}$ iff $d(x, Y\setminus G_{k})>0$, so $f_{k}$ is a finite and continuous function on $X\times X$. Now let 
\begin{align*}
  F_{k}(x, x') & := \frac{f_{k}(x, x')}{1 + f_{k}(x, x')},\\
D(x, x') & := d(x, x') + \sum_{k\in\Nat}2^{-k}\cdot F_{k}(x, x').
\end{align*}
for $x, x'\in X$. Then $D$ is a metric on $X$ (cp. Exercise~\ref{ex-equiv-compl-metric}), and the metrics $d$ and $D$ are equivalent on $X$. Because $d(x, x') \leq D(x, x')$, it is clear that the identity $id: (X, D)\to (X, d)$ is continuous, so it remains to show that $id: (X, d)\to (X, D)$ is continuous. Let $x\in X$ be given, and let $\epsilon>0$, then we find $\ell\in\Nat$ such that 
$
\sum_{k > \ell}2^{-j}\cdot F_{k}(x, x') < \epsilon/3
$ 
for all $x'\in X$. For $k = 1, \dots, \ell$ there exists $\delta_{j}$ such that $F_{j}(x, x')< \epsilon/(3\cdot\ell)$, whenever $d(x, x') < \delta_{j}$, since $x\mapsto d(x, Y\setminus G_{j})$ is positive and continuous. Thus define 
$
\delta := \epsilon/3\wedge\delta_{1}\wedge\dots\wedge\delta_{\ell},
$
then $d(x, x') < \delta$ implies 
\begin{equation*}
  D(x, x') \leq d(x, x') + \sum_{k=1}^{\ell}2^{-j}\cdot F_{j}(x, x') + \frac{\epsilon}{3}
< \frac{\epsilon}{3} + \sum_{k=1}^{\ell}\frac{\epsilon}{3\cdot \ell} + \frac{\epsilon}{3}
= \epsilon.
\end{equation*}
Thus $(X, d)$ and $(X, D)$ have in fact the same open sets. When establishing that $(X, D)$ is complete, we spring Kuratowski's trap. Let $\Folge{x}$ be a $D$-Cauchy sequence. Then this sequence is also a $d$-Cauchy sequence, thus we find $x\in Y$ such that $x_{n}\to  x$, because $(Y, d)$ is complete. We claim that $x\in X$. In fact, if $x\in X$, we find $G_{\ell}$ with $x\not\in G_{\ell}$, so that we can find for each $\epsilon>0$ some index $n_{\epsilon}\in\Nat$ with $F_{\ell}(x_{n}, x_{m})\geq 1 - \epsilon$ for $n, m\geq n_{\epsilon}$. But then $D(x_{n}, x_{m})\geq (1-\epsilon)/2^{\ell}$ for $n, m\geq n_{\epsilon}$, so that $\Folge{x}$ cannot be a $D$-Cauchy sequence. Consequently, $X$ is complete, hence closed. 

Thus we have established:

\BeginTheorem{Polish-is-G-delta}
Let $Y$ be a Polish space. Then the subspace $X \subseteq Y$ is a Polish space
iff $X$ is a $G_\delta$-set.
\QED
\EndTheorem

In particular, open and closed subsets of Polish spaces are Polish
spaces in their subspace topology. Conversely, each Polish space can
be represented as a $G_\delta$-set in the \emph{Hilbert
  cube}\index{Hilbert cube} $[0, 1]^\infty$; this is the famous
characterization of Polish spaces due to
Alexandrov\index{theorem!Alexandrov}~\cite[III.33.VI]{Kuratowski}.

\BeginTheorem{Alexandrov}\textbf{(Alexandrov)}
Let $X$ be a separable metric space, then $X$ is homeomorphic to
a subspace of the Hilbert cube. If $X$ is Polish, this subspace
is a $G_\delta$.
\EndTheorem

\BeginProof
1.
We may and do assume again that the metric $d$ is bounded by 1. Let
$\Folge{x}$ be a countable and dense subset of $X$, and put
\begin{equation*}
f(x) := \langle d(x, x_1), d(x, x_2), \dots \rangle.
\end{equation*}
Then $f$ is injective and continuous. Define
$
g: \Bild{f}{X} \rightarrow X
$
as $f^{-1}$, then $g$ is continuous as well: assume that
$f(y_m) \rightarrow f(y)$ for some $y$, hence
$
\lim_{m \rightarrow \infty} d(y_m, x_n) = d(y, x_n)
$
for each $n \in \Nat$. Since $\Folge{x}$ is dense, we find for
a given $\epsilon > 0$ an index $n$ with $d(y, x_n) < \epsilon$;
by construction we find for $n$ an index $m_0$ with $d(y_m, x_n) < \epsilon$
whenever $m > m_0$. Thus $d(y_m, y) < 2\cdot\epsilon$ for $m > m_0$,
so that $y_m \rightarrow y$. This demonstrates that $g$ is continuous,
thus $f$ is a homeomorphism.

2.
If $X$ is Polish, $\Bild{f}{X} \subseteq [0, 1]^\infty$ is Polish as
well. Thus the second assertion follows from Theorem~\ref{Polish-is-G-delta}.
\EndProof

Recall that a topological Hausdorff space $X$ is \emph{compact} iff
each open cover of $X$ contains a finite cover of
$X$\index{space!topological!compact}\index{topology!compact}.  This
property of compact spaces will be used from time to time. The
Bolzano-Weierstraß Theorem implies that compact metrizable spa\-ces
are Polish. It is inferred from Tihonov's Theorem that the Hilbert
cube $[0, 1]^\infty$ is compact, because the unit interval $[0, 1]$ is
compact, again by the Bolzano-Weierstraß Theorem. Thus Alexandrov's
Theorem~\ref{Alexandrov} embeds a Polish space as a $G_\delta$ into a
compact metric space, the closure of which will be compact.

\Subsubsection{Manipulating Polish Topologies}
\label{sec:properties}

We will show now that a Borel map between Polish spaces can be turned
into a continuous map. Specifically, we will show that, given a measurable map between
Polish spaces, we can find  on the domain a finer Polish topology
with the same Borel sets which renders the map continuous.
This will be established through a sequence of auxiliary statements, each
of which will be of interest and of use in its own right.

We fix for the discussion to follow a Polish space $X$ with topology
$\tau$. Recall that a set is \emph{clopen}\index{set!clopen} in a
topological space iff it is both closed and open.

\BeginLemma{Take-closed-Polish}
Let $F$ be a closed set in $X$. Then there exists a Polish topology $\tau'$ such
that $\tau \subseteq \tau'$ (hence $\tau'$ is finer than $\tau$),
$F$ is clopen in $\tau'$, and
$
\Borel{\tau} = \Borel{\tau'}.
$
\EndLemma

\BeginProof
Both $F$ and $X\setminus F$ are Polish by Theorem~\ref{Polish-is-G-delta}, so the
topological sum of these Polish spaces is Polish again by Lemma~\ref{Products+Coproducts}.
The sum topology is the desired topology.
\EndProof

We will now add a sequence of certain Borel sets to the topology; this will happen step by step, so we should know how to manipulate a sequence of Polish topologies. This is explained now.

\BeginLemma{Take-Sequence-Polish}
Let $\Folge{\tau}$ be a sequence of Polish topologies $\tau_n$
with $\tau \subseteq \tau_n$.
\begin{enumerate}
\item\label{Take-Sequence-Polish-1} The topology $\tau_\infty$
generated by $\bigcup_{n \in \Nat} \tau_n$ is Polish.
\item\label{Take-Sequence-Polish-2} If $\tau_n \subseteq \Borel{\tau}$,
then $\Borel{\tau_\infty} = \Borel{\tau}$.
\end{enumerate}
\EndLemma

\BeginProof
1.
The product $\prod_{n \in \Nat} (X_n, \tau_n)$ is by Lemma~\ref{Products+Coproducts}
a Polish space, where $X_n = X$ for all $n$. Define
the map $f: X \rightarrow \prod_{n \in \Nat} X_n$ through $x \mapsto \langle x, x, \dots\rangle$,
then $f$ is $\tau_\infty$-$\prod_{n \in \Nat} \tau_n$-continuous by construction.
One infers that $\Bild{f}{X}$ is a closed subset of $\prod_{n \in \Nat} X_n$: if
$\Folge{x} \notin \Bild{f}{X}$, take $x_i \not= x_j$ with $i < j$, and let $G_i$ and $G_j$ be disjoint
open neighborhoods of $x_i$ resp. $x_j$. Then
\begin{equation*}
\prod_{\ell < i} X_\ell \times G_i \times \prod_{i < \ell < j} X_\ell \times G_j \times \prod_{\ell > j} X_\ell
\end{equation*}
is an open neighborhood of $\Folge{x}$ that is disjoint from $\Bild{f}{X}$. By Lemma~\ref{ClosedIsPolish},
the latter set is Polish. On the other hand, $f$ is a homeomorphism between $(X, T_\infty)$ and
$\Bild{f}{X}$, which establishes part\emph{~\ref{Take-Sequence-Polish-1}}.

2.
$\tau_n$ has a countable basis
$
\{U_{i, n} \mid i \in \Nat\},
$
with $U_{i, n} \in \Borel{X, \tau}$, since $\tau_{n}\subseteq\Borel{\tau}$. This implies that $\tau_\infty$
has
$
\{U_{i, n} \mid i, n \in \Nat\}
$
as a countable basis, which entails
$
\Borel{X, \tau_\infty} \subseteq \Borel{X, \tau}.
$
The other inclusion is obvious, giving part\emph{~\ref{Take-Sequence-Polish-2}}.
\EndProof

As a consequence, we may add a Borel set to a Polish topology as a clopen set without
destroying the property of the space to be Polish or changing the
Borel sets. This is true as well for sequences of Borel sets,
as we will see now.

\BeginProposition{ClopenSets}
 If
$\Folge{B}$ is a sequence of Borel sets in $X$, then there
  exists a Polish topology $\tau_0$ on $X$ such that
  $\tau_0$ is finer than $\tau$,
  $\tau$ and $\tau_0$ have the same Borel sets, and each $B_n$
  is clopen  in $\tau_0$.
\EndProposition

\BeginProof
1.
We show first that we may add just one Borel set to the topology without
changing the Borel sets. In fact, call a Borel set $B \in \Borel{\tau}$
\emph{neat} if there exists
a Polish topology $\tau_B$ that is finer than $\tau$ such that $B$ is clopen with respect to $tau_B$, and
$\Borel{\tau} = \Borel{\tau_B}$.
\begin{equation*}
\mathcal{H} := \{B \in \Borel{\tau} \mid B\text{ is neat} \}.
\end{equation*}
Then $\tau \subseteq \mathcal{H}$, and each closed set is a
member of $\mathcal{H}$ by Lemma~\ref{Take-closed-Polish}.
Furthermore, $\mathcal{H}$ is closed under
complements by construction, and closed under countable unions by Lemma~\ref{Take-Sequence-Polish}.
Thus we may now infer that
$\mathcal{H} = \Borel{\tau}$, so that each Borel set is neat.

2.
Now construct inductively Polish topologies $\tau_n$ that are
finer than $\tau$ with $\Borel{\tau} = \Borel{\tau_n}$.
Start with $\tau_0 := \tau$.
Adding $B_{n+1}$ to the Polish topology $\tau_n$ according to
the first part yields a finer Polish topology $\tau_{n+1}$ with
the same Borel sets. Thus the assertion follows from Lemma~\ref{Take-Sequence-Polish}.
\EndProof

This permits turning a Borel map into a continuous one, whenever the
domain is Polish and the range is a second countable metric space.

\BeginProposition{Srivastava} Let  $Y$ a separable metric space
with topology $\theta$. If
$f: X \rightarrow Y$ is a $\Borel{\tau}$-$\Borel{\theta}$-Borel
 measurable map, then there
 exists a Polish topology $\tau'$ on $X$ such that
  $\tau'$ is finer than $\tau$, $\tau$ and $\tau'$
  have the same Borel sets, and $f$ is $\tau'$-$\theta$ continuous.
\EndProposition

\BeginProof
The metric topology $\theta$ is generated from the countable basis $\Folge{H}$.
Construct from the Borel sets $\InvBild{f}{H_n}$ and from $\tau$ a Polish
topology $\tau'$ according to Proposition~\ref{ClopenSets}. Because
$\InvBild{f}{H_n} \in \tau'$  for all $n \in \Nat$, the inverse image of each
open set from $\theta$ is $\tau'$-open, hence $f$ is $\tau'$-$\theta$ continuous.
The construction entails $\tau$ and $\tau'$ having the same Borel sets.
\EndProof

This property is most useful, because it permits rendering
measurable maps continuous, when they go into a second countable
metric space (thus in particular into  a Polish space).

As a preparation for dealing with analytic sets, we will show now that
each Borel subset of the Polish space $X$ is the continuous image
of $\Nat^\infty$. We begin with a reduction of the problem space: it
is sufficient to establish this property for closed sets. This is justified by the following observation. 

\BeginLemma{ClosedImageIsEnough}
Assume that each closed set in the Polish space $X$ is a
continuous image of $\Nat^\infty$. Then each Borel set of $X$ is a continuous
image of $\Nat^\infty$.
\EndLemma

\BeginProof (Sketch)
1.
Let
\begin{equation*}
\mathcal{G} := \{B \in \Borel{X} \mid B = \Bild{f}{\Nat^\infty}\text{ for }
f: \Nat^\infty \rightarrow X\text{ continuous}\}
\end{equation*}
be the set of all good guys. Then $\mathcal{G}$ contains by assumption all
closed sets. We show that $\mathcal{G}$ is closed under countable unions and countable
intersections. Then the assertion will follow from Lemma~\ref{simple-Borel}.

2.
Suppose $B_n = \Bild{f_n}{\Nat^\infty}$ for the continuous map $f_n$, then
\begin{equation*}
\mathbb{M} := \{\langle t_1, t_2, \dots\rangle \mid f_1(t_1) = f_2(t_2) = \dots\}
\end{equation*}
is a closed subset of $(\Nat^\infty)^\infty$, and defining
$
f: \langle t_1, t_2, \dots\rangle \mapsto f_1(t_1)
$
yields a continuous map $f: \mathbb{M} \rightarrow X$ with
$
\Bild{f}{\mathbb{M}} = \bigcap_{n \in \Nat} B_n.
$
$\mathbb{M}$ is homeomorphic to $\Nat^\infty$. Thus $\mathcal{G}$ is closed under
countable intersections.

3.
We show that $\mathcal{G}$ is closed also under countable unions. In fact, let $B_{n}\in{\cal G}$ such that $B_{n}=\Bild{f_{n}}{\NatInf}$ with $f_{n}: \NatInf\to X$ continuous. Define
\begin{equation*}
f:
  \begin{cases}
   \NatInf & \to X\\
\langle n, t_{1}, t_{2}, \dots, \rangle & \mapsto f_{n}(t_{1}, t_{2}, \dots). 
  \end{cases}
\end{equation*}
Thus 
\begin{equation*}
\Bild{f}{\NatInf} = \bigcup_{n\in\Nat}\Bild{f_{n}}{\NatInf} = \bigcup_{n\in\Nat}B_{n}.
\end{equation*}
Moreover, $f$ is continuous. If $G\subseteq X$ is open, we have $\InvBild{f}{G} = \bigcup_{n\in\Nat}\{n\}\times\InvBild{f_{n}}{G}$. Since $\InvBild{f_{n}}{G}$ is pen for each $n\in\Nat$, we conclude that $\InvBild{f}{G}$ is open, so that $f$ is indeed continuous. Thus ${\cal G}$ is closed under countable unions, and the assertion follows from Lemma~\ref{simple-Borel}. 
\EndProof

Thus it is sufficient to show that each closed subset of a Polish space
is the continuous image on $\Nat^\infty$. But since a closed subset of a
Polish space is Polish in its own right by Theorem~\ref{Polish-is-G-delta}, we will restrict our attention
to Polish spaces proper.

\BeginProposition{ContimageOfNInfty}
For Polish $X$ there exists a continuous map
$f: \Nat^\infty \rightarrow X$ with $\Bild{f}{\Nat^\infty} = X.$
\EndProposition

\BeginProof
0.
We will define recursively a sequence of closed sets indexed by
elements of $\Nat^*$ that will enable us to define a continuous
map on $\Nat^\infty$.

1.
Let $d$ be a metric that makes $X$ complete.
Represent $X$ as $\bigcup_{n \in \Nat} A_n$ with closed sets $A_n \not= \emptyset$ such
that the diameter $\mathsf{diam}(A_n) < 1$ for each $n \in \Nat$.
Assume that for a word $\alpha \in \Nat^*$ of length $k$
the closed set $A_\alpha \not= \emptyset$
is defined, and write
$
A_\alpha = \bigcup_{n \in \Nat} A_{\alpha n}
$
with closed sets $A_{\alpha n} \not= \emptyset$ such that $\mathsf{diam}(A_{\alpha n}) < 1/(k+1)$
for $n \in \Nat$. This yields for every
$
t = \langle n_1, n_2, \dots\rangle \in \Nat^\infty
$
a sequence of nonempty closed sets  $(A_{n_1n_2 .. n_k})_{k \in \Nat}$ with
diameter $\mathsf{diam}(A_{n_1n_2 .. n_k}) < 1/k$. Because the metric is complete,
$
\bigcap_{k \in \Nat} A_{n_1n_2 .. n_k}
$
contains exactly one point, which is defined to be $f(t)$. This construction renders
$f: \Nat^\infty \rightarrow X$ well defined.

2.
Because we can find for each $x \in X$ an index $n_1' \in \Nat$ with $x \in A_{n_1'}$,
an index $n_2'$ with $x \in A_{n_1'n_2'},$ etc.; the map just defined is onto, so
that $f(\langle n_1', n_2', n_3', \dots\rangle) = x$ for some
$t' := \langle n_1', n_2', n_3', \dots\rangle \in \Nat^\infty$. Suppose $\epsilon > 0$
is given. Since the diameters of the sets $(A_{n_{1}n_{2}\dots n_{k}})_{k\in\Nat}$ tend to $0$, we can find $k_0 \in \Nat$
with $\mathsf{diam}(A_{n_1'n_2' .. n_k'}) < \epsilon$ for all $k > k_0$. Put
$\alpha' := n_1'n_2' .. n_{k_0}'$, then $\Sigma_{\alpha'}$ is an open
neighborhood of $t'$ with
$
\Bild{f}{\Sigma_{\alpha'}} \subseteq B_{\epsilon, d}(f(t')).
$
Thus we find for an arbitrary open neighborhood $V$ of $f(t')$ an open
neighborhood $U$ of $t'$ with $\Bild{f}{U} \subseteq V$, equivalently,
$U \subseteq \InvBild{f}{V}$. Thus $f$
is continuous.
\EndProof

Proposition~\ref{ContimageOfNInfty} permits sometimes the transfer
of arguments pertaining to Polish spaces to arguments using
infinite sequences. Thus a specific space is studied instead of
an abstractly given one, the former permitting some rather special
constructions. This will be capitalized on in the investigation of some
astonishing properties of analytic sets which we will study now.

\Subsection{Analytic Sets and Spaces}
\label{sec:analytic-spaces-1}

We will deal now systematically with analytic sets and spaces. One of
the core results of this section will be the Lusin Separation Theorem,
which permits to separate two disjoint analytic sets through disjoint
Borel sets, and its immediate consequence, the Souslin Theorem, which
says that a set which is both analytic and co-analytic is Borel. These
beautiful results turn out to be very helpful, e.g., in the
investigation of Markov transition systems. In addition, they permit to
state and prove a weak form of Kuratowski's Isomorphism Theorem,
stating that a measurable bijection between two Polish spaces is an
isomorphisms (hence its inverse is measurable as well).

But first the definition of analytic and co-analytic sets for a Polish space $X$.

\BeginDefinition{analytic-and-co-analytic}
An \emph{analytic set}\index{analytic set}\index{set!analytic} in $X$ is
the projection of a Borel subset of $X \times X$. The complement of an analytic set is called a
\emph{co-analytic} set\index{set!co-analytic}.
\EndDefinition

One may wonder whether these projections are Borel sets, but we will
show in a moment that there are strictly more analytic sets than Borel
sets, whenever the underlying Polish space is uncountable. Thus
analytic sets are a proper extension to Borel sets. On the other hand,
analytic sets arise fairly naturally, for example from factoring Polish spaces
through equivalence relations that are generated from a countable
collection of Borel sets. We will see this in
Proposition~\ref{RemainsAnalytic}. Consequently it is sometimes more
adequate to consider analytic sets rather than their Borel cousins,
e.g., when the equivalence of states in a transition system is
at stake.

This is a first characterization of analytic sets (using $\pi_X$ for the projection to $X$).

\BeginProposition{char-analytic}
Let $X$ be a Polish space. Then the following statements are
equivalent for $A \subseteq X$:
\begin{enumerate}
  \item \label{char-analytic-a}$A$ is analytic.
  \item \label{char-analytic-b}There exists a Polish space $Y$ and a
    Borel set $B \subseteq X \times Y$ with $A = \Bild{\pi_X}{B}$.
\item \label{char-analytic-c}There exists a continuous map $f:
  \Nat^\infty \rightarrow X$ with $\Bild{f}{\Nat^\infty} = A$.
  \item \label{char-analytic-d}$A = \Bild{\pi_X}{C}$ for a closed subset $C \subseteq X \times \Nat^\infty$.
\end{enumerate}
\EndProposition

\BeginProof
The implication \emph{\labelImpl{char-analytic-a}{char-analytic-b}} is trivial,
\emph{\labelImpl{char-analytic-b}{char-analytic-c}} follows
from Proposition~\ref{ContimageOfNInfty}: $B = \Bild{g}{\Nat^\infty}$ for some continuous
map $g: \Nat^\infty \rightarrow X \times Y$, so put $f := \pi_X \circ g$. We
obtain \emph{\labelImpl{char-analytic-c}{char-analytic-d}} from the observation that
the graph
$
\{\langle t, f(t)\rangle \mid t \in \Nat^\infty\}
$
of $f$ is a closed subset of $\Nat^\infty \times X$, the first projection of which equals $A$.
Finally, \emph{\labelImpl{char-analytic-d}{char-analytic-a}} is obtained again from
Proposition~\ref{ContimageOfNInfty}.
\EndProof

As an immediate consequence we obtain that a Borel set is
analytic. Just for the record:

\BeginCorollary{BorelIsAnalytic}
Each Borel set in a Polish space is analytic.
\EndCorollary

\BeginProof
Proposition~\ref{char-analytic}
together with Proposition~\ref{ContimageOfNInfty}.
\EndProof

The converse does not hold, as we will show now. This statement is not
only of interest in its own right. Historically it initiated the study
of analytic and co-analytic sets as a separate discipline in set
theory (what is called now \emph{Descriptive Set Theory}).

\BeginProposition{MoreAnalyticThanBorel}
Let $X$ be an uncountable Polish space. Then there exists an analytic
set that is not Borel.
\EndProposition

We show as a preparation for the proof of
Proposition~\ref{MoreAnalyticThanBorel} that analytic sets are closed
under countable unions, intersections, direct and inverse images of
Borel maps. Before doing that, we establish a simple but useful
property of the graphs of measurable maps.

\BeginLemma{GraphIsBorel}
Let $(M, \mathcal{M})$ be a measurable space, $f: M \rightarrow Z$
be a $\mathcal{M}$-$\Borel{Z}$-measurable map, where $Z$ is a separable
metric space. The graph of $f$\index{map!graph}\index{graph!of a map},
\begin{equation*}
\mathsf{graph}(f) := \{\langle m, f(m)\rangle\mid m \in M\},
\end{equation*}
is a member if $\mathcal{M}\otimes\Borel{Z}$.
\EndLemma

\BeginProof
Exercise~\ref{ex-graph-meas}.
\EndProof

Analytic sets have closure properties that are similar to those of
Borel sets, but not quite the same: they are closed under countable
unions and intersections, and under the inverse image of Borel
maps. They are closed under the direct image of Borel maps as well.
Suspiciously missing is the closure under complementation (which
will give rise to Souslin's Theorem). This is different from Borel sets.

\BeginProposition{AnalyticIsStableUnder}
Analytic sets in a Polish space $X$ are closed under countable unions
and countable intersections. If $Y$ is another Polish space, with
analytic sets $A \subseteq X$ and $B \subseteq Y$, and $f: X
\rightarrow Y$ is a Borel map, then $\Bild{f}{A} \subseteq Y$ is
analytic in $Y$, and $\InvBild{f}{B}$ is analytic in $X$.
\EndProposition

\BeginProof
1.  Using the characterization of analytic sets in
Proposition~\ref{char-analytic}, it is shown exactly as in the proof
to Lemma~\ref{ClosedImageIsEnough} that analytic sets are closed under
countable unions and under countable intersections. We trust that the
reader will be able to reproduce those arguments here.

2.  Note first that for $A\subseteq X$ the set $Y \times A$ is analytic in the Polish
space $Y \times X$ by Proposition~\ref{char-analytic}.  In fact, $A = \Bild{\pi_{X}}{B}$ with $B\subseteq X\times X$ Borel by the first part, hence $Y\times A = \Bild{\pi_{Y\times X}}{Y\times B}$ with $Y\times B \subseteq Y\times X \times X$ Borel, which is analytic by the second part. Since $y \in
\Bild{f}{A}$ iff $\langle x, y\rangle \in \mathsf{graph}(f)$ for some
$x \in A$, we write
\begin{equation*}
\Bild{f}{A} = \Bild{\pi_Y}{Y \times A \cap \{\langle y, x\rangle\mid \langle x, y\rangle \in \mathsf{graph}(f)\}}.
\end{equation*}
The set $\{\langle y, x\rangle\mid \langle x, y\rangle \in \mathsf{graph}(f)\}$ is Borel
in $Y \times X$ by Lemma~\ref{GraphIsBorel}, so the assertion follows for the
direct image. The assertion is proved in exactly the same way for
the inverse image.
\EndProof

\BeginProof (of Proposition~\ref{MoreAnalyticThanBorel})
1.
We will deal with the case $X = \Nat^\infty$ first, and apply a diagonal argument.
Let $F \subseteq \Nat^\infty \times (\Nat^\infty \times \Nat^\infty)$ be a universal
closed set according to Proposition~\ref{universal-open}. Thus
each closed set $C \subseteq \Nat^\infty \times \Nat^\infty$ can be
represented as $C = F_t$ for some $t \in \Nat^\infty$.
Taking first projections, we conclude that there exists a universal
analytic set $U \subseteq \Nat^\infty \times \Nat^\infty$ such that
each analytic set $A \subseteq \Nat^\infty$ can be represented as
$U_t$ for some $t \in \Nat^\infty$. In fact, we can write $A = \bigl(\Bild{\pi'_{\NatInf\times\NatInf}}{F}\bigr)_{t}$ with $\pi'_{\NatInf\times\NatInf}$ as the first projection of $(\NatInf\times\NatInf)\times\NatInf$. 

Now set
\begin{equation*}
A := \{\zeta \mid \langle \zeta, \zeta\rangle \in U\}.
\end{equation*}
Because analytic sets are closed under inverse images f Borel maps by
Proposition~\ref{AnalyticIsStableUnder}, $A$ is an analytic
set. Suppose that $A$ is a Borel set, then $\Nat^\infty\setminus A$
is also a Borel set, hence analytic. Thus we find $\xi \in \Nat^\infty$
such that $\Nat^\infty\setminus A = U_\xi$. But now
\begin{equation*}
\xi \in A \Leftrightarrow \langle\xi, \xi\rangle \in U
\Leftrightarrow \xi \in U_\xi \Leftrightarrow \xi \in \Nat^\infty\setminus A.
\end{equation*}
This is a contradiction.

2. The general case is reduced to the one treated above by
observing that an uncountable Polish space contains a homeomorphic
copy on $\Nat^\infty$. But since we are interested mainly in
showing that analytic sets are strictly more general than Borel
sets, we refrain from a very technical discussion of this case and
refer the reader to~\cite[Remark 2.6.5]{Srivastava}.
\EndProof

The representation of an analytic set through a continuous map on
$\Nat^\infty$ has the remarkable consequence that we can separate two
disjoint analytic sets by disjoint Borel sets (Lusin's Separation
Theorem\index{theorem!Lusin}).  This in turn implies a pretty
characterization of Borel sets due to Souslin which says that an
analytic set is Borel iff it is co-analytic as well. Since the latter
characterization will be most valuable to us, we will discuss it in
greater detail now.

We start with Lusin's Separation Theorem.

\BeginTheorem{Lusin}
Given disjoint analytic sets $A$ and $B$ in a Polish space $X$, there
exist disjoint Borel sets $E$ and $F$ with $A \subseteq E$ and $B \subseteq F$.
\EndTheorem

\BeginProof
0.
Call two analytic sets $A$ and $B$ separated by Borel sets iff
$A \subseteq E$ and $B \subseteq F$ for disjoint Borel sets $E$ and $F$.
Observe that if two sequences $\Folge{A}$ and $\Folge{B}$ have the
property that $A_m$ and $B_n$ can be separated by Borel sets for all $m, n \in \Nat$, then
$\bigcup_{n \in \Nat} A_n$ and $\bigcup_{m \in \Nat} B_m$ can also be
separated by Borel sets. In fact, if $E_{m, n}$ and $F_{m, n}$ separate $A_n$
and $B_m$, then
$
E := \bigcap_{m \in \Nat} \bigcup_{n \in \Nat} E_{m, n}
$
and
$
F := \bigcup_{m \in \Nat} \bigcap_{n \in \Nat} F_{m, n}
$
separate $\bigcup_{n \in \Nat} A_n$ and $\bigcup_{m \in \Nat} B_m$.

1.
Now suppose that $A = \Bild{f}{\Nat^\infty}$ and $B = \Bild{g}{\Nat^\infty}$
cannot be separated by Borel sets, where $f, g: \Nat^\infty \rightarrow X$ are continuous and
chosen according to Proposition~\ref{char-analytic}. Because
$
\Nat^\infty = \bigcup_{j \in \Nat} \Sigma_j,
$
($\Sigma_\alpha$ is defined in Example~\ref{Inf-Nat}),
we find indices $k_1$ and $\ell_1$ such that
$\Bild{f}{\Sigma_{j_1}}$ and $\Bild{g}{\Sigma_{\ell_1}}$
cannot be separated by Borel sets. For the same reason, there exist indices
 $k_2$ and $\ell_2$ such that
$\Bild{f}{\Sigma_{j_1j_2}}$
and
$\Bild{g}{\Sigma_{\ell_1\ell_2}}$
cannot be separated by Borel sets. Continuing with this, we define infinite
sequences
$
\kappa := \langle k_1, k_2, \dots\rangle
$
and
$
\lambda := \langle \ell_1, \ell_2, \dots\rangle
$
such that for each $n \in \Nat$ the sets
$\Bild{f}{\Sigma_{j_1  j_2 \dots j_n}}$
and
$\Bild{g}{\Sigma_{\ell_1 \ell_2 \dots \ell_n}}$
cannot be separated by Borel sets. Because $f(\kappa) \in A$ and $g(\lambda) \in B$,
we know $f(\kappa) \not= g(\lambda)$, so we find $\epsilon > 0$ with
$
d(f(\kappa), g(\lambda)) < 2\cdot\epsilon.
$
But we may choose $n$ large enough so that both
$\Bild{f}{\Sigma_{j_1  j_2 \dots j_n}}$
and
$\Bild{g}{\Sigma_{\ell_1 \ell_2 \dots \ell_n}}$
have a diameter smaller than $\epsilon$ each. This is a contradiction
since we now have separated these sets by open balls.
\EndProof

We obtain as a consequence Souslin's Theorem.

\BeginTheorem{Souslin}{\textbf{(Souslin)}} \index{theorem!Souslin}
Let $A$ be an analytic set in a Polish space. If $X\setminus A$ is
analytic, then $A$ is a Borel set.
\EndTheorem

\BeginProof
Let $A$ and $X\setminus A$ be analytic, then they can be separated by
disjoint Borel sets $E$ with $A \subseteq E$ and $F$ with $X\setminus
A \subseteq F$ by Lusin's Theorem~\ref{Lusin}. Thus $A =
E$ is a Borel set.
\EndProof

Souslin's Theorem is important when one wants to show that a
set is a Borel set that is given for example through the image
of another Borel set. A typical scenario for its use is establishing
for a Borel set $A$ and a Borel map
$f: X \rightarrow Y$ that both $C = \Bild{f}{A}$ and
$Y\setminus C = \Bild{f}{X\setminus A}$
hold. Then one  infers from Proposition~\ref{AnalyticIsStableUnder} that
both $C$ and $Y\setminus C$ are analytic, and from Souslin's Theorem
that $A$ is a Borel set. This is a first simple example:

\BeginProposition{invariant-is-Borel}
Let $f: X\to Y$ be surjective and Borel measurable, where $X$ and $Y$
are Polish. Assume that the set $A\in\Borel{X}$ has this
property: $x\in A$ and $f(x) = f(x')$ implies $x'\in A$. Then
$\Bild{f}{A}\in\Borel{Y}$. 
\EndProposition

\BeginProof
Put $C := \Bild{f}{A}$, $D := \Bild{f}{X\setminus A}$, then both $C$
and $D$ are analytic sets by
Proposition~\ref{AnalyticIsStableUnder}. Clearly $Y\setminus
C\subseteq D$. For establishing the other inclusion, let $y\in D$,
hence there exists $x\not\in A$ with $y = f(x)$. But $y\not\in C$, for
otherwise there exists $x'\in A$ with $y = f(x')$, which implies $x\in
A$. Thus $y\in Y\setminus C$. We infer $\Bild{f}{A}\in\Borel{Y}$ now
from Theorem~\ref{Souslin}. 
\EndProof

This yields as an immediate consequence, which will be extended to
analytic spaces in Proposition~\ref{ImageIso} with essentially the
same argument.

\BeginCorollary{isom-polish-bla}
Let $f: X\to Y$ be measurable and bijective with $X$, $Y$
Polish. Then $f$ is a Borel isomorphism. \QED 
\EndCorollary

We state finally \index{theorem!Kuratowski Isomorphism}Kuratowski's Isomorphism Theorem.

\BeginTheorem{kuratowski-isomorphism}
Any two Borel sets of the same cardinality contained in Polish spaces are Borel isomorphic.
\QED
\EndTheorem

The proof requires a reduction to the Cantor ternary set, using the
tools we have discussed here so far. Since giving the proof would lead
us fairly deep into the Wonderland of Descriptive Set Theory, we do
not give it here and refer rather to~\cite[Theorem 3.3.13]{Srivastava},~\cite[Section 15.B]{Kechris}
or~\cite[p. 442]{Kuratowski-Mostowski}.

We make the properties of analytic sets a bit more widely available by
introducing analytic spaces. Roughly, an analytic space is Borel
isomorphic to an analytic set in a Polish space; to be more precise:

\BeginDefinition{AnalyticSpace}
A measurable space $(M, \mathcal{M})$ is called an
\emph{analytic space}\index{space!analytic}
iff there exists a Polish space $X$ and an analytic set $A$
in $X$ such that the measurable spaces $(M, \mathcal{M})$ and $(A, \Borel{X} \cap A)$
are Borel isomorphic. The elements
of  $\mathcal{M}$ are then called the \emph{Borel sets of $M$}.
$\mathcal{M}$ is denoted by $\Borel{M}$.
\EndDefinition

We will omit the $\sigma$-algebra from the notation of an
analytic space.

Analytic spaces share many favorable properties with analytic
sets, and with Polish spaces, but they are a wee bit more general:
whereas an analytic set lives in a Polish space, an analytic space
does only require a Polish space to sit in the background somewhere
and to be Borel isomorphic to it. This makes life considerably
easier, since we are not always obliged to present a Polish space
directly when dealing with properties of analytic spaces.

% An immediate consequence is that the image of an analytic space
% under a Borel map into a Polish space is analytic again.

% \BeginProposition{Remains-Analytic}
% Let $f: X \rightarrow Y$ be a Borel map from the analytic space $X$
% to the Polish space $Y$, then $\Bild{f}{X}$ is an analytic set in $Y$.
% \EndProposition

% \BeginProof
% This is  a mere reformulation from Proposition~\ref{AnalyticIsStableUnder}.
% \EndProof

Take a Borel measurable bijection between two Polish spaces. It is
not a priori clear whether or not this map is an isomorphism. Souslin's
Theorem gives a helpful hand here as well. We will need this property
in a moment for a characterization of countably generated sub-$\sigma$-algebras
of Borel sets, but it appears to be interesting in its own right.

\BeginProposition{ImageIso}
Let $X$ and $Y$ be analytic spaces and $f: X \rightarrow Y$ be
a bijection that is Borel measurable. Then $f$ is a Borel isomorphism.
\EndProposition

\BeginProof
1.
It is no loss of generality to assume that we can find Polish
spaces $P$ and $Q$ such that $X$ and $Y$ are subsets of $P$ resp. $Q$.
We want to show that $\Bild{f}{X \cap B}$ is a Borel set in $Y$,
whenever $B \in \Borel{P}$ is a Borel set. For this we need to
find a Borel set $G \in \Borel{Q}$ such that $\Bild{f}{X \cap B} = G \cap Q$.

2.
Clearly, both $\Bild{f}{X \cap B}$ and   $\Bild{f}{X \setminus B}$ are
analytic sets in $Q$ by Proposition~\ref{AnalyticIsStableUnder}, and because
$f$ is injective, they are disjoint. Thus we can find a Borel set
$G \in \Borel{Q}$ with $\Bild{f}{X \cap B} \subseteq G \cap Y$, and
$\Bild{f}{X \setminus B}\subseteq Q\setminus G \cap Y$. Because $f$
is surjective, we have $\Bild{f}{X \cap B} \cup \Bild{f}{X \setminus B}$,
thus $\Bild{f}{X \cap B} = G \cap Y$
\EndProof

Separable measurable spaces are characterized through
subsets of Polish spaces.

\BeginLemma{Char-separable}
The measurable space $(M, \mathcal{M})$ is separable iff there exists
a Polish space $X$ and a subset $P \subseteq X$ such that
the measurable spaces $(M, \mathcal{M})$ and $(P, \Borel{X}\cap P)$ are
Borel isomorphic.
\EndLemma

\BeginProof
1.
Because $\Borel{X}$ is countably generated for a Polish space $X$ by
Lemma~\ref{DenseWithBorel}, the $\sigma$-algebra $\Borel{X}\cap P$
is countably generated. Since this property is not destroyed by
Borel isomorphisms, the condition above is sufficient.

2.
It is also necessary by Proposition~\ref{sep-is-sep-metr}, because
$(\{0, 1\}^{\Nat}, \bigotimes_{n\in\Nat}\PowerSet{\{0, 1\}})$ is a
Polish space  by Lemma~\ref{Products+Coproducts}. 
% Let $\Folge{A}$ be a generator for $\mathcal{M}$,
% and define
% $$
% f(t) := \langle \chi_{A_1}(t),  \chi_{A_2}(t), \dots\rangle
% $$
% ($\chi_A$ is the indicator function for set $A$), then
% $
% f: M \rightarrow \{0, 1\}^\infty
% $
% is injective, because $\Folge{A}$ separates points. Put
% $
% X := \{0, 1\}^\infty
% $
% and equip $X$ with the product topology, then $X$ is compact, hence
% Polish. Put $P := \Bild{f}{M}$, and let
% $$
% B_n := \{\tau \in P \mid \tau_n = 1\} = P \cap \InvBild{\chi_{A_n}}{\{1\}}.
% $$
% Since $\Borel{X}$ is generated from the sequence
% $\left(\{\tau \in X \mid \tau_n = 1\}\right)_{n \in \Nat}$,
% we infer that $f$ is a Borel isomorphism
% between $(M, \mathcal{M})$ and $(P, \Borel{X}\cap P)$.
\EndProof

Thus analytic spaces are separable.
\BeginCorollary{AnalyticIsSeparable}
An analytic space is a separable measurable space.
\QED
\EndCorollary

% A second consequence is that separable measurable spaces are
% derived from separable metric spaces in a rather straightforward
% way.

% \BeginLemma{SeparableIsMetric}
% For a separable measurable space $(X, \mathcal{A})$ there exists a
% separable metric topology $\tau$ on $X$ such that
% $
% \Borel{\tau} = \mathcal{A}.
% $
% \QED
% \EndLemma

Let us have a brief look at countably generated sub-$\sigma$-algebras
of an analytic space. This will help establishing for example that the
factor space for a particularly interesting and important class of
equivalence relations is an analytic space. The following statement,
which is sometimes referred to as the \emph{\index{theorem!unique
    structure}Unique Structure Theorem}~\cite[Theorem 3.3.5]{Arveson},
says essentially that the Borel sets of an analytic space are uniquely
determined by being countably generated and by separating points. It
comes as a consequence of our discussion of Borel isomorphisms.

\BeginProposition{Unique-Structure}
Let $X$ be an analytic space, $\mathcal{B}_0$ a countably generated sub-$\sigma$-algebra
of $\Borel{X}$ that separates points. Then $\mathcal{B}_0 = \Borel{X}$.
\EndProposition

\BeginProof
1.
$(X, \mathcal{B}_0)$ is a separable measurable space, so there exists a Polish
space $P$ and a subset $Y \subseteq P$ of $P$ such that $(X, \mathcal{B}_0)$
is Borel isomorphic to $P, \Borel{P} \cap Y$ by Lemma~\ref{Char-separable}.
Let $f$ be this isomorphism, then $B_0 = \InvBild{f}{\Borel{P} \cap Y}$.

2.
$f$ is a Borel map from $(X, \Borel{X})$ to $(Y,\Borel{P} \cap Y)$, thus
$Y$ is an analytic set with $\Borel{Y} = \Borel{X} \cap P$ by
Proposition~\ref{RemainsAnalytic}. By Proposition~\ref{AnalyticIsStableUnder},
$f$ is an isomorphism, hence $\Borel{X} =\InvBild{f}{\Borel{P} \cap Y}$.
But this establishes the assertion.
\EndProof

This gives an interesting characterization of measurable spaces to be
analytic, provided they have a separating sequence of sets. Note that
the sequence of sets in the following statement is required to
separate points, but we do not assume that it generates the
$\sigma$-algebra for the underlying space. The statement says that it
does, actually.

\BeginLemma{Sep-Analytic-Map}
Let $X$ be analytic, $f: X \rightarrow Y$ be 
$\Borel{X}$-${\cal B}$-measurable and onto for a measurable space
$(Y, {\cal B})$, which has a sequence of sets in ${\cal B}$ that
separate points. Then $(Y, {\cal B})$ is analytic.
\EndLemma

\BeginProof
1.
The idea is to show that an arbitrary measurable set is contained in the $\sigma$-algebra generated by the sequence in question. Thus let $\Folge{B}$ be the sequence of sets that separates points, take an
arbitrary set $N \in {\cal B}$ and define the $\sigma$-algebra
$
\mathcal{B}_0 := \sigma(\{B_n \mid n \in \Nat\} \cup\{N\}).
$
We want to show that $N\in \sigma(\{B_{n}\mid n\in\Nat\})$, and we show this in a roundabout way by showing that ${\cal B} = \Borel{Y} = {\cal B}_{0}$. Here is, how.

2.
Then $(Y, \mathcal{B}_0)$ is a separable measurable space, so by
Lemma~\ref{Char-separable} we
can find a Polish space $P$ with $Y \subseteq P$ and $\mathcal{B}_0$
as the trace of $\Borel{P}$ on $Y$. Proposition~\ref{AnalyticIsStableUnder}
tells us that $Y = \Bild{f}{X}$ is analytic with
$\mathcal{B}_0 = \Borel{Y}$, and from Proposition~\ref{Unique-Structure} it follows that $\Borel{Y} = \sigma(\{B_n \mid n \in \Nat\})$. Thus $N\in\Borel{Y}$, and since $N\in{\cal B}$ is arbitrary, we conclude $B\subseteq\Borel{Y}$, thus 
${\cal B}\subseteq\Borel{Y}=\sigma(\{B_{n}\mid n\in\Nat\}) \subseteq {\cal B}$. 
\EndProof

We will use Lemma~\ref{Sep-Analytic-Map} for demonstrating that 
factoring an analytic space through a smooth equivalence relation yields an
analytic space again. This class of relations will be defined now and
briefly characterized here. We give a definition in terms of a determining sequence of Borel sets
and relate other characterizations of smoothness in
Lemma~\ref{CharSmooth}.

\BeginDefinition{SmoothRel} Let $X$ be an analytic space and $\rho$
an equivalence relation on $X$. Then $\rho$ is called
\emph{smooth}\index{equivalence relation!smooth} iff there exists a
sequence $\Folge{A}$ of Borel sets such that
\begin{equation*}
\isEquiv{x}{x'}{\rho} \Leftrightarrow \forall n \in \Nat:
[x \in A_n \Leftrightarrow x' \in A_n].
\end{equation*}
$\Folge{A}$ is said to \emph{determine} the relation $\rho$.
\index{equivalence relation!smooth!determining sequence}
\EndDefinition

\BeginExample{modal-yields-smooth}
Given an analytic space $X$, let $M: X\Trans X$ be a transition kernel
which interprets the modal logic presented in
Example~\ref{modal-logic}. Define for a formula $\phi$ and an element
of $x$ the relation $M, x\models \phi$ iff $x\in\Gilt_{M}$, thus $M, x \models \phi$ indicates that 
formula $\phi$ is valid in state $x$. Define the
equivalence relation $\sim$ on $X$ through
\begin{equation*}
  \isEquiv{x}{x'}{\sim} \Longleftrightarrow \forall \phi: [M, x\models \phi \text{ iff } M, x'\models \phi]
\end{equation*}
Thus $x$ and $x'$ cannot be separated through a formula of the
logic. Because the logic has only countably many formulas, the
relation is smooth with the countable set $\{\Gilt_{M}\mid \phi \text{
  is a formula}\}$ as determining relation $\sim$.
\EndExample

We obtain immediately from the
definition that a smooth equivalence relation ---~seen as a subset of
the Cartesian product~--- is a Borel set:

\BeginCorollary{SmoothIsBorel}
Let $\rho$ be a smooth equivalence relation on the analytic space $X$,
then $\rho$ is a Borel subset of $X \times X$.
\EndCorollary

\BeginProof 
Suppose that $\Folge{A}$ determines $\rho$. Since $\isEquiv{x}{x'}{\rho} $ is false iff there exists $n \in \Nat$ with $ \langle
x, x'\rangle \in \left(A_n \times (X\setminus A_n)\right)  \cup \left((X\setminus A_n)
\times A_n\right), $ we obtain
\begin{equation*}
(X \times X) \setminus \rho = \bigcup_{n \in \Nat} \left(A_n \times
(X\setminus A_n)\right)  \cup \left((X\setminus A_n) \times
A_n\right).
\end{equation*}
This is clearly a Borel set in $X \times X$.
\EndProof

The following characterization of smooth equivalence relations is
sometimes helpful and shows that it is not necessary to focus on
sequences of sets. It indicates that the kernels of Borel measurable
maps and smooth relations are intimately related.

\BeginLemma{CharSmooth}
Let $\rho$ be an equivalence relation on an analytic set $X$. Then
these conditions are equivalent:
\begin{enumerate}
 \item \label{char-smmoth-item:1}
  $\rho$ is smooth.
 \item \label{char-smmoth-item:2}
 There exists a sequence $\Folge{f}$ of Borel maps $f_n: X \rightarrow
 Z$ into an analytic space $Z$ such that
 $
 \rho = \bigcap_{n \in \Nat} \Kern{f_n}.
 $
\item \label{char-smmoth-item:3}
 There exists a Borel map $f: X \rightarrow Y$ into an analytic space
 $Y$ with
 $
 \rho = \Kern{f}.
 $
\end{enumerate}
\EndLemma

\BeginProof
\labelImpl{char-smmoth-item:1}{char-smmoth-item:2}:
Let $\Folge{A}$ determine $\rho$, then
\begin{eqnarray*}
  \isEquiv{x}{x'}{\rho}
  & \Leftrightarrow &
  \forall n \in \Nat: [x \in A_n \Leftrightarrow x' \in A_n] \\
  & \Leftrightarrow &
  \forall n \in \Nat: \chi_{A_n}(x) = \chi_{A_n}(x').
\end{eqnarray*}
Thus take $Z = \{0, 1\}$ and $f_n := \chi_{A_n}$.

\labelImpl{char-smmoth-item:2}{char-smmoth-item:3}:
Put
$
Y := Z^\infty.
$
This is an analytic space in the product $\sigma$-algebra, and
\begin{equation*}
f:
\begin{cases}
  X &\rightarrow Y\\
  x &\mapsto (f_n(x))_{n \in \Nat}
\end{cases}
\end{equation*}
is Borel measurable with
$
f(x) = f(x')
$
iff
$
\forall n \in \Nat: f_n(x) = f_n(x').
$

\labelImpl{char-smmoth-item:3}{char-smmoth-item:1}:
Since $Y$ is analytic, it is separable; hence the Borel sets are generated
through a sequence $\Folge{B}$ which separates
points. Put
$
A_n := \InvBild{f}{B_n},
$
then $\Folge{A}$ is a sequence of Borel sets, because the base sets
$B_n$ are Borel in $Y$, and because $f$ is Borel
measurable. We claim that $\Folge{A}$ determines $\rho$:
\begin{eqnarray*}
  f(x) = f(x')
  & \Leftrightarrow &
  \forall n \in \Nat: [f(x) \in B_n \Leftrightarrow f(x') \in B_n]\\
  &&\text{ (since } \Folge{B}\text{ separates points in }Z)\\
  & \Leftrightarrow &
  \forall n \in \Nat: [x \in A_n \Leftrightarrow x' \in A_n].
\end{eqnarray*}
Thus $\langle x, x'\rangle \in\Kern{f}$ is equivalent to being determined by a sequence of measurable sets.
\EndProof

Thus each smooth equivalence relation may be represented as the
kernel of a Borel map, and vice versa.

The interest in analytic spaces comes from the fact that factoring
an analytic space through a smooth equivalence relation will result
in an analytic space again. This requires first and foremost  the
definition of a measurable structure induced by the relation. The
natural choice is the structure imposed by the factor map. The final
$\sigma$-algebra on $\Faktor{X}{\rho}$ with respect to the Borel
sets on $X$ and the natural projection $\fMap{\rho}$ will be chosen;
it is denoted by $\Faktor{\Borel{X}}{\rho}$. Recall that
$\Faktor{\Borel{X}}{\rho}$ is the largest $\sigma$-algebra
$\mathcal{C}$ on $\Faktor{X}{\rho}$ rendering $\fMap{\rho}$ a
$\Borel{X}$-$\mathcal{\mathcal{C}}$-measurable map. Then it turns
out that $ \Borel{\Faktor{X}{\rho}} $ coincides with $
\Faktor{\Borel{X}}{\rho}: $

\BeginProposition{RemainsAnalytic}
Let $X$ be an analytic space, and assume that $\alpha$ is a smooth
equivalence relation on $X$.  Then $\Faktor{X}{\alpha}$ is an analytic
space.
\EndProposition

\BeginProof
In accordance with the characterization of smooth relations
in Lemma~\ref{CharSmooth} we assume that $\alpha$ is
given through a sequence $\Folge{f}$
of measurable maps $f_n: X \rightarrow \Real$.
The factor map is measurable and onto. Put
$
E_{n, r} := \{\Klasse{x}{\alpha} \mid x \in X, f_n(x) < r\},
$
then
$
\mathcal{E} := \{E_{n, r} \mid n \in \Nat, r \in \Rational\}
$
is a countable set of element of the factor $\sigma$-algebra
that separates points. The assertion now follows without difficulties from
Lemma~\ref{Sep-Analytic-Map}.
\EndProof

Let us have a look at invariant sets for an equivalence relation
$\alpha$. 

\BeginDefinition{invarient-meas-set}
Call a subset $A\subseteq X$ \emph{$\alpha$-\index{invariant
    set}\index{equivalence relation!invariant set}invariant} for the
equivalence relation $\alpha$ on $X$ iff $A$ is the union of
$\alpha$-equivalence classes.
\EndDefinition

Thus $A\subseteq X$ is $\alpha$-invariant iff $x\in A$ and
$\isEquiv{x}{x'}{\alpha}$ implies $x'\in A$. Denote by $A^{\nabla} := \bigcup\{\Klasse{x}{\alpha}\mid
x\in A\}$ the smallest $\alpha$-invariant set containing $A$, then we have the representation
$
A^{\nabla} = \Bild{\pi_{2}}{\alpha\cap(X\times A)},
$ 
because $x'\in A^{\nabla}$ iff there exists $x$ with $\langle x', x\rangle\in X\times A$.

An equivalence relation on $X$ is called analytic resp. closed iff it
constitutes an analytic resp. closed subset of the Cartesian product
$X\times X$.

If $X$ is a Polish space, we know that the smooth equivalence relation
$\alpha\subseteq X\times X$ is a Borel subset by
Corollary~\ref{SmoothIsBorel}. We want to show that, conversely, each
closed equivalence relation $\alpha\subseteq X\times X$ is
smooth. This requires the identification of a countable set which
generates the relation, and for this we require the following auxiliary
statement. It may be called separation through  invariant sets.

\BeginLemma{separate-equiv-invariant}
Let $\rho\subseteq X\times X$ be an analytic equivalence relation on the Polish space $X$ with two disjoint analytic sets $A$ and $B$. If $B$ is $\rho$-invariant, then there exists a $\rho$-invariant Borel set $C$ with $A\subseteq C$ and $B\cap C=\emptyset$. 
\EndLemma

\BeginProof
1.  If $D$ is an analytic set, $D^{\nabla}$ is; this follows from the
representation of $D^{\nabla}$ above, and from
Proposition~\ref{AnalyticIsStableUnder}. We construct a sequence
$\Folge{A}$ of invariant analytic sets, and a sequence $\Folge{B}$ of
Borel sets with these properties: $A_{n}\subseteq B_{n}\subseteq
A_{n+1}$, hence $B_{n}$ is sandwiched between consecutive elements
of the first sequence, $A\subseteq A_{1}$, and $B\cap B_{n}=
\emptyset$ for all $n\in \Nat$.

2.
Define $A_{1} := A^{\nabla}$, then $A\subseteq A_{1}$, and $A_{1}$ is $\rho$-invariant. Since $B$ is $\rho$-invariant as well, we conclude $A_{1}\cap B = \emptyset$: if $x\in A_{1}\cap B$, we find $x'\in A$ with $\isEquiv{x}{x'}{\rho}$, hence $x'\in B$, a contradiction. Proceeding inductively, assume that we have already chosen $A_{n}$ and $B_{n}$  with the properties described above, then put $A_{n+1} := B^{\nabla}_{n}$, then $A_{n+1}$ is $\rho$-invariant and analytic, also $A_{n+1}\cap B = \emptyset$ by the argument above. Hence we can find a Borel set $B_{n+1}$ with $A_{n+1}\subseteq B_{n+1}$ and $B_{n+1}\cap B = \emptyset$. 

3.
Now put $C := \bigcup_{n\in\Nat} B_{n}$. Thus $C\in \Borel{X}$ and $C\cap B=\emptyset$, so it remains to show that $C$ is $\rho$-invariant. Let $x\in C$ and $\isEquiv{x}{x'}{\rho}$. Since $x\in B_{n}\subseteq B^{\nabla}_{n}\subseteq B_{n+1}$, we conclude $x'\in B_{n+1}\subseteq C$, and we are done.  
\EndProof

We use this observation now for a closed equivalence relation. Note that the assumption on being analytic in the proof above was made use of in order to establish that the invariant hull of an analytic set is analytic again. 

\BeginProposition{closed-is-smooth}
A closed equivalence relation on a Polish space is smooth.
\EndProposition

\BeginProof
0.
Let $X$ be a Polish space, and $\alpha\subseteq X\times X$ be a closed equivalence relation. We have to find a sequence $\Folge{A}$ of Borel sets which determines $\alpha$. 

1.
Since $X$ is Polish, it has a countable basis ${\cal G}$. Because $\alpha$ is closed, we can write 
\begin{equation*}
  (X\times X)\setminus \alpha = \bigcup \{U_{n}\times U_{m} \mid  U_{n}, U_{m}\in {\cal G}_{0}, U_{n}\cap U_{m} = \emptyset\}
\end{equation*}
for some countable subset ${\cal G}_{0}\subseteq {\cal G}$. Fix $U_{n}$ and $U_{m}$, then also $U^{\nabla}_{n}$ and $U_{m}$ are disjoint. Select the invariant Borel set $A_{n}$ such that $U_{n}\subseteq A_{n}$ and $A_{n}\cap U_{m}=\emptyset$; this is possible by Lemma~\ref{separate-equiv-invariant}. 

2.
We claim that 
\begin{equation*}
  (X\times X)\setminus \alpha = \bigcup_{n\in\Nat}(A_{n}\times (X\setminus A_{n}).
\end{equation*}
In fact, if $\langle x, x'\rangle\not\in \alpha$, select $U_{n}$ and $U_{m}$ with $\langle x, x'\rangle \in U_{n}\times U_{m}\subseteq A_{n}\times (X\setminus A_{n})$. If, conversely, $\langle x, x'\rangle \in A_{n}\times (X\setminus A_{n})$, then $\langle x, x'\rangle \in \alpha$ implies by the invariance of $A_{n}$ that $x'\in A_{n}$, a contradiction. 
\EndProof

The Blackwell-Mackey-Theorem analyzes those Borel sets that are unions of $\mathcal{A}$-atoms for
a sub-$\sigma$-algebra $\mathcal{A} \subseteq \Borel{X}$. If $\mathcal{A}$ is countably generated by, say,
$\Folge{A}$, then it is not difficult to see that an atom in $\mathcal{A}$
can be represented as
\begin{equation*}
\bigcap_{i \in T} A_i \cap \bigcap_{i \in \Nat\setminus T} (X\setminus A_i)
\end{equation*}
for a suitable subset $T \subseteq \Nat$, see Proposition~\ref{probs-char-atom}. It constructs a measurable map $f$ so that the set under consideration is $\Kern{f}$-invariant, which will be helpful in an the application of the Souslin Theorem. But let's see.

\BeginTheorem{Blackwell-Mackey}{\textbf{(Blackwell-Mackey)}}
\index{theorem!Blackwell-Mackey} Let $X$ be an analytic space and
$\mathcal{A} \subseteq \Borel{X}$ be a countably generated
sub-$\sigma$-algebra of the Borel sets of $X$. If $B \subseteq X$ is
a Borel set that is a union of atoms of $\mathcal{A}$, then $B \in
\mathcal{A}$.
\EndTheorem

The idea of the proof is to show that $\Bild{f}{B}$ and $\Bild{f}{X\setminus B}$ are disjoint analytic sets for the measurable map $f$, and to conclude that $B=\InvBild{f}{C}$ for some Borel set $C$, which will be supplied to us through Souslin's Theorem.

\BeginProof
Let $\mathcal{A}$ be generated by $\Folge{A}$, and define
\begin{equation*}
f: X \rightarrow \{0, 1\}^\infty
\end{equation*}
through
\begin{equation*}
x \mapsto \langle\chi_{A_1}(x), \chi_{A_2}(x), \chi_{A_3}(x), \dots\rangle.
\end{equation*}
Then $f$ is $\mathcal{A}$-$\Borel{\{0, 1\}^\infty}$-measurable. We claim that
$\Bild{f}{B}$ and $\Bild{f}{X\setminus B}$ are disjoint. Suppose not, then
we find $t \in \{0, 1\}^\infty$ with $t = f(x) = f(x')$ for some
$x \in B, x' \in X\setminus B$.
Because $B$ is the union of atoms, we find a subset $T \subseteq \Nat$ with
$x \in A_n$, provided $n \in T$, and $x \notin A_n$, provided $n \notin T$.
But since $f(x) = f(x')$, the same  holds for $x'$ as well, which means
that $x' \in B$, contradicting the choice of $x'$.

Because $\Bild{f}{B}$ and $\Bild{f}{X\setminus B}$ are disjoint analytic sets,
we find through Souslin's Theorem~\ref{Souslin} a Borel set $C$ with
\begin{equation*}
  \Bild{f}{B} \subseteq C ,
  \Bild{f}{X\setminus B} \cap C = \emptyset.
\end{equation*}
Thus $\Bild{f}{B} = C$, so that $\InvBild{f}{\Bild{f}{B}} = \InvBild{f}{C}\in{\cal A}$. We show that $\InvBild{f}{\Bild{f}{B}} = B$. It is clear that $B\subseteq\InvBild{f}{\Bild{f}{B}}$, so assume that $f(b)\in \Bild{f}{B}$, so $f(b) = f(b')$ for some $b'\in B$. By construction, this means $b\in B$, since $B$ is an union of atoms, hence $\InvBild{f}{\Bild{f}{B}}\subseteq B$. Consequently, $B = \InvBild{f}{C}\in {\cal A}$. 
\EndProof

When investigating modal logics, one wants to be able to identify the
$\sigma$-algebra which is defined by the validity sets of the
formulas. This can be done through the Blackwell-Mackey-Theorem and is
formulated for generals smooth equivalence relations. 

\BeginProposition{seq-gen-factor-space}
Let $\rho$ be a smooth equivalence relation on the Polish space $X$,
and assume that $\Folge{A}$
generates $\rho$. Then
\begin{enumerate}
\item $\sigma(\{A_{n}\mid n\in\Nat\})$ is the $\sigma$-algebra of
  $\rho$-invariant Borel sets,
\item $\Borel{\Faktor{X}{\rho}} =
\sigma(\{\Bild{\fMap{\rho}}{A_{n}}\mid n\in\Nat\}$.
\end{enumerate}
\EndProposition

\BeginProof
1.
Denote by ${\cal I}$ be the $\sigma$-algebra of $\rho$-invariant Borel sets; we
have to show that ${\cal I} = \sigma(\{A_{n}\mid n\in\Nat\}$. 
\begin{description}
\item[``$\supseteq$''] Each $A_{n}$ is a $\rho$-invariant Borel set.
\item[``$\subseteq$''] Let $B$ be an $\rho$-invariant Borel set, then
  $B = \bigcup_{b\in B} \Klasse{b}{\rho}$. Each class
  $\Klasse{b}{\rho}$ can be written as
  \begin{equation*}
    \Klasse{b}{\rho} = \bigcap_{b\in A_{n}}A_{n}\cap\bigcap_{b\not\in
      A_{n}}(X\setminus A_{n}),
  \end{equation*}
thus $\Klasse{b}{\rho}\in\sigma(\{A_{n}\mid n\in\Nat\}$. Moreover, it
is easy to see that the classes are the atoms of this $\sigma$-algebra
(in fact, we cannot find a proper non-empty $\rho$-invariant subset of an
equivalence class). Thus the Blackwell-Mackey
Theorem~\ref{Blackwell-Mackey} shows that $B\in\sigma(\{A_{n}\mid
n\in\Nat\})$. 
\end{description}

2.
Now let ${\cal E} := \sigma(\{\Bild{\fMap{\rho}}{A_{n}}\mid
n\in\Nat\}$, and let $g: \Faktor{X}{\rho}\to P$ be ${\cal E}$-${\cal
  P}$-measurable for an arbitrary measurable space $(P, {\cal
  P})$. Thus we have for all $C\in{\cal P}$
\begin{align*}
  \InvBild{g}{C}\in{\cal E} & \Leftrightarrow
\InvBild{\fMap{\rho}}{\InvBild{g}{C}}\in\sigma(\{A_{n}\mid
n\in\Nat\})&&\text{ since }A_{n}=\InvBild{\fMap{\rho}}{\Bild{\fMap{\rho}}{A_{n}}}\\
& \Leftrightarrow 
\InvBild{\fMap{\rho}}{\InvBild{g}{C}}\in{\cal I} && \text{ part 1.}\\
& \Leftrightarrow 
\InvBild{\fMap{\rho}}{\InvBild{g}{C}}\in\Borel{X}
\end{align*}
Thus ${\cal E}$ is the final $\sigma$-algebra with respect to
$\fMap{\rho}$, hence equals $\Borel{\Faktor{X}{\rho}}$. 
\EndProof

The following example shows that the equivalence relation generated by
a $\sigma$-algebra need not return the $\sigma$-algebra as its
invariant sets, if the given $\sigma$-algebra is not countably
generated. Proposition~\ref{seq-gen-factor-space} assures us that this
cannot happen in the countably generated case. 

\BeginExample{count-cocount-equiv-rel}
Let $\mathcal{C}$ be the countable-cocountable $\sigma$-algebra on
$\Real$. The equivalence relation $\equiv_{\mathcal{C}}$
generated by $\mathcal{C}$ according to Example~\ref{prob-equiv-gen}
is the identity. Hence is it smooth. The $\sigma$-algebra of
$\equiv_{\mathcal{C}}$-invariant Borel sets equals the Borel set
$\Borel{\Real}$, which is a proper superset of $\mathcal{C}$. 
\EndExample
 
The next example is a somewhat surprising application of the Blackwell-Mackey
Theorem, taken from~\cite[Proposition~57]{Rao-Rao-Borel}. It shows
that the set of countably generated $\sigma$-algebras is
not closed under finite intersections, hence fails to be a lattice under inclusion.

\BeginExample{inters-non-cg}
There exist two countably generated $\sigma$-algebras the intersection
of which is not countably generated. In fact, let $A\subseteq [0, 1]$
be an analytic set which is not Borel, then $\Borel{A}$ is countably
generated by Corollary~\ref{AnalyticIsSeparable}. Let $f: [0, 1]\to A$
be a bijection, and consider ${\cal C} := \InvBild{f}{\Borel{A}}$,
which is countably generated as well. Then ${\cal D} := \Borel{[0,
  1]}\cap{\cal C}$ is a $\sigma$-algebra which has all singletons in
$[0, 1]$ as atoms. Assume that ${\cal D}$ is countably generated, then
${\cal D} = \Borel{[0, 1]}$ by the
Blackwell-Mackey-Theorem~\ref{Blackwell-Mackey}. But this means that
${\cal C} = \Borel{[0, 1]}$, so that $f: [0, 1]\to A$ is a Borel
isomorphism, hence $A$ is a Borel set in $[0, 1]$, contradicting the
assumption.
\EndExample

Among the consequences of Example~\ref{inters-non-cg} is the
observation that the set of smooth equivalence relations of a Polish
space does not form a lattice under inclusion, but is usually only a
$\cap$-semilattice, as the following example shows. Another
consequence is mentioned in Exercise~\ref{non-pushout}.

\BeginExample{inters-non-smooth}
The intersection $\alpha_{1}\cap\alpha_{2}$ of two smooth equivalence
relations $\alpha_{1}$ and $\alpha_{2}$ is smooth again: if
$\alpha_{i}$ is generated by the Borel sets $\{A_{i, n}\mid
n\in\Nat\}$ for $i = 1, 2$, then $\alpha_{1}\cap\alpha_{2}$ is
generated by the Borel sets $\{A_{i, n}\mid i = 1, 2, n\in\Nat\}$. But
now take two countably generated $\sigma$-algebras ${\cal A}_{i}$, and
let $\alpha_{i}$ be the equivalence relations determined by them, see
Example~\ref{prob-equiv-gen}. Then the $\sigma$-algebra
$\alpha_{1}\cup\alpha_{2}$ is generated by ${\cal A}_{1}\cap{\cal
  A}_{2}$, which is by assumption not countably generated. Hence
$\alpha_{1}\cup\alpha_{2}$ is not smooth.
\EndExample

Sometimes one starts not
with a topological space and its Borel sets but rather with a
measurable space: A \emph{standard Borel}\index{space!standard Borel}
space $(X, \mathcal{A})$ is a measurable space such that the
$\sigma$-algebra $\mathcal{A}$ equals $\Borel{\tau}$ for some Polish
topology $\tau$ on $X$. We will not dwell on this distinction. 

\Subsection{The Souslin Operation}
\label{sec:souslin-op}

The collection of analytic sets is closed under Souslin's operation
$\sA$, which we will introduce now. This operation is not only closed
to analytic sets, we will also see that complete measure spaces are
another important class of measurable spaces which are closed under
this operation. Each measurable space can be completed with respect to
its finite measures, so that we do not even need a topology for
carrying out the constructions ahead.

Let $\Nats$ be the set of all finite and non-empty sequences of natural numbers. Denote for  $t = \Folge{x}\in\NatInf$ by $t|k = \langle x_{1}, \dots, x_{k}\rangle$ its first $k$ elements. Given a subset ${\cal C}\subseteq\PowerSet{X}$, denote by 
\begin{equation*}
  \sA({\cal C}) := \{\bigcup_{t\in\NatInf}\bigcap_{k\in\Nat}A_{t|k}\mid A_{v}\in {\cal C}\text{ for all } v\in\Nats\}
\end{equation*}
Note that the outer union may be taken of more than countably many sets. A family $(A_{v})_{v\in\Nats}$ is called a \emph{\index{Souslin scheme}Souslin scheme}, which is called \emph{\index{Souslin scheme!regular}regular} if $A_{w}\subseteq A_{v}$ whenever $v$ is an initial piece of $w$. Because
\begin{equation*}
  \bigcup_{t\in\NatInf}\bigcap_{k\in\Nat}A_{t|k} = \bigcup_{t\in\NatInf}\bigcap_{k\in\Nat}\bigl(\bigcap_{1\leq j \leq k}A_{t|j}\bigr),
\end{equation*}
we can and will restrict our attention to regular Souslin schemes whenever ${\cal C}$ is closed under finite intersections.

We will see now that each analytic set can be represented through a Souslin scheme with a special shape. This has some interesting consequences, among others  that analytic sets are closed under the Souslin operation. 

\BeginProposition{analytic-is-souslin-repr}
Let $X$ be a Polish space and $(A_{v})_{v\in\Nats}$ be a regular Souslin scheme of closed sets such that $diam(A_{v})\to 0$, as the length of $v$ goes to infinity. Then $E := \bigcup_{t\in\NatInf}\bigcap_{k\in\Nat}A_{t|k}$ is an analytic set in $X$. Conversely, each analytic set can be represented in this way. 
\EndProposition

\BeginProof
1.  Assume $E$ is given through a Souslin scheme, then we represent $E
= \Bild{f}{F}$ with $F\subseteq\NatInf$ a closed set and $f: F\to X$
continuous. In fact, put
\begin{equation*}
  F := \{t\in\NatInf \mid A_{t|k} \not=\emptyset\text{ for all }k\}.
\end{equation*}
Then $F$ is a closed subset of $\NatInf$: take $s\in\NatInf\setminus F$, then we can find $k'\in\Nat$ with $A_{s|k'}=\emptyset$, so that $G := \{t\in\NatInf \mid t|k' = s|k'\}$ is open in $\NatInf$, contains $s$ and is disjoint to $F$. Now let $t\in F$, then there exists exactly one point $f(t)\in\bigcap_{k\in\Nat}A_{t|k}$, since $X$ is complete and the diameters of the sets involved tend to zero. Then $E = \Bild{f}{F}$ by construction, and $f$ is continuous.

Let $t\in F$ and  $\epsilon>0$ be given, take $x := f(t)$ and let $B$ be the ball with center $x$ and radius $\epsilon$. Then we can find an index $k$ such that $A_{t|k'}\subseteq S$ for all $k'\geq k$, hence $U := \{s\in F\mid t|k = s|k\}$ is an open neighborhood of $t$ with $\Bild{f}{U}\subseteq B$.

2.
Let $E$ be an analytic set, then $E = \Bild{f}{\NatInf}$ with $f$ continuous by Proposition~\ref{char-analytic}. Define $A_{v}$ as the closure of the set 
$
\Bild{f}{\{t\in\NatInf \mid t|k = v\}},
$
if the length of $v$ is $k$. Then clearly 
\begin{equation*}
  E = \bigcup_{t\in\NatInf}\bigcap_{k\in\Nat}A_{t|k},
\end{equation*}
since $f$ is continuous. It is also clear that $(A_{v})_{v\in\Nats}$ is regular with diameter tending to zero. 
\EndProof

Before we can enter into the demonstration that the Souslin operation is idempotent, we need some auxiliary statements. The first one is readily verified.

\BeginLemma{b-is-bijection}
$b(m, n) := 2^{m-1}(2n-1)$ defines a bijective map $\Nat\times\Nat\to \Nat$. Moreover, $m\leq b(m, n)$ and $n<n'$ implies $b(m, n)<b(m, n')$ for all $n, n', m\in\Nat$.
\QED
\EndLemma

Given $k\in\Nat$, there exists a unique pair $\langle\ell(k), r(k)\rangle\in\Nat\times\Nat$ with $b(\ell(k), r(k)) = k$. We will need the functions $\ell, r: \Nat\to \Nat$ later on. The next function is considerably more complicated, since it caters for a more involved set of parameters.

\BeginLemma{B-is-a-bijection}
Define for and $z=\Folge{z}\in(\Nat^{\Nat})^{\Nat}$ with $z_{n} = (z_{n, m})_{m\in\Nat}$ and  $t\in\Nat^{\Nat}$ 
\begin{equation*}
  B(t, z)_{k} := b(t(k), z_{\ell(k), r(k)}).
\end{equation*}
Then $B: \Nat^{\Nat}\times(\Nat^{\Nat})^{\Nat}\to \Nat^{\Nat}$ is a bijection.
\EndLemma

\BeginProof
1.
We show first that $B$ is injective. Let $\langle t, z\rangle \not= \langle t', z'\rangle$. If $t\not=t'$, we find $k$ with $t(k)\not=t'(k)$, so that $ b(t(k), z_{\ell(k), r(k)}) \not= b(t'(k), z'_{\ell(k), r(k)})$ follows, because $b$ is injective. Now assume that $t=t'$, but $z\not=z'$, so we can find $i, j\in\Nat$ with $z_{i, j}\not=z'_{i, j}$. Let $k := b(i, j)$, so that $\ell(k) = i$ and $r(k) = j$, hence $\langle t(k), z_{\ell(k), r(k)}\rangle \not= \langle t(k), z'_{\ell(k), r(k)}\rangle$, so that $B(t, z)_{k} \not= B(t', z')_{k}$.

2.
Now let $s\in\Nat^{\Nat}$, and define $t\in\Nat^{\Nat}$ and $z\in(\Nat^{\Nat})^{\Nat}$
\begin{align*}
 t_{k} & := \ell(s_{k}),\\
  z_{n, m} & := r(s_{b(n, m)}). 
\end{align*}
Then we have for $k\in\Nat$
\begin{align*}
  B(t, z)_{k} & = b(t_{k}, z_{\ell(k), r(k)})\\
& = b\bigl(\ell(s_{k}), r(s_{b(\ell(k), r(k))})\bigr)\\
& = b\bigl(\ell(s_{k}), r(s_{k})\bigr)\\
& = s_{k}.
\end{align*}
\EndProof

We construct maps $\phi, \psi$ from $b$ and $B$ now with special properties which will be utilized in the proof that the Souslin operation is idempotent.

\BeginLemma{construct-phi-and-psi}
There exist maps $\phi, \psi: \Nats\to \Nats$ with this property: let $w = B(t, z)|b(n, m)$, then $\phi(w) = t|m$ and $\psi(w) = z_{m}|n$. 
\EndLemma

\BeginProof
Fix $v = \langle x_{1}, \dots, x_{k}\rangle$, then define for $m := \ell(k)$ and $n := r(k)$
\begin{align*}
  \phi(v) & := \langle \ell(x_{1}), \dots, \ell(x_{m})\rangle,\\
\psi(v) & := \langle r(x_{b(m, 1}), \dots, r(x_{b(m, n)})\rangle
\end{align*}
We see from Lemma~\ref{b-is-bijection} that these definitions are possible. 

Given $t\in\Nat^{\Nat}$ and $z\in(\Nat^{\Nat})^{\Nat}$, we put $k := b(m, n)$ and $v := B(t, z)|k$, then we obtain from the definition of $\phi$ resp. $\psi$
\begin{equation*}
  \phi(v)  = \langle \ell(v_{1}), \dots, \ell(v_{m})\rangle 
= t|m
\end{equation*}
and
\begin{equation*}
  \psi(v) = \langle r(v_{b(m, 0)}), \dots, r(v_{b(m, n)})\rangle
= z_{m}|n
\end{equation*}
\EndProof

The construction shows that $\sA({\cal C})$ is always closed under
countable unions and countable intersections. We are now in a position
to prove a much more general observation.

\BeginTheorem{souslin-is-idempotent}
$\sA(\sA({\cal C})) = \sA({\cal C})$. 
\EndTheorem

\BeginProof
It is clear that ${\cal C}\subseteq\sA({\cal C})$, so we have to establish the other inclusion. Let 
$
\{B_{v, w}\mid w\in\Nats\}
$
be a Souslin scheme for each $v\in\Nats$, and put 
$
A_{v} := \bigcup_{s\in\NatInf}\bigcap_{m\in\Nat}B_{v, s|m}.
$ 
Then we have
\begin{align*}
A & :=  \bigcup_{t\in\NatInf}\bigcap_{k\in\Nat}A_{t|k}\\
& = \bigcup_{t\in\NatInf}\bigcap_{k\in\Nat}\bigcup_{s\in\NatInf}\bigcap_{m\in\Nat}B_{v, s|m}\\
& = \bigcup_{t\in\NatInf}\bigcup_{\Folge{z}\in(\NatInf)^{\Nat}}\bigcap_{m\in\Nat}\bigcap_{k\in\Nat}B_{t|m, z_{m}|n}\\
& \stackrel{(\ast)}{=}\bigcup_{s\in\NatInf}\bigcap_{k\in\Nat}C_{s|k}
\end{align*}
with 
\begin{equation*}
  C_{v} := B_{\phi(v), \psi(v)}
\end{equation*}
for $v\in\Nats$. So we have to establish the equality marked $(\ast)$. 

``$\subseteq$'': Given $x\in A$, there exists $t\in\NatInf$ and $z\in(\NatInf)^{\Nat}$ such that $x\in B_{t|m, z_{m|n}}$. Put $s := B(t, z)$. Let $k\in\Nat$ be arbitrary, then there exists a pair $\langle m, n\rangle\in\Nat\times\Nat$ with $k = b(m, n)$ by Lemma~\ref{b-is-bijection}. Thus we have 
$
t|m = \phi(s|k)
$
and
$
z_{m}|n = \psi(s|k).
$
by Lemma~\ref{construct-phi-and-psi}, from which $x \in B_{t|m, z_{m}|n} = C_{s|k}$ follows.

``$\supseteq$'': Let $s\in\NatInf$ such that $x\in C_{s|k}$ for all $k\in\Nat$. We can find by Lemma~\ref{B-is-a-bijection} some $t\in\NatInf$ and $z\in(\Nat^{\Nat})^{\Nat}$ with $B(t, z) = s$. Given $k$, there exist $m, n\in\Nat$ with $k = b(m, n)$, hence $C_{s|k} = B_{t|m, z_{m}|n}$. Thus $x\in A$. 
\EndProof

We obtain as an immediate consequence that analytic sets in a Polish
space $X$ are closed under the Souslin operation. This is so because
we have seen that the collection of analytic sets is contained in
$\sA\bigl(\{F\subseteq X\mid F\text{ is closed}\}\bigr)$, so an application of
Theorem~\ref{souslin-is-idempotent} proves the claim. But we can say
even more.

\BeginProposition{compl-closed-souslin}
Assume that the complement of each set in ${\cal C}$ belongs to $\sA({\cal C})$, and $\emptyset\in{\cal C}$. Then $\sigma({\cal C})\subseteq\sA({\cal C})$. 
\EndProposition

\BeginProof
Define 
\begin{equation*}
  {\cal G} := \{A\in\sA({\cal C})\mid  X\setminus A\in\sA({\cal C})\}.
\end{equation*}
Then ${\cal G}$ is closed under complementation. If $\Folge{A}$ is a sequence in ${\cal G}$, then $\bigcap_{n\in\Nat}A_{n}\in{\cal G}$, because $\sA({\cal C})$ is closed under countable unions. Similarly, $\bigcup_{n\in\Nat}A_{n}\in{\cal G}$. Since $\emptyset\in{\cal G}$, we may conclude that ${\cal G}$ is a $\sigma$-algebra, which contains ${\cal C}$ by assumption. Hence $\sigma({\cal C})\subseteq\sigma({\cal G}) = {\cal G}\subseteq\sA({\cal C})$. 
\EndProof

With complete measure spaces we will meet an important class of
measurable spaces, which is closed under the Souslin operation. As a
preparation for this we state and prove an interesting criterion for
being closed. This requires the definition of a particular kind of
cover.

\BeginDefinition{cover-by-alg}
Given a measurable space $(X, {\cal A})$ and a subset $A\subseteq X$, we call $A_{z}\in{\cal A}$ an \emph{\index{${\cal A}$-cover}${\cal A}$-cover} of $A$ iff
\begin{enumerate}
\item $A\subseteq A_{z}$.
\item For every $B\in{\cal A}$ with $A\subseteq B$, $\PowerSet{A_{z}\setminus B}\subseteq {\cal A}$.
\end{enumerate}
\EndDefinition

Thus $A_{z}\in{\cal A}$ covers $A$ in the sense that $A\subseteq
A_{z}$, and if we have another set $B$ is ${\cal A}$ which covers $A$
as well, then \emph{all} the sets which make out the difference
between $A_{z}$ and $B$ are measurable. In addition it follows  that if $A\subseteq A'\subseteq A_{z}$ and $A'\in{\cal A}$, then $A'$ is also an ${\cal A}$-cover. This concept sounds fairly
artificial and somewhat far fetched, but we will see that arises in a
natural way when completing measure spaces. The surprising observation
is that a space is closed under the Souslin operation whenever each
subset has an ${\cal A}$-cover.

\BeginProposition{Marcewski-complete}
Let $(X, {\cal A})$ be a measurable space such that each subset of $X$
has an ${\cal A}$ cover. Then $(X, {\cal A})$ is closed under
the Souslin operation.
\EndProposition

\BeginProof
1.
Let 
\begin{equation*}
  A := \bigcup_{t\in\NatInf}\bigcap_{k\in\Nat}A_{a|k}
\end{equation*}
with $(A_{v})_{v\in\Nats}$ a regular Souslin scheme in ${\cal A}$. Define 
\begin{equation*}
  B_{w} := \bigcup\{\bigcap_{n\in\Nat}A_{t|n}\mid t\in\NatInf, w\text{ is a prefix of } t\}.
\end{equation*}
for $w\in\Nat^{*} = \Nats\cup\{\epsilon\}$. Then $B_{\epsilon} = A$, $B_{w} = \bigcup_{n\in\Nat} B_{wn}$, and $B_{w}\subseteq A_{w}$ if $w\not=\epsilon$.  

By assumption, there exists a minimal ${\cal A}$-cover $C_{w}$ for
$B_{w}$. We may and do assume that $C_{w}\subseteq A_{w}$, and that
$(C_{w})_{w\in\Nat^{*}}$ is regular (we otherwise force this condition
by considering the ${\cal A}$-cover $ \bigl(\bigcap_{v \text{ prefix
    of } w}(C_{v}\cap A_{v})\bigr)_{w\in\Nat^{*}} $ instead). Now put
$D_{w} := C_{w}\setminus\bigcup_{n\in\Nat}C_{wn}$ for
$w\in\Nat^{*}$. We obtain from this construction 
$B_{w}\subseteq C_{w} = \bigcup_{n\in\Nat}C_{wn}\in{\cal A}$, hence
that every subset of $D_{w}$ is in ${\cal A}$, since $C_{w}$ is an
${\cal A}$-cover. Thus every subset of $D :=
\bigcup_{w\in\Nat^{*}}D_{w}$ is in ${\cal A}$.

2.
We claim that $C_{\epsilon}\setminus D\subseteq A$. In fact, let $x\in C_{\epsilon}\setminus D$, then $x\not\in D_{\epsilon}$, so we can find $k_{1}\in\Nat$ with $x\in C_{k_{1}}$, but $x\not\in D_{n_{1}}$. Since $x\not\in D_{k_{1}}$, we find $k_{2}$ with $x\in C_{k_{1},k_{2}}$ such that $x\not\in D_{k_{1}, k_{2}}$. So we inductively define a sequence $t := \Folge{k}$ so that $x\in C_{t|k}$ for all $k\in\Nat$. Because $C_{t|k}\subseteq A_{t|k}$, we conclude that $x\in A$.

3.
Hence we obtain $C_{\epsilon}\setminus A\subseteq D$, and since every subset of $D$ is in ${\cal A}$, we conclude that $C_{\epsilon}\setminus A\in{\cal A}$, which means that $A = C_{\epsilon}\setminus (C_{\epsilon}\setminus A)\in {\cal A}$. 
\EndProof

%%% Local Variables: 
%%% mode: latex
%%% TeX-master: "../Mskr3"
%%% End: 

%\Input{\Folder/UnivMeas}
%spell checked - 24Aug14
\def\cCompl#1#2{\ensuremath{\overline{#1}^{#2}}}
\def\cU#1{\cCompl{#1}{}}

\Subsection{Universally Measurable Sets}
\label{sec:uni-meas}

After this technical preparation we are posed to enter the interesting
world of universally measurable sets with the closure operations that
are associated with them. We define complete measure spaces and show
that an arbitrary ($\sigma$-) finite measure space can be completed,
uniquely extending the measure as we go. This leads also to
completions with respect to families of finite measures, and we show
that the resulting measurable spaces are closed under the Souslin
operation. Two applications are discussed. The first one demonstrates
that a measure defined on a countably generated sub-$\sigma$-algebra
of the Borel sets of an analytic space can be extended to the Borel
sets, albeit not necessarily in a unique way. This result due to Lubin
rests on the important von Neumann Selection Theorem, giving a
universally right inverse to a measurable map from an analytic to a
separable space. Another application of von Neumann's result is the
observation that under suitable topological assumptions for a
surjective map $f$ the lifted map $\FinM{f}$ is surjective as
well. The second application shows that a transition kernel can be
extended to the universal closures of the measurable spaces involved,
provided the target space is separable.

A $\sigma$-finite
measure space $(X, {\cal A}, \mu)$ is called
\emph{\index{complete!measure space}complete} iff $\mu(A) = 0$ with
$A\in{\cal A}$ and $B\subseteq A$ implies $B\in{\cal A}$. Thus if we
have two sets $A, A'\in{\cal A}$ with $A\subseteq A'$ and $\mu(A) = \mu(A')$, then we know
that each set which can be sandwiched between the two will be
measurable as well. We will discuss the completion of a measure space
and investigate some properties. We first note that it is sufficient
to discuss finite measure spaces; in fact, assume that we have a
collection of mutually disjoint sets $\Folge{G}$ with $G_{n}\in{\cal
  A}$ such that $0 < \mu(G_{n})<\infty$ and
$\bigcup_{n\in\Nat}G_{n}=X$, consider the measure
\begin{equation*}
\mu'(B) := \sum_{n\in\Nat}\frac{\mu(B\cap G_{n})}{2^{n}\mu(G_{n})},  
\end{equation*}
then $\mu$ is complete iff $\mu'$ is complete, and $\mu'$ is a probability measure.

We fix for the time being a finite measure $\mu$ on a measurable space $(X, {\cal A})$. The outer measure $\mu^{*}$ is defined through 
\begin{align*}
  \mu^{*}(C) & := \inf\{\sum_{n\in\Nat}\mu(A_{n}) \mid C \subseteq \bigcup_{n\in\Nat}A_{n}, A_{n}\in {\cal A}\text{ for all }n\in\Nat\}\\
& = \inf \{\mu(A) \mid C \subseteq A, A\in{\cal A}\}
\end{align*}
for any subset $C$ of $X$.

\BeginDefinition{the-sigma-ideal}
Call $N\subseteq X$ a \emph{\index{$\mu$-null set}$\mu$-null set} iff $\mu^{*}(N) = 0$. Define ${\cal N}_{\mu}$ as the set of all $\mu$-null sets. 
\EndDefinition

Because $\mu^*$ is countable subadditive, we obtain

\BeginLemma{null-sets-are-sigma-ideal}
${\cal N}_{\mu}$ is a $\sigma$-ideal.
\QED
\EndLemma

Now assume that we have sets $A, A'\in{\cal A}$ and $N, N'\in{\cal N}_{\mu}$ such that their symmetric differences  with $A\Delta N = A'\Delta N'$, so their symmetric differences are the same. Then we may infer $\mu(A) = \mu(A')$, because $A\Delta A' = A\Delta\bigl(A\Delta(N\Delta N')\bigr) = N\Delta N'\subseteq N\cup N'\in{\cal N}_{\mu}$, and $|\mu(A) - \mu(A')| \leq \mu(A\Delta A')$. Thus we may construct an extension of $\mu$ to the $\sigma$-algebra generated by ${\cal A}$ and ${\cal N}_{\mu}$ in an obvious way.

\BeginProposition{completion-of-measure-space}
Define ${\cal A}_{\mu} := \sigma({\cal A}\cup{\cal N}_{\mu})$ and $\cU{\mu}(A\Delta N) := \mu(A)$ for $A\in{\cal A}, N\in{\cal N}_{\mu}$. Then 
\begin{enumerate}
\item ${\cal A}_{\mu}=\{A\Delta N \mid A\in {\cal A}, N\in{\cal N}_{\mu}\}$, and $A\in{\cal A}_{\mu}$ iff there exist sets $A', A''\in{\cal A}$ with $A'\subseteq A \subseteq A''$ and $\mu^*(A''\setminus A') = 0.$
\item $\cU{\mu}$ is a finite measure, and the unique extension of $\mu$ to ${\cal A}_{\mu}$.
\item the measure space $(X, {\cal A}_{\mu}, \cU{\mu})$ is complete. It is called the \emph{\index{completion!$\mu$}$\mu$-completion} of $(X, {\cal A}, \mu)$. 
\end{enumerate}
\EndProposition

\BeginProof
1.
Since ${\cal N}_{\mu}$ is a $\sigma$-ideal, we infer from Lemma~\ref{adjoin-sigma-ideal} that 
$A\in {\cal A}_{\mu}$ iff there exists $B\in{\cal A}$ and $N\in{\cal N}_{\mu}$ with $A = B\Delta N$. Now consider 
\begin{equation*}
  {\cal C} := \{A\in{\cal A}_{\mu}\mid \exists A', A''\in{\cal A}: A'\subseteq A \subseteq A'', \mu^*(A''\setminus A')=0\}.
\end{equation*}
Then ${\cal C}$ is a $\sigma$-algebra which contains ${\cal A}\cup{\cal N}_{\mu}$, thus ${\cal C}={\cal A}_{\mu}$. 

From the observation made just before stating the proposition it becomes clear that $\cU{\mu}$ is well defined on ${\cal A}_{\mu}$. Since $\mu^*$ coincides with $\cU{\mu}$ on ${\cal A}_{\mu}$ and the outer measure is countably subadditive~\SetCite{Lemma 1.107}, we have  to show that $\cU{\mu}$ is additive on ${\cal A}_{\mu}$. This follows immediately from the first part. If $\nu$ is another extension to $\mu$ on ${\cal A}_{\mu}$, ${\cal N}_{\nu} = {\cal N}_{\mu}$ follows, so that 
$
\cU{\mu}(A\Delta N) = \mu(A) = \nu(A) = \nu(A\Delta N)
$
whenever $A\Delta N\in{\cal A}_{\mu}$.

2.
Completeness of $(X, {\cal A}_{\mu}, \cU\mu)$ follows now immediately from the construction. 
\EndProof

Surprisingly, we have received more than we have shopped for, since
complete measure spaces are closed under the Souslin operation. This
is remarkable because the Souslin operation evidently bears no hint at all at
measures which are defined on the base space. In addition, measures
are defined through countable operations, while the Souslin operation
makes use of the uncountable space $\Nat^{\Nat}$.

\BeginProposition{compl-closed-souslin-1}
A complete measure space is closed under the Souslin operation.
\EndProposition

\BeginProof
Let $(X, {\cal A}, \mu)$ be complete, then it is enough to show that each $B\subseteq X$ has an ${\cal A}$-cover (Definition~\ref{cover-by-alg}); then the assertion will follow from Proposition~\ref{Marcewski-complete}. In fact, given $B$, construct $B^*\in{\cal A}$  such that $\mu(B^*) = \mu^{*}(B)$, see~\SetCite{Lemma 1.118}. Whenever $C\in{\cal A}$ with $B\subseteq C$, we evidently have every subset of $B^{*}\setminus C$ in ${\cal A}$ by completeness. 
\EndProof

These constructions work also for $\sigma$-finite measure spaces, as indicated above. Now let $M$ be a non-empty set of $\sigma$-finite measures on the measurable space $(X, {\cal A})$, then define the \emph{\index{completion!$M$}$M$-completion} $\cCompl{{\cal A}}{M}$ and the \emph{\index{completion!universal}universal completion} $\cU{{\cal A}}$ of the $\sigma$-algebra ${\cal A}$ through
\begin{align*}
\cCompl{{\cal A}}{M} & := \bigcap_{\mu\in M}{\cal A}_{\mu},\\
\cU{{\cal A}} & := \bigcap\{{\cal A}_{\mu}\mid \mu \text{ is a $\sigma$-finite measure on } {\cal A}\}.
\end{align*}

As an immediate consequence this yields that the analytic sets in a
Polish space are contained in the universal completion of the Borel
sets, specifically
\BeginCorollary{analytic-are-complete}
Let $X$ be a Polish space and $\mu$ be a finite measure on $\Borel{X}$. Then all analytic sets are contained in $\cU{\Borel{X}}$
\EndCorollary

\BeginProof
Proposition~\ref{compl-closed-souslin-1} together with Proposition~\ref{analytic-is-souslin-repr}.
\EndProof

Just for the record:

\BeginCorollary{univ-clos-is-souslin-closed}
The universal closure of a measurable space is closed under the Souslin operation.
\QED
\EndCorollary

Measurability of maps is preserved when passing to the universal closure.

\BeginLemma{meas-preserv-univ-closure}
Let $f: X\to Y$ be ${\cal A}$-${\cal B}$- measurable, then $f$ is $\cU{{\cal A}}$-$\cU{{\cal B}}$ measurable. 
\EndLemma

\BeginProof
Let $D\in\cU{{\cal B}}$ be a universally measurable subset of $Y$, then we have to show that $E := \InvBild{f}{D}$ is universally measurable in $X$. So we have to show that for every finite measure $\mu$ on ${\cal A}$ there exists $E', E''\in{\cal A}$ with $E'\subseteq E\subseteq E''$ and $\mu(E'\setminus E'') = 0$. Define $\nu$ as the image of $\mu$ under $f$, so that $\nu(B) = \mu(\InvBild{f}{B})$ for each $B\in{\cal B}$, then we know that there exists $D', D''\in{\cal B}$ with $D'\subseteq D \subseteq D''$ such that $\nu(D''\setminus D') = 0$, hence we have for the measurable sets $E' := \InvBild{f}{D'}$, $E'' := \InvBild{f}{D''}$ 
\begin{equation*}
 \mu(E''\setminus E') = \mu(\InvBild{f}{D''\setminus D'}) = \nu(D''\setminus D') = 0. 
\end{equation*}
Thus $\InvBild{f}{D}\in\cU{{\cal A}}$. 
\EndProof

We will give now two applications of this construction. The first will
show that a finite measure on a countably generated
sub-$\sigma$-algebra of the Borel sets of an analytic space has always
an extension to the Borel sets, the second will construct an extension
of a stochastic relation $K: (X, {\cal A})\Trans (Y, {\cal B})$ to a
stochastic relation $\cU{K}: (X, \cU{{\cal A}})\Trans (Y, \cU{{\cal
    B}})$, provided the target space $(Y, {\cal B})$ is
separable. This first application is derived from von Neumann's
Selection Theorem, which is established here as well. It is shown also
that a measurable surjection can be lifted to a measurable map between
finite measure spaces, provided the target space is a separable metric
space.

\Subsubsection{Lubin's Extension through von Neumann's Selectors}
\label{sec:lubin-extension}

Let $X$ be an analytic space, and ${\cal B}$ be a countably generated
sub-$\sigma$-algebra of $\Borel{X}$ we will show that each finite
measure defined on ${\cal B}$ has at least one extension to a measure
on $\Borel{X}$. This is established through a surprising selection argument, as
we will see. 

As a preparation, we require a universally measurable right inverse of
a measurable surjective map $f: X\to Y$.  We know from the Axiom of
Choice that we can find for each $y\in Y$ some $x\in X$ with $f(x) =
y$, because $\{\InvBild{f}{\{y\}}\mid y\in Y\}$ is a partition of $X$
into non-empty sets. Set $g(y) := x$. Selecting an inverse image in
this way will not guarantee, however, that $g$ has any favorable
properties, even if, say, both $X$ and $Y$ are compact metric and $f$
is continuous. Hence we will have to proceed in a more systematic way.

We will use the observation that each analytic set in a Polish space
can be represented as the continuous image of $\NatInf$, as discussed
in Proposition~\ref{char-analytic}.

We will first formulate a sequence of auxiliary statements that deal
with finding for a given surjective map $f: X \rightarrow Y$ a map $g: Y
\rightarrow X$ such that $f \circ g = id_Y$. This map $g$ should have
some sufficiently pleasant properties.

Thus in order to make the first step it turns out to be helpful
focusing the attention to analytic sets being the continuous images
of $\NatInf$. This looks a bit far fetched, because we want to deal
with universally measurable sets, but remember that analytic sets are
universally measurable.

We can lexicographically order
$\NatInf$ by saying that $\Folge{t} \preceq \Folge{t'}$ iff there
exists $k \in \Nat$ such that $t_k \leq t'_k$, and $t_j = t'_j$ for
all $\ell$ with $1 \leq j < k$. Then $\preceq$ defines a total order on $\NatInf$.
We will capitalize on this order, to be more precise, on the interplay
between the order and the topology. Let us briefly look into the order
structure of $\NatInf$.

\BeginLemma{closed-has-minimal}
Each nonempty closed set $F \subseteq \Nat^\infty$ has a minimal element
in the lexicographic order.
\EndLemma

\BeginProof
Let $n_1$ be the minimal first component of all elements of $F$, $n_2$
be the minimal second component of those elements of $F$ that start
with $n_1$, etc. This defines an element $t := \langle n_1, n_2,
\ldots \rangle$. We claim that $t\in F$. Let $U$ be an open
neighborhood of $t$, then there exists $k\in\Nat$ such that $t\in
\Sigma_{n_{1}\dots n_{k}}\subseteq U$. By construction,
$\Sigma_{n_{1}\dots n_{k}}\cap F\not= \emptyset$, thus each open
neighborhood of $t$ contains an element of $F$. Hence $t$ is an
accumulation point of $F$, and since $f$ is closed, it contains all
its accumulation points. Thus $t\in F$.
\EndProof

Now if we have $f: \NatInf\to X$ continuous, we know that the inverse
images $\InvBild{f}{\{y\}}$ for $y\in \Bild{f}{\NatInf}$ are
closed. Thus me may pick for each $y\in \Bild{f}{\NatInf}$ this
smallest element. This turns out to be a suitable choice, as the
following statement shows.

\BeginLemma{first-Borel-section}
Let $X$ be Polish, $Y \subseteq X$ analytic with $Y = \Bild{f}{\NatInf}$
for some continuous $f: \NatInf \rightarrow X$. Then there exists
$g: Y \rightarrow\NatInf$  such that
\begin{enumerate}
  \item $f \circ g = id_Y$,
  \item $g$ is $\cU{\Borel{Y}}$-$\cU{\Borel{\NatInf}}$-measurable.
\end{enumerate}
\EndLemma

\BeginProof
1. Since $f$ is continuous, the inverse image $\InvBild{f}{\{y\}}$ for
each $y \in Y$ is a closed and nonempty set in $\Nat^\infty$. Thus
this set contains a minimal element $g(y)$ in the lexicographic order
$\preceq$ by Lemma~\ref{closed-has-minimal}. It is clear that $f(g(y))
= y$ holds for all $y \in Y$.

2. Denote by $ A(t') := \{t \in \Nat^\infty \mid t \prec
t'\}, $ then $A(t')$ is open: let $\Folge{\ell} = t \prec
t'$ and $k$ be the first component in which $t$ differs from
$t'$, then $ \Sigma_{\ell_1\ldots \ell_{k-1}} $ is an open
neighborhood of $t$ that is entirely contained in $A(t')$. It is
easy to see that $ \{A(t') \mid t' \in \Nat^\infty\} $ is a
generator for the Borel sets of $\Nat^\infty$.

3. We claim that $ \InvBild{g}{A(t')} = \Bild{f}{A(t')} $ holds.
In fact, let $y \in \InvBild{g}{A(t')}$, so that $ g(y) \in
A(t')$, then $y = f(g(y)) \in \Bild{f}{A(t')}. $ If, on the
other hand, $y = f(t)$ with $t \prec t'$, then by
construction $t \in \InvBild{f}{\{y\}}$, thus $g(y) \preceq t
\prec t'$, settling the other inclusion.

This equality implies that $\InvBild{g}{A(t')}$ is an analytic
set, because it is the image of an open set under a continuous
map. Consequently, $\InvBild{g}{A(t')}$ is universally measurable
for each $A(t')$ by Corollary~\ref{analytic-are-complete}.
Thus  $g$ is a universally measurable map.
\EndProof

This statement is the work horse for establishing that a right inverse
exists for surjective Borel maps between an analytic space and a
separable measurable space. All we need to do now is to massage things
into a shape that will render this result applicable in the desired
context.  The following theorem is usually attributed to von Neumann\index{theorem!von Neumann Selection}.

\BeginTheorem{von-neumann-selection}
Let $X$ be an analytic space, $(Y, \mathcal{B})$ a separable measurable space and
$f: X \rightarrow Y$ a surjective measurable map. Then there exists
$g: Y \rightarrow X$ with these properties:
\begin{enumerate}
  \item $f \circ g = id_Y$,
  \item $g$ is $\cU{\mathcal{B}}$-$\cU{\Borel{X}}$-measurable.
\end{enumerate}
\EndTheorem

\BeginProof
1.
We may and do assume by Lemma~\ref{Sep-Analytic-Map} that $Y$ is an analytic
subset of a Polish space $Q$, and that $X$ is an analytic subset of a Polish
space $P$. $x \mapsto \langle x, f(x)\rangle$ is a bijective Borel map
from $X$ to the graph of $f$, so $\mathsf{graph}(f)$ is an analytic set by
Proposition~\ref{AnalyticIsStableUnder}. Thus we can find a continuous
map $F: \NatInf \rightarrow P \times Q$ with
$
\Bild{F}{\NatInf} = \mathsf{graph}(f).
$
Consequently, $\pi_Q \circ F$ is a continuous map from $\NatInf$ to
$Q$ with
$$
\Bild{(\pi_Q \circ F)}{\NatInf} = \Bild{\pi_Q}{\mathsf{graph}(f)} = Y.
$$
Now let $G: Y \rightarrow \NatInf$ be chosen according to
Lemma~\ref{first-Borel-section} for $\pi_Q \circ F$. Then
$
g := \pi_P \circ F \circ G: Y \rightarrow X
$
is the map we are looking for:
\begin{itemize}
  \item $g$ is universally measurable, because $G$ is, and because
  $\pi_P \circ F$ are continuous, hence universally measurable as well,
  \item $f \circ g = f \circ (\pi_P \circ F \circ G) = (f \circ \pi_P) \circ F \circ G = \pi_Q \circ F \circ G
  = id_Y$, so $g$ is right inverse to $f$.
\end{itemize}
\EndProof

Due to its generality, the von Neumann Selection Theorem has many
applications in diverse areas, many of them surprising. The art is
plainly to reformulate the problem so that an application of this
selection theorem is possible. We pick two applications, viz., showing
that the image $\FinM{f}$ of a surjective Borel map $f$ yields a
surjective Borel map again, and Lubin's measure extension.

\BeginProposition{SubProbIsOnto}
Let $X$ be an analytic space,  $Y$ a second countable metric space. If $f: X
\rightarrow Y$ is a surjective Borel map, so is
$
\FinM{f}: \FinM{X} \rightarrow \FinM{Y}.
$
\EndProposition

\BeginProof
1. 
From Theorem~\ref{von-neumann-selection} we find a map
$
g: Y \rightarrow X
$
such that
$
f \circ g = id_Y
$
and $g$ is $\cU{\Borel{Y}} -
\cU{\Borel{X}}$-measurable.

2. Let $\nu \in \FinM{Y}$, and define
$
\mu := \FinM{g}(\nu),
$
then
$
\mu \in \FinM{X, \cU{\Borel{X}}}
$
by construction. Restrict $\mu$ to the Borel sets on $X$, obtaining
$
\mu_0 \in \FinM{X, \Borel{X}}.
$
Since we have for each set $B \subseteq Y$ the equality
$
g^{-1}\bigl[f^{-1}[B]\bigr] = B,
$
we see that for each $B \in \Borel{Y}$
\begin{equation*}
\FinM{f}(\mu_0)(B) = \mu_0(\InvBild{f}{B}) =
\mu(\InvBild{f}{B}) =
\nu(g^{-1}\bigl[\InvBild{f}{B}\bigr]) = \nu(B)
\end{equation*}
holds.
\EndProof

This has as a consequence that $\FinSenza$ is an endofunctor
on the category of Polish or analytic spaces with surjective Borel maps
as morphisms; it displays a pretty interaction of reasoning in
measurable spaces and arguing in categories. 

The following extension theorem due to \index{theorem!Lubin} Lubin
shows that one can extend a finite measure from a countably generated
sub-$\sigma$-algebra to the Borel sets of an analytic space. In
contrast to classical extension theorems it does not permit to conclude that
the extension is uniquely determined.

\BeginTheorem{lubin-extension}
Let $X$ be an analytic space, and $\mu$ be a finite measure on a countably generated sub-$\sigma$-algebra ${\cal A}\subseteq\Borel{X}$. Then there exists an extension of $\mu$ to a finite measure $\nu$ on $\Borel{X}$. 
\EndTheorem

\BeginProof
Let $\Folge{A}$ be the generator of ${\cal A}$, and define the map $f: X\to \{0, 1\}^{\Nat}$ through $x \mapsto (\chi_{A_{n}})_{n\in\Nat}$. Then $M := \Bild{f}{X}$ is an analytic space, and $f$ is $\Borel{X}$-$\Borel{M}$ measurable by Proposition~\ref{sep-is-sep-metr} and Proposition~\ref{AnalyticIsStableUnder}. Moreover, 
\begin{equation}
\label{this-sigma-alg}
{\cal A} = \{\InvBild{f}{C}\mid C\in\Borel{M}\}.
\end{equation}
By von Neumann's Selection Theorem~\ref{von-neumann-selection} there exists $g: M\to X$ with $f\circ g = id_{M}$ which is $\cU{\Borel{M}}$-$\cU{\Borel{X}}$-measurable. Define 
\begin{equation*}
\nu(B) := \cU{\mu}\bigl(\InvBild{(g\circ f)}{B}\bigr)
\end{equation*}
for $B\in\Borel{X}$ with $\cU{\mu}$ as the completion of $\mu$ on $\cU{{\cal A}}$. Since we have for $B\in\Borel{X}$ that $\InvBild{g}{B}\in\cU{\Borel{M}}$, we may conclude from (\ref{this-sigma-alg}) that $\InvBild{f}{\InvBild{g}{B}}\in\cU{{\cal A}}$. $\nu$ is an extension to $\mu$. In fact, given $A\in{\cal A}$, we know that $A=\InvBild{f}{C}$ for some $C\in\Borel{M}$, so that we obtain
\begin{align*}
  \nu(A) & = \cU{\mu}\bigl(\InvBild{(g\circ f)}{\InvBild{f}{C}}\bigr)\\
& = \cU{\mu}\bigl(f^{-1}\circ g^{-1}\circ f^{-1}[C]\bigr)\\
& = \cU{\mu}\bigl(\InvBild{f}{C}\bigr), && \text{ since }f\circ g = id_{M}\\
& = \cU{\mu}(A)\\
& = \mu(A).
\end{align*}
\EndProof

This can be rephrased in a slightly different way. The identity $id_{{\cal A}}: (X, \Borel{X})\to (X, {\cal A})$ is measurable, because ${\cal A}$ is a sub-$\sigma$-algebra of $\Borel{X}$. Hence it induces a measurable map $\SubProbSenza{id_{{\cal A}}}: \SubProbSenza(X, \Borel{X})\to \SubProbSenza(X, {\cal A})$. Lubin's Theorem then implies that $\SubProbSenza{id_{{\cal A}}}$ is surjective. This is so since for a given $\mu\in \SubProbSenza(X, \Borel{X})$, $\SubProbSenza{id_{{\cal A}}}(\mu)$ is just the restriction of $\mu$ to the sub-$\sigma$-algebra ${\cal A}$. 

\Subsubsection{Completing a Transition Kernel}
\label{sec:ext-stoch-rel-to-compl}

{
% Beware of local macros
%
\def\measb{\Borel}
\newcommand{\theCompl}[2]{\cCompl{#1}{#2}}
\newcommand{\tee}[1]{\ensuremath{{\sigmaFont S}_{#1}}}
\def\seqFont{\mathfrak}
\def\sigmaFont{\mathcal}
\newcommand{\abbrev}[1]{\ensuremath{\seqFont{#1}}}
\def\g{\abbrev{g}}
\def\cpl{\cU}

In some probabilistic models for modal logics it becomes sometimes necessary to assume that the state space is closed under Souslin's operation, see for example~\cite{EED-PDL-TR}, on the other hand one may not always assume that a complete measure space is given. Hence one wants to complete it, but it is then also mandatory to complete the transition law as well. This means that an extension of the transition law to the completion becomes necessary. This problem will be studied now. 

The completion of a measure space is described in terms of null sets and using inner and outer approximations, see Proposition~\ref{completion-of-measure-space}. We will use the latter here, fixing measurable spaces $(X; {\cal A})$ and $(Y, {\cal B})$. Denote by $\tee{X}$ the smallest $\sigma$-algebra on $X$ which contains ${\cal A}$ and which is closed under the Souslin operation, hence $\tee{X}\subseteq \cU{{\cal A}}$ by Corollary~\ref{univ-clos-is-souslin-closed}. 

Fix $K: (X, {\cal A}) \Trans (Y, {\cal B})$ as a transition kernel, and assume first that ${\cal B}$ is the $\sigma$-algebra of Borel sets for a \emph{second countable metric space}. This means that the topology $\tau$ of $Y$ has a countable base $\tau_0$, which in turn means that
$
G = \bigcup\{H\in \tau_0 \mid H \subseteq G\}
$
for each open set $G\in\tau$. 

For each $x\in X$ we have through the transition kernel $K$ a finite measure $K(x)$, to which we may associate an out measure $\bigl(K(x)\bigr)^{*}$ on the power set of $X$. We want to show that the map 
\begin{equation*}
x \mapsto \bigl(K(x)\bigr)^*(A)
\end{equation*}
is $\tee{X}$-measurable for each $A\subseteq Y$; define for convenience
\begin{equation*}
K^*(x) := \bigl(K(x)\bigr)^*.
\end{equation*}
Establishing measurability is broken into a sequence of steps. 

We need the following regularity argument (but compare Exercise~\ref{ex-meas-non-regular} for the non-metric case)

\BeginLemma{are-regular}
Let $\mu$ be a finite measure on $(Y, \Borel{Y})$,
$B\in\Borel{Y}$. Then we can find for each $\epsilon>0$ an open set
$G\subseteq Y$ with $B\subseteq G$ and a closed set $F\supseteq B$
such that $\mu(G\setminus F) < \epsilon$.
\EndLemma

\BeginProof
Let 
\begin{equation*}
  {\cal G} := \{B\in\Borel{Y} \mid \text{ the assertion is true for $B$}\}.
\end{equation*}
Then plainly ${\cal G}$ is closed under complementation and contains
the open as well as the closed sets. If $F\subseteq Y$ is closed, we
can represent $F = \bigcap_{n\in\Nat} G_{n}$ with $\Folge{G}$ as a
decreasing sequence of open sets, hence $ \mu(F) =
\inf_{n\in\Nat}\mu(F_{n}) = \lim_{n\to\infty}\mu(F_{n}), $ so that
${\cal G}$ also contains the closed sets; one arguments similarly for
the open sets as increasing unions of open sets.

Now let $\Folge{B}$ be a sequence of mutually disjoint sets in ${\cal
  G}$, select $G_{n}$ open for $B_{n}$ and $\epsilon/2^{-(n+1)}$, then
$G := \bigcup_{n\in\Nat}G_{n}$ is open with $B := \bigcup_{n\in\Nat}
B_{n} \subseteq G$ and $\mu(G\setminus B) \leq \epsilon$. Similarly,
select the sequence $\Folge{F}$ with $F_{n}\subseteq B_{n}$ and
$\mu(B_{n}\setminus F_{n}) < \epsilon/2^{-(n+1)}$ for all $n\in\Nat$,
put $F := \bigcup_{n\in\Nat}F_{n}$ and select $m\in\Nat$ with
$\mu(F\setminus \bigcup_{n=1}^{m}F_{n})< \epsilon/2$, then $F' :=
\bigcup_{n=1}^{m}F_{n}$ is closed, $F'\subseteq B$, and
$\mu(B\setminus F') < \epsilon$.

Hence ${\cal G}$ is closed under complementation as well as countable
disjoint unions; this implies ${\cal G} = \Borel{Y}$ by the
$\pi$-$\lambda$ Theorem~\ref{Pi-Lambda}.
\EndProof

Fix $A\subset Y$ for the moment. We claim that 
\begin{equation*}
K^*(x)(A) = \inf\{K(x)(G) \mid A \subseteq G \text{ open}\}
\end{equation*}
holds for each $x \in X$. In fact, given $\epsilon > 0$, there exists $A \subseteq A_0\in\Borel{Y}$ with 
$
K(x)(A_0) - K^*(x)(A) < \epsilon/2.
$
Applying Lemma~\ref{are-regular} to $K(x)$, we find an open set $G \supseteq A_0$ with 
$
K(x)(G) - K(x)(A_0) < \epsilon/2,
$
thus 
$
K(x)(G) - K^*(x)(A) < \epsilon.
$

$\tau_0$ is a countable base for the open sets, which we may assume to be closed under finite unions (because otherwise
$
\{G_1\cup\dots\cup G_k \mid k \in \Nat, G_1, \dots, G_k\in\tau_0\}
$
is a countable base which has this property). Hence we obtain 
\begin{equation}
\label{ext-B}
K^*(x)(A) = \inf\{\sup_{n\in\Nat}K(x)(G_n) \mid A \subseteq \bigcup_{n\in\Nat}G_n, \Folge{G}\subseteq \tau_0\text{ increases}\}.
\end{equation}
Let 
\begin{equation*}
%\label{ }
\mathcal{G}_A := \{\Folge{G} \subseteq \tau_0\mid \Folge{G}\text{ increases and } A \subseteq \bigcup_{n\in\Nat} G_n\}
\end{equation*}
be the set of all increasing sequences from base $\tau_0$ which cover $A$. Partition $\mathcal{G}_A$ into the sets
\begin{align*}
%\label{ }
\mathcal{N}_A & := \{\g \in \mathcal{G}_A \mid \g \text{ contains only a finite number of sets}\},\\
\mathcal{M}_A & := \mathcal{G}_A\setminus\mathcal{N}_A.
\end{align*}
Because $\tau_0$ is countable, $\mathcal{N}_A$ is. %Denote for a sequence $\Folge{x}$ its first $k$ elements $x_1, \dots, x_k$ by $\Folge{x}\mid k$. 

\BeginLemma{construct-f}
There exists an injective map 
$
\Phi: \mathcal{M}_A\to \Nat^\Nat
$
such that $\g\mid k = \g'\mid k$ implies $\Phi(\g)\mid k = \Phi(\g')\mid k$ for all $k\in\Nat$. 
\EndLemma

\BeginProof
1.
Build an infinite tree in this way: the root is the empty set, a node $G$ at level $k$ has all elements $G'$ from $\tau_0$ with $G \subseteq G'$ as offsprings. Remove from the tree all paths $H_1, H_2, \dots$ such that 
$
A \not\subseteq \bigcup_{n\in\Nat} H_n.
$
Call the resulting tree $\mathcal{T}$. 

2.
Put $G_0 := \emptyset$, and let $\mathcal{T}_{1, G_0}$ be the set of nodes of $\mathcal{T}$ on level 1 (hence just the offsprings of the root $G_0$), then there exists an injective map 
$
\Phi_{1, G_0}: \mathcal{T}_{1, G_0}\to \Nat.
$
If $G_1, \dots, G_k$ is a finite path to inner node $G_k$ in $\mathcal{T}$, denote by 
$
\mathcal{T}_{k+1, G_1, \dots, G_k}
$
the set of all offsprings of $G_k$, and let 
%$
\begin{equation*}
%\label{ }
\Phi_{k+1, G_1, \dots, G_k}:  \mathcal{T}_{k+1, G_1, \dots, G_k}\to\Nat
\end{equation*}
%$
be an injective map. Define %for $\g := \Folge{G}$
\begin{equation*}
%\label{ }
\Phi: \begin{cases}
      \mathcal{M}_A& \to\Nat^\Nat, \\
      \Folge{G} &\mapsto \bigl(\Phi_{n, G_1, \dots, G_{n-1}}(G_n)\bigr)_{n\in\Nat}.
\end{cases}
\end{equation*}

3.
Assume $\Phi(\g) = \Phi(\g')$, then an inductive reasoning shows that $\g = \g'$. In fact, $G_1 = G_1'$, since $\Phi_{1, \emptyset}$ is injective. If  $\g\mid k = \g'\mid k$ has already been established, we know that 
$
\Phi_{k + 1, G_1, \dots, G_k} = \Phi_{k+1, G_1', \dots, G_k'}
$
is injective, so that $G_{k+1} = G'_{k+1}$ follows. A similar inductive argument shows that  $\Phi(\g)\mid k = \Phi(\g')\mid k$, provided $\g \mid k = \g'\mid k$ for each $k \in \Nat$ holds. 
\EndProof

The following lemmata collect some helpful properties.

\BeginLemma{prop-phi-1}
$\g=\g'$ iff $\Phi(\g)\mid k = \Phi(\g')\mid k$ for all $k \in \Nat$, whenever $\g, \g'\in \mathcal{M}_A$.  
\QED
\EndLemma

\BeginLemma{prop-phi-2}
Denote by 
$
J_k := \{\alpha \mid k \mid \alpha\in\Bild{\Phi}{\mathcal{M}_A}\}
$
all initial pieces of sequences in the image of $\Phi$. Then $\alpha \in \Bild{\Phi}{\mathcal{M}_A}$ iff $\alpha\mid k \in J_k$ for all $k\in\Nat$.  
\EndLemma

\BeginProof
Assume that $\alpha =\Phi(\g)\in \Bild{\Phi}{\mathcal{M}_A}$  with $\g = \Folge{C}\in\mathcal{M}_A$ and $\alpha\mid k \in J_k$  for all $k\in\Nat$, so for given $k$ there exists $\g^{(k)} = \Folge{C^{(k)}}\in \mathcal{M}_A$ with $\alpha\mid k = \Phi(\g^{(k)})\mid k$. Because $\Phi_1$ is injective, we obtain $C_1 = C^{(1)}_1. $ Assume for the induction step that 
$
G_i = G_i^{(j)}
$
has been shown for $ 1 \leq i, j \leq k$. Then we obtain from 
$
\Phi(\g)\mid k+1 = \Phi(\g^{(k+1)})\mid k+1 
$
that 
$
G_1 = G_1^{(k+1)}, \dots, G_k = G_k^{(k+1)}.
$
Since 
$
\Phi_{k+1, G_1, \dots, G_k}
$
is injective, the equality above implies $G_{k+1} = G_{k+1}^{(k+1)}$. Hence $\g = \g^{(k)}$ for all $k\in\Nat$, and $\alpha \in \Bild{\Phi}{\mathcal{M}_A}$ is established. The reverse implication is trivial. 
\EndProof

\BeginLemma{er-meas}
$
E_r := \{x \in X \mid K^*(x)(A) \leq r\}\in\tee{X}
$
for $r\in\pReal$. 
\EndLemma

\BeginProof
The set $E_r$ can be written as
\begin{equation*}
%\label{ }
E_r = \bigcup_{\g\in\mathcal{N}_A}\{x \in X \mid K(x)\bigl(\bigcup\g\bigr) \leq r\}
\cup
\bigcup_{\g\in\mathcal{M}_A}\{x \in X \mid K(x)\bigl(\bigcup\g\bigr) \leq r\}
\end{equation*}
Because $\mathcal{N}_A$ is countable, and $K: X \Trans Y$ is a transition kernel, we infer 
\begin{equation*}
%\label{ }
\bigcup_{\g\in\mathcal{N}_A}\{x \in X \mid K(x)\bigl(\bigcup\g\bigr) \leq r\}\in \measb{X}
\end{equation*}
Put for $v\in\Nats$ 
\begin{equation*}
D_v := 
%\label{ }
\begin{cases}
      \emptyset,& \text{if } v \notin\bigcup_{k\in\Nat}J_k, \\
      \{x \in X \mid K(x)(G_n) \leq r\},& \text{if } v = \Phi\bigl(\Folge{G}\bigr)\mid n.
\end{cases}
\end{equation*}

Lemma~\ref{prop-phi-1} and Lemma~\ref{prop-phi-2} show that $D_v \in \measb{X}$ is well defined. Because
\begin{equation}
\label{ect-C}
\bigcup_{\g\in\mathcal{M}_A}\{x \in X \mid K(x)\bigl(\bigcup\g\bigr) \leq r\}
=
\bigcup_{\alpha\in\Nat^\Nat}\bigcap_{n\in\Nat} D_{\alpha\mid n},
\end{equation}
and because $\tee{X}$ is closed under the Souslin operation and contains $\measb{X}$, we conclude that $E_r\in\tee{X}$. 
\EndProof

\BeginProposition{the-extension}
Let $K: (X; {\cal A})\Trans (Y, {\cal B})$ be a transition kernel, and assume that $Y$ is a separable metric space. Let $\tee{X}$ be the smallest $\sigma$-algebra which contains ${\cal A}$ and which is closed under the Souslin operation. Then there exists a unique transition kernel 
\begin{equation*}
%\label{ }
\overline{K}: (X, \tee{X}) \Trans (Y, \theCompl{\Borel{Y}}{\{K(x)\mid x\in X\}})
\end{equation*}
extending $K$.
\EndProposition

\BeginProof
1.
Put 
$
\overline{K}(x)(A) := K^*(x)(A)
$
for $x\in X$ and $A\in\theCompl{\Borel{Y}}{\{K(x)\mid x\in X\}}$. Because $A$ is an element of the $K(x)$-completion of $\Borel{Y}$, we know that $\cU{K}(x) = \cU{K(x)}$ defines a sub probability on 
$
\theCompl{\Borel{Y}}{\{K(x)\mid x\in X\}}.
$
It is clear that $\overline{K}(x)$ is the unique extension of $K(x)$ to the latter $\sigma$-algebra. It remains to be shown that $\overline{K}$ is a transition kernel. 

2.
Fix $A\in\theCompl{\Borel{Y}}{\{K(x)\mid x\in X\}}$ and $q\in[0, 1]$, then
\begin{equation*}
%\label{ }
\{x \in X \mid K^*(x)(A) < q\} = 
\bigcup_{\ell\in\Nat}\bigcup_{\g\in\mathcal{G}_A}\{x \in X \mid K(x)\bigl(\bigcup\g\bigr) \leq q - \frac{1}{\ell}\}
\end{equation*}
The latter set is a member of $\tee{X}$ by Lemma~\ref{er-meas}. 
\EndProof

Separability of the target space is required because it is this property which makes sure that the measure for each Borel set can be approximated arbitrarily well from within by closed sets, and from the outside by open sets~\cite[Lemma 3.4.14]{Srivastava}. 

Before discussing consequences, a mild generalization to separable measurable spaces should be mentioned. Proposition~\ref{the-extension} yields as an immediate consequence:
\BeginCorollary{ext-to-completion}
Let $K:(X; {\cal A}) \Trans (Y, {\cal B})$ be a transition kernel such that $(Y, {\cal B})$ is a separable measurable space. Assume that $\sigmaFont{X}$  is a $\sigma$-algebra on $X$ which is closed under the Souslin operation with $\tee{X}\subseteq \sigmaFont{X}$, and that  
$\sigmaFont{Y}$ is a $\sigma$-algebra on $X$ with
$
{\cal B}\subseteq\sigmaFont{Y}\subseteq \theCompl{{\cal B}}{\{K(x)\mid x\in X\}}. 
$
Then there exists a unique extension 
$
(X, \sigmaFont{X}) \Trans (Y, \sigmaFont{Y})
$
to $K$. In particular $K$ has a unique extension to a transition kernel $\overline{K}:(X, \cU{\cal A})\Trans(Y, \cU{\cal B})$. 
\EndCorollary

\BeginProof
This follows from Proposition~\ref{the-extension} and the characterization of separable measurable spaces in Proposition~\ref{sep-is-sep-metr}.  
\EndProof 
}
%%% Local Variables: 
%%% mode: latex
%%% TeX-master: "../Mskr3"
%%% End: 

%\Input{\Folder/Selections}
%spell checked - 24Aug14
\def\weak#1{\ensuremath{{#1}^w}}
\Subsection{Measurable Selections}
\label{sec:meas-selections}
Looking again at von Neumann's Selection
Theorem~\ref{von-neumann-selection}, we have found for a given
surjection $f: X\to Y$ a universally measurable map $g: Y \to X$ with
$f\circ g = id_{Y}$. This can be rephrased: we have $g(y) \in
\InvBild{f}{\{y\}}$ for each $y\in Y$, so $g$ may be considered a
universal measurable selection for the set valued map $y \mapsto
\InvBild{f}{\{y\}}$. We will consider this problem from a slightly
different angle by assuming that $(X, {\cal A})$ is a measurable, $Y$
is a Polish space, and that we have a set valued map $F: X\to
\PowerSet{Y}\setminus\{\emptyset\}$ for which a measurable selection
is to be constructed, i.e., a measurable (not merely universally
measurable) map $g: X\to Y$ such that $g(y)\in F(y)$ for all $y\in
Y$. Clearly, the Axiom of Choice guarantees the existence of a map
which picks an element from $F(y)$ for each $y$, but this is not
enough. 

We assume that $F(y)$ is always a closed subset of $Y$, and
that it is measurable. Since $F$ does not necessarily take single
values only, we have to define measurability in this case. Denote by
$\Closed(Y)$ the set of all closed and non-empty subsets of $Y$.

\BeginDefinition{set-valued-measurability}
A map $F: X\to \Closed(Y)$ from a measurable space $(X, {\cal A})$ to
the closed non-empty subsets of a Polish space $Y$ is called
\emph{\index{measurable!set-valued map}measurable} (or a \emph{measurable \index{measurable!relation}relation}) iff
\begin{equation*}
  \weak{F}(G) := \{x\in X\mid F(x)\cap G\not=\emptyset\} \in{\cal A}
\end{equation*}
for every open subset $G\subseteq Y$. The map $s: X\to Y$ is called a
\emph{\index{measurable!selector}measurable selector} for $F$ iff $s$ is ${\cal
  A}$-$\Borel{Y}$-measurable such that $s(x)\in F(x)$ for all $x\in
X$.
\EndDefinition

Since $\{f(x)\}\cap G \not=\emptyset$ iff $f(x)\in G$, measurability
as defined in this definition is a generalization of measurability for
point valued maps $f: X\to Y$.

The selection theorem  due to \index{theorem!Kuratowski and Ryll-Nardzewski}Kuratowski and Ryll-Nardzewski tell us that a measurable selection exists for a measurable closed valued map, provided $Y$ is Polish. To be specific:

\BeginTheorem{meas-selections-exist}
Given a measurable space $(X, {\cal A})$ and a Polish space $Y$, 
a measurable map $F: X \to \Closed(Y)$ has a measurable selector. 
\EndTheorem

\BeginProof
Fix a complete metric $d$ on $Y$.  Denote by $B(y, r)$ the open ball
around $y\in Y$ with radius $r>0$; $d$ is the metric on $Y$ such that
the metric space $(Y, d)$ is complete. Recall that the distance of an
element $y$ to a closed set $C$ is $d(y, C) := \inf\{d(y, y')\mid
y'\in C\},$ hence $d(y, C) = 0$ iff $y\in C$.

Let $\Folge{y}$ be dense, and define $f_{1}(x) := y_{n}$, if $n$ is
the smallest index $k$ so that $F(x)\cap B(y_{k}, 1)\not=\emptyset$. Then $f_{1}: X \to Y$ is ${\cal A}$-$\Borel{Y}$ measurable, because the map takes only a countable number of values and 
\begin{equation*}
  \{x\in X\mid f_{1}(x) = y_{n}\} = \weak{F}(B(y_{n}, 1))\setminus\bigcup_{k=1}^{n-1}\weak{F}(B(y_{k}, 1)).
\end{equation*}
Proceeding inductively, assume that we have defined measurable maps $f_{1}, \dots, f_{n}$ such that
\begin{align*}
d(f_{j}(x), f_{j+1}(x)) & < 2^{-(j-1)},&& 1 \leq j < n\\
d(f_{j}(x), F(x)) & < 2^{-j}, && 1 \leq j \leq n
\end{align*}
Put $X_{k} := \{x\in X \mid f_{n}(x) = y_{k}\}$, and define $f_{k+1}(x) := y_{\ell}$ for $x\in X_{k}$, where $\ell$ is the smallest index $m$ such that $F(x)\cap B(y_{k}, 2^{-n})\cap B(y_{m}, 2^{-(n+1)}) \not= \emptyset.$ Moreover, there exists $y'\in  B(y_{k}, 2^{-n})\cap B(y_{m}, 2^{-(n+1})$, thus 
\begin{equation*}
d(f_{n}(x), f_{n+1}(x)) \leq d(f_{n}(x), y') + f(f_{n+1}(x), y') < 2^{-n} + 2^{-(n+1)}
\end{equation*}
The argumentation from above shows that $f_{n+1}$ takes only countably many values, and we know that $d(f_{n+1}(x), F(x)) < 2^{-(n+1)}.$

Thus $(f_{n}(x))_{n\in\Nat})$ is a Cauchy sequence for each $x\in
X$. Since $(Y, d)$ is complete, the limit $f(x)
:= \lim_{n\to \infty} f_{n}(x)$ exists with $d(f(x), F(x)) = 0$, hence
$f(x)\in F(x)$, because $F(x)$ is closed. Moreover as a pointwise
limit of a sequence of measurable functions $f$ is measurable, so $f$
is the desired measurable selector.
\EndProof

It is possible to weaken the conditions on $F$ and on ${\cal A}$, see
Exercise~\ref{ex-selection-algebra}. This theorem has an interesting
consequence, viz., that we can find a sequence of dense selectors for
$F$.

\BeginCorollary{castaing-representation}
Under the assumptions of Theorem~\ref{meas-selections-exist}, a
measurable map $F: X \to \Closed(Y)$ has a sequence $\Folge{f}$ of
measurable selectors such that $\{f_{n}(x)\mid n\in\Nat\}$ is dense in
$F(x)$ for each $x\in X$.
\EndCorollary

\BeginProof
1.
We use notations from above. Let again $\Folge{y}$ be a dense sequence in $Y$, and define for $n, m\in\Nat$ the map
\begin{equation*}
F_{n, m}(x) := 
\begin{cases}
  F(x)\cap B(y_{n}, 2^{-m}), & \text{ if } x\in\weak{F}(B(y_{n}, 2^{-m}))\\
F(x), & \text{ otherwise}
\end{cases}
\end{equation*}
Denote by $H_{n, m}(x)$ the closure of $F_{n, m}(x)$. 

2.
$H_{n, m}: X \to \Closed(Y)$ is measurable. In fact, put $A_{1} := \weak{F}(B(y_{n}, 2^{-m})), A_{2} := X\setminus A_{1}$, then $A_{1}, A_{2}\in {\cal A}$, because $F$ is measurable and $B(y_{n}, 2^{-m})$ is open. But then we have for an open set $G\subseteq Y$
\begin{align*}
  \{x \in X \mid H_{n, m}\cap G \not= \emptyset\} 
  & =
  \{x \in X \mid F_{n, m} \cap G \not= \emptyset\}\\
  & = \{x \in A_1 \mid F(x)\cap G \cap B(y_n, 2^{-m}) \not= \emptyset\} \cup 
  \{x \in A_2\mid F(x)\cap G\not= \emptyset\},
\end{align*}
thus $\weak{H}_{n, m}(G)\in \mathcal{A}$.

3.
We can find a measurable selector $s_{n, m}$ for $H_{n, m}$ by  Theorem~\ref{meas-selections-exist}, so we have to show that 
$
\{s_{n, m}(x) \mid n,m\in \Nat\}
$ 
is dense in $F(x)$ for each $x\in X$. Let $y\in F(x)$. Given $\epsilon>0$, select $m$ with $2^{-m} < \epsilon/2$; there exists $y_n$ with $d(y, y_n) < 2^{-m}$. Thus $x\in \weak{H}_{n, m}(B(y_n, 2^{-m}))$, and $s_{n, m}(x)$ is a member of the closure of $B(y_n, 2^{-m})$, which means $d(y, s_{n, m}(x)) < \epsilon$. Now arrange $\{s_{n, m}(x) \mid n,m\in \Nat\}$ as a sequence, then the assertion follows. 
\EndProof

This is a first application of measurable selections. 

\BeginExample{hit-meas-selections}
% Let $Y$ be a Polish space, then $\Closed(Y)$ is a Polish space with
% the Vietoris topology\Rand{Quelle}, and the hit-$\sigma$-algebra
% ${\cal H}_{{\cal G}}(\Closed(Y))$ are Borel sets for this topology
Call a map $h: X\to \Borel{Y}$ for the Polish space $Y$
\emph{\index{hit!measurable}hit-measurable} iff $h$ is measurable with
respect to ${\cal A}$ and ${\cal H}_{{\cal G}}(\Borel{Y})$, where
${\cal G}$ is the set of all open sets in $Y$, see
Example~\ref{hit-sigma-algebra}.  Thus $h$ is hit-measurable iff
$\{x\in X \mid h(x)\cap U\not=\emptyset\}\in {\cal A}$ for each open
set $U\subseteq Y$. If $h$ is image finite (i.e., $h(x)$ is always
non-empty and finite), then there exists a sequence $\Folge{f}$ of
measurable maps $f_{n}: X\to Y$ such that $h(x) = \{f_{n}(x)\mid
n\in\Nat \}$ for each $x\in X$. This is so because $h: X \to
\Closed(Y)$ is measurable, hence
Corollary~\ref{castaing-representation} is applicable.
\EndExample

Transition kernels into Polish spaces induce a measurable closed
valued map, for which selectors exist.

\BeginExample{transitins-have-selections}
Let under the assumptions of Theorem~\ref{meas-selections-exist} $K: (X, {\cal A})\Trans (Y, \Borel{Y})$ be a transition kernel with $K(x)(Y) > 0$ for all
$x\in X$. Then there exists a measurable map $f: X\to Y$ such that
$K(x)(U) > 0$, whenever $U$ is an open neighborhood of $f(x)$. 

In fact, $\Gamma: x \mapsto \supp(K(x))$ takes non-empty and closed
values by Lemma~\ref{support-for-regular}. If $G\subseteq Y$ is open,
then 
\begin{equation*}
  \weak{\Gamma}(G) = \{x\in X\mid \supp(K(x))\cap G \not=\emptyset\} =
  \{x\in X \mid K(x)(G) > 0\}\in {\cal A}.
\end{equation*}
Thus $\Gamma$ has a measurable selector $f$ by
Theorem~\ref{meas-selections-exist}. The assertion now follows from Corollary~\ref{member-of-support}
\EndExample

Perceiving a stochastic relation $K: (X, {\cal A})\Trans (Y,
\Borel{Y})$ as a probabilistic model for transitions such that
$K(x)(B)$ is the probability for making a transition from $x$ to $B$
(with $K(x)(Y)\leq 1$), we may interpret the selection $f$ as one possible
deterministic version for a transition: the state $f(x)$ is possible,
since $f(x)\in\supp(K(x))$, and it may even be undertaken with
positive probability.
%%% Local Variables: 
%%% mode: latex
%%% TeX-master: "../Mskr3"
%%% End: 

%\Input{\Folder/Integration}
%spell checked - 24Aug14
\Subsection{Integration}
\label{sec:integration}

After having studied the structure of measurable sets under various
conditions on the underlying space with an occasional side glance at
real-valued measurable functions, we will discuss now
integration. This is a fundamental operation associated with
measures. The integral of a function with respect to a measure will
be what you expect it to be, viz., for non-negative functions the area
between the curve and the $x$-axis. This view will be confirmed later
on, when Fubini's Theorem will be available for computing measures in
Cartesian products. For the time being, we build up the integral in a
fairly straightforward way through an approximation by step functions,
obtaining a linear map with some favorable properties, for example the
Lebesgue Dominated Convergence Theorem. All the necessary
constructions are given in this section, offering more than one
occasion to exercise the well-known $\epsilon$-$\delta$-arguments,
which are necessary, but not particularly entertaining. But that's
life.

The second part of this section offers a complementary view --- it
starts from a positive linear map with some additional continuity
properties and develops a measure from it. This is Daniell's approach,
suggesting that measure and integral are really most of the time two
sides of the same coin. We show that this duality comes to life
especially when we are dealing with a compact metric space: Here the
celebrated Riesz Representation Theorem gives a bijection between
probability measures on the Borel sets and positive linear functions
mapping $1$ to $1$ on the continuous real-valued functions. We
formulate and prove this theorem here; it should be mentioned that
this is not the most general version available, as with most other
results discussed here (but probably there is no such thing as a
\emph{most general version}, since the development did branch out into
wildly different directions).

This section will be fundamental for the discussions and results
later in this chapter. Most results are formulated for finite or
$\sigma$-finite measures, and usually no attempt has been made to find
the boundary delineating a development.

\Subsubsection{From Measure to Integral}
\label{sec:step-functions}

We fix a measure space $(X, {\cal A},\mu)$. Denote for the moment by
$\StepFnct{X, \mathcal{A}}$ the set of all measurable step functions,
and by $\StepFnctP{X, \mathcal{A}}$ the non-negative step functions;
similarly, $\MeasbFnctP{X, \mathcal{A}}$ are the non-negative
measurable functions. Note that $\StepFnct{X, \mathcal{A}}$ is a
vector space under the usual operations, and that it is a lattice
under finite or countable pointwise suprema and infima. We know from
Proposition~\ref{ApproxStepFncts} that we can approximate each bounded
measurable function by a sequence of step functions from below.

Define
\begin{equation}
\label{def-int-step-fnct}
  \int_{X}\sum_{i=1}^{n}\alpha_{i}\cdot\chi_{A_{i}}\ d\mu := \sum_{i=1}^{n}\alpha_{i}\cdot\mu(A_{i})
\end{equation}
as the \emph{integral with respect to $\mu$} for the step function $\textstyle{\sum_{i=1}^{n}\alpha_{i}\cdot\chi_{A_{i}}}\in\StepFnct{X, \mathcal{A}}$. Exercise~\ref{step-fnct-val} tells us that the integral is well defined: if $f, g\in\StepFnct{X, \mathcal{A}}$ with $f = g$, then 
\begin{equation*}
  \sum_{\alpha\in\Real}\alpha\cdot\mu(\{x\in X\mid f(x) = \alpha\}) = \sum_{\beta\in\Real}\beta\cdot\mu(\{x\in X\mid g(x) = \beta\}).
\end{equation*}
Thus the definition (\ref{def-int-step-fnct}) yields the same value for the integral. These are some elementary properties of the integral for step functions.

\BeginLemma{prop-int-step-fncts}
Let $f, g\in\StepFnct{X, \mathcal{A}}$ be step functions, $\alpha\in\Real$. Then
\begin{enumerate}
\item\label{item:1} $\int_{X}\alpha\cdot f\ d\mu = \alpha\cdot\int_{X}f\ d\mu$,
\item\label{item:2} $\int_{X}(f+g)\ d\mu = \int_{X}f\ d\mu + \int_{X}g\ d\mu$,
\item\label{item:3} if $f\geq 0$, then $\int_{X}f\ d\mu \geq 0$, in particular, $f\mapsto \int_{X}f\ d\mu$ is monotone,
\item\label{item:4} $\int_{X}\chi_{A}\ d\mu = \mu(A)$ for $A\in{\cal A}$,
\item\label{item:5} $|\int_{X}f\ d\mu| \leq \int_{X}|f|\ d\mu$.
\end{enumerate}
Moreover the map $A\mapsto \int_{A}f\ d\mu := \int_{X}f\cdot\chi_{A}\ d\mu$ is additive on ${\cal A}$ whenever $f\in \StepFnctP{X, \mathcal{A}}$. \QED
\EndLemma

We know from Proposition~\ref{ApproxStepFncts} that we can find for $f\in\MeasbFnctP{X, \mathcal{A}}$ a sequence $\Folge{f}$ in $\StepFnctP{X, \mathcal{A}}$ such that $f_{1}\leq f_{2}\leq \dots $ and $sup_{n\in\Nat}f_{n} = f$. This observation is used for the definition of the integral for $f$. We define
\begin{equation*}
  \int_{X}f\ d\mu := \sup\bigl\{\int_{X}g\ d\mu\mid g\leq f\text{ and } g\in\StepFnctP{X, \mathcal{A}}\bigr\}
\end{equation*}
Note that the right hand side may be infinite; we will discuss this shortly. 

The central observation is formulated in Levi's Theorem:

\BeginTheorem{beppo-levi}
Let $\Folge{f}$ be an increasing sequence of functions in $\MeasbFnctP{X, \mathcal{A}}$ with limit $f$, then the limit $\bigl(\int_{X}f_{n}\ d\mu\bigr)_{n\in\Nat}$ exists and equals $\int_{X}f\ d\mu$. 
\EndTheorem

\BeginProof
1.
Because the integral is monotone in the integrand by Lemma~\ref{prop-int-step-fncts}, the limit
\begin{equation*}
  \ell := \lim_{n\to \infty}\Int{X}{f_{n}}{\mu}
\end{equation*}
exists (possibly in $\Real\cup\{\infty\}$), and we know from monotonicity that  $\ell\leq \dInt{f}{\mu}$.

2. 
Let $f = c > 0$ be a constant, and let $0 < d < c$. Then 
$
\sup_{n\in\Nat}d\cdot\chi_{\{x\in X \mid f_{n}(x) \geq d\}} = d,
$ hence 
we obtain
\begin{equation*}
  \dInt{f}{\mu} \geq \dInt{f_{n}}{\mu} \geq \dInt[\{x\in X \mid f_{n}(x) \geq d\}]{f_{n}}{\mu}\geq d\cdot\mu(\{x\in X \mid f_{n}(x) \geq d\})
\end{equation*}
for every $n\in\Nat$, thus
\begin{equation*}
  \dInt{f}{\mu} \geq d\cdot\mu(X).
\end{equation*}
Letting $d$ approaching $c$, we see that 
\begin{equation*}
  \dInt{f}{\mu} \geq \lim_{n\to \infty}\dInt{f_{n}}{\mu} \geq c\cdot\mu(X) = \dInt{f}{\mu}.
\end{equation*}
This gives the desired equality.

3.
If $f=c\cdot\chi_{A}$ with $A\in{\cal A}$, we restrict the measure space to $(A, {\cal A}\cap A, \mu)$, so the result is true also for step functions based on one single set.

4.
Let $f = \sum_{i=1}^{n}\alpha_{i}\cdot\chi_{A_{i}}$ be a step function, then we may assume that the sets $A_{1}, \dots, A_{n}$ are mutually disjoint. Consider $f_{i} := f\cdot\chi_{A_{i}} = \alpha_{i}\cdot\chi_{A_{i}}$ and apply the previous step to $f_{i}$, taking additivity from Lemma~\ref{prop-int-step-fncts}, part~\ref{item:2}. into account.

5. Now consider the general case. Select step functions $\Folge{g}$ with $g_{n}\in\StepFnctP{X, \mathcal{A}}$ such that $g_{n}\leq f_{n}$ and $|\dInt{f_{n}}{\mu} -\dInt{g_{n}}{\mu}| < 1/n$. We may and do assume that $g_{1}\leq g_{2}\leq \dots $, for we otherwise may pass to the step function $h_{n} := \sup\{g_{1}, \dots, g_{n}\}$. Let $0\leq g\leq f$ be a step function, then $\lim_{n\to \infty} (g_{n}\wedge g) = g$, so that we obtain from the previous step
\begin{align*}
  \dInt{g}{\mu} & = \lim_{n\to \infty}\dInt{g_{n}\wedge g}{\mu}\\
& \leq \lim_{n\to \infty}\dInt{g_{n}}{\mu}\\
& \leq \lim_{n\to \infty}\dInt{f_{n}}{\mu}
\end{align*}
Because $\dInt{g}{\mu}$ may be chosen arbitrarily close to $\ell$, we finally obtain
\begin{equation*}
  \lim_{n\to \infty}\dInt{f_{n}}{\mu}\leq\dInt{f}{\mu}\leq \lim_{n\to \infty}\dInt{f_{n}}{\mu},
\end{equation*}
which implies the assertion for arbitrary $f\in\MeasbFnctP{X, \mathcal{A}}$. 
\EndProof

Since we can approximate each non-negative measurable function from below and from above by step functions (Proposition~\ref{ApproxStepFncts} and Exercise~\ref{ex-meas-approx-from-below}), we obtain from Levi's Theorem for $f\in\MeasbFnctP{X, \mathcal{A}}$ the representation
\begin{equation*}
  \sup\bigl\{\dInt{g}{\mu} \mid \StepFnctP{X, \mathcal{A}}\ni g\leq f\bigr\} 
= \dInt{f}{\mu} 
= \inf\bigl\{\dInt{g}{\mu} \mid  f\leq g\in\StepFnctP{X, \mathcal{A}}\bigr\}.
\end{equation*}
This strongly resembles ---~and generalizes~--- the familiar construction of the Riemann integral for a continuous function $f$ over a bounded interval by sandwiching it between lower and upper sums of step functions.

Compatibility of the integral with scalar
multiplication and with addition is now an easy consequence of Levi's
Theorem:

\BeginCorollary{integral-is-linear}
Let $a\geq0$ and $b\geq0$ be non-negative real numbers, then
\begin{equation*}
  \dInt{a\cdot f + b\cdot g}{\mu} = a\cdot\dInt{f}{\mu} + b\cdot\dInt{g}{\mu}
\end{equation*}
for $f, g\in\MeasbFnctP{X, \mathcal{A}}$. 
\EndCorollary

\BeginProof
Let $\Folge{f}$ and $\Folge{g}$ be sequences of step functions which converge monotonically to $f$ resp. $g$. Then $(a\cdot f_{n} + b\cdot g_{n})_{n\in\Nat}$ is a sequence of step functions converging monotonically to $a\cdot f+b\cdot g$. Apply Levi's Theorem~\ref{beppo-levi} and the linearity of the integral on step functions from Lemma~\ref{prop-int-step-fncts} to obtain the assertion.
\EndProof

Given an arbitrary $f\in\MeasbFnct{X, \mathcal{A}}$, we can decompose $f$ into a positive and a negative part $f^{+} := f\vee 0$ resp. $f^{-} := (-f)\vee 0$, so that $f = f^{+}- f^{-}$ and $|f| = f^{+}+f^{-}$. 

A function $f\in\MeasbFnct{X, \mathcal{A}}$ is called \emph{integrable} (with respect to $\mu$) iff 
\begin{equation*}
\dInt{|f|, \mu}<\infty,
\end{equation*}
in this case we set 
\begin{equation*}
  \dInt{f}{\mu} := \dInt{f^{+}}{\mu} - \dInt{f^{-}}{\mu}.
\end{equation*}

In fact, because $f^{+}\leq f$, we obtain from Lemma~\ref{prop-int-step-fncts} that $\dInt{f^{+}}{\mu} < \infty$, similarly we see that $\dInt{f^{-}}{\mu}<\infty$. The integral is well defined, because if $f = f_{1} - f_{2}$ with $f_{1}, f_{2}\geq 0$, we conclude $f_{1}\leq f \leq |f|$, hence $\dInt{f_{1}}{\mu}<\infty$, and $f_{2}\leq |f|$, so that $\dInt{f_{2}}{\mu}<\infty$, which implies $\dInt{f^{+}}{\mu} + \dInt{f_{2}}{\mu} = \dInt{f^{-}}{\mu} + \dInt{f_{1}}{\mu}$ by Corollary~\ref{integral-is-linear}. This we obtain in fact $\dInt{f^{+}}{\mu} - \dInt{f^{-}}{\mu} = \dInt{f_{1}}{\mu} - \dInt{f_{2}}{\mu}$. 

This special case is also of interest: let $A\in{\cal A}$, define for $f$ integrable 
\begin{equation*}
  \dInt[A]{f}{\mu} := \dInt{f\cdot\chi_{A}}{\mu}
\end{equation*}
(note that $|f\cdot\chi_{A}|\leq |f|$). We emphasize occasionally the integration variable by writing $\dInt{f(x)}{\mu(x)}$ instead of $\dInt{f}{\mu}$. 

Collecting some useful and a.e. used properties, we state

\BeginProposition{collect-some-poperties}
Let $f, g\in\MeasbFnct{X, \mathcal{A}}$ be measurable functions, then
\begin{enumerate}
\item If $f\geq_{\mu}0$, then $\dInt{f}{\mu} = 0$ iff
  $f=_{\mu}0$.
\item If $f$ is integrable, and $|g|\leq_{\mu}|f|$, then $g$ is
  integrable.
\item If $f$ and $g$ are integrable, then so are $a\cdot f + b\cdot g$
  for all $a, b\in\Real$, and $\dInt{a\cdot f+b\cdot g}{\mu} = a\cdot\dInt{f}{\mu} + b\cdot\dInt{g}{\mu}$. 
\item If $f$, and $g$ are integrable, and $f\leq_{\mu}g$, then
  $\dInt{g}{\mu}\leq \dInt{f}{\mu}$.
\item If $f$ is integrable, then $|\dInt{f}{\mu}| \leq
  \dInt{|f|}{\mu}$.
\end{enumerate}
\QED
\EndProposition

We now state and prove some statements which relate sequences of functions to their integrals. The first one is traditionally called Fatou's Lemma. 

\BeginProposition{fatous-lemma}
Let $\Folge{f}$ be a sequence in $\MeasbFnctP{X, \mathcal{A}}$. Then 
\begin{equation*}
  \dInt{\liminf_{n\to \infty} f_{n}}{\mu} \leq \liminf_{n\to \infty}\dInt{f_{n}}{\mu}
\end{equation*}
\EndProposition

\BeginProof
Since $(\inf_{m\geq n}f_{m})_{n\in\Nat}$ is an increasing sequence of
measurable functions in $\MeasbFnctP{X, \mathcal{A}}$, we obtain from
Levi's Theorem~\ref{beppo-levi}
\begin{equation*}
  \dInt{f}{\mu} = \lim_{n\to \infty}\dInt{\inf_{m\geq n}f_{m}}{\mu} = \sup_{n\in\Nat}\dInt{\inf_{m\geq n}f_{m}}{\mu}.
\end{equation*}
Because we plainly have by monotonicity
$
\dInt{\inf_{m\geq n}f_{m}}{\mu}\leq\inf_{m\geq n}\dInt{f_{m}}{\mu},  
$
the assertion follows. 
\EndProof

The \emph{\index{theorem!Lebesgue Dominated Convergence}Lebesgue Dominated Convergence Theorem} is a very important and much used tool; it can be derived now easily from Fatou's Lemma.

\BeginTheorem{lebesgue-dominated-convergence}
Let $\Folge{f}$ be a sequence of measurable functions with 
$f_{n}\aeC f$ for some measurable function $f$, and $|f_{n}|\leq_{\mu} g$ for all $n\in\Nat$ and an integrable function $g$. Then $f_{n}$ and $f$ are integrable, and
\begin{equation*}
  \lim_{n\to \infty}\dInt{f_{n}}{\mu} = \dInt{f}{\mu}\text{ and } \lim_{n\to \infty}\dInt{|f_{n}-f|}{\mu} = 0.
\end{equation*}
\EndTheorem

\BeginProof
1.
It is no loss of generality to assume that $f_{n}\to f$ and $\forall n\in\Nat: f_{n}\leq g$ pointwise (otherwise modify the $f_{n}$, $f$ and $g$ on a set of $\mu$-measure zero). Because $|f_{n}|\leq g$, we conclude from Proposition~\ref{collect-some-poperties} that $f_{n}$ is integrable, and since $f\leq g$ holds as well, we infer that $f$ is integrable as well. 

2.
Put $g_{n} := |f| + g - |f_{n}-f|$, then $g_{n}\geq 0$, and $g_{n}$ is integrable for all $n\in\Nat$. We obtain from Fatou's Lemma
\begin{align*}
  \dInt{|f| + g}{\mu} & = \dInt{\liminf_{n\to \infty} g_{n}}{\mu}\\
& \leq \liminf_{n\to \infty}\dInt{g_{n}}{\mu}\\
& = \dInt{|f|+g}{\mu} - \limsup_{n\to \infty}\dInt{|f_{n}-f|}{\mu}.
\end{align*}
Hence we obtain $\limsup_{n\to \infty}\dInt{|f_{n}-f|}{\mu}=0$, thus $\lim_{n\to \infty}\dInt{|f_{n}-f|}{\mu} = 0$. 

3.
We finally note that 
\begin{equation*}
   \bigl|\dInt{f_{n}}{\mu}-\dInt{f}{\mu}\bigr| = \bigl|\dInt{(f_{n}-f)}{\mu}\bigr| \leq \dInt{|f_{n}-f|}{\mu},
\end{equation*}
which completes the proof.
\EndProof

These are immediate consequences of the Lebesgue Theorem:

\BeginCorollary{summation-is-intbar}
Let $\Folge{f}$ be a sequence of measurable functions, $g$ integrable, such that $|\sum_{k=1}^{n}f_{k}|\leq_{\mu}g$ for all $n\in\Nat$. Then all $f_{n}$ as well as $f := \sum_{n\in\Nat}f_{n}$ are integrable, and $\dInt{f}{\mu} = \sum_{n\in\Nat}\dInt{f_{n}}{\mu}$. \QED
\EndCorollary

\BeginCorollary{integral-is-measure}
Let $f\geq_{\mu}0$ be an integrable function, then $A\mapsto \dInt[A]{f}{\mu}$ defines a finite measure on ${\cal A}$.
\EndCorollary

\BeginProof
All the properties of a measure are immediate, $\sigma$-additivity follows from Corollary~\ref{summation-is-intbar}. 
\EndProof

Integration with respect to an image measure is also available right away. It yields the fairly helpful \emph{\index{change of variables!image measure}change of variables formula} for image measures. 

\BeginCorollary{int-image-measure}
Let $(Y, {\cal B})$ a measurable space and $g: X\to Y$ be ${\cal A}$-${\cal B}$-measurable. Then $h\in\MeasbFnct{Y, \mathcal{B}}$ is $\FinM{g}(\mu)$ integrable iff $g\circ h$ is $\mu$-integrable, and in this case we have
\begin{equation}
\label{change-of-variables}
  \dInt[Y]{h}{\FinM{g}(\mu)} = \dInt{h\circ g}{\mu}.
\end{equation}
\EndCorollary

\BeginProof
We show first that formula~(\ref{change-of-variables}) is true for step functions. In fact, if $h = \chi_{B}$ with a measurable set $B$, then we obtain from the definition
\begin{equation*}
  \dInt[Y]{\chi_{B}}{\FinM{g}(\mu)} = \FinM{g}(\mu)(B) = \mu(\InvBild{g}{B}) = \dInt{\chi_{B}\circ g}{\mu}
\end{equation*}
(since $\chi_{B}(g(x)) = 1 $ iff $x \in\InvBild{g}{B}$). This observation extends by linearity to step functions, so that  we obtain for $h = \sum_{i=1}^{n}b_{i}\cdot\chi_{B_{i}}$
\begin{equation*}
  \dInt[Y]{h}{\FinM{g}(\mu)} = \sum_{i=1}^{n}b_{i}\cdot\dInt{\chi_{B_{i}}\circ g}{\mu} = \dInt{h\circ g}{\mu}
\end{equation*}
Thus the assertion now follows from Levi's Theorem~\ref{beppo-levi}. 
\EndProof

The reader is probably familiar with the change of variables formula in
classical calculus. It deals with $k$-dimensional Lebesgue measure
$\lambda^{k}$, and a differentiable and injective map $T: V\to W$ from
an open set $V\subseteq\Real^{k}$ to a bounded set
$W\subseteq\Real^{k}$. T is assumed to have a continuous inverse. Then the integral of a
measurable and bounded function $f: \Bild{T}{V}\to \Real$ can be expressed in terms of the integral over $V$ of $f\circ T$ and the Jacobian $J_{T}$ of $T$. To be specific\index{change of variables!calculus}\label{change of variables_calculus}
\begin{equation*}
  \dInt[\Bild{T}{V}]{f}{\lambda^{k}} = \dInt[V]{(f\circ T)\cdot |J_{T}|}{\lambda^{k}}.
\end{equation*}
Recall that the Jacobian $J_{T}$ of $T$ is the determinant of the partial derivatives of $T$, i.e., 
\begin{equation*}
  J_{T}(x) = \mathsf{det}\bigl((\frac{\partial T_{i}(x)}{\partial x_{j}})\bigr).
\end{equation*}
This representation can be derived from the representation for the
integral with respect to the image measure from
Corollary~\ref{int-image-measure} and from the Radon-Nikodym
Theorem~\ref{radon-nikodym} through a somewhat lengthy application of
results from fairly elementary linear algebra. We do not want to
develop this apparatus in the present presentation, we will, however,
provide a glimpse at the one-dimensional situation in
Proposition~\ref{integration-by-substitution}. The reader is referred
rather to Rudin's exposition~\cite[p. 181 - 188]{Rudin} or to
Stromberg's more elementary discussion in~\cite[p. 385 -
392]{Stromberg}; if you read German, Elstrodt's derivation~\cite[§
V.4]{Elstrodt} should not be missed.

\Subsubsection{The Daniell Integral and Riesz's Representation Theorem}
\label{sec:daniell-integral}

The previous section developed the integral from a finite or
$\sigma$-finite measure; the result was a linear functional on a
subspace of measurable functions, which will be investigated in
greater detail later on. This section will demonstrate that it is
possible to obtain a measure from a linear functional on a well
behaved space of functions. This approach was proposed by
P.~J. Daniell ca. 1920, it is called in his honor the
\emph{Daniell integral}. It is useful when a linear functional is
given, and one wants to show that this functional is actually defined
by a measure, which then permits putting the machinery of
measure theory into action. We will encounter such a situation, e.g., when
studying linear functionals on spaces of integrable
functions. Specifically, we derive the Riesz Representation Theorem,
which shows that there is a one-to-one correspondence between
probability measures and normed positive linear functionals on the vector
lattice of continuous real valued functions on a compact metric space.

Let us fix a set $X$ throughout. We will also fix a set ${\cal F}$ of functions $X\to \Real$ which is assumed to be a vector space (as always, over the reals) with a special property.

\BeginDefinition{vector-lattice}
A vector space ${\cal F} \subseteq\Real^{X}$ is called a \emph{\index{vector lattice}vector lattice} iff $|f|\in{\cal F}$ whenever $f\in{\cal F}$.  
\EndDefinition

Now fix the vector lattice ${\cal F}$. Each vector lattice is indeed a lattice: define 
\begin{align*}
  f\vee g & := (|f - g| + f + g)/2,\\
f\wedge g & := -\bigl((-f)\vee(-g)\bigr)\\
f \leq g & \Leftrightarrow f\vee g = g\\
& \Leftrightarrow f\wedge g = f
\end{align*}
Thus ${\cal F}$ contains with $f$ and $g$ also $f\wedge g$ and $f\vee g$, and it is easy to see that $\leq$ defines a partial order on ${\cal F}$ such that $\sup\{f, g\} = f\vee g$ and $\inf\{f, g\} = f\wedge g$, see, e.g.,~\SetCite{2.5.5}. Note that we have $\max\{\alpha, \beta\} = (|\alpha - \beta| + \alpha+ \beta)/2$ for $\alpha, \beta\in\Real$, thus we conclude that $f\leq g$ iff $f(x) \leq g(x)$ for all $x\in\Real$. 

We will find these properties helpful; they will be used silently below.

\BeginLemma{helpful-vector-lattice}
If $0\leq\alpha\leq \beta\in\Real$ and $f\in{\cal F}$ with $f\geq 0$, then $\alpha\cdot f \leq \beta\cdot f$. If $f, g\in{\cal F}$ with $f\leq g$, then $f+h\leq g+h$ for all $h\in{\cal F}$. Also, $f\wedge g + f\vee g = f+g$. 
\EndLemma

\BeginProof
Because $f\geq 0$, we obtain
\begin{equation*}
2\cdot\bigl((\alpha\cdot f)\vee(\beta\cdot f)\bigr) = (|\alpha-\beta|+\alpha+\beta)\cdot f = 2\cdot\alpha\vee\beta\cdot f = 2\cdot\beta\cdot f.
\end{equation*}
This establishes the first claim. The second one follows from 
\begin{equation*}
2\cdot\bigl((f+h)\vee(g+h)\bigr) = |f-g|+f+g+2\cdot h = 2\cdot(g + h).
\end{equation*}
The third one is estabished through the observation that it holds pointwise, and from the observation that $f\leq g$ iff $f(x) \leq g(x)$ for all $x\in X$. 
\EndProof

We assume that $1\in {\cal F}$, and that a function $L: {\cal F}\to \Real$ is given, which has these properties:
\begin{itemize}
\item $L(\alpha\cdot f + \beta\cdot g) = \alpha\cdot L(f) + \beta\cdot L(g)$, so that $L$ is linear,
\item if $f\geq 0$, then $L(f)\geq 0$, so that $L$ is positive,
\item $L(1) = 1$, so that $L$ is normed,
\item If $\Folge{f}$ is a sequence in ${\cal F}$ which decreases to $0$, then $\lim_{n\to \infty}L(f_{n}) = 0$, so that $L$ is continuous from above at $0$. 
\end{itemize}

These are some immediate consequences from the properties of $L$.

\BeginLemma{first-props-daniell}
If $f, g\in{\cal F}$, then $L(f\wedge g) + L(f\vee g) = L(f)+L(g)$. If $\Folge{f}$ and $\Folge{g}$ are increasing sequences of non-negative functions in ${\cal F}$ with $\lim_{n\to \infty}f_{n} \leq \lim_{n\to \infty}g_{n}$, then $\lim_{n\to \infty}L(f_{n}) \leq \lim_{n\to \infty}L(g_{n}).$
\EndLemma

\BeginProof
The first property follows from the linearity of $L$. For the second one, we observe that $\lim_{k\to \infty}(f_{n}\wedge g_{k}) = f_{n}\in{\cal F}$, the latter sequence being increasing. Consequently, we have 
\begin{equation*}
  L(f_{n}) \leq \lim_{k\to \infty}L(f_{n}\wedge g_{k}) \leq \lim_{k\to \infty}L(g_{k})
\end{equation*}
for all $n\in\Nat$, which implies the assertion. 
\EndProof

${\cal F}$ determines a $\sigma$-algebra ${\cal A}$ on $X$, viz., the smallest $\sigma$-algebra which renders each $f\in{\cal F}$ measurable. We will show now that $L$ determines a unique probability measure on ${\cal A}$ such that 
\begin{equation*}
  L(f) = \dInt{f}{\mu}
\end{equation*}
holds for all $f\in{\cal F}$. 

This will be done in a sequence of steps. A brief outline looks like
this: We will first show that $L$ can be extended to the set ${\cal
  L}^{+}$ of all bounded monotone limits from the non-negative
elements of ${\cal F}$, and that the extension respects monotone
limits. From ${\cal L}^{+}$ we extract via indicator functions 
an algebra of sets, and from the extension to $L$ an outer
measure. This will then turn out to yield the desired probability.

Define
\begin{equation*}
  {\cal L}^{+} := \{f: X\to \Real\mid f\text{ is bounded, there exists } 0\leq f_{n}\in{\cal F}\text{ increasing with }f = \lim_{n\to \infty} f_{n}\}.
\end{equation*}
Define $L(f) := \lim_{n\to \infty}L(f_{n})$ for $f\in{\cal L}^{+}$, whenever $f = \lim_{n\to \infty}f_{n}$ with the increasing sequence $\Folge{f}\subseteq{\cal F}$. Then we obtain from Lemma~\ref{first-props-daniell} that this extension $L$ on ${\cal L}^{+}$ is well defined, and it is clear that $L(f)\geq 0$, and that  $L(\alpha\cdot f + \beta\cdot g) = \alpha\cdot L(f) + \beta\cdot L(g)$, whenever $f, g\in{\cal L}^{+}$ and $\alpha, \beta\in\pReal$. We see also that $f, g\in{\cal L}^{+}$ implies that $f\wedge g, f\vee g\in {\cal L}^{+}$ with $L(f\wedge g) + L(f\vee g) = L(f) + L(g)$.  It turns out that $L$ also respects the limits of increasing sequences.

\BeginLemma{daniell-cont-above}
Let $\Folge{f}\subseteq{\cal L}^{+}$ be an increasing and uniformly bounded sequence, then $L(\lim_{n\to \infty}f_{n}) = \lim_{n\to \infty}L(f_{n})$. 
\EndLemma

\BeginProof
Because $f_{n}\in{\cal L}^{+}$, we know that there exists for each $n\in\Nat$ an increasing sequence $(f_{m, n})_{m\in\Nat}$ of elements $f_{m, n}\in{\cal F}$ such that $f_{n} = \lim_{m\to \infty}f_{m, n}$. Define
\begin{equation*}
  g_{m} := \sup_{n\leq m}f_{m, n}.
\end{equation*}
Then $\Folge[m]{g}$ is an increasing sequence in ${\cal F}$ with $f_{m,
  n}\leq g_{m}$, and $g_{m}\leq f_{1}\vee f_{2}\vee \dots\vee f_{m} =
f_{m}$, so that $g_{m}$ is sandwiched between $f_{m, n}$ and $f_{m}$
for all $m\in\Nat$ and $n\leq m$. This yields $L(f_{m, n})\leq L(g_{m}) \leq L(f_{m})$ for these $n, m$. Thus 
$
\lim_{n\to \infty}f_{n} = \lim_{m\to \infty}g_{m},
$
and hence
\begin{equation*}
  \lim_{n\to \infty}L(f_{n}) = \lim_{m\to \infty}L(g_{m}) = L(\lim_{m\to \infty}g_{m}) = L(\lim_{n\to \infty}f_{n}). 
\end{equation*}
Thus we have shown that $\lim_{n\to \infty}f_{n}$ can be obtained as the limit of an increasing sequence of functions from ${\cal F}$; because $\Folge{f}$ is uniformly bounded, this limit is an element of ${\cal L}^{+}$.  
\EndProof

Now define
\begin{align*}
  {\cal G} & := \{G\subseteq X \mid  \chi_{G}\in {\cal L}^{+}\},\\
\mu(G) & := L(\chi_{G})\text{ for }G\in{\cal G}.
\end{align*}
Then ${\cal G}$ is closed under finite intersections and finite unions
by the remarks made before Lemma~\ref{daniell-cont-above}. Moreover, 
${\cal G}$ is closed under countable unions with
$\mu(\bigcup_{n\in\Nat}G_{n}) = \lim_{n\to \infty}\mu(G_{n})$, if
$\Folge{G}$ is an increasing sequence in ${\cal G}$. Also $\mu(X) =
1$. Now define, as in the Carathéodory process, see~\SetCite{2.6.3}
\begin{align*}
  \mu^{*}(A) & := \inf\{\mu(G) \mid G\in{\cal G}, A\subseteq G\},\\
{\cal B} & := \{B\subseteq X \mid \mu^{*}(B) + \mu^{*}(X\setminus B) = 1\}.
\end{align*}
 
We obtain from the Carathéodory extension process

\BeginProposition{b-does-it}
${\cal B}$ is a $\sigma$-algebra, and $\mu^{*}$ is countably additive on ${\cal B}$.
\EndProposition

\BeginProof
\SetCite{Proposition 1.127}
\EndProof

Put $\mu(B) := \mu^{*}(B)$ for $B\in{\cal B}$, then $(X, {\cal B}, \mu)$ is a measure space, and $\mu$ is a probability measure on $(X, {\cal B})$.

In order to carry out the programme sketched above, we need a
$\sigma$-algebra. We have on one hand the $\sigma$-algebra ${\cal A}$
generated by ${\cal F}$, and on the other hand ${\cal B}$, gleaned
from the Carathéodory extension. It is not immediately clear how these
$\sigma$-algebras are related to each other. And then we also have ${\cal
  G}$ as an intermediate family of sets, obtained from ${\cal
  L}^{+}$. This diagram shows the objects we will to discuss, together with a short hand indication of the respective relationships:
\begin{equation*}
\xymatrix{
{\cal F}\ar[d]\ar[rrrr]^{\sigma(\cdot )}&&&&{\cal A}\ar@{-}[d]^{=}
\\
{\cal L}^{+}\ar[rr]_{\chi}
&&{\cal G}\ar[rr]^{\subseteq}_{\mathrm{Carath\acute{e}odory}}\ar[urr]^{\sigma(\cdot ) =}
&&{\cal B}
}
\end{equation*}
We investigate the relationship of ${\cal A}$ and ${\cal
  G}$ first.

\BeginLemma{sigma-g-is-a}
${\cal A} = \sigma({\cal G})$.
\EndLemma

\BeginProof
1.
Because ${\cal A}$ is the smallest $\sigma$-algebra rendering all elements of ${\cal F}$ measurable, and because each element of ${\cal L}^{+}$ is the limit of a sequence of elements of ${\cal F}$, we obtain ${\cal A}$-measurability for each element of ${\cal L}^{+}$. Thus ${\cal G}\subseteq {\cal A}$. 

2.
Let $f\in{\cal L}^{+}$ and $c\in\pReal$, then $f_{n} := 1\wedge n\cdot\sup\{f-c, 0\}\in{\cal L}^{+}$, and
$
\chi_{\{x\in X\mid f(x) > c\}} = \lim_{n\to \infty}f_{n}.
$
This is a monotone limit. Hence $\{x\in X\mid f(x) > c\}\in {\cal G}$, thus in particular each element of ${\cal F}$ is $\sigma({\cal G})$-measurable. This implies that ${\cal A}\subseteq\sigma({\cal G})$ holds. 
\EndProof

The relationship between ${\cal B}$ and ${\cal G}$ is a bit more difficult to establish.

\BeginLemma{g-subset-b}
${\cal G}\subseteq{\cal B}$.
\EndLemma

\BeginProof
We have to show that $\mu^{*}(G) + \mu^{*}(X\setminus G) = 1$ for all  $G\in{\cal G}$. Fix $G\in{\cal G}$. We obtain from additivity that $\mu(G) + \mu(H) = \mu(G\cap H) + \mu(G\cup H) \geq \mu(X) = 1$ holds for any $H\in{\cal G}$ with $X\setminus G\subseteq H$, so that $\mu^{*}(G) + \mu^{*}(X\setminus G) \leq 1$ remains to be shown. 

Because $G\in{\cal G}$, there exists an increasing sequence $\Folge{f}$ of elements in ${\cal F}$ such that $\chi_{G}=\sup_{n\in\Nat}f_{n}$, consequently, $\chi_{X\setminus G} = \inf_{n\in\Nat}(1-f_{n})$. Now let $n\in\Nat$, and $0 < c \leq 1$, then 
$
X\setminus G \subseteq U_{n, c} := \{x\in X \mid 1-f_{n}(x) > c\}$ with $U_{n, c}\in{\cal G}$. Because $\chi_{U_{n, c}}\leq (1-f_{n})/c$, we obtain $\mu^{*}(X\setminus G) \leq L(1-f_{n})/c$; this inequality holds for all $c$ and all $n\in\Nat$. Letting $c\to 1$ and $n\to \infty$, this yields $\mu^{*}(X\setminus G) \leq 1 - \mu^{*}(G)$.  

Consequently, $\mu^{*}(G)+\mu^{*}(X\setminus G) = 1$ for all $G\in{\cal G}$, which establishes the claim. 
\EndProof
 
This yields the desired relationship of ${\cal A}$, the $\sigma$-algebra generated by the functions in ${\cal F}$, and ${\cal B}$, the $\sigma$-algebra obtained from the extension process.

\BeginCorollary{each-l-is-measb}
${\cal A} \subseteq {\cal B}$, and each element of ${\cal L}^{+}$ is ${\cal B}$-measurable.
\EndCorollary

\BeginProof
We have seen that ${\cal A} = \sigma({\cal G})$ and that ${\cal G}\subseteq {\cal B}$, so the first assertion follows from Proposition~\ref{b-does-it}. The second assertion is immediate from the first one.
\EndProof

Because $\mu$ is countably additive, hence a probability measure on ${\cal B}$, and because each element of ${\cal F}$ is ${\cal B}$-measurable, the integral $\dInt{f}{\mu}$ is defined, and we are done.

\BeginTheorem{daniell-integration}
Let ${\cal F}$ be a vector lattice of functions $X\to \Real$ with $1\in{\cal F}$, $L: {\cal F}\to \Real$ be a linear and monotone functional on ${\cal F}$ such that $L(1) = 1$, and $L(f_{n})\to 0$, whenever $\Folge{f}\subseteq{\cal F}$ decreases to $0$. Then there exists a unique probability measure $\mu$ on the $\sigma$-algebra ${\cal A}$ generated by ${\cal F}$ such that
\begin{equation*}
  L(f) = \dInt{f}{\mu}
\end{equation*}
holds for all $f\in{\cal F}$.
\EndTheorem

\BeginProof
Let ${\cal G}$ and ${\cal B}$ be constructed as above.

\emph{Existence}: Because ${\cal A}\subseteq{\cal B}$, we may restrict $\mu$ to ${\cal A}$, obtaining a probability measure. Fix $f\in{\cal F}$, then $f$ is ${\cal B}$-measurable, hence $\dInt{f}{\mu}$ is defined.  Assume first that $0\leq f \leq 1$, hence $f\in{\cal L}^{+}$. We can write $f = \lim_{n\to \infty}f_{n}$ with step functions $f_{n}$, the contributing sets being members of ${\cal G}$. Hence $L(f_{n}) = \dInt{f_{n}}{\mu}$, since $L(\chi_{G}) = \mu(G)$ by construction. Consequently, we obtain from Lemma~\ref{daniell-cont-above} and Lebesgue's Dominated Convergence Theorem~\ref{lebesgue-dominated-convergence}
\begin{equation*}
  L(f) = L(\lim_{n\to \infty}f_{n}) = \lim_{n\to \infty}L(f_{n}) = \lim_{n\to \infty}\dInt{f_{n}}{\mu} = \dInt{\lim_{n\to \infty}f_{n}}{\mu} = \dInt{f}{\mu}.
\end{equation*}
This implies the assertion also for bounded $f\in{\cal F}$ with $f\geq 0$. If $0\leq f$ is unbounded, write $f = \sup_{n\in\Nat}(f\wedge n)$ and apply Levi's Theorem~\ref{beppo-levi}. In the general case, decompose $f = f^{+}-f^{-}$ with $f^{+} := f\vee 0$ and $f^{-} := (-f)\vee 0$, and apply the foregoing. 

\emph{Uniqueness}: Assume that there exists a probability measure $\nu$ on ${\cal A}$ with $L(f) = \dInt{f}{\nu}$ for all $f\in{\cal F}$, then the construction shows that $\mu(G) = L(\chi_{G}) = \nu(G)$ for all $G\in{\cal G}$. Since ${\cal G}$ is closed under finite intersections, and since ${\cal A} =\sigma({\cal G})$, we conclude that $\nu(A) = \mu(A)$ for all $A\in{\cal A}$. 
\EndProof

We obtain as a consequence the famous \emph{\index{theorem!Riesz
    Representation}Riesz Representation Theorem}, which we state and
formulate for the metric case. Recall that ${\cal C}(X)$ is the linear
space of all continuous functions $X\to \Real$ on a topological $X$, and ${\cal C}_{b}(X)$ is the subspace of all bounded functions. We state the result first for metric spaces and for bounded continuous functions, specializing the result then to the compact metric case.

\BeginCorollary{pre-riesz-repr}
Let $X$ be a metric space, and let $L: {\cal C}_{b}(X)\to \Real$ be a positive linear function with $\lim_{n\to \infty}L(f_{n}) =0$ for each sequence $\Folge{f}\subseteq{\cal C}_{b}(X)$ which decreases monotonically to $0$. Then there exists a unique finite Borel measure $\mu$ such that 
\begin{equation*}
  L(f) = \dInt{f}{\mu}
\end{equation*}
holds for all $f\in{\cal C}_{b}(X)$. 
\EndCorollary

\BeginProof
It is clear that ${\cal C}_{b}(X)$ is a vector lattice with $1\in{\cal C}_{b}(X)$. We may and do assume that $L(1) = 1$. The result follows immediately from Theorem~\ref{daniell-integration} now.
\EndProof

If we take a compact metric space, then each continuous map $X\to \Real$ is bounded. We show that the assumption on $L$'s continuity follows from compactness (which is usually referred to as \emph{Dini's Theorem}). 

\BeginTheorem{riesz-representation}
Let $X$ be a compact metric space. Given a positive linear functional $L: {\cal C}(X)\to \Real$, there exists a unique finite Borel measure $\mu$ such that 
\begin{equation*}
  L(f) = \dInt{f}{\mu}
\end{equation*}
holds for all $f\in{\cal C}(X)$. 
\EndTheorem

\BeginProof
It is clear that ${\cal C}(X)$ is a vector lattice which contains $1$. Again, we assume that $L(1) = 1$ holds. In order to apply Theorem~\ref{daniell-integration}, we have to show that $\lim_{n\to \infty}L(f_{n})= 0$, whenever $\Folge{f}\subseteq{\cal C}(X)$ decreases monotonically to $0$. But since $X$ is compact, we claim that $\sup_{x\in X}f_{n}(x)\to 0$, as $n\to \infty$. This is so because $\{x\in X\mid f_{n}\geq c\}$ is a family of closed sets with empty intersection for any $c > 0$, so we find by compactness a finite subfamily with empty intersection. Hence the assumption that $\sup_{x\in X}f_{n}(x) \geq c > 0$ for all $n\in\Nat$ would lead to a contradiction. Thus the assertion follows from Theorem~\ref{daniell-integration}. 
\EndProof

Because $f\mapsto \dInt{f}{\mu}$ defines for each Borel measure $\mu$
a positive linear functional on ${\cal C}(X)$, and because a measure
on a metric space is uniquely determined by its integral on the
bounded continuous functions, we obtain:

\BeginCorollary{unique-riesz-corr}
For a compact metric space $X$ there is a bijection between positive
linear functionals on ${\cal C}(X)$ and finite Borel measures. \QED
\EndCorollary

The reason for not formulating the Riesz Representation Theorem
immediately for general topological spaces is that
Theorem~\ref{daniell-integration} works with the $\sigma$-algebra
generated ---~in this case~--- by ${\cal C}(X)$; this is in general
the $\sigma$-algebra of Baire sets, which in turn may be properly contained in
the Borel sets. Thus one obtains in the general case a Baire measure
which then would have to be extended uniquely to a Borel measure. This
is discussed in detail in\cite[Sec. 7.3]{Bogachev}.

A typical scenario for the application of the Riesz Theorem runs like
this: one starts with a probability measure on a metric space
$X$. This space can be embedded into a compact metric space $X'$, say,
and one knows that the integral on the bounded continuous functions on
$X$ extends to a positive linear map on the continuous functions on
$X'$. Then the Riesz Representation Theorem kicks in and gives a
probability measure on $X'$. We will see a situation like this when
investigating the weak topology on the space of all finite measures on
a Polish space in Section~\ref{sec:weak-top}.
%%% Local Variables:  
%%% mode: latex
%%% TeX-master: "../Mskr3"
%%% End: 

%\Input{\Folder/ProductMeasures}
%spell checked - 24Aug14
\Subsection{Product Measures}
\label{sec:prod-measures}

As a first application of integration we show that the product of two
finite measures yields a measure again. This will lead to the Fubini's
Theorem on product integration, which evaluates a product integrable
function on a product along its vertical or its horizontal cuts (in
this sense it may be compared to a line sweeping algorithm --- you
traverse the Cartesian product, and in each instance you measure the
cut). We apply this then to infinite products, first with a countable
index set, then for an arbitrary one. Infinite products are a special
case of projective systems, which may be described as sequences of
probabilities which are related through projections. We show that such
a projective system has a projective limit,i.e., a measure on the set
of all sequences such that the projective system proper is obtained
through a projection. This construction is, however, only feasible in
a Polish space, since here a compactness argument is available which
ascertains that the measure we are looking for is $\sigma$-additive. A
small step leads to projective limits for stochastic relations. We
demonstrate an application for projective limits through the
interpretation for the logic \texttt{CSL}.

Fix for the time being two finite measure spaces $(X, {\cal A}, \mu)$
and $(Y, {\cal B}, \nu)$. The Cartesian product $X\times Y$ is endowed
with the product $\sigma$-algebra ${\cal A}\otimes{\cal B}$ which is
the smallest $\sigma$-algebra containing all measurable rectangles
$A\times B$ with $A\in{\cal A}$ and $B\in{\cal B}$, see
Section~\ref{sec:measurable-sets}.

Recall that for $Q\subseteq X\times Y$ the cuts $Q_{x} := \{y\in Y\mid
\langle x, y\rangle \in Q\}$ and $Q^{y} := \{x\in X \mid \langle x,
y\rangle \in Q\}$ are defined. It is clear that $Q_{x}\in {\cal B}$
and $Q^{y}\in {\cal A}$ holds for $x\in X, y\in Y$, whenever
$Q\in{\cal A}\otimes{\cal B}$. In fact, take for example the vertical cut $Q_{x}$ and consider the set 
\begin{equation*}
  {\cal Q} := \{Q\in {\cal A}\otimes{\cal B}\mid Q_{x}\in{\cal B}\}.
\end{equation*}
Then $A\times B\in {\cal Q}$, whenever $A\in{\cal A}, B\in{\cal
  B}$;this is so since the set of all measurable rectangles forms a
generator for the product $\sigma$-algebra which is closed under
finite intersections. Because $(X\times Y)\setminus Q)_{x} =
Y\setminus Q_{x}$, we infer that ${\cal Q}$ is closed under
complementation, and because
$(\bigcup_{n\in\Nat}Q_{n})_{x}=\bigcup_{n\in\Nat}Q_{n, x}$, we
conclude that ${\cal Q}$ is closed under disjoint countable
unions. Hence ${\cal Q} = {\cal A}\otimes{\cal B}$ by the
$\pi$-$\lambda$-Theorem~\ref{Pi-Lambda}.

\BeginLemma{cuta-give-equal}
Let $Q\in{\cal A}\otimes{\cal B}$ be a measurable set, then both 
$
\phi(x) := \nu(Q_{x})
$
and
$
\psi(y) := \mu(Q^{y})
$
define bounded measurable functions with
\begin{equation*}
  \dInt{\nu(Q_{x})}{\mu(x)} = \dInt[Y]{\mu(Q^{y})}{\nu(y)}.
\end{equation*}
\EndLemma

\BeginProof
We use the same argument as above to establish that both $\phi$ and
$\psi$ are measurable functions, noting that $\nu((A\times B)_{x}) =
\chi_{A}(x)\cdot\nu(B)$, similarly, $\mu((A\times B)^{y}) =
\chi_{B}(y)\cdot\mu(A)$; in the next step the set of all $Q\in{\cal
  A}\otimes{\cal B}$ is shown to satisfy the assumptions of the
$\pi$-$\lambda$-Theorem~\ref{Pi-Lambda}. 

In the same way, the equality
of the integrals is established, noting that
\begin{equation*}
  \dInt{\nu((A\times B)_{x})}{\mu(x)} = \mu(A)\cdot\nu(B) = \dInt{\mu((A\times B)^{y})}{\nu(y)}.
\end{equation*}
\EndProof

This yields without much ado

\BeginTheorem{product-measure-exists}
Given the finite measure spaces $(X, {\cal A}, \mu)$ and $(Y, {\cal B}, \nu)$, there exists a unique finite measure $\mu\otimes\nu$ on ${\cal A}\otimes{\cal B}$ such that $(\mu\otimes\nu)(A\times B) = \mu(A)\cdot\nu(B)$ for $A\in{\cal A}, B\in{\cal B}$. Moreover,
\begin{equation*}
  (\mu\otimes\nu)(Q) = \dInt{\nu(Q_{x})}{\mu(x)} = \dInt[Y]{\mu(Q^{y})}{\nu(y)}
\end{equation*}
holds for all $Q\in{\cal A}\otimes{\cal B}$. 
\EndTheorem

\BeginProof
1.
We establish the existence of $\mu\otimes\nu$ by an appeal to Lemma~\ref{cuta-give-equal} and to the properties of the integral according to Proposition~\ref{collect-some-poperties}. Define 
\begin{equation*}
  (\mu\otimes\nu)(Q) := \dInt{\nu(Q_{x})}{\mu(x)},
\end{equation*}
then this defines a finite measure on ${\cal A}\otimes{\cal B}$.
\begin{itemize}
\item Let $Q\subseteq Q'$, then $Q_{x}\subseteq Q'_{x}$ for all $x\in X$, hence $\dInt{\nu(Q_{x})}{\mu(x)} \leq \dInt{\nu(Q'_{x})}{\mu(x)}$. Thus $\mu\otimes\nu$ is monotone.
\item If $Q$ and $Q'$ are disjoint, then $Q_{x}\cap Q'_{x} = (Q\cap Q')_{x} =\emptyset$ for all $x\in X$. Thus $\mu\otimes\nu$ is additive.
\item Let $\Folge{Q}$ be a sequence of disjoint measurable sets, then $(Q_{n, x})_{n\in\Nat}$ is disjoint for all $x\in X$, and \begin{equation*}
\dInt{\nu(\bigcup_{n\in\Nat}Q_{n, x})}{\mu(x)} = \dInt{\sum_{n\in\Nat}\nu(Q_{n, x})}{\mu(x)} = \sum_{n\in\Nat}\dInt{\nu(Q_{n, x})}{\mu(x)}
\end{equation*}
by Corollary~\ref{summation-is-intbar}. Thus $\mu\otimes\nu$ is $\sigma$-additive.
\end{itemize}

2.
Suppose that $\rho$ is a finite measure on ${\cal A}\otimes{\cal B}$ with $\rho(A\times B) = \mu(A)\cdot\nu(B)$ for all $A\in{\cal A}$ and all $B\in{\cal B}$. Then
\begin{equation*}
  {\cal G} := \{Q\in {\cal A}\otimes{\cal B}\mid \rho(Q) = (\mu\otimes\nu)(Q)\}
\end{equation*}
contains the generator $\{A\times B\mid A\in{\cal A}, B\in{\cal B}\}$ of ${\cal A}\otimes{\cal B}$, which is closed under finite intersections. Because both $\rho$ and $\mu\otimes\nu$ are measures, ${\cal G}$ is closed under countable disjoint unions, because both contenders are finite, ${\cal G}$ is also closed under complementation. The $\pi$-$\lambda$-Theorem~\ref{Pi-Lambda} shows that ${\cal G} = {\cal A}\otimes{\cal B}$. Thus $\mu\otimes\nu$ is uniquely determined. 
\EndProof

Theorem~\ref{product-measure-exists} holds also for $\sigma$-finite measures. In fact, assume that the contributing measure spaces are $\sigma$-finite, and let $\Folge{X}$ resp. $\Folge{Y}$ be increasing sequences in ${\cal A}$ resp. ${\cal B}$ such that $\mu(X_{n})<\infty$ and $\nu(Y_{n})<\infty$ for all $n\in\Nat$, and $\bigcup_{n\in\Nat}X_{n} = X$ and $\bigcup_{n\in\Nat}Y_{n}=Y$. Localize $\mu$ and $\nu$ to $X_{n}$ resp. $Y_{n}$ by defining $\mu_{n}(A) := \mu(A\cap X_{n})$, similarly, $\nu_{n}(B) := \nu(B\cap Y_{n})$; since these measures are finite, we can extend them uniquely to a measure $\mu_{n}\otimes\nu_{n}$ on ${\cal A}\otimes{\cal B}$. Since $\bigcup_{n\in\Nat}X_{n}\times Y_{n} = X\times Y$ with the increasing sequence $(X_{n}\times Y_{n})_{n\in\Nat}$, we set 
\begin{equation*}
  (\mu\otimes\nu)(Q) := \sup_{n\in\Nat}(\mu_{n}\otimes\nu_{n})(Q).
\end{equation*}
Then $\mu\otimes\nu$ is a $\sigma$-finite measure on ${\cal A}\otimes{\cal B}$. Now assume that we have another $\sigma$-finite measure $\rho$ on ${\cal A}\otimes{\cal B}$ with $\rho(A\times B) = \mu(A)\cdot\nu(B)$ for all $A\in{\cal A}$ and $B\in{\cal B}$. Define $\rho_{n}(Q) := \rho(Q\cap (X_{n}\times Y_{n}))$, hence $\rho_{n} = \mu_{n}\otimes\nu_{n}$ by uniqueness of the extension to $\mu_{n}$ and $\nu_{n}$, so that we obtain
\begin{equation*}
  \rho(Q) = \sup_{n\in\Nat}\rho_{n}(Q) = \sup_{n\in\Nat}(\mu_{n}\otimes\nu_{n})(Q) = (\mu\otimes\nu)(Q)
\end{equation*}
for all $Q\in {\cal A}\otimes{\cal B}$. Thus we have shown

\BeginCorollary{sigma-finite-product-measure-exists}
Given two $\sigma$-finite measure spaces $(X, {\cal A}, \mu)$ and $(Y, {\cal B}, \nu)$, there exists a unique $\sigma$-finite measure $\mu\otimes\nu$ on ${\cal A}\otimes{\cal B}$ such that $(\mu\otimes\nu)(A\times B) = \mu(A)\cdot\nu(B)$. We have
\begin{equation*}
  (\mu\otimes\nu)(Q) = \dInt{\nu(Q_{x})}{\mu(x)} = \dInt[Y]{\mu(Q^{y})}{\nu(y)}
\end{equation*}
\QED
\EndCorollary

The construction of the product measure has been done here through integration of cuts. An alternative would have been the canonical approach. This approach would have investigated the map $\langle A, B\rangle \mapsto \mu(A)\cdot \nu(B)$ on the set of all rectangles, and then put the extension machinery developed through the Carathéodory approach into action. It is a matter of taste which approach to prefer. --- 

The following example displays a slight generalization (a finite measure is but a constant transition kernel). 

\BeginExample{measure-times-kernel}
Let $K: (X, {\cal A})\Trans (Y, {\cal B})$ a transition kernel (see
Definition~\ref{transition-kernels}) such that the map $x\mapsto
K(x)(Y)$ is integrable with respect to the finite measure $\mu$. Then
\begin{equation*}
  (\mu\otimes K)(Q) := \dInt{K(x)(Q_{x})}{\mu(x)}
\end{equation*}
defines a finite measure on $(X\times Y, {\cal A}\otimes{\cal
  B})$. The $\pi$-$\lambda$-Theorem~\ref{Pi-Lambda} tell us that this
measure is uniquely determined by the condition $(\mu\otimes
K)(A\times B) = \dInt[A]{K(x)(B)}{\mu(x)}$ for $A\in{\cal A},
B\in{\cal B}$.

Interpret in a probabilistic setting $K(x)(B)$ as the
probability that an input $x\in X$ yields an output in $B\in{\cal B}$,
and assume that $\mu$ gives the initial probability with which the
system starts, then $\mu\otimes K$ gives the probability of all
pairings, i.e., $(\mu\otimes K)(Q)$ is the probability that a pair
$\langle x, y\rangle$ consisting of an input value $x\in X$ and an
output value $y\in Y$ will be a member of $Q\in{\cal A}\otimes{\cal
  B}$.
\EndExample

This may be further extended, replacing the measure on $K$'s domain by a transition kernel as well.

\BeginExample{kernel-otimes-kernel}
Consider the scenario of Example~\ref{measure-times-kernel} again, but
take a third measurable space $(Z, {\cal C})$ with a transition kernel
$L: (Z, {\cal C})\Trans (X, {\cal A})$ into account; assume
furthermore that $x\mapsto K(x)(Y)$ is integrable for each
$L(z)$. Then $L(z)\otimes K$ defines a finite measure on $(X\times Y,
{\cal A}\otimes{\cal B})$ for each $z\in Z$ according to
Example~\ref{measure-times-kernel}. We claim that this defines a
transition kernel $(Z, {\cal C})\Trans (X\times Y, {\cal
  A}\otimes{\cal B})$. For this to be true, we have to show that 
$
z\mapsto \dInt{K(x)(Q_{x})}{L(z)(x)}
$
is measurable for each $Q\in {\cal A}\otimes{\cal B}$. 

Consider 
\begin{equation*}
  {\cal Q} := \{Q\in {\cal A}\otimes{\cal B} \mid \text{ the assertion is true for } Q\}.
\end{equation*}
Then ${\cal Q}$ is closed under complementation. It is also closed under countable disjoint unions by Corollary~\ref{summation-is-intbar}. If $Q = A\times B$ is a measurable rectangle, we have
$
\dInt{K(x)(Q_{x})}{L(z)(x)} = \dInt[A]{K(x)(B)}{L(z)(x)}.
$
Then Exercise~\ref{ex-kernel-yields} shows that this is a measurable function $Z\to \Real$. Thus ${\cal Q}$ contains all measurable rectangles, so ${\cal Q} = {\cal A}\otimes{\cal B}$ by the $\pi$-$\lambda$-Theorem~\ref{Pi-Lambda}. This establishes measurability of $z\mapsto \dInt{K(x)(Q_{x})}{L(z)(x)}$ and shows that it defines a transition kernel. 
\EndExample

As a slight modification, the next example shows the composition of
transition kernels, usually called \emph{convolution}.

\BeginExample{convolution-kernels}
Let $K: (X, {\cal A})\Trans (Y, {\cal B})$ and $L: (Y, {\cal B})\Trans
(Z, {\cal C})$ be transition kernels, and assume that the map $y
\mapsto L(y)(Z)$ is integrable with respect to measures $K(x)$ for an
arbitrary $x\in X$. Define for $x\in X$ and $C\in{\cal C}$
\begin{equation*}
  (L*K)(x)(C) := \dInt{L(y)(C)}{K(x)(y)}.
\end{equation*}
Then $L*K: (X, {\cal A})\Trans (Z, {\cal C})$ is a transition
kernel. In fact, $(L*K)(x)$ is for $x\in X$ fixed a finite measure on
${\cal C}$ according to Corollary~\ref{summation-is-intbar}. From
Exercise~\ref{ex-kernel-yields} we infer that $x\mapsto
\dInt{L(y)(C)}{K(x)(y)}$ is a measurable function, since $y\mapsto
L(y)(C)$ is measurable for all $C\in{\cal C}$.

Because transition kernels are the Kleisli morphisms for the
endofunctor $\FinSenza$ on the category of measurable spaces~\CategCite{Example 1.99}, it is
not difficult to see that this defines Kleisli composition; in
particular it follows that this composition is associative.
\EndExample

\BeginExample{choquet-example}
Let $f\in\MeasbFnctP{X, \mathcal{A}}$, then we know that ``the area under the graph'', viz., 
\begin{equation*}
  C_{\leq}(f) := \{\langle x, r\rangle \mid x\in X, 0\leq r \leq f(x)\}
\end{equation*}
is a member of ${\cal A}\otimes\Borel{\Real}$. This was shown in
Corollary~\ref{pre-choquet}. Then
Corollary~\ref{sigma-finite-product-measure-exists} tells us that
\begin{equation*}
  (\mu\otimes\lambda)(C_{\leq}(f)) = \dInt{\lambda\bigl((C_{\leq}(f))_{x}\bigr)}{\mu(x)},
\end{equation*}
where $\lambda$ is Lebesgue measure on $\Borel{\Real}$. Because
\begin{equation*}
\lambda\bigl((C_{\leq}(f))_{x}\bigr) = \lambda(\{r\mid 0\leq r \leq f(x)\}) = f(x),
\end{equation*}
we obtain
\begin{equation*}
  (\mu\otimes\lambda)(C_{\leq}(f)) = \dInt{f}{\mu}.
\end{equation*}
On the other hand, 
\begin{equation*}
  (\mu\otimes\lambda)(C_{\leq}(f)) = \dInt[\pReal]{\mu\bigl((C_{\leq}(f)_{r}\bigr)}{\lambda(r)},
\end{equation*}
and this gives the integration formula
\begin{equation}
\label{eq-choquet}
  \dInt{f}{\mu} = \int_{0}^{\infty}\mu(\{x\in X \mid f(x) \geq r\})\ dr.
\end{equation}
In this way, the integral of a non-negative function may be
interpreted as measuring the area under its graph.
\EndExample

%%% Local Variables: 
%%% mode: latex
%%% TeX-master: "../Mskr3"
%%% End: 

%\Input{\Folder/Fubini}
%spell checked - 24Aug14
\Subsubsection{Fubini's Theorem}
\label{sec:fubini}
In order to discuss integration with respect to a product measure, we
introduce the cuts of a function $f: X\times Y\to \Real$, defining
$f_{x} := \lambda y.f(x, y)$ and $f^{y} := \lambda x.f(x, y)$.  Thus
we have $f(x, y) = f_{x}(y) = f^{y}(x)$, the first equality resembling
currying.

For the discussion to follow, we will admit also the values
$\{-\infty, +\infty\}$ as function values. So define $\RealT :=
\Real\cup\{-\infty, +\infty\}$, and let $B\subseteq \RealT$ be a Borel
set iff $B\cap\Real\in\Borel{\Real}$. Measurability of functions
extends accordingly: if $f: X\to \RealT$ is measurable, then in
particular $\{x\in X\mid f(x)\in\Real\}\in{\cal A}$, and the set of
values on which $f$ takes the values $+\infty$ or $-\infty$ is a
member of ${\cal A}$. Denote by $\MeasbFnctT{X, \mathcal{A}}$ the set
of measurable functions with values in $\RealT$, and by
$\MeasbFnctTP{X, \mathcal{A}}$ those which take non-negative
values. The integral $\dInt{f}{\mu}$ and integrability is defined in
the same way as above for $f\in\MeasbFnctTP{X, \mathcal{A}}$. Then it
is clear that $f\in\MeasbFnctTP{X, \mathcal{A}}$ is integrable iff
$f\cdot\chi_{\{x\in X\mid f(x)\in\Real\}}$ is integrable and
$\mu(\{x\in X\mid f(x)=\infty\}) = 0.$

With this in mind, we tackle the integration of a measurable function
$f: X\times Y\to \RealT$ for the finite measure spaces $(X, {\cal A},
\mu)$ and $(Y, {\cal B}, \nu)$.

\BeginProposition{product-integration-positive}
Let $f\in\MeasbFnctTP{X\times Y, \mathcal{A}\otimes{\cal B}}$, then 
\begin{enumerate}
\item $\lambda x. \dInt[Y]{f_{x}}{\nu}$ and $\lambda y.\dInt{f^{y}}{\mu}$ are measurable functions $X\to \RealT$ resp. $Y\to \RealT$.
\item we have
  \begin{equation*}
    \dInt[X\times Y]{f}{\mu\otimes\nu} = \dInt{\bigl(\dInt[Y]{f_{x}}{\nu}\bigr)}{\mu(x)} = \dInt[Y]{\bigl(\dInt{f^{y}}{\mu}\bigr)}{\nu(y)}
  \end{equation*}
\end{enumerate}
\EndProposition

\BeginProof
1.
Let $f = \sum_{i=1}^{n}a_{i}\cdot\chi_{Q_{i}}$ be a step function with $a_{i}\geq 0$ and $Q_{i}\in{\cal A}\otimes{\cal B}$ for $i = 1, \dots, n$. Then 
\begin{equation*}
  \dInt[Y]{f_{x}}{\nu} = \sum_{i=1}^{n}a_{i}\cdot\nu(Q_{i, x}).
\end{equation*}
This is a measurable function $X\to \Real$ by Lemma~\ref{cuta-give-equal}. We obtain
\begin{align*}
  \dInt[X\times Y]{f}{\mu\otimes\nu} & = \sum_{i=1}^{n}a_{i}\cdot(\mu\otimes\nu)(Q_{i})\\
& = \sum_{i=1}^{n}a_{i}\cdot\dInt{\nu(Q_{i, x})}{\mu(x)}\\
& = \dInt{\sum_{i=1}^{n}a_{i}\cdot\nu(Q_{i, x})}{\mu(x)}\\
& = \dInt{\bigl(\dInt[Y]{f_{x}}{\nu}\bigr)}{\mu(x)}
\end{align*}
Interchanging the r\^oles of $\mu$ and $\nu$, we obtain the representation of $\lambda y.\dInt[X\times Y]{f}{\mu\otimes\nu}$ in terms of $\dInt{f^{y}}{\mu}$ and $\nu$. Thus the assertion is true for step functions. 

2.
In the general case we know that we can find an increasing sequence $\Folge{f}$ of step functions with $f = \sup_{n\in\Nat}f_{n}$. Given $x\in X$, we infer that $f_{x} = \sup_{n\in\Nat}f_{x, n}$, so that 
\begin{equation*}
  \dInt[Y]{f_{x}}{\nu} = \sup_{n\in\Nat}\dInt{f_{n, x}}{\nu}
\end{equation*}
by Levi's Theorem~\ref{beppo-levi}. This implies
measurability. Applying Levi's Theorem again to the results from part
1., we have
\begin{align*}
  \dInt[X\times Y]{f}{\mu\otimes\nu} & = \sup_{n\in\Nat}\dInt[X\times Y]{f_{n}}{\mu\otimes\nu}\\
& = \sup_{n\in\Nat}\dInt{\bigl(\dInt[Y]{f_{n, x}}{\nu}\bigr)}{\mu(x)}\\
& = \dInt{\bigl(\sup_{n\in\Nat}\dInt[Y]{f_{n, x}}{\nu}\bigr)}{\mu(x)}\\
& = \dInt[X\times Y]{\bigl(\dInt[Y]{f_{x}}{\nu}\bigr)}{\mu(x)}
\end{align*}
Again, interchanging r\^oles yields the symmetric equality. 
\EndProof

This yields as an immediate consequence that the cuts of a product integrable function are almost everywhere integrable, to be specific:

\BeginCorollary{inf-only-on-zero-set}
Let $f: X\times Y\to \Real$ be $\mu\otimes\nu$ integrable, and put 
\begin{align*}
  A & := \{x\in X \mid f_{x}\text{ is not $\nu$-integrable}\},\\
B & := \{y\in Y \mid f^{y}\text{ is not $\mu$-integrable}\}.
\end{align*}
Then $A\in{\cal A}$, $B\in{\cal B}$, and $\mu(A) = \nu(B) = 0$. 
\EndCorollary

\BeginProof
Because $A = \{x\in X\mid \dInt[Y]{|f_{x}|}{\nu} = \infty\}$, we see that $A\in{\cal A}$. By the additivity of the integral, we have
\begin{equation*}
  \dInt[X\times Y]{|f|}{\mu\otimes\nu} = \dInt[X\setminus A]{\bigl(\dInt[Y]{|f_{x}|}{\nu}\bigr)}{\mu(x)} + 
\dInt[A]{\bigl(\dInt[Y]{|f_{x}|}{\nu}\bigr)}{\mu(x)} < \infty,
\end{equation*}
hence $\mu(A) = 0$. $B$ is treated in the same way. 
\EndProof

It is helpful to extend our integral in a minor way. Assume that
$\dInt{|f|}{\mu} < \infty$ for $f: X\to \RealT$ measurable, and that
$\mu(A) = 0$ with $A := \{x\in X\mid |f(x)| = \infty\}$. Change $f$ on
$A$ to a finite value, obtaining a measurable function $f_{*}: X\to
\Real$, and define
\begin{equation*}
  \dInt{f}{\mu} := \dInt{f_{*}}{\mu}.
\end{equation*}
Thus $f\mapsto \dInt{f}{\mu}$ \emph{does not notice} this change on a
set of measure zero. In this way, we assume always that an integrable
function takes finite values, even if we have to convince it to
do so on a set of measure zero.

With this in mind, we obtain

\BeginCorollary{fubinito-1}
Let $f: X\times Y\to \Real$ be integrable, then $\lambda x.\dInt[Y]{f_{x}}{\nu}$ and $\lambda y.\dInt{f^{y}}{\nu}$ are integrable with respect to $\mu$ resp. $\nu$, and
\begin{equation*}
  \dInt[X\times Y]{f}{\mu\otimes\nu} = \dInt{\bigl(\dInt[Y]{f_{x}}{\nu}\bigr)}{\mu(x)} = \dInt[Y]{\bigl(\dInt{f^{y}}{\mu}\bigr)}{\nu(y)}.
\end{equation*}
\EndCorollary

\BeginProof
After the modification on a set of $\mu$-measure zero, we know that 
\begin{equation*}
  \bigl|\dInt{f_{x}}{\nu}\bigr| \leq \dInt[Y]{|f_{x}|}{\nu} < \infty
\end{equation*}
for all $x\in X$, so that $\lambda x.\dInt[Y]{f_{x}}{\nu}$ is integrable with respect to $\mu$; similarly,  $\lambda y.\dInt{f^{y}}{\nu}$ is integrable with respect to $\nu$ for all $y\in Y$. We obtain from Proposition~\ref{product-integration-positive} and the linearity of the integral
\begin{align*}
  \dInt[X\times Y]{f}{\mu\otimes\nu} & = \dInt[X\times Y]{f^{+}}{\mu\otimes\nu}-\dInt[X\times Y]{f^{-}}{\mu\otimes\nu}\\
& = \dInt{\bigl(\dInt[Y]{f_{x}^{+}}{\nu}\bigr)}{\mu(x)} - \dInt{\bigl(\dInt[Y]{f_{x}^{-}}{\nu}\bigr)}{\mu(x)}\\
& = \dInt{\bigl(\dInt[Y]{f_{x}^{+}}{\nu} - \dInt[Y]{f_{x}^{-}}{\nu}\bigr)}{\mu(x)}\\
& = \dInt{\bigl(\dInt[Y]{f_{x}}{\nu}\bigr)}{\mu(x)}
\end{align*}
The second equation is treated in exactly the same way.
\EndProof

Now we now how to treat a function which is integrable, but we do not
yet have a criterion for integrability. The elegance of Fubini's
Theorem shines through the observation that the existence of the
iterated integrals yields integrability for the product integral. To
be specific:

\BeginTheorem{fubini-full}
Let $f: X\times Y\to \Real$ be measurable. Then these statements are equivalent
\begin{enumerate}
\item\label{fubini-full:1} $\dInt[X\times Y]{|f|}{\mu\otimes\nu} < \infty$.
\item\label{fubini-full:2} $\dInt{\bigl(\dInt[Y]{|f_{x}|}{\nu}\bigr)}{\mu(x)} < \infty$.
\item\label{fubini-full:3} $\dInt[Y]{\bigl(\dInt{|f^{y}|}{\mu}\bigr)}{\nu(y)} < \infty$.
\end{enumerate}
Under one of these conditions, $f$ is $\mu\otimes\nu$-integrable, and
\begin{equation}
\label{fubini-full:4}
  \dInt[X\times Y]{f}{\mu\otimes\nu} = \dInt{\bigl(\dInt[Y]{f_{x}}{\nu}\bigr)}{\mu(x)} = \dInt[Y]{\bigl(\dInt{f^{y}}{\mu}\bigr)}{\nu(y)}.
\end{equation}
\EndTheorem

\BeginProof
We discuss only \labelImpl{fubini-full:1}{fubini-full:2}, the other
implications are proved similarly.  From
Proposition~\ref{product-integration-positive} it is inferred that
$|f|$ is integrable, so \ref{fubini-full:2}.~holds by
Corollary~\ref{fubinito-1}, from which we also obtain
representation~(\ref{fubini-full:4}).
\EndProof

%%% Local Variables: 
%%% mode: latex
%%% TeX-master: "../Mskr3"
%%% End: 

%\Input{\Folder/InfiniteProducts}
%spell checked - 24Aug14
\Subsubsection{Infinite Products and Projective Limits}
\label{sec:inf-prod-proj-limit}
Corollary~\ref{sigma-finite-product-measure-exists} extends to a finite number of $\sigma$-finite measure spaces in a natural way. Let $(X_{i}, {\cal A}_{i}, \mu_{i})$ be $\sigma$-finite measure spaces for $1 \leq i \leq n$, the uniquely determined product measure on ${\cal A}_{1}\otimes\dots\otimes{\cal A}_{n}$ is denoted by $\mu_{1}\otimes\dots\otimes\mu_{n}$, and we infer from Corollary~\ref{sigma-finite-product-measure-exists} that we may write 
\begin{align*}
  (\mu_{1}\otimes\dots\otimes\mu_{n})(Q) 
 & = \dInt[X_{2}\times\dots\times X_{n}]{\mu_{1}(Q^{x_{2},\dots, x_{n}})}{(\mu_{2}\otimes\dots\otimes\mu_{n})(x_{2},\dots, x_{n})},\\
& = \dInt[X_{1}\times\dots\times X_{n-1}]{\mu_{n}(Q_{x_{1},\dots, x_{n-1}})}{(\mu_{1}\otimes\dots\otimes\mu_{n-1})(x_{1},\dots, x_{n-1})}
\end{align*}
whenever $Q\in {\cal A}_{1}\otimes\dots\otimes{\cal A}_{n}$. 

We will have a closer look now at infinite products, where we restrict ourselves to probability measures, and here we consider the countable case first. So let $(X_{n}, {\cal A}_{n}, \varpi_{n})$ be a measure space with a probability measure $\varpi_{n}$ on ${\cal A}_{n}$ for $n\in\Nat$. 

Let us fix some notations first. Put
\begin{align*}
  X^{(n)} & := \prod_{k\geq n}X_{k},\\
{\cal A}^{(n)} & := \{A\times X^{^{(n+\ell)}}\mid A\in{\cal A}_{n}\otimes\dots\otimes{\cal A}_{n+\ell-1}\text{ for some }\ell\in\Nat\}
\end{align*}
The elements of ${\cal A}^{(n)}$ are the \emph{\index{cylinder sets}cylinder sets} for $X^{(n)}$. Thus $X^{(1)} = \prod_{n\in\Nat}X_{n}$, and $\bigotimes_{n\in\Nat}{\cal A}_{n} = \sigma({\cal A}^{(1)})$. Given $A\in{\cal A}^{(n)}$, we can write $A$ as $A = C\times X^{n+\ell}$ with $C\in{\cal A}_{n}\otimes\dots\otimes{\cal A}_{n+\ell-1}$. So if we set 
\begin{equation*}
  \varpi^{(n)}(A) := \varpi_{n}\otimes\dots\otimes\varpi_{n+\ell-1}(C),
\end{equation*}
then $\varpi^{(n)}$ is well defined on ${\cal A}^{(n)}$, and it is readily verified that it is monotone and additive with $\varpi^{(n)}(\emptyset) = 0$ and $\varpi^{(n)}(X^{(n)}) = 1$. Moreover, we infer from Theorem~\ref{product-measure-exists} that 
\begin{equation*}
  \varpi^{(n)}(C) = \dInt[X_{n+1}\times\dots\times X_{n+m}]{\varpi_{n}(C^{x_{n+1},\dots,x_{n+m}})}{(\varpi_{n+1}\otimes\dots\otimes\varpi_{n+m})(x_{n+1}\dots x_{n+m})}
\end{equation*}
for all $C\in{\cal A}^{(n)}$. 

The goal is to show that there exists a unique probability measure
$\varpi$ on $\bigotimes_{n\in\Nat}(X_{n}, {\cal A}_{n})$ such that
$\varpi\bigl(A\times X^{(n+1)}\bigr) =
(\varpi_{1}\otimes\dots\otimes\varpi_{n})(A)$ whenever $A\in{\cal
  A}_{n}\otimes{\cal A}^{(n+1)}$. If we can show that $\varpi^{(1)}$ is
$\sigma$-additive on ${\cal A}^{(1)}$, then we can extend
$\varpi^{(1)}$ to the desired $\sigma$-algebra
by~\SetCite{Theorem~2.112}. For this it is sufficient to show that
$\inf_{n\in\Nat}\varpi^{(1)}(A_{n}) > \epsilon >0$ implies
$\bigcap_{n\in\Nat}A_{n}\not=\emptyset$ for any decreasing sequence
$\Folge{A}$ in ${\cal A}^{(1)}$.

The basic idea is to construct a sequence
$\Folge{x}\in\bigcap_{n\in\Nat}A_{n}$. We do this step by step. First
we determine an element $x_{1}\in X_{1}$ such that we can expand the
---~admittedly very short~--- partial sequence $x_{1}$ to a sequence
which is contained in all $A_{n}$; this means that we have to have
$A_{n}^{x_{1}}\not=\emptyset$ for all $n\in\Nat$, because
$A_{n}^{x_{1}}$ contains all possible continuations of $x_{1}$ into
$A_{n}$. We conclude that these sets are non-empty, because their
measure is strictly positive. If we have such an $x_{1}$, we start
working on the second element of the sequence, so we have a look at
some $x_{2}\in X_{2}$ such that we can expand $x_{1}, x_{2}$ to a
sequence which is contained in all $A_{n}$ so we have to have
$A_{n}^{x_{1}, x_{2}}\not=\emptyset$ for all $n\in\Nat$. Again, we
look for $x_{2}$ so that the measure of $A_{n}^{x_{1}, x_{2}}$ is
strictly positive for each $n$. Continuing in this fashion, we obtain
the desired sequence, which then has to be an element of
$\bigcap_{n\in\Nat}A_{n}$ by construction.

This is the plan. Let us have a look at how to find $x_{1}$. Put $E_{1}^{(n)} := \{x_{1}\in X_{1}\mid \varpi^{(2)}(A_{n}^{x_{1}})>\epsilon/2\}$. Because 
\begin{equation*}
  \varpi^{(1)}(A_{n}) = \dInt[X_{1}]{\varpi^{(2)}(A_{n}^{x_{1}})}{\varpi_{1}(x_{1})}
\end{equation*}
we have 
\begin{align*}
  0 < \epsilon < \varpi^{(1)}(A_{n}) & = \dInt[E^{(n)}_{1}]{\varpi^{(2)}(A_{n}^{x_{1}})}{\varpi_{1}(x_{1})}
+ \dInt[X_{1}\setminus E^{(n)}_{1}]{\varpi^{(2)}(A_{n}^{x_{1}})}{\varpi_{1}(x_{1})}\\
& \leq \varpi_{1}(E^{(n)}_{1}) + \epsilon/2\cdot\varpi^{(1)}(X_{1}\setminus E^{(n)}_{1})\\
& \leq \varpi_{1}(E^{(n)}_{1}) + \epsilon/2.
\end{align*}
Thus $\varpi_{1}(E^{(n)}_{1}) \geq \epsilon/2$ for all $n\in\Nat.$
Since $A_{1}\supseteq A_{2}\supseteq \dots$, we have also
$E_{1}^{(1)}\supseteq E_{1}^{(2)}\supseteq\dots$, so let $E_{1} :=
\bigcap_{n\in\Nat}E_{1}^{(n)}$, then $E_{1}\in{\cal A}_{1}$ with
$\varpi_{1}(E_{1})\geq\epsilon/2 > 0$. In particular,
$E_{1}\not=\emptyset$. Pick and fix $x_{1}\in E_{1}$. Then $A_{n}^{x_{1}}\in{\cal A}^{(2)}$, and $\varpi^{(2)}(A_{n}^{x_{1}}) > \epsilon/2$ for all $n\in\Nat$. 

Let us have a look at how to find the second element; this is but a small variation of the idea just presented. Put 
$
E_{2}^{(n)} := \{x_{2}\in X_{2}\mid \varpi^{(3)}(A_{n}^{x_{1}, x_{2}})>\epsilon/4\}
$
for $n\in\Nat$. Because 
\begin{equation*}
  \varpi^{(2)}(A_{n}^{x_{1}}) = \dInt[X_{2}]{\varpi^{(3)}(A_{n}^{x_{1}, x_{2}})}{\varpi_{2}(x_{2})},
\end{equation*}
we obtain similarly $\varpi_{2}(E_{2}^{(n)}) \geq \epsilon/4$ for all
$n\in\Nat$. Again, we have a decreasing sequence, and putting $E_{2}
:= \bigcap_{n\in\Nat}E_{2}^{(n)}$, we have
$\varpi_{2}(E_{2})\geq\epsilon/4$, so that $E_{2}\not=\emptyset$. Pick
$x_{2}\in E_{2}$, then $A_{n}^{x_{1}, x_{2}}\in{\cal A}^{(3)}$ and
$\varpi^{(3)}(A_{n}^{x_{1}, x_{2}})>\epsilon/4$ for all $n\in\Nat$. In
this manner we determine inductively for each $k\in\Nat$ the finite
sequence $\langle x_{1}, \dots, x_{k}\rangle\in X_{1}\times\dots\times
X_{k}$ such that $ \varpi^{(k+1)}(A_{n}^{x_{1}, \dots, x_{k}}) >
\epsilon/2^{k} $ for all $n\in\Nat$. Consider now the sequence
$\Folge{x}$. From the construction it is clear that $\langle x_{1},
x_{2}, \dots, x_{k}, \dots\rangle \in\bigcap_{n\in\Nat}A_{n}$. This
shows that $\bigcap_{n\in\Nat}A_{n}\not=\emptyset$, and it implies
that $\varpi^{(1)}$ is a premeasure on the algebra ${\cal A}^{(1)}$.

Hence we have established~\SetCite{Theorem~2.112}

\BeginTheorem{countably-infinite-product}
Let $(X_{n}, {\cal A}_{n}, \varpi_{n})$ be probability spaces for all $n\in\Nat$. Then there exists a unique probability measure $\varpi$ on $\bigotimes_{n\in\Nat}(X_{n}, {\cal A}_{n})$ such that 
\begin{equation*}\textstyle
  \varpi(A\times\prod_{k > n}X_{k}) = (\varpi_{1}\otimes\dots\otimes\varpi_{n})(A)
\end{equation*}
for all $A\in\bigotimes_{i=1}^{n}{\cal A}_{i}$. 
\QED
\EndTheorem

Define the projection $\pi_{n}^{\infty}: \Folge{x}\mapsto \langle x_{1}, \dots,
x_{n}\rangle$ from $\prod_{n\in\Nat}X_{n}$ to $\prod_{i=1}^{n}X_{i}$. In terms of image measures, 
the theorem states that there exists a unique probability measure
$\varpi$ on the infinite product such that
$\SubProb{\pi_{n}^{\infty}}(\varpi) =
\varpi_{1}\otimes\dots\otimes\varpi_{n}$.

Now let us have a look at the general case, in which the index set is
not necessarily countable. Let $(X_{i}, {\cal A}_{i}, \mu_{i})$ be a
family of probability spaces for $i \in I$, put $X := \prod_{i\in
  I}X_{i}$ and ${\cal A} := \bigotimes_{i\in I}{\cal A}_{i}$. Given
$J\subseteq I$, define $\pi_{J}: (x_{i})_{i\in I}\mapsto (x_{i})_{i\in
  J}$ as the projection $X\to \prod_{i\in J}X_{i}$. Put ${\cal A}_{J} := \InvBild{\pi_{J}}{\bigotimes_{j\in J}{\cal A}_{j}}$.

Although the index set $I$ may be large, the measurable sets in ${\cal
  A}$ are always determined by a countable subset of the index set:

\BeginLemma{coutab-determined}
Given $A\in {\cal A}$, there exists a countable subset $J\subseteq I$
such that $\chi_{A}(x) = \chi_{A}(x')$, whenever $\pi_{J}(x) =
\pi_{J}(x')$.
\EndLemma

\BeginProof
Let ${\cal G}$ be the set of all $A\in {\cal A}$ for which the assertion is true. Then ${\cal G}$ is a $\sigma$-algebra which contains $\InvBild{\pi_{\{i\}}}{{\cal A}_{i}}$ for every $i\in I$, hence ${\cal G} = {\cal A}$. 
\EndProof

This yields as an immediate consequence

\BeginCorollary{coutab-determined-cor}
${\cal A} = \bigcup\{{\cal A}_{J}\mid J\subseteq I\text{ is countable}\}$.
\EndCorollary

\BeginProof
It is enough to show that the set on the right hand side is a $\sigma$-algebra. This follows easily from Lemma~\ref{coutab-determined}.
\EndProof

We obtain from this observation, and from and the previous result for the countable case that arbitrary products exist.

\BeginTheorem{arbitrary-infinite-product}
Let  $(X_{i}, {\cal A}_{i}, \mu_{i})$ be a family of probability spaces for $i \in I$. Then there exists a unique probability measure $\mu$ on $\bigotimes_{i\in I}(X_{i}, {\cal A}_{i})$ such that
\begin{equation}
\label{arbitrary-infinite-product-ext-prop}
  \mu(\InvBild{\pi_{\{i_{1}, \dots, i_{k}\}}}{C}) = (\mu_{i_{1}}\otimes\dots\otimes\mu_{i_{k}})(C)
\end{equation}
for all $C\in\bigotimes_{j=1}^{k}{\cal A}_{i_{j}}$ and all $i_{1}, \dots, i_{k}\in I$.
\EndTheorem

\BeginProof
Let $A\in{\cal A}$, then there exists a countable subset $J\subseteq I$ such that $A\in{\cal A}_{J}$. Let $\mu_{J}$ be the corresponding product measure on ${\cal A}_{J}$. Define $\mu(A) := \mu_{J}(A)$, then it it easy to see that $\mu$ is a well defined measure on ${\cal A}$, since the extension to countable products is unique. From the construction it follows also that the desired property~(\ref{arbitrary-infinite-product-ext-prop}) is satisfied. 
\EndProof

For the interpretation of some logics the projective limit of a
projective family of stochastic relations is helpful; this is the
natural extension of a product. It will be discussed now. Denote by $ X^\infty := \prod_{k \in
  \Nat} X$ the infinite product of $X$ with itself; recall that $\ProbSenza$ is the probability functor, assigning to each measurable space its probability measures.

\BeginDefinition{Proj-Limit-Meas} Let $X$ be a Polish space, and
$\Folge{\mu}$ a sequence of probability measures $\mu_n \in
\Prob{X^n}$. This sequence is called a \emph{projective
system}\index{measure!projective system} iff
$
\mu_n(A) = \mu_{n+1}(A \times X)
$
for all $n \in \Nat$ and all Borel sets $A \in \Borel{X^n}$. A
probability measure $\mu_\infty \in \Prob{X^\infty}$ is called the
\emph{projective limit}\index{measure!projective limit}
of the projective system $\Folge{\mu}$ iff
\begin{equation*}\textstyle
\mu_n(A) = \mu_\infty(A \times \prod_{j > n} X)
\end{equation*}
for all $n \in \Nat$ and $A \in \Borel{X^n}$.
\EndDefinition

Thus a sequence of measures is a projective system iff each measure is
the projection of the next one; its projective limit is characterized
through the property that its values on cylinder sets coincides with
the value of a member of the sequence, after taking projections. A
special case is given by product measures. Assume that $\mu_{n} =
\nu_{1}\otimes\dots\otimes\nu_{n}$, where $\Folge{\nu}$ is a sequence
of probability measures on $X$. Then the condition on projectivity is
satisfied, and the projective limit is the infinite product
constructed above. It should be noted, however, that the projectivity condition
does not express $\mu_{n+1}(A\times B)$ in terms of $\mu_{n}(A)$ for
an arbitrary measurable set $B\subseteq X$, as the product measure
does.

It is not immediately obvious that a projective limit exists in
general, given the rather weak dependency of the measures. In general,
it will not, and this is why. The basic idea for the construction of
the infinite product has been to define the limit on the cylinder sets
and then to extend this premeasure --- but it has to be established
that it is indeed a premeasure, and this is difficult in general. The
crucial property in the proof above has been that $ \mu_{n_k}(A_k)
\rightarrow 0 $ whenever $\Folge{A}$ is a sequence of cylinder sets
$A_k$ (with at most $n_k$ components that do not equal $X$) that
decreases to $\emptyset$. This property has been established above for
the case of the infinite product through Fubini's Theorem, but this is
not available in the general setting considered here. We will see,
however, that a topological argument will be helpful. This is why we
did postulate the base space $X$ to be Polish.

We start with an even stronger topological condition, viz., that the
space under consideration is compact and metric.  The central statement
is

\BeginProposition{Projective-1}
Let $X$ be a compact metric space.  Then the projective system $\Folge{\mu}$ has a unique projective limit
$\mu_\infty$.
\EndProposition

\BeginProof
1.
Let
$
A = A'_k \times \prod_{j > k} X
$
be a cylinder set with $A'_k \in \Borel{X^k}$. Define
$
\mu^*(A) := \mu_k(A'_k).
$
Then $\mu^*$ is well defined on the cylinder sets, since the sequence forms a projective
system. In order to show that $\mu^*$ is a premeasure on the cylinder
sets, we take a decreasing sequence $\Folge{A}$ of cylinder
sets with
$
\bigcap_{n \in \Nat} A_n = \emptyset
$
and show that
$
\inf_{n \in \Nat} \mu^*(A_n) = 0.
$
In fact, suppose that $\Folge{A}$ is decreasing with
$
\mu^*(A_n) \geq \delta
$
for all $n \in \Nat$, then we show that
$
\bigcap_{n \in \Nat} A_n \not= \emptyset.
$

We can write
$
A_n = A'_n \times \prod_{j > k_n} X
$
for some $ A'_n \in \Borel{X^{k_n}}. $ From Lemma~\ref{are-regular}
we get for each $n$ a closed, hence compact set $ K'_n \subseteq A'_n $ such that
$ \mu_{k_n}(A'_n\setminus K'_n) < \delta/2^{n}. $ Because
$X^\infty$ is compact by Tichonov's Theorem, 
\begin{equation*}
K''_n := K'_n \times \prod_{j > k_n} X
\end{equation*}
is a compact set, and $ K_n := \bigcap_{j =
1}^n K''_j \subseteq A_n $ is compact as well, with
\begin{equation*}
  \mu^*(A_n\setminus K_n)
    \leq
  \mu^*(\bigcup_{j=1}^n A''_n\setminus K''_j)
    \leq 
  \sum_{j = i}^n \mu^*(A''_j \setminus K''_j)
    =
  \sum_{j = 1}^n \mu_{k_j}(A'_j \setminus K'_j)
    \leq
  \sum_{j = 1}^\infty \delta/2^{j}
    =
  \delta.
\end{equation*}
Thus $\Folge{K}$ is a decreasing sequence of nonempty compact sets;
consequently,
\begin{equation*}
\emptyset \not= \bigcap_{n \in \Nat} K_n \subseteq \bigcap_{n \in
  \Nat} A_n.
\end{equation*}

2. Since the cylinder sets generate the Borel sets of $X^\infty$,
and since $\mu^*$ is a premeasure, we know that there exists a
unique extension $\mu_\infty \in \Prob{X^\infty}$ to it. Clearly, if
$A \subseteq X^n$ is a Borel set, then
\begin{equation*}\textstyle
\mu_\infty(A \times \prod_{j > n} X)
=
\mu^*(A \times \prod_{j > n} X)
=
\mu_n(A),
\end{equation*}
so we have constructed a projective limit.

3. Suppose that $\mu'$ is another probability measure in
$\Prob{X^\infty}$ that has the desired property. Consider
\begin{equation*}
\mathcal{D} := \{D \in \Borel{X^\infty} \mid  \mu_\infty(D) =
\mu'(D)\}.
\end{equation*}
It is clear the $\mathcal{D}$ contains all cylinder sets, that it
is closed under complements, and under countable disjoint unions.
By the $\pi$-$\lambda$-Theorem~\ref{Pi-Lambda} $\mathcal{D}$
contains the $\sigma$-algebra generated by the cylinder sets,
hence all Borel subset of $X^\infty$. This establishes uniqueness
of the extension. 
\EndProof

The proof makes critical use of the observation that we can approximate the measure of a Borel set arbitrarily well by compact sets from within; see Lemma~\ref{are-regular}. It is also important that compact sets have the finite intersection property: if each finite intersection of a family of compact sets is nonempty, the intersection of the entire family cannot be empty. Consequently the proof given above works in general Hausdorff spaces, provided the measures under consideration have the approximation property mentioned above.

We free ourselves from the restrictive assumption of having a compact metric space using the Alexandrov embedding\index{theorem!Alexandrov} of a Polish space into a compact metric space.

\BeginProposition{Projective-2} 
Let $X$ be a Polish space, $\Folge{\mu}$ be a projective system on $X$. Then there exists a unique projective limit\index{measure!projective limit} $ \mu_\infty \in \Prob{X^\infty} $ for $\Folge{\mu}$. 
\EndProposition

\BeginProof 
$X$ is a dense measurable subset of a compact metric space $\Dach{X}$
by Alexandrov's Theorem~\ref{Alexandrov}. Defining $ \Dach{\mu}_n(B)
:= \mu_n(B \cap X^n) $ for the Borel set $B \subseteq \Dach{X}^n$
yields a projective system $ \left(\Dach{\mu}_n\right)_{n \in \Nat} $
on $\Dach{X}$ with a projective limit $ \Dach{\mu}_\infty $ by
Proposition~\ref{Projective-1}.  Since by construction $
\Dach{\mu}_\infty(X^\infty) = 1, $ restrict $\Dach{\mu}_\infty$ to the
Borel sets of $X^\infty$, then the assertion follows.
\EndProof

An interesting application of this construction arises through stochastic relations that form a projective system. We will show now that there exists a kernel which may be perceived as a (pointwise) projective limit.

\BeginCorollary{Projective-3} Let $X$ and $Y$ be Polish spaces, and assume that $ J^{(n)}$ is a stochastic relation on $X$ and $Y^n$ for each $n \in \Nat$ such that the sequence $ \left(J^{(n)}(x)\right)_{n \in \Nat} $ forms a projective system on $Y$ for each $x \in X$, in particular $J^{(n)}(x)(Y^n) = 1$ for all $x \in X$. Then there exists a unique sub-Markov kernel $ J_\infty$ on $X$ and $Y^\infty$ such that $ J_\infty(x) $ is the projective limit of $\left(J^{(n)}(x)\right)_{n \in \Nat}$ for each $x \in X$.  
\EndCorollary

\BeginProof 
0. 
Let for $x$ fixed $J_\infty(x)$ be the projective limit of the projective system $ \left(J^{(n)}(x)\right)_{n \in \Nat} $. By the definition of a stochastic relation we need to show that the map $ x \mapsto J_\infty(x)(B) $ is measurable for every $B \in \Borel{Y^\infty}$.

1.
In fact, consider
\begin{equation*}
\mathcal{D} := \{B \in \Borel{Y^\infty} \mid x \mapsto
J_\infty(x)(B)\text{ is measurable}\}
\end{equation*}
then the general properties of measurable functions imply that
$\mathcal{D}$ is a $\sigma$-algebra on $Y^\infty$. Take a cylinder set
$
B = B_0 \times \prod_{j > k} Y
$
with
$
B_0 \in \Borel{Y^k}
$
for some $k \in \Nat$, then, by the properties of the projective
limit, we have
$
J_\infty(x)(B) = J^{(k)}(x)(B_0).
$
But
$
x \mapsto J^{(k)}(x)(B_0)
$
constitutes a measurable function on $X$. Consequently,  $B \in
\mathcal{D}$, and so $\mathcal{D}$ contains
the cylinder sets which generate $\Borel{Y^\infty}$. Thus
measurability is established for each Borel set $B \subseteq Y^\infty,$
arguing with the $\pi$-$\lambda$-Theorem~\ref{Pi-Lambda} as in the last part of the proof for
Proposition~\ref{Projective-1}.
\EndProof

\Subsubsection{Case Study: Continuous Time Stochastic Logic}
\label{sec:case-study-csl}

We illustrate this construction through the interpretation of a path
logic over infinite paths; the logic is called \texttt{CSL} ---
\emph{\index{logic!continuous time stochastic}continuous time
  stochastic logic}. Since the discussion of this application requires
some preparations, some of which are of independent interest, we
develop the example in a series of steps.

We introduce \texttt{CSL} now and describe it informally first.

Fix $P$ as a countable set of atomic propositions. We define
recursively state formulas\index{logic!CSL}\index{CSL} and path
formulas for \texttt{CSL}: 
\begin{description}
 \item[State formulas] are defined through the syntax
$$
\varphi ::= \top \mid  a \mid  \neg \varphi \mid  \varphi \wedge
\varphi' \mid  \Steady{\varphi} \mid  \PathQuant{\psi}.
$$
Here $a \in P$ is an atomic proposition, $\psi$ is a path formula,
$\RelOp$ is one of the relational operators $<, \leq, \geq, >$, and
$p \in [0, 1]$ is a rational number.
\item[Path formulas] are defined through
$$
\psi ::= \Next{\varphi} \mid  \Until{\varphi}{\varphi'}
$$
with $\varphi, \varphi'$ as state formulas, $I \subseteq \pReal$ a
closed interval of the real numbers with rational bounds (including $I = \pReal$).
\end{description}
We denote the set of all state formulas by $\AllFormulas$.

The operator $\Steady{\varphi}$ gives the \emph{steady-state
probability} for $\varphi$ to hold with the boundary condition
$\RelOp p$\index{CSL!steady-state}\index{steady-state!CSL};
the formula $\PathQuantSenza$ replaces quantification:
the \emph{path-quantifier} formula $\PathQuant{\psi}$ holds in a
state $s$ iff the probability of all paths starting in $s$ and
satisfying $\psi$ is specified by
$\RelOp p$\index{CSL!path quantifier}\index{path quantifier!CSL}.
Thus $\psi$ holds on almost
all paths starting from $s$ iff $s$ satisfies $\PathQuant[\geq
1]{\psi}$, a path being an alternating infinite sequence $\sigma =
\langle s_0, t_0, s_1, t_1, \dots\rangle$ of states $x_i$ and of
times $t_i$.  Note that the time is being made explicit here. The
\emph{next-operator}
$\Next{\varphi}$\index{CSL!next operator}\index{next operator!CSL}
is assumed to hold on path
$\sigma$ iff $s_1$ satisfies $\varphi$, and $t_0 \in I$ holds.
Finally, the \emph{until-operator}
$\Until{\varphi_1}{\varphi_2}$\index{CSL!until operator}\index{until operator!CSL}
holds on path $\sigma$ iff we can find a point in time $t \in I$
such that the state $\At{t}$ which $\sigma$ occupies at time $t$
satisfies $\varphi_2$, and for all times $t'$ before that, $\At{t'}$
satisfies $\varphi_1$.

A Polish state space $S$ is fixed; this space is used for modelling a
transition system takes also time into account. We are not only
interested in the next state of a transition but also in the time
after which to make a transition. So the basic probabilistic datum
will be a stochastic relation $M: S\Trans \pReal\times S$; if we are
in state $s$, we will do a transition to a new state $s'$ after we did
wait some specified time $t$; $M(s)(D)$ will give the probability that
the pair $\langle t, s' \rangle\in D$.  We assume that
$M(s)(\pReal\times S) = 1$ holds for all $s\in S$.

A \emph{path} $\sigma$ is an element of the set $ (S \times \pReal)^\infty. $
Path $\sigma = \langle s_0, t_0, s_1, t_1, \dots \rangle$ may be
written as $ s_0 \Pfeil[t_0] s_1 \Pfeil[t_1] \dots $ with the
interpretation that $t_i$ is the time spent in state $s_i$. Given $i
\in \Nat$, denote $s_i$ by $\sigma[i]$ as the $(i+1)$-st state of
$\sigma$, and let $\delta(\sigma, i) := t_i$. Let for $t \in \pReal$ the
index $j$ be the smallest index $k$ such that $ t < \sum_{i=0}^k t_i,
$ and put $\At{t} := \sigma[j]$, if $j$ is defined; set $\At{t} :=
\#$, otherwise (here $\#$ is a new symbol not in $S \cup \pReal$).
$S_{\#}$ denotes $S \cup \{\#\}$; this is a Polish space when endowed
with the sum $\sigma$-algebra. The definition of $\At{t}$ makes sure
that for any time $t$ we can find a rational time $t'$ with $ \At{t} =
\At{t'}. $

We will deal only with infinite paths. This is no loss of generality
because events that happen at a certain time with probability 0 will
have the effect that the corresponding infinite paths occur only with
probability 0. Thus we do not prune the path; this makes the notation
somewhat easier to handle.

The Borel sets $\Borel{\Paths}$ are the smallest $\sigma$-algebra
which contains all the cylinder sets
\begin{equation*}
\{\prod_{j = 1}^n (B_j \times I_j) \times \prod_{j > n} (S \times
\pReal) \mid  n \in \Nat, I_1, \dots, I_n\text{ rational
intervals, } B_1, \dots, B_n \in \Borel{S}\}.
\end{equation*}
Thus a cylinder set is an infinite product that is determined
through the finite product of an interval with a Borel set in $S$. It will be helpful to
remember that the intersection of two cylinder sets is again a
cylinder set.

Given $M: S\Trans \pReal\times S$ with Polish $S$, define inductively $M_{1} := M$, and 
\begin{equation*}
  M_{n+1}(s_{0})(D) := 
\dInt[(\pReal\times S)^{n}]{M(s_{n})(D_{t_{0}, s_{1}, \dots, t_{n-1}, s_{n}})}{M_{n}(s_{0})(t_{0}, s_{1}, \dots, t_{n-1}, s_{n})}
\end{equation*}
for the Borel set $D\subseteq (\pReal\times S)^{n+1}$. Let us
illustrate this for $n=1$. Given $D\in\Borel{(\pReal\times S)^{2}}$
and $s_{0}\in S$ as a state to start from, we want to calculate the
probability $M_{2}(s_{0})(D)$ that $\langle t_{0}, s_{1}, t_{1},
s_{2}\rangle \in D$. This is the probability for the initial path
$\langle s_{0}, t_{0}, s_{1}, t_{1}, s_{2}\rangle$ (a path\emph{let}),
given the initial state $s_{0}$.  Since $\langle t_{0}, s_{1}\rangle$
is taken care of in the first step, we fix it and calculate $
M(s_{1})(\{\langle t_{1}, s_{2}\rangle \mid \langle t_{0}, s_{1},
t_{1}, s_{2}\rangle \in D\}) = M(s_{1}(D_{t_{0}, s_{1}}), $ by
averaging, using the probability provided by $M(s_{0})$, so that we
obtain
\begin{equation*}
M_{2}(s_{0})(D) 
=
\dInt[\pReal\times S]{M(s_{1})(D_{t_{0}, s_{1}})}{M(s_{0})(t_{0}, s_{1})}
\end{equation*}

Consequently, for the general case we obtain $M_{n+1}(s_{0})(D)$ as the probability for $\langle s_{0}, t_{0},\dots, s_{n}, t_{n}, s_{n+1}\rangle$ as the initial piece of an infinite path to be a member of $D$. This probability indicates that we start in $s_{0}$, remain in this state for $t_{0}$ units of time, then enter state $s_{1}$, remain there for $t_{1}$ time units, etc., and finally leave state $s_{n}$ after $t_{n}$ time units, entering $s_{n+1}$, all this happening within $D$. 

We claim that $(M_{n}(s))_{n\in\Nat}$ is a projective system. We first see from Example~\ref{kernel-otimes-kernel} that $M_{n}: S\Trans (\pReal\times S)^{n}$ defines a transition kernel for each $n\in\Nat$. Now let $D = A\times (\pReal\times S)$ with $A\in\Borel{(\pReal\times S)^{n}}$, then $M(s_{n})(D_{t_{0}, s_{1}, \dots, t_{n-1}, s_{n}}) = M(s_{n})(\pReal\times S) = 1$
 for all $\langle t_{0}, s_{1}, \dots, t_{n-1}, s_{n}\rangle\in A$, so that we obtain 
$
M_{n+1}(s)(A\times (\pReal\times S)) = M_{n}(s)(A).
$
The condition on projectivity is satisfied. Hence there exists a unique projective limit, hence a transition kernel 
\begin{equation*}
M_{\infty}: S\Trans (\pReal\times S)^{\infty}
\end{equation*}
with
\begin{equation*}
M_{n}(s)(A) = M_{\infty}(s)\bigl(A\times \prod_{k>n}(\pReal\times S)\bigr)
\end{equation*}
for all $s\in S$ and for all $A\in\Borel{(\pReal\times S)^{n}}$. 

The projective limit displays indeed limiting behavior: suppose
$B$ is an infinite measurable cube
$\prod_{n \in \Nat} B_n$ with $B_n \in \Borel{\pReal \times S}$
as Borel sets. Because
\begin{equation*}\textstyle
B = \bigcap_{n \in \Nat} \left(\prod_{1 \leq j \leq n} B_j
\times \prod_{j > n} (\pReal \times S)\right),
\end{equation*}
is represented as the intersection of a monotonically decreasing sequence,
we have for all $s \in S$
\begin{align*}\textstyle
  \Inff{M}(s)(B) & = \textstyle\lim_{n\rightarrow\infty} \Inff{M}(s)\bigl(\prod_{1 \leq j \leq n} B_j
\times \prod_{j > n} (\pReal \times S)\bigr) \\
    & = \textstyle\lim_{n\rightarrow\infty} M_n(s)\bigl(\prod_{1 \leq j \leq n} B_j\bigr).
  \end{align*}
Hence $\Inff{M}(s)(B)$ is the limit of the probabilities $M_n(s)(B_n)$
at step $n$.

In this way models based on a Polish state space $S$ yield stochastic
relations $ S \Trans (\pReal \times S)^\infty $ through projective
limits.  Without this limit it would be difficult to model the
transition behavior on infinite paths; the assumption that we work in
Polish spaces makes sure that these limits in fact do exist. To get
started, we need to assume that given a state $s \in S$, there is
always a state to change into after a finite amount of time.

We obtain as a first consequence of the construction for the projective
limit a recursive formulation for the
transition law $M: X \Trans (\pReal \times S)^\infty$. Interestingly, it reflects
the domain equation
$
(\pReal \times S)^\infty = (\pReal \times S) \times (\pReal \times
X)^\infty. 
$

\BeginLemma{int-rec-form}
If $D \in \Borel{(\pReal \times S)^\infty }$, then
\begin{equation*}
\Inff{M}(s)(D) = \int_{\pReal \times S} \Inff{M}(s')(D_{\langle t, s'\rangle})\ M_{1}(s)(d\langle t,
s'\rangle)
\end{equation*}
holds for all $s \in S$.
\EndLemma

\BeginProof 
Recall that $ D_{\langle t, s'\rangle} = \{\tau \mid
\langle t, s', \tau\rangle \in D\}. $ Let 
\begin{equation*}
D = \left(H_1 \times
\ldots \times H_{n+1}\right) \times \prod_{j > n} (\pReal \times S)
\end{equation*}
be a cylinder set with $H_i \in \Borel{\pReal \times S}, 1 \leq i \leq
n+1$. The equation in question in this case boils down to
\begin{equation*}
M_{n+1}(s)(H_1 \times \ldots \times H_{n+1})
=
\int_{H_1} M_{n}(s')(H_2 \times \ldots \times H_{n+1})
M_{1}(s)(d\langle t, s'\rangle).
\end{equation*}
This may easily be derived from the definition of the projective
sequence. Consequently, the equation in question holds for all
cylinder sets, thus the $\pi$-$\lambda$-Theorem~\ref{Pi-Lambda}
implies that it holds for all Borel subsets of $(\pReal \times
S)^\infty$.
\EndProof

This decomposition indicates that we may first select in state $s$ a
new state and a transition time; with these data the system then works
just as if the selected new state would have been the initial
state. The system does not have a memory but reacts depending on its
current state, no matter how it arrived there.
Lemma~\ref{int-rec-form} may accordingly be interpreted as a Markov
property\index{Markov property} for a process the behavior of which is
independent of the specific step that is undertaken.

We need some information about the $@$-operator before continuing.

\BeginLemma{ein-hilfssatz}
$\langle \sigma, t\rangle \mapsto \At{t}$ is
a Borel measurable map from $\Paths \times \pReal$ to $S_{\#}$. In
particular, the set $ \{\langle \sigma, t\rangle \mid  \At{t} \in
S\} $ is a measurable subset of $\Paths \times \pReal$.
\EndLemma

\BeginProof

0. Note that we claim joint measurability in both components (which is strictly stronger than measurability in each component). Thus we have to show that $ \{\langle \sigma, t\rangle \mid \At{t} \in A\} $ is a measurable subset of $\Paths \times \pReal$, whenever $A \subseteq S_{\#}$ is Borel.

1. Because for fixed $i \in \Nat$ the map $ \sigma \mapsto
\delta(\sigma, i) $ is a projection, $\delta(\cdot, i)$ is
measurable, hence $ \sigma \mapsto \sum_{i=0}^j \delta(\sigma, i) $
is. Consequently,
\begin{equation*}\textstyle
  \{\langle \sigma, t\rangle \mid  \At{t} = \#\}
  = 
  \{\langle \sigma, t\rangle \mid  \forall j: t \geq \sum_{i=0}^j
  \delta(\sigma, i)\}
  = 
  \bigcap_{j \geq 0} \{\langle \sigma, t\rangle \mid  t \geq \sum_{i=0}^j
  \delta(\sigma, i)\}.
\end{equation*}
This is clearly a measurable set.

2. Put
$
stop(\sigma, t) := \inf\{k\geq 0 \mid  t < \sum_{i=0}^k
\delta(\sigma, i)\},
$
thus $stop(\sigma, t)$ is the smallest index for which the accumulated waiting in $\sigma$ times exceed $t$. 
\begin{equation*}\textstyle
X_k := \{\langle \sigma, t\rangle \mid  stop(\sigma, t) = k\} =
\{\langle \sigma, t\rangle \mid
   \sum_{i=0}^{k-1} \delta(\sigma, i) \leq t <
   \sum_{i=0}^k \delta(\sigma, i)\}
 \end{equation*}
is a measurable set by Corollary~\ref{pre-choquet}. Now let $B \in \Borel{S}$
be a Borel set, then
\begin{align*}
  \{\langle \sigma, t\rangle \mid  \At{t} \in B\}
  & =  
  \bigcup_{k \geq 0} \{\langle \sigma, t\rangle \mid  \At{t} \in B,
  stop(\sigma, t) = k\} \\
  & =  
  \bigcup_{k \geq 0} \{\langle \sigma, t\rangle \mid  \sigma[k] \in B,
  stop(\sigma, t) = k\} \\
  & =  
  \bigcup_{k \in \Nat}\left(X_k \cap \bigl(\prod_{i<k}(S \times \pReal) \times (B \times \pReal)
  \times \prod_{i > k} (S \times \pReal)\bigr)\right).
\end{align*}
Because $X_k$ is measurable, the latter set is measurable. This
establishes measurability of the $@$-map. 
\EndProof

As a consequence, we establish that some sets and maps, which will be important for the later development, are actually measurable. A notational convention for improving readability is proposed: the letter $\sigma$ will always denote a generic element of $\Paths$, and the letter $\tau$ always a generic element of $\pReal \times \Paths$.

\BeginProposition{SomeMeas} 
We observe the following properties:
\begin{enumerate}
 % \item \label{p-item:1}
 %  The set of \emph{Zeno paths}\index{logic!{CSL}!Zeno paths}
 %  $\{\sigma \mid  \sum_{i \geq 0} \delta(\sigma, i) \text{ exists and is finite}\}$ is a
 %  measurable subset of $\Paths$,
 \item \label{p-item:2}
  $\{\langle \sigma, t\rangle \mid
  \lim_{i \rightarrow \infty} \delta(\sigma, i) = t\}$
  is a measurable subset of $\Paths \times \pReal$,
 \item \label{p-item:3}
  let $N_\infty: S \Trans (\pReal \times S)^\infty$ be a stochastic relation, then
\begin{align*}
  s & \mapsto \liminf_{t \rightarrow \infty}
  N_\infty(s)(\{\tau \mid  \At[\langle s, \tau\rangle]{t} \in A\})\\
  s &\mapsto \limsup_{t \rightarrow \infty}
  N_\infty(s)(\{\tau \mid  \At[\langle s, \tau\rangle]{t} \in A\})
\end{align*}
constitute measurable maps $X \rightarrow \pReal$ for each Borel set $A \subseteq S$.
\end{enumerate}
\EndProposition

\BeginProof 0. The proof makes crucial use of the fact that the real
line is a complete metric space (so each Cauchy sequence converges),
and that the rational numbers are a dense and countable set.

% 1. Since $ \sum_{i \geq 0} \delta(\sigma, i) $ exists and is finite
% iff given $\epsilon > 0$ there exists $n \in \Nat$ such that $ \mid
% \sum_{i=n_1}^{n_2} \delta(\sigma, i)\mid  < \epsilon $ whenever
% $n_1, n_2 \geq n$, we see that
% \begin{multline*}
%   \{\sigma \mid  \sum_{i \geq 0} \delta(\sigma, i) \text{ exists and is
%     finite}\}
% =\\ \bigcap_{\Rational\ni\epsilon>0} \bigcup_{n \in \Nat}
% \bigcap_{n_1, n_2 \geq n} \{\sigma \mid  \left| \sum_{i=n_1}^{n_2}
% \delta(\sigma, i)\right|  < \epsilon\}.
% \end{multline*}
% Measurability of $ \sigma \mapsto \delta(\sigma, i) $ for each $i$
% follows from Lemma~\ref{Hilfssatz}. This implies measurability of
% the set in part$~\ref{p-item:1}$, since the union and the
% intersections are defined over countable index sets.

1. In order to establish part$~\ref{p-item:2}$, write
\begin{equation*}
\{\langle \sigma, t\rangle \mid  \lim_{i \rightarrow \infty}
\delta(\sigma, i) = t\} = \bigcap_{\Rational\ni\epsilon>0}
\bigcup_{n \in \Nat} \bigcap_{m \geq n} \{\langle \sigma, t\rangle
\mid  \mid  \delta(\sigma, m) - t \mid  < \epsilon\}.
\end{equation*}
By Lemma~\ref{ein-hilfssatz}, the set
\begin{equation*}
\{\langle \sigma, t\rangle \mid  \mid  \delta(\sigma, m) - t \mid  <
\epsilon\} = \{\langle \sigma, t\rangle \mid  \delta(\sigma, m) > t
- \epsilon\} \cap \{\langle \sigma, t\rangle \mid   \delta(\sigma,
m) < t + \epsilon\}
\end{equation*}
is a measurable subset of $\Paths \times \pReal$, and since the
union and the intersections are countable, measurability is
inferred.

2. From the definition of the $@$-operator it is immediate that
given an infinite path $\sigma$ and a time $t\in \pReal$, there
exists a rational $t'$ with $ \At{t} = \At{t'}. $ Thus we obtain for
an arbitrary real number $x$, an arbitrary Borel set $A \subseteq S$
and $s \in S$
\begin{align*}
  \liminf_{t \rightarrow \infty}
  N_\infty(s)(\{\tau \mid  \At[\langle s, \tau\rangle]{t} \in
  A\}) \leq x   
& \Leftrightarrow
  \sup_{t\geq 0} \inf_{r\geq t}
  N_\infty(s)(\{\tau \mid  \At[\langle s, \tau\rangle]{r} \in
  A\}) \leq x\\
& \Leftrightarrow
  \sup_{\Rational\ni t\geq 0} \inf_{\Rational\ni r\geq t}
  N_\infty(s)(\{\tau \mid  \At[\langle s, \tau\rangle]{r} \in
  A\}) \leq x\\
& \Leftrightarrow
  s \in
  \bigcap_{\Rational\ni t\geq 0}\bigcup_{\Rational \ni r \geq t} A_{r,
  x}
\end{align*}
with
\begin{equation*}
A_{r, x} := \{s' \mid
  N_\infty(s')(\{\tau \mid  \At[\langle s', \tau\rangle]{r} \in
  A\}) \leq x\}.
\end{equation*}
We infer that $A_{r, x}$ is a measurable subset of $S$ from the fact
that $N_\infty$ is a stochastic relation and from
Exercise~\ref{ex-joint-measurability}. Since a map $f: W \rightarrow
\Real$ is measurable iff each of the sets $ \{w \in W\mid f(w) \leq
s\} $ is a measurable subset of $W$, the assertion follows for the
first map. The second part is established in exactly the same way,
using that $f: W \rightarrow \Real$ is measurable iff $ \{w \in W \mid
f(w) \geq s\} $ is a measurable subset of $W$, and observing
\begin{equation*}
  \limsup_{t \rightarrow \infty}
  N_\infty(s)(\{\tau\mid  \At[\langle x, \tau\rangle]{t} \in
  A\}) \geq x
  \Leftrightarrow
  \inf_{\Rational\ni t\geq 0} \sup_{\Rational\ni r\geq t}
  N_\infty(x)(\{\tau\mid  \At[\langle s, \tau\rangle]{r} \in
  A\}) \geq x.
\end{equation*}
\EndProof

This has some consequences which will come in useful for the
interpretation of \texttt{CSL}. Before stating them, it is noted that
the statement above (and the consequences below) do not make use of
$N_{\infty}$ being a projective limit; in fact, we assume $N_{\infty}:
S\Trans (\pReal\times S)^{\infty}$ to be an arbitrary stochastic
relation. A glimpse at the proof shows that these statements even hold
for finite transition kernels, but since we will use it for the
probabilistic case, we stick to stochastic relations.

Now for the consequences. As a first consequence we obtain that the
set on which the asymptotic behavior of the transition times is
reasonable (in the sense that it tends probabilistically to a limit) is
well behaved in terms of measurability:

\BeginCorollary{LimitBehavior} 
Let $A \subseteq X$ be a Borel set,
and assume that $N_\infty: S \Trans (\pReal \times S)^\infty$ is a
stochastic relation. Then
\begin{enumerate}
 \item \label{p-item:5}
  the set
  $
  Q_A := \{s \in S \mid  \lim_{t \rightarrow \infty}
  N_\infty(s)(\{\tau \mid  \At[\langle s,
  \tau\rangle]{t} \in A\})
  \text{ exists}\}
  $
  on which the limit exists is a Borel subset of $S$,
 \item \label{p-item:6}
  $
  s \mapsto \lim_{t \rightarrow \infty}
  N_\infty(s)(\{\tau\mid  \At[\langle s,
  \tau\rangle]{t} \in A\}
  $
  is a measurable map $Q_A \rightarrow \pReal$.
\end{enumerate}
\EndCorollary

\BeginProof 
Since $ s \in Q_A $ iff
\begin{equation*}
\liminf_{t \rightarrow \infty}
  N_\infty(x)(\{\tau\mid  \At[\langle s,
  \tau\rangle]{t} \in A\})
= \limsup_{t \rightarrow \infty}
  N_\infty(x)(\{\tau\mid  \At[\langle s,
  \tau\rangle]{t} \in A\}),
\end{equation*}
and since the set on which two Borel measurable maps coincide is a
Borel set itself, the first assertion follows from
Proposition~\ref{SomeMeas}, part$~\ref{p-item:3}$. This  implies the second assertion as well.
\EndProof

When dealing with the semantics of the until operator later, we will also need to establish measurability of certain sets. Preparing for that, we state:

\BeginLemma{UntilMeas} 
Assume that $A_1$ and
$A_2$ are Borel subsets of $S$, and let $I \subseteq \pReal$ be an
interval, then
\begin{equation*}
U(I, A_1, A_2) := \{\sigma\mid  \exists t \in I: \At{t} \in A_2
\wedge \forall t'\in [0, t[: \At{t'} \in A_1\}
\end{equation*}
is a measurable set of paths, thus $ U(I, A_1, A_2) \in
\Borel{\Paths}. $
\EndLemma

\BeginProof 
0. Remember that, given a path $\sigma$ and a time $t
\in \pReal$, there exists a rational time $t_r \leq t$ with $ \At{t}
= \At{t_r}. $ Consequently,
\begin{equation*}
U(I, A_1, A_2) =  \bigcup_{t\in\Rational \cap I} \bigl( \{\sigma
\mid  \At{t} \in A_2\} \cap \bigcap_{t' \in \Rational \cap [0, t]}
\{\sigma \mid  \At{t'} \in A_1\} \bigr).
\end{equation*}
The inner intersection is countable and is performed over measurable
sets by Lemma~\ref{ein-hilfssatz}, thus forming a measurable set of paths.
Intersecting it with a measurable set and forming a countable union
yields a measurable set again. 
\EndProof

Now that we know how to probabilistically describe the behavior of
paths, we are ready for a probabilistic interpretation of
\texttt{CSL}. We have started from the assumption that the one-step
behavior is governed through a stochastic relation $M: S \Trans \pReal
\times S$ with $M(s)(\pReal \times S) = 1$ for all $s \in S$ from
which the stochastic relation $ \Inff{M}: S \Trans \pReal \times
\Paths $ has been constructed. The interpretations for the formulas
can be established now, and we show that the sets of states resp.
paths on which formulas are valid are Borel measurable.

To get started on the formal definition of the semantics, we assume
that we know for each atomic proposition which state it is satisfied
in. Thus we fix a map $\ell$ that maps $P$ to $\Borel{S}$, assigning
each atomic proposition a Borel set of states.

The semantics is described as usual recursively through relation
$\models$ between states resp. paths, and formulas. Hence $s\models\phi$ means that state formula $\phi$ holds in state $s$, and $\sigma\models\psi$ means that path formula $\psi$ is true on path $\sigma$. 

Here we go:
\begin{enumerate}
 \item
  $s \models \top$ is true for all $s \in S$.
 \item
  $s \models a$ iff $s \in \ell(a)$.
\item
 $s \models \varphi_1 \wedge \varphi_2$ iff $s \models \varphi_1$ and $s
 \models \varphi_2$.
\item
 $s \models \neg \varphi$ iff $s \models \varphi$ is false.
\item
 $s \models \Steady{\varphi}$ iff
 $
 \lim_{t \rightarrow \infty}\Inff{M}(s)(\{\tau \mid
 \At[\langle s, \tau\rangle]{t} \models \varphi \})
 $
 exists and is $\RelOp p$.
\item
 $s \models \PathQuant{\psi}$ iff
 $
 \Inff{M}(s)(\{\tau \mid  \langle s, \tau\rangle \models \psi\}) \RelOp p.
 $
\item
 $\sigma \models \Next{\varphi}$ iff $
 \sigma[1] \models \varphi \text{ and } \delta(\sigma, 0) \in I.
 $
\item
 $\sigma \models \Until{\varphi_1}{\varphi_2}$ iff
 $
 \exists t \in I: \At{t} \models \varphi_2 \text{ and }
 \forall t' \in [0, t[: \At{t'} \models \varphi_1.
 $
\end{enumerate}

Most interpretations should be obvious. Given a state $s$, we say that
$s\models\Steady{\phi}$ iff the asymptotic behavior of the paths
starting at $s$ gets eventually stable with a limiting probability
given by $\RelOp p$. Similarly, $s \models \PathQuant{\psi}$ holds iff
the probability that path formula $\psi$ holds for all $s$-paths is
specified through $\RelOp p$. For $\langle s_{0}, t_{0}, s_{1}, \dots,
\rangle \models \Next{\varphi}$ to hold we require $s_{1}\models\phi$
after a waiting time $t_{0}$ for the transition to be a member of
interval $I$. Finally, $\sigma \models \Until{\varphi_1}{\varphi_2}$
holds iff we can find a time point $t$ in the interval $I$ such that
the corresponding state $\At{t}$ satisfies $\phi_{2}$, and for all
states on that path before $t$, formula $\phi_{1}$ is assumed to hold.
The kinship to \texttt{CTL} is obvious.

Denote by $ \Gilt $ and $\Gilt[\psi]$ the set of all
states for which the state formula $\varphi$ holds, resp. the set of
all paths for which the path formula $\varphi$ is valid. We do not
distinguish notationally between these sets, as far as the basic
domains are concerned, since it should always be clear whether we
describe a state formula or a path formula.

We show that we are dealing with measurable sets. Most of the work
for establishing this has been done already. What remains to be done is to fit in the
patterns that we have set up in Proposition~\ref{SomeMeas} and its Corollaries.

\BeginProposition{AllMeasurable} 
The set $\Gilt[\xi]$ is  Borel, whenever
$\xi$ is a state formula or a path formula. 
\EndProposition

\BeginProof 0. The proof proceeds by induction on the structure of
the formula $\xi$. The induction starts with the formula $\top$,
for which the assertion is true, and with the atomic propositions,
for which the assertion follows from the assumption on $\ell$: $
\Gilt[a] = \ell(a) \in \Borel{S}. $ We assume for the induction step
that we have established  that $\Gilt, \Gilt[\varphi_1] \text{ and }
\Gilt[\varphi_2]$ are Borel measurable.

1. For the next-operator we write
\begin{equation*}
\Gilt[\Next{\varphi}] = \{\sigma \mid  \sigma[1] \in \Gilt \text{
and } \delta(\sigma, 0) \in I\}.
\end{equation*}
This is the cylinder set $ \left(S \times I \times \Gilt \times
\pReal\right) \times \Paths, $ hence is a Borel set.

2. The until-operator may be represented through
\begin{equation*}
\Gilt[\Until{\varphi_1}{\varphi_2}] = U(I, \Gilt[\varphi_1],
\Gilt[\varphi_2]),
\end{equation*}
which is a Borel set by Lemma~\ref{UntilMeas}.

3. Since $\Inff{M}: S \Trans (\pReal \times S)^\infty$ is a
stochastic relation, we know that
\begin{equation*}
\Gilt[\PathQuant{\psi}] = \{s \in S \mid  \Inff{M}(s)(\{\tau \mid
\langle s, \tau\rangle \in \Gilt\}) \RelOp p\}
\end{equation*}
is a Borel set.

4. We know from Corollary~\ref{LimitBehavior} that the set
\begin{equation*}
  Q_{\Gilt} := \{s \in S \mid  \lim_{t \rightarrow \infty}
  \Inff{M}(s)(\{\tau \mid  \At[\langle s,
  \tau\rangle]{t} \in \Gilt\})
  \text{ exists}\}
\end{equation*}
is a Borel set, and that
\begin{equation*}
  J_\varphi: Q_{\Gilt} \ni s \mapsto \lim_{t \rightarrow \infty}
  \Inff{M}(x)\left(\{\tau \mid  \At[\langle s,
  \tau\rangle]{t} \in \Gilt\}\right) \in [0, 1]
\end{equation*}
is a Borel measurable function. Consequently,
\begin{equation*}
\Gilt[\Steady{\varphi}] = \{s \in Q_{\Gilt}\mid  J_\varphi(s)
\RelOp p\}
\end{equation*}
is a Borel set. 
\EndProof

Measurability of the sets on which a given formula is valid constitutes of
course a prerequisite for computing interesting properties. So we
can compute, e.g.,
$$
\PathQuant[\geq 0.5]{\Until[{[10, 20]}]{(\neg down)}{\Steady[\geq
    0.8]{up_2 \vee up_3}}})
$$
as the set of all states that with probability at least $0.5$ will
reach a state between $10$ and $20$ time units so that the system is
operational ($up_2, up_3 \in P$) in a steady state with a
probability of at least $0.8$; prior to reaching this state, the
system must be operational continuously ($down \in P$).

The description of the semantics is just the basis for entering into the investigation of expressivity of the models associated with $M$ and with $\ell$. We leave CSL here, however, and note that the construction of the projective limit is the basic ingredient for further investigations.

%%% Local Variables: 
%%% mode: latex
%%% TeX-master: "../Mskr3"
%%% End: 

%\Input{\Folder/WeakTop}
%spell checked - 24Aug14
\Subsection{The Weak Topology}
\label{sec:weak-top}

Now that we have integration at our disposal, we will look again at topological issues for the space of finite measures. 
We fix in this section $(X, d)$ as a metric space; recall that ${\cal C}_{b}(X)$ is the space of all bounded continuous functions $X\to \Real$. This space induces the weak  topology on the space $\FinM{X} = \FinM{X, \Borel{X}}$ of all finite Borel measures on $(X, \Borel{X})$. This is the smallest topology which renders the evaluation map
\begin{equation*}
  \mu\mapsto \dInt{f}{\mu}
\end{equation*}
continuous for every continuous and bounded map $f: X\to \Real$. This topology is fairly natural, and it is related to the topologies on $\FinM{X}$ considered so far, the Alexandrov topology and the topology given by the Levy-Prohorov metric, which are discussed in Section~\ref{sec:topology-on-measures}. We will show that these topologies are the same, provided the underlying space is Polish, and we will show in this case show that $\FinM{X}$ is itself a Polish space. Somewhat weaker results may be obtained if the base space is only separable metric, and it turns out that tightness, i.e., inner approximability through compact sets, the the property which sets Polish spaces apart for our purposes. We introduce also a very handy metric for the weak topology due to Hutchinson. Two case studies on bisimulations of Markov transition systems and on quotients for stochastic relations demonstrate the interplay of topological considerations with selection arguments, which become available on $\FinM{X}$ once this space is identified as Polish.
 
Define as the basis for the topology the sets
\begin{equation*}\textstyle
  U_{f_{1}, \dots, f_{n}, \epsilon}(\mu) := \{\nu\in\FinM{X}\mid \bigl|\dInt{f_{i}}{\nu} - \dInt{f_{i}}{\mu}\bigr|<\epsilon \text{ for } 1 \leq i \leq n\}
\end{equation*}
with $\epsilon>0$ and $f_{1}, \dots, f_{n}\in {\cal C}_{b}(X)$. Call the topology the \emph{\index{weak topology}weak \index{topology!weak}topology} on $\FinM{X}$.

With respect to convergence, we have this characterization, which indicates the relationship between the weak topology and the Alexandrov-topology investigated in Section~\ref{sec:topology-on-measures}.

\BeginTheorem{portmanteau-weak}
The following statements are equivalent for a sequence $\Folge{\mu}\subseteq\FinM{X}$.
\begin{enumerate}
\item\label{port:1}
$\mu_{n}\to \mu$ in the weak topology.
\item\label{port:2} $\dInt{f}{\mu_{n}}\to \dInt{f}{\mu}$ for all $f\in{\cal C}_{b}(X)$.
\item\label{port:3} $\dInt{f}{\mu_{n}}\to \dInt{f}{\mu}$ for all bounded and uniformly continuous $f: X\to \Real$.
\item\label{port:4} $\mu_{n}\to \mu$ in the A-topology. 
\end{enumerate}
\EndTheorem

\BeginProof
The implications \labelImpl{port:1}{port:2} and \labelImpl{port:2}{port:3} are trivial.

\labelImpl{port:3}{port:4}: Let $G\subseteq X$ be open, then $f_{k}(x) := 1\wedge k\cdot d(x, X\setminus G)$ defines a uniformly continuous map (because $|d(x, X\setminus G) - d(x', X\setminus G)| \leq d(x, x')$), and $0\leq f_{1}\leq f_{2}\leq \dots$ with $\lim_{k\to \infty}f_{k} = \chi_{G}$. Hence $\dInt{f_{k}}{\mu} \leq \dInt{\chi_{G}}{\mu} = \mu(G)$, and  by monotone convergence $\dInt{f_{k}}{\mu}\to \mu(G)$.  From the assumption we know that $\dInt{f_{k}}{\mu_{n}}\to \dInt{f_{k}}{\mu}$, as $n\to \infty$, so that we obtain for all $k\in \Nat$
\begin{equation*}
  \lim_{n\to \infty}\dInt{f_{k}}{\mu_{n}}\leq \liminf_{n\to \infty}\mu_{n}(G),
\end{equation*}
with in turn implies $\mu(G) \leq \liminf_{n\to \infty}\mu_{n}(G).$

\labelImpl{port:4}{port:2}
We may assume that $f \geq 0$, because the integral is linear.
Then we can represent the
integral through (see Example~\ref{choquet-example}, equation~(\ref{eq-choquet}))
\begin{equation*}
\int_X f\ d\nu = \int_0^\infty \nu(\{x \in X \mid f(x) > t\})\ dt.
\end{equation*}
Since $f$ is continuous, the set
$
\{x \in X \mid f(x) > t\}
$
is open. By Fatou's Lemma (Proposition~\ref{fatous-lemma})
we obtain from the assumption
\begin{align*}
\liminf_{n\rightarrow\infty} \int_X f\ d\mu_n
&= \liminf_{n\rightarrow\infty} \int_0^\infty \mu_n(\{x \in X \mid f(x) > t\})\ dt \\
&\geq
\int_0^\infty \liminf_{n\rightarrow\infty} \mu_n(\{x \in X \mid f(x) > t\})\ dt \\
&\geq
\int_0^\infty \mu(\{x \in X \mid f(x) > t\})\ dt \\
&=
\int_X f\ d\mu.
\end{align*}
Because $f \geq 0$ is bounded, we find $T \in \Real$ such that $f(x) \leq T$
for all $x \in X$, hence $g(x) := T - f(x)$ defines a non-negative and bounded
function. Then by the preceding argument
$
\liminf_{n\rightarrow\infty} \int_X g\ d\mu_n \geq \int_X g\ d\mu.
$
Since $\mu_n(X) \rightarrow \mu(X)$, we infer
\begin{equation*}
\limsup_{n\rightarrow\infty} \int_X f\ d\mu_n \leq \int_X f\ d\mu,
\end{equation*}
which implies the desired equality.
\EndProof

Let $X$ be separable, then the A-topology is metrized by the Prohorov
metric (Theorem~\ref{a-top-is-metrizable}). Thus we have established
that the metric topology and the topology of weak convergence are the
same for separable metric spaces. Just for the record:

\BeginTheorem{Metric-Equiv-Weak}
Let $X$ be a separable metric space, then
the Prohorov metric is a metric for  the  topology of weak convergence.
\QED
\EndTheorem

It is now easy to find a dense subset in $\FinM{X}$. As one might expect,
the measures living on discrete subsets are dense. Before stating and proving the corresponding statement, we have a brief look at the embedding of $X$ into $\FinM{X}$.

\BeginExample{embed-x-into-all-finite}
The base space $X$ is embedded into $\FinM{X}$ as a closed subset
through $x\mapsto\delta_{x}$. In fact, let
$(\delta_{x_{n}})_{n\in\Nat}$ be a sequence which converges weakly to
$\mu\in\FinM{X}$. We have in particular $\mu(X) = \lim_{n\to
  \infty}\delta_{x_{n}}(X) = 1$, hence $\mu\in\Prob{X}$. Now assume
that $\Folge{x}$ does not converge, hence it does not have a
convergent subsequence in $X$. Then the set $S := \{x_{n}\mid
n\in\Nat\}$ is closed in $X$, so are all subsets of $S$. Take an
infinite subset $C\subseteq S$ with an infinite complement $S\setminus
C$, then $\mu(C) \geq \limsup_{n\to \infty}\delta_{x_{n}}(C) = 1$, and
with the same argument $\mu(S\setminus C) = 1$. This contradicts
$\mu(X) = 1$. Thus we can find $x\in X$ with $x_{n}\to x$, hence
$\delta_{x_{n}}\to \delta_{x}$, so that the image of $X$ in $\FinM{X}$
is closed.
\EndExample

\BeginProposition{Discrete-is-dense}
Let $X$ be a separable metric space.
The set
\begin{equation*}\textstyle
\bigl\{\sum_{k \in \Nat} r_k\cdot \delta_{x_k} \mid x_k \in X, r_k \geq 0\bigr\}
\end{equation*}
of discrete measures  is dense in the topology of weak convergence.
\EndProposition

\BeginProof
Fix $\mu \in \FinM{X}$. Cover $X$ for each $k \in  \Nat$ with mutually disjoint Borel sets
$(A_{n, k})_{n \in \Nat}$, each of which has a diameter not less that
$1/k$. Select an arbitrary $x_{n, k} \in A_{n, k}$. We claim that
$
\mu_n := \sum_{k \in \Nat} \mu(A_{n, k})\cdot \delta_{x_{n, k}}
$
converges weakly to $\mu$. In fact, let $f: X \rightarrow \Real$ be a uniformly continuous and bounded map.
Since $f$ is uniformly continuous,
\begin{equation*}
\eta_n := \sup_{k \in \Nat}\bigl(\sup_{x \in A_{n, k}} f(x) - \inf_{x \in A_{n, k}} f(x)\bigr)
\end{equation*}
tends to $0$, as $n \rightarrow\infty$. Thus
\begin{align*}
  \left| \int_X f\ d\mu_n - \int_X f\ d\mu\right|
  &= 
  \bigl| \sum_{k \in \Nat} \bigl(\int_{A_{n, k}} f\ d\mu_n - \int_{A_{n, k}} f\ d\mu\bigr)\bigr| \\
    &\leq  \eta_n\cdot\sum_{k \in \Nat} \mu(A_{n, k})\\
    &\leq  \eta_n\\
    &\rightarrow  0.
\end{align*}
\EndProof

This yields immediately

\BeginCorollary{Weak-is-separable}
If $X$ is a separable metric space, then $\FinM{X}$ is a separable
metric space in the topology of weak convergence.
\EndCorollary

\BeginProof
Because
$
\sum_{k=1}^{n} r_k\cdot\delta_{x_k}\to  \sum_{k \in \Nat} r_k\cdot\delta_{x_k},
$
as $n\to  \infty$ in the weak topology, and because the rationals $\Rational$ are
dense in the reals, we obtain from
Proposition~\ref{Discrete-is-dense} that
$%\begin{equation*}
\bigl\{\sum_{k=1}^{n} r_k\cdot\delta_{x_k}\mid x_k \in D, 0 \leq r_k \in \Rational, n \in \Nat\bigr\}
$ %\end{equation*}
is a countable and dense subset of $\FinM{X}$,
whenever $D \subseteq X$ is a countable and dense subset of $X$.
\EndProof

Another immediate consequence refers to the
weak-*-$\sigma$-algebra. We obtain from
Lemma~\ref{weak-star-vs-a-borel} together with
Corollary~\ref{Weak-is-separable}

\BeginCorollary{weak-star-is-borel}
Let $X$ be a metric space, then the weak-*-$\sigma$-algebra are the Borel sets of the A-topology.
\QED
\EndCorollary

We will show now that $\FinM{X}$ is a Polish space,
provided $X$ is one; thus applying the $\FinSenza$-functor
to a Polish space does not leave the realm of Polish
spaces. 

We know by Alexandrov's Theorem~\ref{Alexandrov} that a separable
metrizable space is Polish iff it can be embedded as a $G_\delta$-set
into the Hilbert cube. We show first that for compact metric $X$ the
space $\SubProb{X}$ of all subprobability measures with the topology
of weak convergence is itself a compact metric space. This is
established by embedding it as a closed subspace into $[-1,
+1]^\infty$. But there is nothing special about taking
$\SubProbSenza$; the important property is that all measures are
uniformly bounded (by $1$, in this case). Any other bound would also
do.

We require for this the Stone-Weierstraß Theorem which states (in the
form needed here) that the unit ball in the space of all bounded
continuous functions on a compact metric space is separable
itself~\cite[Chapter 7, Problem S~(e), p.~245]{Kelley}. The idea of
the embedding is to take a countable dense sequence $\Folge{g}$ of this
unit ball. Since we are dealing with probability measures, and since
we know that each $g_{n}$ maps $X$ into the interval $[-1, 1]$, we
know that $-1 \leq \dInt{g_{n}}{\mu} \leq 1$ for each $\mu$. This then
spawns the desired map, which together with its inverse is shown
through the Riesz Representation Theorem to be continuous.

Well, this is the plan of attack for establishing

\BeginProposition{compact-embeds}
Let $X$ be a compact metric space. Then $\SubProb{X}$ is a
compact metric space.
\EndProposition

\BeginProof
1.  The space ${\cal C}(X)$ of continuous maps
into the reals is for compact metric $X$ a separable Banach space
under the $\sup$-norm $\|\cdot \|_{\infty}$. The closed unit ball
\begin{equation*}
\mathbf{C}_1 := \{f \in {\cal C}(X)\mid \|f\|_\infty \leq 1\}
\end{equation*}
is, as mentioned above, a separable metric space in its own right. Let
$\Folge{g}$ be a countable sense subset in $\mathbf{C}_1$, and define
\begin{equation*}\textstyle
\Theta: \SubProb{X} \ni \nu
\mapsto
\langle\dInt{g_{1}}{\nu}, \dInt{g_{2}}{\nu}, \dots\rangle \in [-1, 1]^\infty.
\end{equation*}
Then $\Theta$ is injective, because the sequence $\Folge{g}$ is dense. 

2.
Also, $\Theta^{-1}$ is continuous. In fact, let $\Folge{\mu}$ be a sequence in $\SubProb{X}$ such that
$\bigl(\Theta(\mu_n)\bigr)_{n \in \Nat}$ converges in $[-1, 1]^\infty,$
put
$
\alpha_i := \lim_{n  \rightarrow\infty} \int_X g_i\ d\mu_n.
$
For each $f \in \mathbf{C}_1$ there exists a subsequence $(g_{n_k})_{k \in \Nat}$
such that
$
\| f - g_{n_k}\|_\infty \rightarrow 0
$
as
$
k \rightarrow \infty,
$
because $\Folge{g}$ is dense in $\mathbf{C}_{1}$. Thus
\begin{equation*}
L(f) := \lim_{n  \rightarrow\infty} \int_X f\ d\mu_n
\end{equation*}
exists. Define 
$
L(\alpha\cdot f) := \alpha\cdot L(f),
$
for $\alpha \in \Real$, then it is immediate that $L:{\cal C}(X)\to  \Real$
is linear and that $L(f) \geq 0$, provided $f \geq 0$. 
The Riesz Representation Theorem~\ref{riesz-representation} now
gives a unique $\mu \in \SubProb{X}$ with
\begin{equation*}
L(f) = \int_X f\ d\mu,
\end{equation*}
and the construction shows that
\begin{equation*}
\lim_{n  \rightarrow\infty}\Theta(\mu_n) = \langle\int_X g_1\ d\mu, \int_X g_2\ d\mu, \dots\rangle.
\end{equation*}

3.
Consequently, $\Theta: \SubProb{X}\to \Bild{\Theta}{\SubProb{X}}$ is a homeomorphism, and $\Bild{\Theta}{\SubProb{X}}$ is closed, hence compact. Thus $\SubProb{X}$ is compact.
\EndProof

We obtain as a first consequence
\BeginProposition{compact-iff}
$X$ is compact iff $\SubProb{X}$ is, whenever $X$ is a Polish space.
\EndProposition

\BeginProof
It remains to show that $X$ is compact, provided $\SubProb{X}$
is. Choose a complete metric $d$ for $X$. Thus $X$ is isometrically
embedded into $\SubProb{X}$ by $x\mapsto \delta_{x}$ with $A :=
\{\delta_{x}\mid x\in X\}$ being closed. We could appeal to
Example~\ref{embed-x-into-all-finite}, but a direct argument is
available as well. In fact, if $\delta_{x_{n}}\to \mu$ in the weak
topology, then $\Folge{x}$ is a Cauchy sequence in $X$ on account of
the isometry. Since $(X, d)$ is complete, $x_{n}\to x$ for some $x\in
X$, hence $\mu=\delta_{x}$, thus $A$ is closed, hence compact.
\EndProof

The next step for showing that $\FinM{X}$ is Polish is nearly canonical. If $X$ is a
Polish space, it may be embedded as a $G_\delta$-set into a compact
space $\widetilde{X}$, the subprobabilities of which are topologically
a closed subset of $[-1, +1]^\infty$, as we have just seen. We will
show now that $\FinM{X}$ is a $G_\delta$ in $\FinM{\widetilde{X}}$ as
well.

\BeginProposition{it-is-itself-polish}
Let $X$ be a Polish space. Then $\FinM{X}$ is a Polish
space in the topology of weak convergence.
\EndProposition

\BeginProof
1.
Embed $X$ as a $G_\delta$-subset into a compact metric space
$\widetilde{X}$, hence $X \in \Borel{\widetilde{X}}$. Put
\begin{equation*}
\FinSenza_0 :=
\{\mu \in \FinM{\widetilde{X}} \mid \mu(\widetilde{X}\setminus X) = 0\},
\end{equation*}
so $\FinSenza_0$ contains exactly those finite measures on $\widetilde{X}$ that
are concentrated on $X$. Then $\FinSenza_0$ is homeomorphic to $\FinM{X}$.

2.  
Write $X$ as $X =\bigcap_{n \in \Nat} G_n$, where $\Folge{G}$ is a
sequence of open sets in $\widetilde{X}$. Given $r > 0$, the set
\begin{equation*}
\Gamma_{k, r} := \{\mu \in \FinM{\widetilde{X}} \mid \mu(\widetilde{X}\setminus G_k) < r\}
\end{equation*}
is open in $\FinM{\widetilde{X}}$. In fact, if $\mu_n \notin
\Gamma_{k, r}$ converges to $\mu_{0}$ in the weak topology, then
\begin{equation*}
\mu_0(\widetilde{X}\setminus G_k) \geq
\limsup_{n \rightarrow\infty} \mu_n(\widetilde{X}\setminus G_k)
\geq r
\end{equation*}
by Theorem~\ref{portmanteau-weak}, since $\widetilde{X}\setminus G_k$
is closed. Consequently,
$
\mu_0 \notin \Gamma_{k, r}.
$
This shows that $\Gamma_{k, r}$ is open, because its complement is closed. Thus
\begin{equation*}
\FinSenza_0 = \bigcap_{n \in \Nat} \bigcap_{k \in \Nat} \Gamma_{n, 1/k}
\end{equation*}
is a $G_\delta$-set, and the assertion follows.
\EndProof

Thus we obtain as a consequence

\BeginProposition{polish-iff-polish}
$\FinM{X}$ is a Polish space in the topology of weak convergence iff $X$ is. 
\EndProposition

\BeginProof
Let $\FinM{X}$ be Polish. The base space $X$ is embedded into
$\FinM{X}$ as a closed subset by
Example~\ref{embed-x-into-all-finite}, hence is a Polish space by Theorem~\ref{Polish-is-G-delta}.
\EndProof

Let $\mu\in\FinM{X}$ with $X$ Polish. Since $X$ has a countable basis,
we know from Lemma~\ref{support-for-regular} that $\mu$ is supported
by a closed set, since $\mu$ is $\tau$-regular. But in the presence of
a complete metric we can say a bit more, viz., that the value of
$\mu(A)$ may be approximated from within by compact sets to arbitrary
precision.

\BeginDefinition{tightness}
A finite Borel measure $\mu$ is called \emph{\index{measure!tight}tight} iff 
\begin{equation*}
  \mu(A) = \sup\{\mu(K) \mid  K\subseteq A \text{ compact}\}
\end{equation*}
holds for all $A\in\Borel{X}$. 
\EndDefinition

Thus tightness means for $\mu$ that we can find for any $\epsilon>0$ and for any Borel set $A\subseteq X$ a compact set $K\subseteq A$ with $\mu(A\setminus K) < \epsilon$. Because a finite measure on a separable metric space is regular, i.e., $\mu(A)$ can be approximated from within $A$ by closed sets (Lemma~\ref{are-regular}), it suffices in this case to consider tightness at $X$, hence to postulate that there exists for any $\epsilon>0$ a compact set $K\subseteq X$ with $\mu(X\setminus K) < \epsilon$. We know in addition that each finite measure is $\tau$-regular by Lemma~\ref{regular-tau-in-separable}. Capitalizing on this and on completeness, we find

\BeginProposition{tight-on-polish}
Each finite Borel measure on Polish space $X$ is tight.
\EndProposition

\BeginProof
1.
We show first that we can find for each $\epsilon>0$ a compact set $K\subseteq X$ with $\mu(X\setminus K) < \epsilon$. In fact, given a complete metric $d$, consider 
\begin{equation*}
  {\cal G} := \bigl\{\{x\in X \mid  d(x, M) < 1/n\}\mid M\subseteq X\text{ is finite}\bigr\}.
\end{equation*}
Then ${\cal G}$ is a directed collection of open sets with $\bigcup{\cal G} = X$, thus we know from $\tau$-regularity of $\mu$ that $\mu(X) = \sup\{\mu(G)\mid G\in{\cal G}\}$. Consequently, given $\epsilon>0$ there exists for each $n\in\Nat$ a finite set $M_{n}\subseteq X$ with $\mu(\{x\in X\mid d(x, M_{n})<1/n\})>\mu(X)-\epsilon/2^{n}$. Now define
\begin{equation*}
  K := \bigcap_{n\in\Nat}\{x\in X\mid d(x, M_{n}\leq 1/n\}.
\end{equation*}

Then $K$ is closed, and complete (since $(X, d)$ is complete). Because each $M_{n}$ is finite, $K$ is totally bounded. Thus $K$ is compact. We obtain
\begin{equation*}
  \mu(X\setminus K) \leq \sum_{n\in\Nat}\mu(\{x\in X\mid d(x, M_{n})\geq 1/n\}) \leq \sum_{n\in\Nat}\epsilon\cdot2^{-n} = \epsilon.
\end{equation*}

2.
Now let $A\in\Borel{X}$, then for $\epsilon>0$ there exists $F\subseteq A$ closed with $\mu(A\setminus F) < \epsilon/2$, and chose $K\subseteq X$ compact with $\mu(X\setminus K) < \epsilon/2$. Then $K\cap F\subseteq A$ is compact with $\mu(A\setminus(F\cap K)) < \epsilon$. 
\EndProof

Tightness is sometimes an essential ingredient when arguing about
measures on a Polish space. The discussion on the Hutchinson metric in
the next section provides an example, it shows that at crucial point
tightness kicks in and saves the day.

%%% Local Variables:  
%%% mode: latex
%%% TeX-master: "../Mskr3"
%%% End: 

%\Input{\Folder/HutchMetric}
%spell checked - 24Aug14
\Subsubsection{The Hutchinson Metric}
\label{sec:hutchinson-metric}

We will explore now another approach to the weak topology for Polish spaces through the Hutchinson metric. Given a fixed metric $d$ on $X$, define 
\begin{equation*}\textstyle
  V_{\gamma}:= \{f: X\to \Real\mid |f(x)-f(y)|\leq d(x, y)\text{ and }|f(x)|\leq \gamma\text{ for all }x, y\in X\},
\end{equation*}
Thus $f$ is a member of $V_{\gamma}$ iff $f$ is non-expanding (hence has a Lipschitz constant $1$), and iff its supremum norm $\|f\|_{\infty}$ is bounded by $\gamma$. Trivally, all elements of $V_{\gamma}$ are uniformly continuous. Note the explicit dependence on the metric $d$. 
The \emph{Hutchinson distance} $H_{\gamma}(\mu, \nu)$ between $\mu, \nu\in\FinM{X}$ is defined as 
\begin{equation*}
H_{\gamma}(\mu, \nu) := \sup_{f\in V_{\gamma}}\bigl(\dInt{f}{\mu} - \dInt{f}{\nu}\bigr).
\end{equation*}
Then $H_{\gamma}$ is easily seen to be a metric on $\FinM{X}$. \index{$H_{\gamma}$}$H_{\gamma}$ is called the \emph{\index{Hutchinson metric}Hutchinson metric} (sometimes also Hutchinson-Monge-Kantorovicz metric). 

The relationship between this metric and the topology of weak convergence is stated now (\cite[Theorem 2.5.17]{Edgar}):

\BeginProposition{weak-is-hutchinson}
Let $X$ be a Polish space. Then $H_{\gamma}$ is a metric for the topology of weak convergence on $\FinM{X}$ for any $\gamma>0$. 
\EndProposition

\BeginProof
1.
We may and do assume that $\gamma = 1$, otherwise we scale accordingly. Now let $H_{1}(\mu_{n}, \mu)\to 0$ as $n\to \infty$, then $\lim_{n\to \infty}\mu_{n}(X) = \mu(X)$. Let $F\subseteq X$ be closed, then we can find for given $\epsilon>0$ a function $f\in V_{1}$ such that $f(x) = 1$ for $x\in F$, and $\dInt{f}{m}\leq \mu(F)+\epsilon$. This gives 
\begin{equation*}
  \limsup_{n\to \infty}\mu_{n}(F) \leq \lim_{n\to \infty}\dInt{f}{\mu_{n}} = \dInt{f}{\mu} \leq \mu(F) + \epsilon
\end{equation*}
Thus convergence in the Hutchinson metric implies convergence in the A-topology, hence in the topology of weak convergence, by Proposition~\ref{convergence-a-topology}. 

2.
Now assume that $\mu_{n}\to \mu$ in the topology of weak convergence, thus $\mu_{n}(A)\to \mu(A)$ for all $A\in\Borel{X}$ with $\mu(\partial A) = 0$ by Corollary~\ref{conv-a-top-boundary}; we assume that $\mu_{n}$ and $\mu$ are probability measures, otherwise we scale again. Because $X$ is Polish, $\mu$ is tight by Proposition~\ref{tight-on-polish}. 

Fix $\epsilon>0$, then there exists a compact set $K\subseteq X$ with 
\begin{equation*}
\mu(X\setminus K) < \frac{\epsilon}{5\cdot\gamma}.
\end{equation*}
Given $x\in K$, there exists an open ball $B_{r}(x)$ with center $x$ and radius $r$ such that $0< r < \epsilon/10$ such that $\mu(\partial B_{r}(x))=0$, see Corollary~\ref{countably-many-point-values}. Because $K$ is compact, a finite number of these balls will suffice, thus $K \subseteq B_{r_{1}}(x_{1})\cup\dots\cup B_{r_{p}}(x_{p})$. Transform this cover into a disjoint cover by setting
\begin{align*}
  E_{1} & := B_{r_{1}}(x_{1}),\\
E_{2} & := B_{r_{2}}(x_{2})\setminus E_{1}, \\
\dots\\
E_{p} & := B_{r_{p}}(x_{p})\setminus (E_{1}\cup\dots\cup E_{p-1})\\
E_{0} & := S\setminus (E_{1}\cup\dots\cup E_{p})
\end{align*}
We observe these properties:
\begin{enumerate}
\item For $i = 1, \dots, p$, the diameter of each $E_{i}$ is not greater than $2\cdot r_{i}$, hence smaller that $\epsilon/5$,
\item For $i = 1, \dots, p$, $\partial E_{i}\subseteq \partial \bigl(B_{r_{1}}(x_{1})\cup\dots\cup B_{r_{p}}(x_{p})\bigr)$, thus $\partial E_{i}\subseteq (\partial B_{r_{1}}(x_{1}))\cup\dots\cup(\partial B_{r_{p}}(x_{p}))$, hence $\mu(\partial E_{i}) = 0$. 
\item Because the boundary of a set is also the boundary of its complement, we conclude $\mu(\partial E_{0}) = 0$ as well. Moreover, $\mu(E_{0})<\epsilon/(5\cdot\gamma)$, since $E_{0}\subseteq X\setminus K$. 
\end{enumerate}

Eliminate all $E_{i}$ which are empty. Select $\eta>0$ such that $p\cdot\eta < \epsilon/5$, and determine $n_{0}\in\Nat$ so that $|\mu_{n}(E_{i})-\mu(E_{i})|<\eta$ for $i = 0, \dots, p$ and $n\geq n_{0}$. 

We have to show that 
\begin{equation*}
  \sup_{f\in V_{\gamma}}\bigl(\dInt{f}{\mu_{n}}-\dInt{f}{\mu}\bigr)\to 0, \text{ as }n\to \infty.
\end{equation*}
So take $f\in V_{\gamma}$ and fix $n\geq n_{0}$. Let $i = 1,\dots, p$, pick an arbitrary $e_{i}\in E_{i}$; because each $E_{i}$ has a diameter not greater than $\epsilon/5$, we know that 
$
|f(x) - f(e_{i})| < \epsilon/5
$
for each $x\in E_{i}$. If $x\in E_{0}$, we have $|f(x)|\leq \gamma$. Now we are getting somewhere: let $n\geq n_{0}$, then we obtain
\begin{align*}
  \dInt{f}{\mu_{n}} 
& = \sum_{i=0}^{p}\dInt[E_{i}]{f}{\mu_{n}}\\
& \leq \gamma\cdot\mu_{n}(E_{0}) + \sum_{i=1}^{p}\bigl(f(t_{i}) + \frac{\epsilon}{5}\bigr)\cdot\mu_{n}(E_{i})\\
& \leq \gamma\cdot(\mu(E_{0}) + \eta) + \sum_{i=1}^{p}\bigl(f(t_{i}) + \frac{\epsilon}{5}\bigr)\cdot(\mu(E_{i})+\eta)\\
& \leq \gamma\cdot(\frac{\epsilon}{5\cdot\gamma} + \eta) + \sum_{i=1}^{p}(f(t_{i})-\frac{\epsilon}{5})\cdot\mu(E_{i}) + 
\frac{2\cdot\epsilon}{5}\sum_{i=1}^{p}\mu(E_{i}) + \frac{p\cdot\epsilon\cdot\eta}{5}\\
& \leq \dInt{f}{\mu} + \epsilon
\end{align*}
Recall that 
\begin{equation*}
\sum_{i=1}^{p}\mu(E_{i}) \leq \sum_{i=0}^{p}\mu(E_{i}) = \mu(X) = 1,
\end{equation*}
and that 
\begin{equation*}
\dInt[E_{i}]{f}{\mu} \geq \mu(E_{i})\cdot(f(t_{i})-\epsilon/5).
\end{equation*}
In a similar fashion, we obtain
$
\dInt{f}{\mu_{n}} \geq \dInt{f}{\mu} - \epsilon$, so that we have established 
\begin{equation*}
  |\dInt{f}{\mu}-\dInt{f}{\mu_{n}}|<\epsilon
\end{equation*}
for $n\geq n_{0}$. Since $f\in V_{\gamma}$ was arbitrary, we have shown that $H_{\gamma}(\mu_{n}, \mu)\to 0$. 
\EndProof

The Hutchinson metric is sometimes easier to use that the Prohorov
metric, because integrals may sometimes easier manipulated in
convergence arguments than $\epsilon$-neighborhoods of sets.

%%% Local Variables: 
%%% mode: latex
%%% TeX-master: "../Mskr3"
%%% End: 

%\Input{\Folder/CaseBisim}
%spell checked - 24Aug14
\Subsubsection{Case Study: Bisimulation}
\label{sec:case-bisim}

Bisimilarity is an important notion in the theory of concurrent
systems, introduced originally by Milner for transition systems, see~\CategCite{Section 1.6.1} for a general discussion. We
will show in this section that the methods developed so far may be
used in the investigation of bisimilarity for stochastic systems. We
will first show that the category of stochastic relations has
semi-pullbacks and use this information for a construction of
bisimulations for these systems. 

If we are in a general category $\catK$, then the semi-pullback for two morphisms $f: a\to c$ and $g: b\to c$ with common range $c$ consists of an object $x$ and of morphisms $p_{a}: x\to a$ and $p_{b}: x\to b$ such that $f\circ p_{a} = g\circ p_{b}$, i.e., such that this diagram commutes in $\catK$:
  \begin{equation*}
\xymatrix{
x\ar[d]_{p_{a}}\ar[rr]^{p_{b}}&&b\ar[d]^{g}\\
a\ar[rr]_{f}&&c
}
\end{equation*}
We want to show that semi-pullbacks exist for stochastic relations over Polish spaces. This requires some preparations, provided through selection arguments. The next statement appears to be interesting in its own right; it shows that a measurable selection for weakly continuous stochastic relations exist.  

\BeginProposition{6-selection}
Let $X_i$, $Y_i$ be Polish spaces, $K_i: X_i\Trans Y_i$ be a weakly continuous stochastic relation, $i = 1, 2$. Let $A\subseteq X_1\times X_2$ and $B\subseteq Y_1 \times Y_2$ be closed subsets of the respective Cartesian products with projections equal to the base spaces, and assume that for $\langle x_1, x_2\rangle\in A$ the set 
\begin{equation*}
%\label{ }
\Gamma(x_1, x_2) := \{\,\mu\in\SubProb{B} \mid \SubProb{\beta_i}(\mu)
= K_i(x_i), \,i = 1, 2\,\}
\end{equation*}
is not empty, $\beta_i: B\to Y_i$ denoting the projections. Then there
exists a stochastic relation $M: A\Trans B$ such that $M(x_1,
x_2)\in\Gamma(x_1, x_2)$ for all $\langle x_1, x_2\rangle\in A$.
\EndProposition

\BeginProof
1.
Let $\comp{Y_i}$ for $i = 1, 2$ be the Alexandrov compactification of $Y_i$ and $\comp{B}$ the closure of $B$ in $\comp{Y}_1\times\comp{Y}_2$. Then $\comp{B}$ is compact and contains the embedding of $B$ into $\comp{Y}_1\times\comp{Y}_2$, which we identify with $B$, as a Borel subset. This is so since $Y_i$ is a Borel subset in its compactification. The projections $\comp{\beta}_i: \comp{B}\to \comp{Y}_i$ are the continuous extensions to the projections $\beta_i: B \to Y_i$. 

2.
The map 
$
r_i:  \SubProb{Y_i}\to \SubProb{\comp{Y}_i} 
$
with 
$
r_i(\mu)(G) := \mu(G\cap Y_i)
$
for $G\in\Borel{\comp{Y}_i}$ is continuous; in fact, it is an isometry with respect to the respective Hutchinson metrics, once we have fixed metrics for the underlying spaces. Define for $\langle x_1, x_2\rangle\in A$ the set 
\begin{equation*}
\Gamma_0(x_1, x_2) := \{\,\mu\in\SubProb{\comp{B}} \mid \SubProb{\comp{\beta}_i}(\mu) = (r_i\circ K_i)(x_i),\, i = 1, 2\,\}.
\end{equation*}
Thus $\Gamma_0$ maps $A$ to the nonempty closed subsets of $\SubProb{\comp{B}}$, since $\SubProb{\comp{\beta}_i}$ and $r_i\circ K_i$ are continuous for $i = 1, 2$. If $\mu\in\Gamma_0(x_1, x_2)$, then 
\begin{align*}
\mu(\comp{B}\setminus B) 
& \leq
\mu\bigl(\comp{B}\cap(\comp{Y}_1\setminus Y_1\times \comp{Y}_2)\cup(\comp{Y}_1\times \comp{Y}_2\setminus Y_2)\bigr)\\
& =
\SubProb{\comp{\beta}_1}(\mu)(\comp{Y}_1\setminus Y_1) + \SubProb{\comp{\beta}_2}(\mu)(\comp{Y}_2\setminus Y_2)\\
& = 
\bigl(r_1\circ K_1\bigr)(x_1)(\comp{Y}_1\setminus Y_1) + \bigl(r_2\circ K_2\bigr)(x_2)(\comp{Y}_2\setminus Y_2)\\
& =
0.
\end{align*}
Hence all members of $\Gamma_0(x_1, x_2)$ are concentrated on $B$.

3.
Let $C\subseteq \SubProb{\comp{B}}$ be compact, and assume that $(t_n)_{n\in \Nat}$ is a converging sequence in $A$ with $t_n\in\weak{\Gamma}_{0}(C)$ for all $n\in\Nat$ such that $t_n\to t_0\in A$. Then there exists some $\mu_n\in C\cap\Gamma_0(t_n)$ for each $n\in\Nat$. Since $C$ is compact, there exists a convergent subsequence, which we assume to be the sequence itself, so $\mu_n\to\mu$ for some $\mu\in C$ in the topology of weak convergence. Continuity of $\SubProb{\comp{\beta}_i}$ and of $K_i(x_i)$ for $i = 1, 2$ implies  $\mu\in\Gamma_0$. Consequently, $\weak{\Gamma}_{0}(C)$ is a closed subset of $A$.  

4.
Since $\SubProb{\comp{B}}$ is compact, we may represent each open set $G$ as a countable union of compact sets $\Folge{C}$, so that 
\begin{equation*}
%\label{ }
\weak{\Gamma}_{0}(G) = \bigcup_{n\in\Nat}\weak{\Gamma}_{0}(C_n),
\end{equation*} 
hence $\weak{\Gamma}_{0}(G)$ is a Borel set in $A$. The Kuratowski\&Ryll-Nardzewski Selection Theorem~\ref{meas-selections-exist} together with Lemma~\ref{char-trans-kernel} gives us a stochastic relation $M_0: A\Trans \comp{B}$ with $M_0(x_1, x_2)\in\Gamma_0(x_1, x_2)$ for all $\langle x_1, x_2\rangle\in A$. Define $M(x_1, x_2)$ as the restriction of $M_0(x_1, x_2)$ to the Borel sets of $B$, then $M: A \Trans B$ is the desired relation, because $M_0(x_1, x_2)(\comp{B}\setminus B) = 0$. 
\EndProof

% Given two measure spaces $(X, {\cal A}, \mu)$ and $(Y, {\cal B}, \nu)$, a measurable map $f: X\to Y$ is called a \emph{\index{measure space!morphism}morphism} of measure spaces  $f: (X, {\cal A}, \mu) \to (Y, {\cal B}, \nu)$ iff $\nu = \FinM{f}(\mu)$ holds; thus a morphism of measure spaces carries the measure $\mu$ on one space to the measure $\nu$ on the other one in the sense that $\nu(B) = \mu(\InvBild{f}{B})$ for all $B\in{\cal B}$. Note that this relationship relates the measures' values on the $\sigma$-algebras ${\cal B}$ and $\InvBild{f}{{\cal B}}\subseteq {\cal A}$; we do not say anything about the comportment of $\mu$ on ${\cal A}\setminus\InvBild{f}{{\cal B}}$. 

For the construction we are about to undertake we will put to work the selection machinery just developed; this requires us to show that the set from which we want to select is non-empty. The following technical argument will be of assistance. 

Assume that we have Polish spaces $X_{1}, X_{2}$ and a separable measure space $(Z, {\cal C})$ with surjective and measurable maps $f_{i}: X_{i}\to Z$ for $i = 1, 2$. We also have subprobability measures $\mu_{i}\in\SubProb{X_{i}}$. Since $(Z, {\cal C})$ is separable, we may assume that ${\cal C}$ constitutes the Borel sets for some metric space $(Z, d)$ so that $d$ has a countable dense subset, see~Proposition~\ref{sep-is-sep-metr}. Proposition~\ref{Srivastava} then tells us that we may assume that $f_{1}$ and $f_{2}$ are continuous. Now define
\begin{align*}
  S & := \{\langle x_{1}, x_{2}\rangle \in X_{1}\times X_{2}\mid f_{1}(x_{1}) = f_{2}(x_{2})\}\\
{\cal A} & := S\cap \InvBild{(f_{1}\times f_{2})}{{\cal C}\otimes{\cal C}}.
\end{align*}
Since $\Delta_{Z} := \{\langle z, z\rangle \mid  z\in Z\}$ is a closed subset of $Z\times Z$, and since $f_{1}$ and $f_{2}$ are continuous, $S = \InvBild{(f_{1}\times f_{2})}{\Delta_{Z}}$ is a closed subset of the Polish space $X_{1}\times X_{2}$, hence a Polish space itself by Lemma~\ref{ClosedIsPolish}. Now assume that we have a finite measure $\theta$ on ${\cal A}$ such that 
$
\SubProb{\pi_{i}}(\theta)(E_{i}) = \mu_{i}(E_{i})
$
for all $E_{i}\in\InvBild{f_{i}}{{\cal C}}$, $i = 1, 2$ with $\pi_{1}: X_{1}\to  Z$ and $\pi_{2}: X_{2}\to Z$ as the projections. Now ${\cal A}\subseteq\Borel{S}$ is usually not the $\sigma$-algebra of Borel sets for some Polish topology on $S$, which, however, will be needed. Here Lubin's construction steps in.

\BeginLemma{application-lubin}
In the notation above, there exists a measure $\theta^{+}$ on the Borel sets of $S$ extending $\theta$ such that 
$
\SubProb{\pi_{i}}(\theta^{+})(E_{i}) = \mu_{i}(E_{i})
$
holds for all $E_{i}\in{\cal B}(S)$. 
\EndLemma

\BeginProof
Because ${\cal C}$ is countably generated, ${\cal C}\otimes{\cal C}$
is, so ${\cal A}$ is a countably generated $\sigma$-algebra. By
Lubin's Theorem~\ref{lubin-extension} there exists an extension
$\theta^{+}$ to $\theta$.%\footnote{Der Beweis ist noch unvollständig.}
\EndProof

So much for the technical preparations; we will now turn to bisimulations. A bisimulation relates two transition systems which are connected through a mediating system. In order to define this, we need morphisms. In the case of stochastic systems, recall that a \emph{morphism} $m = (f, g): K_{1}\to K_{2}$ for \index{morphism!stochastic relations} stochastic relations $K_{i}: (X_{i}, {\cal A}_{i})\Trans (Y_{i}, {\cal B}_{i})$ ($i = 1, 2$) over general measurable spaces is given through the measurable maps $f: X_{1}\to X_{2}$ and $g: Y_{1}\to Y_{2}$ such that this diagram of measurable maps commutes
\begin{equation*}
\xymatrix{
(X_{1}, {\cal A}_{1})\ar[d]_{K_{1}}\ar[rr]^{f}&& (X_{2}, {\cal A}_{2})\ar[d]^{K_{2}}\\
\SubProb{Y_{1}, {\cal B}_{1}}\ar[rr]_{\SubProb{g}}&&\SubProb{Y_{2}, {\cal B}_{2}}
}
\end{equation*}
Equivalently, $K_{2}(f(x_{1})) = \SubProb{g}(K_{1}(x_{1}))$, which translates to $K_{2}(f(x_{1}))(B) = K_{1}(x_{1})(\InvBild{g}{B})$ for all $B\in{\cal B}_{2}$. 

\BeginDefinition{6-Bisimilar} 
The stochastic relations $K_{i}: (X_{i}, {\cal A}_{i})\Trans (Y_{i},
{\cal B}_{i})$ ($i = 1, 2$), are called
\emph{\index{bisimilar}bisimilar} iff there exist a stochastic
relation $M: (A, {\cal X})\Trans (B, {\cal Y})$ and surjective
morphisms $m_{i} = (f_{i}, g_{i}): M\to K_{i}$ such that the $\sigma$-algebra
$ \InvBild{g_{1}}{{\cal B}_{1}} \cap \InvBild{g_{2}}{{\cal B}_{2}} $ is nontrivial,
i.e., contains not only $\emptyset$ and $B$.  The relation $M$ is
called \emph{mediating}.
\EndDefinition

The first condition on bisimilarity is in accordance with the general
definition of bisimilarity of coalgebras; it requests that $ m_{1} $
and $ m_{2} $ form a span of morphisms
\begin{equation*}
\xymatrix{
K_{1} & M\ar[l]_{m_{1}}\ar[r]^{m_{2}} & K_{2}.
}
\end{equation*}
Hence, the following diagram of measurable maps is supposed to commute with $m_{i} = (f_{i}, g_{i})$ for $i = 1, 2$
\begin{equation*}
\xymatrix{
(X_{1}, {\cal A}_{1})\ar[d]_{K_{1}} 
&& (A, {\cal X})\ar[ll]_{f_{1}}\ar[rr]^{f_{2}}\ar[d]^{M} 
&& (X_{2}, {\cal A}_{2})\ar[d]^{K_{2}}\\
\SubProb{Y_{1}, {\cal B}_{1}} 
&& \SubProb{B, {\cal Y}}\ar[ll]^{\SubProb{g_{1}}}\ar[rr]_{\SubProb{g_{2}}} 
&& \SubProb{Y_{2}, {\cal B}_{2}}
}
\end{equation*}
Thus, for each $a \in A, D \in{\cal B}_{1}, E \in {\cal B}_{2}$
the equalities 
\begin{alignat*}{2}
K_{1}\bigl(f(a)\bigr)(D) & = \bigl(\SubProb{g_{1}}\circ  M\bigr)(a)(D) && = M(a)\bigl(\InvBild{g_{1}}{D}\bigr)\\
K_{2}\bigl(f_{2}(a)\bigr)(E) & = \bigl(\SubProb{g_{2}}\circ  M\bigr)(a)(E) && = M(a)\bigl(\InvBild{g_{2}}{E}\bigr)
\end{alignat*}
should be satisfied.
The second condition, however, is special; it states that we can find an event $ C^* \in{\cal Y} $ which is common to both $K_{1}$ and $K_{2}$ in the sense that
\begin{equation*}
%\label{ }
\InvBild{g_{1}}{B_{1}} = C^* = \InvBild{g_{2}}{B_{2}}
\end{equation*}
for some
$
B_{1} \in {\cal B}_{1}
$
and
$
B_{2} \in {\cal B}_{2}
$
such that both $C^* \not= \emptyset$ and $C^* \not= B$ hold (note that
for $C^* =\emptyset$ or $C^* = B$ we can always take the empty and the full set,
respectively). Given such a $C^*$ with $B_{1}, B_{2}$ from above we get for each $a \in
A$
\begin{equation*}
  K_{1}(f_{1}(a))(B_{1})
= 
  M(a)(\InvBild{g_{1}}{B_{1}}) 
=  
  M(a)(C^*)
=  
  M(a)(\InvBild{g_{2}}{B_{2}})
=  
  K_{2}(g_{2}(a))(B_{2});
\end{equation*}
thus the event $C^*$ ties $K_{1}$ and $K_{2}$ together. Loosely
speaking, $ \InvBild{g_{1}}{{\cal B}_{1}} \cap \InvBild{g_{2}}{{\cal
    B}_{2}} $ can be described as the $\sigma$-algebra of common
events, which is required to be nontrivial.  

Note that without the second condition two relations $K_{1}$ and
$K_{2}$ which are strictly probabilistic (i.e., for which the entire
space is always assigned probability $1$) would always be bisimilar:
Put $ A := X_{1}\times X_{2}$, $B :=Y_{1}\times Y_{2} $ and set for
$\langle x_{1}, x_{2}\rangle \in A$ as the mediating relation $
M(x_{1}, x_{2}) := K_{1}(x_{1}) \otimes K_{2}(x_{2});$ that is, define
$M$ pointwise to be the product measure of $K_{1}$ and $K_{2}$. Then
the projections will make the diagram commutative. But this is
way too weak, because bisimulations relate transition systems, and
it does not offer particularly interesting insights when two arbitrary systems can
be related. It is also clear that using products for mediation does not work for the
subprobabilistic case.

We will show now that we can construct a bismulation for stochastic
relations which are linked through a co-span $ \xymatrix{ K_{1}&&
  K\ar[ll]\ar[rr] && K_{2}.  } $ The center $K$ of this co-span should
be defined over second countable metric spaces, $K_{1}$ and $K_{2}$
over Polish spaces. This situation is sometimes easy to obtain, e.g.,
when factoring Kripke models over Polish spaces through a suitable
logic; then $K$ is defined over analytic spaces, which are separable
metric. This is described in greater detail in
Example~\ref{modal-logics}.

\BeginProposition{6-HasSemiPullbacks} 
Let $K_{i}: X_{i}\Trans Y_{i}$ be stochastic relations over Polish spaces,
and assume that $K: X\Trans Y$ is a stochastic relation,
where $X, Y$ are second countable metric spaces. Assume that we have a cospan of morphisms $m_{i}: K_{i}\to K, i = 1, 2$, then there exists a stochastic relation $M$ and morphisms $m^{+}_{i}: M\Trans K_{i}, i = 1, 2$ rendering this diagram commutative. \begin{equation*}
%\label{ }
\xymatrix{
M\ar[rr]^{m^{+}_1}\ar[d]_{m^{+}_2}&&K_2\ar[d]^{m_2}\\
K_1\ar[rr]_{m_1}&&K
}
\end{equation*}
The stochastic relation $M$ is defined over Polish spaces. 
\EndProposition

\BeginProof
1. Assume $K_i = (X_i, Y_i, K_i)$ with
$m_i = (f_i, g_i)$, $i = 1, 2$.
Because of
Proposition~\ref{Srivastava} we may assume that the respective
$\sigma$-algebras on $X_1$ and $X_2$ are obtained from Polish
topologies which render $f_1$ and $K_1$ as well as
$f_2$ and $K_2$ continuous. These topologies are fixed for
the proof. Put
\begin{align*}
  A & := \{\langle x_{1}, x_{2}\rangle \in X_1 \times X_2 \mid
  f_1(x_1) = f_2(x_2)\},\\
  B & :=  \{\langle y_{1}, y_{2}\rangle\in Y_1 \times Y_2 \mid  g_1(y_1) =
  g_2(y_2)\},
\end{align*}
then both $A$ and $B$ are closed, hence Polish. $\alpha_i: A \rightarrow X_i$ and
$\beta_i: B \rightarrow Y_i$ are the projections, $i = 1, 2$.
The diagrams
\begin{equation*}
%\label{ }
\xymatrix{
X_1\ar[rr]^{f_1}\ar[d]_{K_1} && X\ar[d]_{K}  && X_2\ar[ll]_{f_2}\ar[d]^{K_2} \\
\SubProb{Y_1}\ar[rr]_{\SubProb{g_1}} && \SubProb{Y} && \SubProb{Y_2}\ar[ll]^{\SubProb{g_2}}
}
\end{equation*}
are commutative by assumption, thus we know that for $x_i \in X_i$
\begin{equation*}
  K(f_1(x_1)) = \SubProb{g_1}(K_1(x_1)) \text{ and }
  K(f_2(x_2)) = \SubProb{g_2}(K_2(x_2))
\end{equation*}
holds. The construction implies that
$
(g_1\circ  \beta_1)(y_1, y_2) = (g_2\circ  \beta_2)(y_1, y_2)
$
is true for $\langle y_1, y_2\rangle\in B$, and
$
g_1\circ  \beta_1: B \rightarrow Y
$
is surjective.

2.  Fix $\langle x_1, x_2\rangle \in A$.  Separability of the target
spaces now enters: We know that the image of a surjective map under
$\SubProbSenza$ is onto again by Proposition~\ref{SubProbIsOnto}, so
that there exists $ \mu_0 \in \SubProb{B} $ with 
$ \SubProb{g_1\circ
  \beta_1}(\mu_0) = K(f_1(x_1)), $ consequently, $ \SubProb{g_i \circ
  \beta_i}(\mu_0) = \SubProb{g_i}(K_i(x_i))\ (i = 1, 2).  
$ 
But this
means for $i = 1, 2$
\begin{equation*}
\forall E_i \in \InvBild{g_i}{\Borel{Y}}:
\SubProb{\beta_i}(\mu_0)(E_i) = K_i(x_i)(E_i).
\end{equation*}
Put
\begin{equation*}
%\label{ }
\Gamma(x_1, x_2):= \{\mu \in \SubProb{B} \mid
\SubProb{\beta_1}(\mu) = K_1(x_1)\text{ and } \SubProb{\beta_2}(\mu) =
K_2(x_2)\}, 
\end{equation*}
then Lemma~\ref{application-lubin} shows that
$
\Gamma(x_1, x_2) \not= \emptyset.
$

3. 
The set
\begin{equation*}
\weak{\Gamma}(C) = \{\langle x_1, x_2\rangle \in A \mid
\Gamma(x_1, x_2) \cap C \not= \emptyset\} 
\end{equation*}
is closed in $A$ for compact
$
C \subseteq \SubProb{B}.
$
This is shown exactly as in the second part of the proof for Proposition~\ref{6-selection}, from which now is inferred that there exists a measurable map
$
M: A \rightarrow \SubProb{B}
$
such that
$
M(x_1, x_2) \in \Gamma(x_1, x_2)
$
holds for every
$
\langle x_1, x_2\rangle \in A.
$
Thus $M: A \Trans B$ is a stochastic relation
with
\begin{equation*}
  K_1\circ  \alpha_1 = \SubProb{\beta_1}\circ  M \text{ and }
  K_2 \circ \alpha_2 = \SubProb{\beta_2}\circ  M.
\end{equation*}
Thus $M$ with $m^{+}_{1} := (\alpha_{1}, \beta_{1})$ and $m^{+}_{2} :=
(\alpha_{2}, \beta_{2})$ is the desired semi-pullback.
\EndProof

Now we know that we may construct from a co-span of stochastic
relations a span. Let us have a look at a typical situation in which
such a co-span may occur.

\BeginExample{modal-logics}
Consider the modal logic from Example~\ref{modal-logic} again, and
interpret the logic through stochastic relations $K: S\Trans S$ and
$L: T\Trans T$ over the Polish spaces $S$ and $T$. The equivalence
relations $\sim_{K}$ and $\sim_{L}$ are defined as in
Example~\ref{modal-yields-smooth}. Because we have only countably many
formulas, these relations are smooth. For readability, denote the
equivalence class associated with $\sim_{K}$ by $[\cdot]_{K}$, similar
for $[\cdot]_{L}$. Because $\sim_{K}$ and $\sim_{L}$ are smooth, the
factor spaces $\Faktor{S}{K}$ resp. $\Faktor{T}{L}$ are analytic
spaces when equipped with the final $\sigma$-algebra with respect to
$\fMap{K}$ resp. $\fMap{L}$ by Proposition~\ref{RemainsAnalytic}. The
factor relation $K_{F}: \Faktor{S}{K}\Trans \Faktor{S}{K}$ is then the
unique relation which makes this diagram commutative
\begin{equation*}
\xymatrix{
S\ar[d]_{K}\ar[rr]^{\fMap{K}} && \Faktor{S}{K}\ar[d]^{K_{F}}\\
\SubProb{S}\ar[rr]_{\SubProb{\fMap{K}}}&&\SubProb{\Faktor{S}{K}}
}
\end{equation*}
This translates to $K(s)(\InvBild{\fMap{K}}{B}) =
K_{F}(\Klasse{s}{K})(B)$ for all $B\in\Borel{\Faktor{S}{K}}$ and all
$s\in X$. 

Associate with each formula
$\phi$ its validity sets $\Gilt_{K}$ resp. $\Gilt_{L}$, and call $s\in
S$ \emph{logically equivalent} to $t\in T$ iff we have for each formula
$\phi$
\begin{equation*}
  s\in\Gilt_{K}\Leftrightarrow t\in\Gilt_{L}
\end{equation*}
Hence $s$ and $t$ are logically equivalent iff no formula is able to
distinguish between states $s$ and $t$; call the stochastic relations
$K$ and $L$ \emph{logically equivalent} iff given $s\in S$ there
exists $t\in T$ such that $s$ and $t$ are logically equivalent, and
vice versa.

Now assume that $K$ and $L$ are logically equivalent, and consider 
\begin{equation*}
  \Phi := \bigl\{\langle\Klasse{s}{K}, \Klasse{t}{L}\rangle \mid s\in S\text{
    and }t\in T\text{ are logically equivalent}\bigr\}.
\end{equation*}
Then $\Phi$ is the graph of a bijective map; this is easy to
see. Denote the map by $\Phi$ as well. 
Since
$\InvBild{\Phi}{\Bild{\fMap{L}}{\Gilt_{L}}} =
\Bild{\fMap{K}}{\Gilt_{K}}$, and since the set
$
\{\Bild{\fMap{L}}{\Gilt_{L}}\mid \phi\text{ is a formula}\}
$
generates $\Borel{\Faktor{T}{L}}$  by Proposition~\ref{seq-gen-factor-space}, $\Phi: \Faktor{S}{K}\to
\Faktor{T}{L}$ is Borel measurable; interchanging the r\^oles of $K$
and $L$ yields measurability of $\Phi^{-1}$. 

Hence we have this picture for logical equivalent $K$ and $L$:
\begin{equation*}
\xymatrix{
&&L\ar[d]^{\fMap{L}}\\
K\ar[rr]_{\Phi\circ \fMap{K}}&&L_{F}
}
\end{equation*}
\EndExample

This example can be generalized to the case that the relations operate
on two spaces rather than only on one. Let $K: X\Trans Y$ be a
transition kernel over the Polish spaces $X$ and $Y$. Then the pair
$(\kappa, \lambda)$ of smooth equivalence relations $\kappa$ on $X$
and $\lambda$ on $Y$ is called a \emph{\index{congruence}congruence}
for $K$ iff there exists a transition kernel $K_{\kappa, \lambda}:
\Faktor{X}{\kappa}\Trans\Faktor{Y}{\lambda}$ rendering the diagram
commutative:
\begin{equation*}
\xymatrix{
X\ar[d]_{K}\ar[rr]^{\fMap{\kappa}} 
&& \Faktor{X}{\kappa}\ar[d]^{K_{\kappa, \lambda}}\\
Y\ar[rr]_{\SubProb{\fMap{\lambda}}} 
&& \Faktor{Y}{\lambda}
}
\end{equation*}
Because $\fMap{\kappa}$ is an epimorphism, $K_{\kappa, \lambda}$ is
uniquely determined, if it exists (for a discussion of congruences for
stochastic coalgebras, see~\CategCite{Section~1.6.2}). Commutativity of the diagram translates to
\begin{equation*}
  K(x)(\InvBild{\fMap{\lambda}}{B}) = K_{\kappa, \lambda}(\Klasse{x}{\kappa})(B)
\end{equation*}
for all $x\in X$ and all $B\in\Borel{\Faktor{Y}{\lambda}}$. Call in
analogy to Example~\ref{modal-logics} the transition kernels $K_{1}:
X_{1}\Trans Y_{1}$ and $K_{2}: X_{2}\Trans Y_{2}$ \emph{logically equivalent}
iff there exist congruences $(\kappa_{1}, \lambda_{1})$ for $K_{1}$
and $(\kappa_{2}, \lambda_{2})$ for $K_{2}$ such that the factor
relations $K_{\kappa_{1}, \lambda_{1}}$ and $K_{\kappa_{2},
  \lambda_{2}}$ are isomorphic. 

In the spirit of this discussion, we obtain from
Proposition~\ref{6-HasSemiPullbacks} 

\BeginTheorem{log-equiv-are-bisim}
Logically equivalent stochastic relations over Polish spaces are bisimilar.
\EndTheorem

\BeginProof
1.
The proof applies Proposition~\ref{6-HasSemiPullbacks}; first it
has to show how to satisfy the assumptions of that statement. Let
$K_{i}: X_{i}\Trans Y_{i}$ be stochastic relations over Polish spaces for
$i = 1, 2$. We assume that $K_{1}$ is logically equivalent to $K_{2}$,
hence there exist congruences $(\kappa_{i}, \lambda_{i})$ for $K_{i}$
such that the associated stochastic relations $K_{\kappa_{i},
  \lambda_{i}}: \Faktor{X_{i}}{\kappa_{i}}\Trans
\Faktor{Y_{i}}{\lambda_{i}}$ are isomorphic. Denote this isomorphism by
$(\phi, \psi)$, so $\phi: \Faktor{X_{1}}{\kappa_{1}}\to
\Faktor{X_{2}}{\kappa_{2}}$ and $\psi: \Faktor{Y_{1}}{\lambda_{1}}\to
\Faktor{Y_{2}}{\lambda_{2}}$ are in particular measurable bijections, so are their
inverses. 

2.
Let $\eta_{2} := (\fMap{\kappa_{2}}, \fMap{\lambda_{2}})$ be the
factor morphisms $\eta_{2}: K_{2}\to K_{\kappa_{2}, \lambda_{2}}$, and
put $\eta_{1} := (\phi\circ \fMap{\kappa_{1}}, \psi\circ
\fMap{\lambda_{1}})$, thus we obtain this co-span of morphisms
\begin{equation*}
\xymatrix{
K_{1}\ar[rr]^{\eta_{1}}&&K_{\kappa_{2}, \lambda_{2}}&&K_{2}\ar[ll]_{\eta_{2}}
}
\end{equation*}
Because both $\Faktor{X_{2}}{\kappa_{2}}$ and $\Faktor{Y_{2}}
  {\lambda_{2}}$ are analytic spaces on account of $\kappa_{2}$ and
$\lambda_{2}$ being smooth, see Proposition~\ref{RemainsAnalytic}, we
apply Proposition~\ref{6-HasSemiPullbacks} and obtain a mediating
relation $M: A\Trans B$ with Polish $A$ and $B$ such that the
projections $\alpha_{i}: A\to X_{i}$ and $\beta_{i}: B\to Y_{i}$ are
morphisms for $i = 1, 2$. Here
\begin{align*}
  A & := \{\langle x_{1}, x_{2}\rangle \mid
  \phi(\Klasse{x_{1}}{\kappa_{1}}) = \Klasse{x_{2}}{\kappa_{2}}\}\\
 B & := \{\langle y_{1}, y_{2}\rangle \mid
  \phi(\Klasse{y_{1}}{\lambda_{1}}) = \Klasse{y_{2}}{\lambda_{2}}\}
\end{align*}
It remains to be demonstrated that the $\sigma$-algebra of common
events, viz., the intersection
$\InvBild{\beta_{1}}{\Borel{Y_{1}}}\cap\InvBild{\beta_{2}}{\Borel{Y_{2}}}$ 
is not trivial. 

3.
Let $U_{2}\in\Borel{Y_{2}}$ be $\lambda_{2}$-invariant. Then
$\Bild{\fMap{\lambda_{2}}}{U_{2}}
\in\Borel{\Faktor{Y_{2}}{\lambda_{2}}}$,
because
$U_{2}=\InvBild{\fMap{\lambda_{2}}}{\Bild{\fMap{\lambda_{2}}}{U_{2}}}$
on account of $U_{2}$ being $\lambda_{2}$-invariant. Thus 
$
U_{1} := \InvBild{\fMap{\lambda_{1}}}{%
\InvBild{\psi}{%
\Bild{\fMap{\lambda_{2}}}{U_{2}}
}
}
$
is an $\lambda_{1}$-invariant Borel set in $Y_{1}$ with
\begin{align*}
  \langle y_{1}, y_{2}\rangle \in (Y_{1}\times U_{2})\cap B
& \Leftrightarrow
y_{2}\in U_{2}\text{ and } \psi(\Klasse{y_{1}}{\lambda_{1}}) =
\Klasse{y_{2}}{\lambda_{2}}\\
& \Leftrightarrow
\langle y_{1}, y_{2}\rangle \in (U_{1}\times U_{2})\cap B.
\end{align*}
One shows in exactly the same way 
\begin{equation*}
\langle y_{1}, y_{2}\rangle \in (U_{1}\times Y_{2})\cap B
\Leftrightarrow
\langle y_{1}, y_{2}\rangle \in (U_{1}\times U_{2})\cap B.
\end{equation*}
Consequently, $(U_{1}\times U_{2})\cap B$ belongs to both 
$
\InvBild{\beta_{1}}{\Borel{Y_{1}}}
$
and
$
\InvBild{\beta_{2}}{\Borel{Y_{2}}},
$
so that this intersection is not trivial. 
\EndProof

Call a class $\mathfrak{A}$ of spaces closed under bisimulations if
the mediating relation for stochastic relations over spaces from
$\mathfrak{A}$ is again defined over spaces from $\mathfrak{A}$. Then
the result above shows that Polish spaces are closed under
bisimulations. This generalizes a result by Desharnais, Edalat and
Panangaden~\cite{Edalat,Desharnais-Edalat-Panangaden} which
demonstrates ---~through a completely different approach~--- that
analytic spaces are closed under bisimulations; Sánchez
Terraf~\cite{Sanchez} has shown that general measurable spaces are not
closed under bisimulations. In view of von Neumann's Selection
Theorem~\ref{von-neumann-selection} it might be interesting to see
whether complete measurable spaces are closed.

We have finally a look at a situation in which no semi-pullback
exists. A first example in this direction was presented
in~\cite[Theorem 12]{Sanchez}. It is based on the extension of
Lebesgue measure to a $\sigma$-algebra which does contain the Borel
sets of $[0, 1]$ augmented by a non-measurable set, and it shows that
one can construct Markov transition systems which do not have a
semi-pullback. The example below extends this by showing that one does
not have to consider transition systems, but that a look at the
measures on which they are based suffices.

\BeginExample{no-semi-pullback-exists}
A morphism $f: (X, {\cal A}, \mu)\to (Y, {\cal B}, \nu)$ of \index{morphism!measure spaces} measure spaces is an ${\cal A}$-${\cal B}$-measurable map $f: X\to Y$ such that $\nu = \FinM{f}(\mu)$. Since each finite measure can be viewed as a transition kernel, this is a special case of morphisms for transition kernels. If ${\cal B}$ is a sub-$\sigma$-algebra of ${\cal A}$ with $\mu$ an extension to $\nu$, then the identity is a morphisms $(X, {\cal A}, \mu)\to (X, {\cal B}, \nu)$. 

Denote Lebesgue measure on $([0, 1], \Borel{[0, 1]})$ by $\lambda$. Assuming the Axiom of Choice, we know that there exists $W\subseteq [0, 1]$ with $\lambda_{*}(W) = 0$ and $\lambda^{*}(W) = 1$ by~\SetCite{Lemma 1.7.7}. Denote by ${\cal A}_{W} := \sigma(\Borel{[0, 1]}\cup\{W\}$ the smallest $\sigma$-algebra containing the Borel sets of $[0, 1]$ and $W$. Then we know from Exercise~\ref{ex-extend-measure} that we can find for each $\alpha\in[0, 1]$ a measure $\mu_{\alpha}$ on ${\cal A}_{W}$ which extends $\lambda$ such that $\mu_{\alpha}(W) = \alpha$. 

Hence by the remark just made, the identity yields a morphism $f_{\alpha}: ([0, 1], {\cal A}_{W}, \mu_{\alpha}) \to ([0, 1], \Borel{[0, 1]}, \lambda)$. Now let $\alpha\not=\beta$, then 
\begin{equation*}
\xymatrix{
([0, 1], {\cal A}_{W}, \mu_{\alpha})\ar[rr]^{f_{\alpha}} &&
 ([0, 1], \Borel{[0, 1]}, \lambda) &&
 ([0, 1], {\cal A}_{W}, \mu_{\beta})\ar[ll]_{f_{\beta}}
}
\end{equation*}
is a co-span of morphisms. 

We claim that this co-span does not have a semi-pullpack. In fact, assume that $(P, {\cal P}, \rho)$ with morphisms $\pi_{\alpha}$ and $\pi_{\beta}$ is a semi-pullback, then $f_{\alpha}\circ \pi_{\alpha} = f_{\beta}\circ \pi_{\beta}$, so that $\pi_{\alpha} = \pi_{\beta}$, and $\InvBild{\pi_{\alpha}}{W} = \InvBild{\pi_{\beta}}{W}\in{\cal P}$. But then
\begin{equation*}
  \alpha = \mu_{\alpha}(W) = \rho(\InvBild{\pi_{\alpha}}{W}) = \rho(\InvBild{\pi_{\beta}}{W}) = \mu_{\beta}(W) = \beta.
\end{equation*}
This contradicts the assumption that $\alpha\not=\beta$. 
\EndExample

This example shows that the topological assumptions imposed above are
indeed necessary. It assumes the Axiom of Choice, so one might ask
what happens if this axiom is replaced by the Axiom of Determinacy. We
know that the latter one implies that each subset of the unit interval
is $\lambda$-measurable by~\SetCite{Theorem 1.7.14}, so
$\lambda_{*}(W) = \lambda^{*}(W)$ holds for each $W\subseteq[0,
1]$. Then at least the construction above does not work (on the other
hand, we made use of Tihonov's Theorem, which is known to be
equivalent to the Axiom of Choice~\cite[Theorem
4.68]{Herrlich-Choice}, so there is probably no escape from the Axiom
of Choice).

%\Input{\Folder/CaseQuotient}
%spell checked - 24Aug14
\Subsubsection{Case Study: Quotients for Stochastic Relations}
\label{sec:case-quot-stoch-rel}
{
%beware of local macros!
\def\tau{H}
\def\theta{J}
As Monty Python used to say, ``And now for something completely
different!'' We will deal now with quotients for stochastic relations,
perceived as morphisms in the Kleisli category over the monad which is
given by the subprobability functor (which is sometimes called the
Giry monad). We will first have a look at surjective maps as
epimorphisms in the category of sets, explaining the problem there,
show that a straightforward approach gleaned from the category of sets
does not appear promising, and show then that measurable selections
are the appropriate tool for tackling the problem.

For motivation, we start with surjective maps on a fixed set
$M$, serving as a domain. Let $f: M\to X$ and $g: M\to Y$ be onto, and define the partial
order $f\leq g$ iff $f = \zeta\circ g$ for some $\zeta: Y\to
X$. Clearly, $\leq$ is reflexive and transitive; the equivalence
relation $\sim$ defines through $f\sim g$ iff $f \leq g$ and $g\leq f$
is of interest here. Thus $f = \zeta\circ g$ and $g = \xi\circ f$ for
suitable $\zeta: Y\to  X$ and $\xi: X\to Y$. Because surjective
maps are epimorphisms in the category of sets with maps as morphisms,
we obtain $\zeta\circ \xi = id_{X}$ and $\xi\circ \zeta=id_{Y}$. Hence
$\zeta$ and $\xi$ are bijections. The surjections $f$ and $g$, both
with domain $M$, are 
equivalent iff there exists a bijection $\beta$ with $f = \beta\circ
g$. This is called a \emph{\index{quotient object}quotient object} for $M$ We know that the surjection $f: M\to Y$ can be factored as $f =
\widetilde{f}\circ \fMap{\Kern{f}}$ with $\widetilde{f}:
\Klasse{x}{\Kern{f}}\mapsto f(x)$ as the bijection. Thus for maps, the
quotient objects for $M$ may be identified through the quotient maps
$\fMap{\Kern{f}}$, in a similar way, the quotient objects in the
category of groups can be identified through normal subgroups;
see~\cite[V.7]{MacLane} for a discussion. Thus quotients seem to
be interesting. 

We turn to stochastic relations. The subprobability functor on the
category of measurable spaces is the functorial part of the Giry monad, and
the stochastic relations are just the Kleisli morphism for this monad,
see~\CategCite{Example 1.99}. Let $K: (X, {\cal A})\Trans (Y, {\cal
  B})$ be a stochastic relation, then
Exercise~\ref{ex-kernel-yields} shows that
\begin{equation*}
  \comp{K}(\mu): B \mapsto \dInt{K(x)(B)}{\mu(x)}
\end{equation*}
defines a $\schwach{X, {\cal A}}$-$\schwach{Y, {\cal B}}$-measurable
map $\SubProb{X, {\cal A}}\to \SubProb{Y, {\cal B}}$; $\comp{K}$ is
sometimes called the \emph{Kleisli map} associated with the Kleisli morphism $K$(it should not be confused with the completion of $K$ as discussed in Section~\ref{sec:ext-stoch-rel-to-compl}). It is clear that
$K\mapsto \comp{K}$ is injective, because $\comp{K}(\delta_{x}) =
K(x)$. 

It will helpful to evaluate the integral with respect to
$\comp{K}(\mu)$: let $g: Y\to \Real$ be bounded and measurable, then 
\begin{equation}
\label{kleisli-integral-form}
  \dInt[Y]{g}{K(\mu)} = \dInt{\dInt[Y]{g(y)}{K(x)(y)}}{\mu(x)}.
\end{equation}
In order to show this, assume first that $g=\chi_{B}$ for
$B\in{\cal B}$, then both sides evaluate to $K(\mu)(B)$, so the
representation is valid for indicator functions. Linearity of the
integral yields the representation for step functions. Since we
may find for general $g$ a sequence $\Folge{g}$ of step functions with
$\lim_{n\to \infty}g_{n}(y) = g(y)$ for all $y\in Y$, and since $g$ is bounded, hence
integrable with respect to all finite measures, we obtain from
Lebesgue's Dominated Convergence
Theorem~\ref{lebesgue-dominated-convergence} that
\begin{align*}
  \dInt[Y]{g}{\comp{K}(\mu)} 
& = \lim_{n\to \infty}\dInt[Y]{g_{n}}{\comp{K}K(\mu)}\\
& = \lim_{n\to \infty}\dInt{\dInt[Y]{g_{n}(y)}{K(x)(y)}}{\mu(x)}\\
& = \dInt{\lim_{n\to \infty}\dInt[Y]{g_{n}(y)}{K(x)(y)}}{\mu(x)}\\
& = \dInt{\dInt[Y]{g(y)}{K(x)(y)}}{\mu(x)}
\end{align*}
This gives the desired representation.

The Kleisli map is related to the convolution operation defined in
Example~\ref{convolution-kernels}:

\BeginLemma{bar-vs-convolution}
Let $K: (X, {\cal A})\Trans (Y, {\cal
  B})$ and $L: (Y, {\cal B})\Trans (Z, {\cal
  C})$, then $\comp{L*K} = \comp{L}\circ \comp{K}$. 
\EndLemma

\BeginProof
Evaluate both the left and the right hand side for $\mu\in\SubProb{X,
  {\cal A}}$ and $C\in {\cal C}$:
\begin{align*}
  \comp{L*K}(\mu)(C) 
& = \dInt{\dInt[Y]{L(y)(C)}{K(x)(y)}}{\mu(x)}\\
& = \dInt[Y]{L(y)(C)}{\comp{K}(\mu)(y)}&&\text{ by }(\ref{kleisli-integral-form})\\
& = \comp{L}(\comp{K})(\mu)(C)
\end{align*}
This implies the desired equality. 
\EndProof

Associate with each measurable $f: Y \to Z$ a stochastic relation $\delta_f: Y \Trans Z$ through 
$
\delta_f(y)(C) := \delta_y(\InvBild{f}{C}),
$
then 
$
\delta_f = \SubProb{f}\circ\delta,
$
and a direct computation shows
$
\delta_f *  K  = \SubProb{f}\circ K.
$
In fact,
\begin{align*}
\bigl(\delta_f *  K\bigr)(x)(C)
& =
\int_Y \delta_f(y)(C)\ K(x)(dy)\\
& =
\int_Y \chi_{\InvBild{f}{C}}(y)\ K(x)(dy)\\
& =
K(x)(\InvBild{f}{C})\\
& =
\bigl(\SubProb{f}\circ K\bigr)(x)(C).
\end{align*}

On the other hand, if $f: W \to X$ is measurable, then
\begin{equation*}
%\label{ }
\bigl(K * \delta_f\bigr)(w)(B) = \int_X K(x)(B)\ \delta_f(w)(dx) = (K\circ f)(w)(B).
\end{equation*}

In particular, it follows that $e_{X} := \SubProb{id_{X}}$ is the neutral element:
$K = e_{X}*K =K*e_{X} = K$. 
Recall that $K$ is an epimorphism in the Kleisli category iff $L_{1} * K = L_{2} * K$ implies
$L_{1} = L_{2}$ for any stochastic relations $L_{1}, L_{2}: (Y, {\cal
  B})\Trans (Z, {\cal C})$. Lemma~\ref{bar-vs-convolution} tells us
that if the Kleisli map $\comp{K}$  is onto, then $K$ is an
epimorphism. Now let $K: (X, {\cal A})\Trans (Y, {\cal B})$ and $L:
(X, {\cal A})\Trans (Z, {\cal C})$ be stochastic relations, and assume
that both $K$ and $L$ are epis. Define as above
\begin{equation*}
  K\leq L \Leftrightarrow \exists \theta: (Z, {\cal C})\Trans (Y,
  {\cal B}): K = \theta*L
\end{equation*}
Hence we can find in case  $K\leq L$ a stochastic relation $\theta$
such that $K(x)(B) = \dInt[Z]{\theta(z)(B)}{L(x)(z)}$ for $x\in X$ and
$B\in{\cal B}$.  

We will deal for the rest of this section with Polish spaces. Fix $X$
as a Polish spaces. For identifying the quotients with respect to
Kleisli morphisms, one could be tempted to mimic the approach observed
for the sets as outlined above. This is studied in the next example.

\BeginExample{quotient-ex3}
{\renewcommand{\widetilde}[1]{{#1}_{\sharp}}
Let $K: X \Trans Y$ be a stochastic relation with Polish $Y$ which is
an epi. $\Faktor{X}{\Kern{K}}$ is an analytic space, since $K: X\to
\SubProb{Y}$ is a measurable
map into the  Polish space $\SubProb{Y}$
by~Proposition~\ref{polish-iff-polish}, so that $\Kern{K}$ is smooth. Define the map $E_K: X\to\SubProb{\Faktor{X}{\Kern{K}}}$ through 
$
E_K(x) := \delta_{\Klasse{x}{\Kern{K}}},
$
hence we obtain for each $x\in X$, and each Borel set $G\in\Borel{\Faktor{X}{\Kern{K}}}$ 
\begin{equation*}
%\label{ }
E_K(x)(B) = \delta_{\Klasse{x}{\Kern{K}}}(G) = \delta_x(\InvBild{\fMap{\Kern{K}}}{G}) = \SubProb{\fMap{\Kern{K}}}(\delta_x)(G).
\end{equation*}
Thus $E_K$ is an epi as well: take $\mu\in\SubProb{X}$ and $G \in
\Borel{\Faktor{X}{\Kern{K}}}$, then
\begin{align*}
\comp{E}_K(\mu)(G)
& =
\dInt{E_K(x)(G)}{\mu(x)}\\
& =
\dInt{\delta_x(\InvBild{\fMap{\Kern{K}}}{G})}{\mu(x)}\\
& =
\mu(\InvBild{\fMap{\Kern{K}}}{G}) \\
& =
\SubProb{\fMap{\Kern{K}}}(\mu)(G),
\end{align*}
so that $\comp{E}_K = \SubProb{\fMap{\Kern{K}}}$; since the image of a
surjective map under $\SubProbSenza$ is surjective again by Proposition~\ref{SubProbIsOnto}, we conclude that $E_K$ is an epi. Now define for $x\in X$ the map 
\begin{equation*}
\widetilde{K}(\Klasse{x}{\Kern{K}}) := K(x),
\end{equation*}
then standard arguments show that $\widetilde{K}$ is well defined and
constitutes a stochastic relation $\widetilde{K}:
\Faktor{X}{\Kern{K}}\Trans Y$. Moreover we obtain for $x\in X, H\in
\Borel{Y}$ by the change of variables formula in Corollary~\ref{int-image-measure}
\begin{align*}
(\widetilde{K}\klComp E_K)(x)(H) 
& =
\dInt[\Faktor{X}{\Kern{K}}]{\widetilde{K}(t)(H)}{E_K(x)(t)}\\
& = 
\dInt[\Faktor{X}{\Kern{K}}]{\widetilde{K}(t)(H)}{\SubProb{\fMap{\Kern{K}}}(\delta_x)(t)}\\
& =
\dInt{\widetilde{K}(\Klasse{w}{\Kern{K}})(H)}{\delta_x(w)}\\
& =
\dInt{K(w)(H)}{\delta_x(w)}\\
& =
K(x)(H).
\end{align*}
Consequently, $K$ can be factored as $K = \widetilde{K}*E_K$ with the
epi $E_K$. But there is no reason why in general $\widetilde{K}$
should be invertible; for this to hold, the map $\comp{\widetilde{K}}:
\SubProb{\Faktor{X}{\Kern{K}}}\to\SubProb{Y}$ is required to be
injective. Hence $\isEquiv{K}{E_K}{\approx}$ holds only in special cases. 
} 
\EndExample

This last example indicates that a characterization of quotients
for the Kleisli category at least for the Giry monad cannot be derived
directly by carrying over a characterization for the underlying
category. 

For the rest of the section we discuss the Kleisli category
for the Giry monad over Polish spaces, hence we deal with stochastic
relations. Let $X$, $Y$ and $Z$ be Polish, and fix $K: X \Trans Y$ and
$L: X \Trans Z$ so that $\isEquiv{K}{L}{\approx}$. Hence there exists
$\theta: Y \Trans Z$ with inverse $\tau: Z \Trans Y$ and
$ L = \theta *  K$ and $K = \tau * L$. Because both $K$ and $L$ are
epis, we obtain these simultaneous
equations
$
\tau *  \theta = e_Y,
\theta * \tau = e_Z.
$
They entail
$
\label{a-1}\dInt[Z]{\tau(z)(B)}{\theta(y)(z)} = \delta_y(B)
$
and
$
\dInt[Y]{\theta(y)(C)}{\tau(z)(y)} = \delta_z(C)
$
for all $y\in Y, z\in Z$ and $B\in\Borel{Y}, C\in\Borel{Z}$. Because
singletons are Borel sets, these equalities imply
$
\dInt[Z]{\tau(z)(\{y\})}{\theta(y)(z)} = 1
$
and
$
\dInt[Y]{\theta(y)(\{z\})}{\tau(z)(y)} = 1.
$
Consequently, we obtain
\begin{align*}
\forall y \in Y: \theta(y)(\{z\in Z \mid \tau(z)(\{y\}) = 1\}) = 1,\\
\forall z \in Z: \tau(z)(\{y\in Y \mid \theta(y)(\{z\}) = 1\}) = 1.
\end{align*}

\BeginProposition{Borel-graph}
There exist Borel maps $f: Y\to Z$ and $g: Z\to Y$ such that
$ 
\tau\bigl(f(y)\bigr)(\{y\})  = 1
$
and
$
\label{b-4}\theta\bigl(g(z)\bigr)(\{z\}) = 1
$
for all $y\in Y, z \in Z$ 
\EndProposition

\BeginProof
1.
Define
$
P := \{\langle y, z\rangle \in Y\times Z\mid \tau(z)(\{y\}) = 1\},
$
and
$
Q := \{\langle z, y\rangle \in Z\times Y\mid \theta(y)(\{z\}) = 1\}, 
$
then $P$ and $Q$ are Borel sets. We establish this for $P$, the argumentation for $Q$ is very similar.

2. 
With a view towards Proposition~\ref{Srivastava} we may and do
assume that $\tau: Z\to\SubProb{Y}$ is continuous. Let $\bigl(\langle
y_n, z_n\rangle\bigr)_{n\in \Nat}$ be a sequence in $P$ with $\langle
y_n, z_n\rangle\to\langle y, z\rangle$, hence the sequence
$\bigl(\tau(z_n)\bigr)_{n\in\Nat}$ converges weakly $\tau(z)$. Given
$m\in \Nat$ there exists $n_0\in\Nat$ such that $y_n\in V_{1/m}(y)$
for all $n_0\geq n$, where $V_{1/m}(y)$ is the closed ball of radius
$1/m$ around $y$. Since $\tau$ is weakly continuous, we obtain 
$
\limsup_{n\to\infty}\tau(z_n)\bigl(V_{1/m}(y)\bigr) \leq \tau(z)\bigl(V_{1/m}(y)\bigr)
$
from Proposition~\ref{convergence-a-topology}, hence
$
%\begin{equation*}
%\label{ }
\tau(z)\bigl(V_{1/m}(y)\bigr) = 1.
%\end{equation*}
$
Because 
$ %\begin{equation*}
%\label{ }
\bigcap_{m\in\Nat}V_{1/m}(y) = \{y\},
$ %\end{equation*}
we conclude $\tau(z)(\{y\}) = 1$, thus $\langle y, z\rangle\in
P$. Consequently, $P$ is a closed subset of $Y \times Z$, hence a
Borel set. 

3.
Since $P$ is closed, the cut $P_{y}$ at $y$ is closed as well, and we
have 
$
\theta(y)(P_{y}) = \theta(y)(\{z\in Z \mid \tau(z)(\{y\}) = 1\} = 1,
$
thus we obtain $\supp(\theta(y))\subseteq P_{y}$, because the
support $\supp(\theta(y))$ is the smallest closed set $C$ with
$\theta(y)(C) = 1$. Since $y\mapsto \supp(\theta(y))$ is measurable
(cp. Example~\ref{transitins-have-selections}), we obtain from
Theorem~\ref{meas-selections-exist} a measurable map $f: Y\to Z$ with
$f(y)\in\supp(\theta(y))\subseteq P_{y}$ for all $y\in Y$, thus
$\tau(f(y))(\{y\}) = 1$ for all $y\in Y$. 

4.
In the same way we obtain measurable $g: Z\to Y$ with the desired properties.
\EndProof

Discussing the maps $f, g$ obtained above from $\tau$ and $\theta$, we
see that $\tau\circ f = e_Y, \theta\circ g = e_Z$, and we calculate
through change of variables formula in
Corollary~\ref{int-image-measure} for each $z_0\in Y$ and each
$H\in\Borel{Z}$
\begin{align*}
\bigl(\tau\klComp(\SubProb{f}\circ\tau)\bigr)(z_0)(H)
& =
\int_Z\tau(z)(H)\ (\SubProb{f}\circ\tau)(z_0)(dz)\\
& =
\int_Y \tau(f(y))(H)\ \tau(z_0)(dy)\\
& =
\int_Y \delta_y(H)\ \tau(z_0)(dy)\\
& =
\tau(z_0)(H).
\end{align*}
Thus 
$ 
\tau\klComp(\SubProb{f}\circ\tau) = \tau,
$
and because $\tau$ is a mono, we infer that 
$
\SubProb{f}\circ\tau = e_Z.
$
Since 
\begin{equation*}
%\label{ }
\SubProb{f}\circ\tau = (e_Z\circ f)\klComp \tau = \theta\klComp \tau
\end{equation*}
we infer on account of $\tau$ being an epi that $\theta = e_Z\circ
f$. Similarly we see that $\tau = e_Y\circ g$. 

\BeginLemma{char-unit-prod}
Given stochastic relations $\theta: Y \Trans Z$ and $\tau: Z \Trans Y$ with 
$ 
\tau * \theta  = e_Y 
$
and 
$
\theta *\tau = e_Z,
$
there exist Borel isomorphisms $f: Y\to Z$ and $g: Z\to Y$ with 
$
\theta  = e_Z\circ f,
$
and
$
\tau = e_Y\circ g.
$
\EndLemma

\BeginProof
We infer for $y\in Y$ from
\begin{align*}
%\label{ }
\delta_y(G) & = e_Y(y)(G) \\
& = 
(\tau * \theta)(y)(G)\\
& =
\dInt[Z]{\tau(z)(G)}{\theta(y)(z)}\\
& =
\delta_{f(y)}(\InvBild{g}{G})\\
& =
\delta_y(\InvBild{f}{\InvBild{g}{G}})
\end{align*}
for all Borel sets $G\in\Borel{Y}$ that 
$
g\circ f = id_Y,
$
similarly,
$
f \circ g = id_Z
$
is inferred. Hence the Borel maps $f$ and $g$ are bijections, thus Borel isomorphisms.
\EndProof

This yields a characterization of the quotient equivalence relation in the Kleisli category for the Giry monad.
\BeginProposition{char-equiv-giry}
Assume the stochastic relations $K: X\Trans Y$ and $L: X \Trans Z$ are
both epimorphisms with respect to Kleisli composition, then these conditions are equivalent
\begin{enumerate}
\item \label{char-equiv-giry-1}  $\isEquiv{K}{L}{\approx}.$
\item \label{char-equiv-giry-2} $L = \SubProb{f}\circ K$ for a Borel isomorphism $f: Y \to Z$. 
\end{enumerate}
 \EndProposition

\BeginProof
\labelImpl{char-equiv-giry-1}{char-equiv-giry-2}: Because $\isEquiv{K}{L}{\approx}$, there exists an invertible $\theta: Y \Trans Z$ with inverse $\tau: Z \Trans Y$ and
$
L = \theta * K.
$
We infer from Lemma~\ref{char-unit-prod} the existence of a Borel isomorphism $f: Y\to Z$ such that $\theta = \unit_Z\circ f$. Consequently, we have for $x\in X$ and the Borel set $H\in \Borel{Z}$
\begin{align*}
L(x)(H) 
&=
\dInt[Y]{\theta(y)(H)}{K(x)(y)}\\
& =
\dInt[Y]{\delta_{f(y)}(H)}{K(x)(y)}\\
& =
K(x)(\InvBild{f}{H})\\
& = 
\bigl(\SubProb{f}\circ K\bigr)(x)(H)
\end{align*}

\labelImpl{char-equiv-giry-2}{char-equiv-giry-1}: If $L =
\SubProb{f}\circ K = (\unit_Z\circ f) * K$ for the Borel
isomorphism $f: Y \to Z$, then $K = (\unit_Y\circ g) * L$ with
$g: Z \to Y$ as the inverse to $f$.
\EndProof

Consequently, given the epimorphisms $K: X \Trans Y$ and $L: X \Trans Z$, the relation
$\isEquiv{K}{L}{\approx}$ entails their 
their base spaces $Y$ and $Z$ being Borel isomorphic, and vice versa. Hence the Borel isomorphism classes are the quotient objects for this relation. 

This classification should be complemented by a characterization of epimorphic Kleisli morphisms for this monad. This seems to be an open question.  
} 
%%% Local Variables: 
%%% mode: latex
%%% TeX-master: "../Mskr3"
%%% End: 

%\Input{\Folder/LpSpaces}
%spell checked - 24Aug14
\Subsection{$\rLpS$-Spaces}
\label{lp-spaces}

We will construct now for a measure space $(X, {\cal A}, \mu)$ a family $\bigl(\rLp{\mu}\bigr)_{1\leq p \leq \infty}$ of Banach spaces. Some properties of these spaces are discussed now, in particular we will identify their dual spaces. The case $p=2$ gives the particularly interesting space $\rLp[2]{\mu}$, which is a Hilbert space under the inner product 
$
\langle f, g\rangle \mapsto \dInt{f\cdot g}{\mu}.
$
Hilbert spaces have some properties which will turn out to be helpful, and which will be exploited for the underlying measure spaces. For example, von Neumann obtained from a representation of their continuous linear maps both the Lebesgue decomposition and the Radon-Nikodym Theorem derivative in one step! We join Rudin's exposition~\cite[Section 6]{Rudin} in giving the truly ravishing proof here. But we are jumping ahead. After investigating the basic properties of Hilbert spaces including the closest approximation property and the identification of continuous linear functions we mode to a discussion of the more general $\rLpS$-spaces and investigate the positive linear functionals on them. 

Some important developments like the definition of signed measures are briefly touched, some are not. The topics which had to be omitted here include the weak topology induced by $\rLpS[q]$ on $\rLpS$ for conjugate pairs $p, q$; this would have required some investigations into convexity, which would have led into a wonderful, wondrous but far-away country.

The last section deals with disintegration as an application of both
the Radon-Nikodym derivative and the Hahn Extension Theorem. It deals
with the problem of decomposing a finite measure on a product into its
projection onto the first component and an associated transition
kernel. This corresponds to reversing a Markov transition system with
a given initial distribution akin the converse of a relation in a
set-oriented setting.

\Subsubsection{A Spoonful Hilbert Space Theory}
\label{hilbert-spaces}

Let $H$ be a real vector space. A map $(\cdot, \cdot): H\times H\to \Real$ is said to be an \emph{\index{inner product}inner product} iff these conditions hold for all $x, y, z\in H$ and all $\alpha, \beta\in\Real$:
\begin{enumerate}
\item $(x, y) = (y, x)$,so the inner product is commutative.
\item $(\alpha\cdot x + \beta\cdot z, y) = \alpha\cdot(x, y) + \beta\cdot(z, y)$, so the inner product is linear in the first, hence also in the second component.
\item $(x, x)\geq 0$, and $(x, x) = 0$ iff $x = 0$.  
\end{enumerate}

We confine ourselves to real vector spaces. Hence the laws for the
inner product are somewhat simplified in comparison to vector spaces
over the complex number. There one would, e.g. postulate that $(y, x)$
is the complex conjugate for $(x, y)$.

The inner product is the natural generalization of the scalar product in Euclidean spaces
\begin{equation*}
  (\langle x_{1}, \dots, x_{n}\rangle, \langle y_{1}, \dots, y_{n}\rangle) := \sum_{i=1}^{n}x_{i}\cdot y_{i}, 
\end{equation*}
which satisfies these laws, as one verifies readily.

We fix an inner product $(\cdot , \cdot )$ on $H$. Define the norm of $x\in H$ through 
\begin{equation*}
  \|x\| := \sqrt{(x, x)},
\end{equation*}
this is possible because $(x, x)\geq 0$. Before investigating $\|\cdot\|$ in detail, we need the \emph{\index{Schwarz inequality}Schwarz \index{inequality!Schwarz}inequality} as a very helpful tool. It relates the norm to the inner product of two elements.

\BeginLemma{schwarz-inequality}
$|(x, y)| \leq \|x\|\cdot\|y\|$.
\EndLemma

\BeginProof
Let $a := \|x\|^{2}$, $b := \|y\|^{2}$, and $c := |(x, y)|$. Then $c = t\cdot(x, y)$ with $t\in\{-1, +1\}$. We have for each real $r$
\begin{equation*}
  0 \leq (x - r\cdot t\cdot y, x - r\cdot t \cdot y) = (x, x) - 2\cdot r \cdot t\cdot (x, y) + r^{2}\cdot (y, y),
\end{equation*}
thus $a - 2\cdot r\cdot c + r^{2}\cdot b\geq 0$. If $b = 0$, we must also have $c = 0$, otherwise the inequality would be false for large positive $r$. Hence the inequality is true in this case. So we may assume that $b \not= 0$. Put $r := c/b$, so that $a \geq c^{2}/b$, so that $a\cdot b \geq c^{2}$, from which the desired inequality follows.
\EndProof

Schwarz's inequality will help in establishing that a vector space with an inner product is a normed space.

\BeginProposition{is-a-norm}
Let $H$ be a real vector space with an inner product, then $(H, \|\cdot \|)$ is a normed space.
\EndProposition

\BeginProof
It is clear from the definition of the inner product that $\|\alpha\cdot x\| = |\alpha|\cdot \|x\|$, and that $\|x\| = 0$ iff $x = 0$; the crucial point is the triangle inequality. We have
\begin{align*}
  \|x + y\|^{2} & = (x + y, x + y) = \|x\|^{2} + \|y\|^{2} + 2\cdot (x, y)\\
& \leq \|x\|^{2} + 2\cdot \|x\|\cdot \|y\| + \|y\|^{2}&&\text{by Lemma~\ref{schwarz-inequality}}\\
& = (\|x\|+\|y\|)^{2}.
\end{align*}
\EndProof

Thus each inner product space yields a normed space,consequently it spawns a metric space through $\langle x, y\rangle \mapsto \|x-y\|$. Finite dimensional vector spaces $\Real^{n}$ are Hilbert spaces under the inner product mentioned above. It produces for $\Real^{n}$ the familiar Euclidean distance 
\begin{equation*}
  \|x - y\| = \sqrt{\sum_{i=1}^{n} (x_{i}-y_{i})^{2}}. 
\end{equation*}
We will meet square integrable functions as another class of Hilbert spaces, but before discussing them, we need some preparations.

\BeginCorollary{continuous-inner}
The maps $x \mapsto \|x\|$ and $x\mapsto (x, y)$ with $y\in H$ fixed are continuous. 
\EndCorollary

\BeginProof
We obtain from $\|x\|\leq \|y\| + \|x-y\|$ and $\|y\|\leq \|x\|+\|x-y\|$ that $\bigl|\|x\|-\|y\|\bigr|\leq \|x-y\|$, hence the norm is continuous. From Schwarz's inequality we see that $|(x, y) - (x', y)| = |(x - x', y)| \leq \|x-x'\|\cdot \|y\|$, which shows that $(\cdot , y)$ is continuous. 
\EndProof

From the properties of the inner product it is apparent that $x \mapsto (x, y)$ is a continuous linear functional:

\BeginDefinition{cont-lin-fnct}
Let $H$ be an inner product space with norm $\|\cdot \|$. A linear map $L: H\to \Real$ which is continuous in the norm topology is called a \emph{continuous linear functional} on $H$. 
\EndDefinition

If $L: H\to \Real$ is a continuous linear functional, then its \emph{\index{kernel!linear functional}kernel}
\begin{equation*}
  Kern(L) := \{x\in H \mid L(x) = 0\}
\end{equation*}
is a closed linear subspace of $H$, i.e., is a real vector space in its own right. Say that $x\in H$ is \index{orthogonal!vector}\emph{orthogonal} to $y\in H$ iff $(x, y) = 0$, and denote this by $x\bot y$. This is the generalization of the familiar concept of orthogonality in Euclidean spaces, which is formulated also in terms of the inner product. 
Given a linear subspace $M\subseteq H$, define the \emph{\index{orthogonal!complement}orthogonal complement} $M^{\bot}$ of $M$ as 
\begin{equation*}
  M^{\bot} := \{y \in H\mid x \bot y\text{ for all }x\in M\}.
\end{equation*}
The orthogonal complement is a linear subspace as well, and it is closed by Corollary~\ref{continuous-inner}, since $M = \bigcap_{x\in M}\{y\in H\mid (x, y) = 0\}$. Then $M\cap M^{\bot} = \{0\}$, since a vector $z\in M\cap M^{\bot}$ is orthogonal to itself, hence $(z, z) = 0$, which implies $z = 0$. 

Hilbert spaces are introduced now as those linear spaces for which
this metric is complete. Our goal is to show that continuous linear
functionals on a Hilbert space $H$ are given exactly through the inner
product.

\BeginDefinition{hilbert-space}
A \emph{\index{Hilbert space}Hilbert space} is a real vector space which is a complete metric space under the induced metric.  
\EndDefinition

Note that we fix the metric for which the space is to be complete,
noting that completeness is not a property of the underlying
topological space but rather of a specific metric.

Recall that a subset $C\subseteq H$ is called \emph{\index{convex}convex} iff it contains with two points also the straight line between them, thus iff $\alpha\cdot x + (1-\alpha)\cdot y\in C$, whenever $x, y\in C$ and $0 \leq \alpha\leq 1$.

A key tool for our development is the observation that a closed convex
subset of a Hilbert space has a unique element of smallest norm. This
property is familiar from Euclidean spaces. Visualize a compact convex
set in $\Real^{3}$, then this set has a unique point which is closest
to the origin. The statement below is more general, because it refers
to closed and convex sets.

\BeginProposition{convex-minimal}
Let $C\subseteq H$ be a closed and convex subset of the Hilbert space $H$. Then there exists a unique $y\in C$ such that $\|y\| = \inf_{z\in C}\|z\|$. 
\EndProposition

\BeginProof
1.
Put $r := \inf_{z\in C}\|z\|$, and let $x, y\in C$, hence by convexity $(x+y)/2\in C$ as well. The parallelogram law from Exercise~\ref{ex-parallelogram-law} gives
\begin{equation*}
  \|x-y\|^{2} = 2\cdot \|x\|^{2}+2\cdot \|y\|^{2} - 4\cdot \|(x+y)/2\|^{2} \leq 2\cdot \|x\|^{2}+2\cdot \|y\|^{2}-4\cdot r^{2}.
\end{equation*}
Hence if we have two vectors $x\in C$ and $y\in C$ of minimal norm, we obtain $x = y$. Thus, if such a vector exists, it must be unique.

2.
Let $\Folge{x}$ be a sequence in $C$ such that $\lim_{n\to \infty}\|x_{n}\|= r$. At this point, we have only informations about the sequence $\bigl(\|x_{x}\|\bigr)_{n\in\Nat}$ of real numbers, but we can actually show that the sequence proper is a Cauchy sequence. It works like this. We obtain, again from the parallelogram law, the estimate
\begin{equation*}
  \|x_{n}-x_{m}\| \leq 2\cdot (\|x_{n}\|^{2}+\|x_{m}\|^{2}-2\cdot r^{2}),
\end{equation*}
so that for each $\epsilon>0$ we can find $n_{0}$ such that $\|x_{n}-x_{m}\|<\epsilon$ if $n, m\geq n_{0}$. Hence $\Folge{x}$ is actually a Cauchy sequence, and since $H$ is complete, we find some $x$ such that $\lim_{n\to \infty}x_{n} = x$. Clearly, $\|x\| = r$, and since $C$ is closed, we infer that $x\in C$. 
\EndProof

Note how the geometric properties of an inner product space,
formulated through the parallelogram law, and the metric properties of
being complete cooperate. 

This unique approximation property has two remarkable
consequences. The first one establishes for each element $x\in H$ a
unique representation as $x = x_{1}+ x_{2}$ with $x_{1}\in M$ and
$x_{2}\in M^{\bot}$ for a closed linear subspace $M$ of $H$, and the
second one shows that the only continuous linear maps on the Hilbert
space $H$ are given by $\lambda x.(x, y)$ for $y\in H$. We need the
first one for establishing the second one, so both find their place in
this somewhat minimal discussion of Hilbert spaces.

\BeginProposition{orthogonal-projection}
Let $H$ be a Hilbert space, $M\subseteq H$ a closed linear subspace. Each $x\in H$ has a unique representation $x = x_{1} + x_{2}$ with $x_{1}\in M$ and $x_{2}\in M^{\bot}$.
\EndProposition

\BeginProof
1.
If such a representation exists, it must be unique. In fact, assume that $x_{1} + x_{2} = x = y_{1} + y_{2}$ with $x_{1}, y_{1}\in M$ and $x_{2}, y_{2}\in M^{\bot}$, then $x_{1}-y_{1} = y_{2}-x_{2}\in M\cap M^{\bot}$, which implies $x_{1} = y_{1}$ and $x_{2}= y_{2}$ by the remark above. 

2.
Fix $x\in H$, we may and do assume that $x\not\in M$, and define $C := \{x-y\mid y\in M\}$, then $C$ is convex, and, because $M$ is closed, it is closed as well. Thus we can find an element in $C$ which is of smallest norm, say, $x - x_{1}$ with $x_{1}\in M$. Put $x_{2} := x-x_{1}$, and we have to show that $x_{2}\in M^{\bot}$, hence that $(x_{2}, y) = 0$ for any $y\in M$. Let $y\in M, y\not=0$ and choose $\alpha\in\Real$ arbitrarily (for the moment, we'll fix it later). Then $x_{2} -\alpha\cdot y = x - (x_{1}+\alpha\cdot y)\in C$, thus $\|x_{2}-\alpha\cdot y\|^{2}\geq \|x_{2}\|^{2}$. Expanding, we obtain
\begin{equation*}
  (x_{2}-\alpha\cdot y, x_{2}-\alpha\cdot y) = (x_{2}, x_{2}) - 2\cdot \alpha\cdot (x_{2}, y) + \alpha^{2}\cdot (y, y) \geq (x_{2}, x_{2}).
\end{equation*}
Now put $\alpha := (x_{2}, y)/(y, y)$, then the above inequality yields
\begin{equation*}
  -2\cdot \frac{(x_{2}, y)^{2}}{(y, y)} + \frac{(x_{2}, y)^{2}}{(y, y)} \geq 0,
\end{equation*}
which implies $-(x_{2}, y)^{2}\geq 0$, hence $(x_{2}, y) = 0$. Thus $x_{2}\in M^{\bot}$. 
\EndProof

Thus $H$ is decomposed into $M$ and $M^{\bot}$ for any closed linear subspace $M$ of $H$ in the sense that each element of $H$ can be written as a sum of elements of $M$ and of $M^{\bot}$, and, even better, this decomposition is unique. These elements are perceived as the projections to the subspaces. In the case that we can represent $M$ as the kernel $\{x\in H\mid  L(x) = 0\}$ of a continuous linear map $L: H\to \Real$ with $L \not= 0$ we can say actually more.

\BeginLemma{decomp-factor-linear}
Let for Hilbert space $H$, $L: H\to \Real$ be a continuous linear functional with $L\not= 0$. Then $Kern(L)^{\bot}$ is isomorphic to $\Real$
\EndLemma

\BeginProof 
Define $\phi(y) := L(y)$ for $y\in Kern{L}^{\bot}$. Then
$\phi(\alpha\cdot y+\beta\cdot y') = \alpha\cdot \phi(y)+\beta\cdot
\phi(y')$ follows from the linearity of $L$. If $\phi(y) = \phi(y')$,
then $y-y'\in Kern(L)\cap Kern(L)^{\bot}$, so that $y = y'$, so that $\phi$ is one-to-one. Given $t\in\Real$, we can find $x\in H$ with $L(x) = t$; decompose $x$ as $x_{1}+x_{2}$ with $x_{1}\in Kern(L)$ and $x_{2}\in Kern(L)^{\bot}$, then $\phi(x_{2}) = L(x - x_{1}) = t$. Thus $\phi$ is onto. Hence we have found a linear and bijective map $Kern(L)^{\bot}\to \Real$. 
\EndProof

Returning to the decomposition of an element $x\in H$, we fix an arbitrary $y\in Kern(L)\setminus\{0\}$. Then we may write $x = x_{1} +\alpha\cdot y$, where $\alpha\in\Real$. This follows immediately from Lemma~\ref{decomp-factor-linear}, and it has the consequence we are aiming at.

\BeginTheorem{repr-for-lin-op}
Let $H$ be a Hilbert space, $L: H\to \Real$ be a continuous linear functional. Then there exists $y\in H$ with $L(x) = (x, y)$ for all $x\in H$. 
\EndTheorem

\BeginProof
If $L = 0$, this is trivial. Hence we assume that $L \not= 0$. Thus we can find $z\in Kern(L)^{\bot}$ with $L(z) = 1$; put $y = \gamma\cdot z$ so that $L(y) = (y, y)$. Each $x\in H$ can be written as $x = x_{1} + \alpha\cdot y$ with $x_{1}\in Kern(L)$. Hence 
\begin{equation*}
L(x) = L(x_{1} + \alpha\cdot y) = \alpha\cdot L(y) = \alpha\cdot (y, y) = (x_{1}+\alpha\cdot y, y) = (x, y).
\end{equation*}
Thus $L = \lambda x.(x, y)$ is established. 
\EndProof

This rather abstract view of Hilbert spaces will be put to use now to the more specific case of integrable functions. 

\Subsubsection{The $\rLpS$-Spaces are Banach Spaces}
\label{sec:are-banach-spaces}

We will investigate now the structure of integrable functions for a
fixed $\sigma$-finite measure space $(X, {\cal A}, \mu)$. We will
obtain a family of Banach spaces, which have some interesting
properties. In the course of investigations, we will usually not
distinguish between functions which differ only on a set of measure
zero (because the measure will not be aware of the differences). For
this, we introduced above the equivalence relation $=_{\mu}$ (``equal
$\mu$.almost everywhere'') in with $f =_{\mu} g$ iff $\mu(\{x\in X\mid
f(x)\not= g(x)\}) = 0$, see Section~\ref{sec:ess-bounded-fncts} on
page~\pageref{sec:ess-bounded-fncts}. In those cases where we will
need to look at the value of a function at certain points, we will
make sure that we will point out the difference.

Let us see how this works in practice. Define
\begin{equation*}\textstyle
  \cLp[1]{\mu} := \{f\in \MeasbFnct{X, {\cal A}} \mid \dInt{|f|}{\mu}<\infty\},
\end{equation*}
thus $f\in\cLp[1]{\mu}$ iff $f: X\to \Real$ is measurable and has a finite $\mu$-integral. 

Then this space defines a vector space and closed with respect to $|\cdot|$, hence we have immediately

\BeginProposition{lp1-is-vector-lattice}
$\cLp[1]{\mu}$ is a vector lattice.\QED
\EndProposition

Now put
\begin{equation*}
  \rLp[1]{\mu} := \{[f]\mid f\in\rLp[1]{\mu}\},
\end{equation*}
then we have to explain how to perform the algebraic operations on the equivalence classes (note that we write $[f]$ rather than $[f]_{\mu}$, which we will do when more than one measure has to be involved). Since the set of all nullsets is a $\sigma$-ideal, these operations are easily shown to be well-defined:
\begin{align*}
  [f] + [g] & := [f + g],\\
  [f] \cdot  [g] & := [f\cdot  g],\\
\alpha\cdot [f] & := [\alpha\cdot f]
\end{align*}

Thus we obtain
\BeginProposition{lp1-is-vector-lattice-1}
$\rLp[1]{\mu}$ is a vector lattice.\QED
\EndProposition

Let $f\in\rLp[1]{\mu}$, then we define
\begin{equation*}
  \|f\|_{1}:= \dInt{|f|}{\mu}
\end{equation*}
as the $\rLpS[1]$-norm for $f$. Let us have a look at the properties which a decent norm should have. First, we have $\|f\|_{1}\geq 0$, and $\|\alpha\cdot f\|_{1} = \alpha\cdot \|f\|_{1}$, this is immediate. Because $|f + g| \leq |f| + |g|$, the triangle inequality holds. Finally, let $\|f\|_{1} = 0$, thus $\dInt{|f|}{\mu} = 0$, consequently, $f=_{\mu}0$, which means $f = [0]$. 

This will be a basis for the definition of a whole family of linear spaces of integrable functions. Call the positive real numbers $p$ and $q$ \emph{\index{conjugate numbers}conjugate} iff they satisfy
\begin{equation*}
  \frac{1}{p} + \frac{1}{q} = 1
\end{equation*}
(for example, $2$ is conjugate to itself). This may be extended to $p = 0$, so that we also consider $0$ and $\infty$ as conjugate numbers, but using this pair will be made explicit. 

The first step for extending the definition of $\rLpS[1]$ will be \emph{Hölder's \index{inequality!H\"older}inequality}, which is based on this simple geometric fact:

\BeginLemma{basis-for-hoelder}
Let $a, b$ be positive real numbers, $p>0$ conjugate to $q$, then 
\begin{equation*}
  a\cdot b \leq \frac{a^{p}}{p} + \frac{b^{q}}{q},
\end{equation*}
equality holding iff $b = a^{p-1}$.
\EndLemma

\BeginProof
The exponential function is convex, i.e., we have 
\begin{equation*}
  e^{(1-\alpha)\cdot x - \alpha\cdot y}\leq (1-\alpha)\cdot e^{x} + \alpha\cdot e^{y}
\end{equation*}
for all $x, y\in\Real$ and $0\leq \alpha\leq 1$. Because both $a> 0$ and $b>0$, we find $r, s$ such that $a = e^{r/p}$ and $b = e^{s/q}$. Since $p$ and $q$ are conjugate, we obtain from $1/p = 1-1/q$
\begin{equation*}
  a\cdot b = e^{r/p+s/q} \leq \frac{e^{s}}{p} + \frac{e^{q}}{q} = \frac{a^{p}}{p} + \frac{b^{q}}{q}.
\end{equation*}
\EndProof

This betrays one of the secrets of conjugate $p$ and $q$, viz., that they give rise to a convex combination. 

We are ready to formulate and prove \emph{\index{H\"older's inequality}Hölder's \index{inequality!H\"older}inequality}, arguably one of the most frequently used inequalities in integration (as we will see as well); the proof follows the one given for~\cite[Theorem 3.5]{Rudin}.

\BeginProposition{hoelder-inequality}
Let $p>0$ and $q>0$ be conjugate, $f$ and $g$ be non-negative measurable functions on $X$. Then
\begin{equation*}
  \dInt{f\cdot g}{\mu} \leq \bigl(\dInt{f^{p}}{\mu}\bigr)^{1/p}\cdot \bigl(\dInt{g^{q}}{\mu}\bigr)^{1/q}.
\end{equation*}
\EndProposition

\BeginProof
Put for simplicity
\begin{equation*}
  A := \bigl(\dInt{f^{p}}{\mu}\bigr)^{1/p}\text{ and }
B := \bigl(\dInt{g^{q}}{\mu}\bigr)^{1/q}.
\end{equation*}
If $A = 0$, we may conclude from $f=_{\mu}0$ that $f\cdot g=_{\mu}0$, so there is nothing to prove. If $A > 0$ but $b=\infty$, the inequality is trivial, so we assume that $0 < A < \infty, 0 < B < \infty$. Put 
\begin{equation*}
  F := \frac{f}{A}, G := \frac{g}{B},
\end{equation*}
thus we obtain
\begin{equation*}
  \dInt{F^{p}}{\mu} = \dInt{G^{q}}{\mu} = 1.
\end{equation*}
We obtain $F(x)\cdot G(x)\leq F(x)^{p}/p + G(x)^{q}/q$ for every $x\in X$ from Lemma~\ref{basis-for-hoelder}, hence
\begin{equation*}
  \dInt{F\cdot G}{\mu} \leq \frac{1}{p}\cdot \dInt{F^{p}}{\mu} + \frac{1}{q}\cdot \dInt{G^{q}}{\mu} \leq \frac{1}{p}+\frac{1}{q} = 1.
\end{equation*}
Multiplying both sides with $A\cdot B>0$ now yields the desired result.
\EndProof

This gives \emph{\index{Minkowski's inequality}Minkowski's \index{inequality!Minkowski}inequality} as a consequence. Put for $f: X\to \Real$ measurable, and for $p\geq 1$
\begin{equation*}
  \|f\|_{p} := \bigl(\dInt{|f|^{p}}{\mu}\bigr)^{1/p}.
\end{equation*}

\BeginProposition{minkowski-inequality}
Let $1 \leq p < \infty$ and let $f$ and $g$ be non-negative measurable functions on $X$. Then
\begin{equation*}
  \|f+g\|_{p}\leq \|f\|_{p}+\|g\|_{p}
\end{equation*}
\EndProposition

\BeginProof
The inequality follows for $p = 1$ from the triangle inequality for $|\cdot |$, so we may assume that $p > 1$. We may also assume that $f, g\geq 0$. Then we obtain from Hölder's inequality with $q$ conjugate to $p$ 
\begin{align*}
  \|f+g\|_{p}^{p} & = \dInt{(f+g)^{p-1}\cdot f}{\mu}+\dInt{(f+g)^{p-1}\cdot g}{\mu}\\
& \leq \|f+g\|_{p}^{p/q}\cdot \bigl(\|f\|_{p}+\|g\|_{p}\bigr)
\end{align*}
Now assume that $\|f+g\|_{p}=\infty$, we may divide by the factor $\|f+g\|_{p}^{p/q}$, and we obtain the desired inequality from $p - p/q = p\cdot (1-1/q) = 1$. If, however, the left hand side is infinite, then the inequality
\begin{equation*}
  (f+g)^{p}\leq 2^{p}\cdot max\{f^{p}, g^{p}\}\leq 2^{p}\cdot (f^{p}+g^{p})
\end{equation*}
shows that the right hand side is infinite as well. 
\EndProof

Given $1 \leq p < \infty$, define
\begin{equation*}
  \cLp{\mu} := \{f\in\MeasbFnct{X, \mathcal{A}} \mid \|f\|_{p}<\infty\}
\end{equation*}
with $\rLp{\mu}$ as the corresponding set of $=_{\mu}$-equivalence classes. An immediate consequence from Minkowski's inequality is

\BeginProposition{lp-are-vector-spaces}
$\cLp{\mu}$ is a linear space over $\Real$, and $\|\cdot \|_{p}$ is a pseudo-norm on it. $\rLp{\mu}$ is a normed space. 
\EndProposition

\BeginProof
It is immediate from Proposition~\ref{minkowski-inequality} that $f+g\in\cLp{\mu}$ whenever $f, g\in\cLp{\mu}$, and $\cLp{\mu}$ is closed under scalar multiplication as well. That $\|\cdot \|_{p}$ is a pseudo-norm is also immediate. Because scalar multiplication and addition are compatible with forming equivalence classes, the set $\rLp{\mu}$ of classes is a real vector space as well. As usual, we will identify $f$ with its class, unless otherwise stated. Now $f\in\rLp{\mu}$ with $\|f\|_{p}=0$, then $|f|=_{\mu} 0$, hence $f=_{\mu}0$, thus $f = 0$. So $\|\cdot \|_{p}$ is a norm on $\rLp{\mu}$. 
\EndProof

In Section~\ref{sec:ess-bounded-fncts} the vector spaces ${\cal L}_{\infty}(\mu)$ and $L_{\infty}(\mu)$ are introduced, so we have now a family  $\bigl({\cal L}_{p}(\mu)\bigr)_{1\leq p \leq \infty}$ of vector spaces together with their associated spaces $\bigl({L}_{p}(\mu)\bigr)_{1\leq p \leq \infty}$ of $\mu$-equivalence classes, which are normed spaces. They share the property of being Banach spaces.

\BeginProposition{all-lp-are-banach}
$\rLp{\mu}$ is a Banach space for $1\leq p \leq \infty$.
\EndProposition

\BeginProof
1.
Let us first assume that the measure is finite. 
We know already from Proposition~\ref{linf-is-banach-space} that ${\cal L}_{\infty}(\mu)$ is a Banach space, so we may assume that $p< \infty$. 

Let $\Folge{f}$ be a Cauchy sequence in $\rLp{\mu}$, then we obtain
\begin{equation*}
  \epsilon^{p}\cdot \mu(\{x\in X\mid |f_{n}-f_{m}|\geq \epsilon\}) \leq \dInt{|f_{n}-f_{m}|^{p}}{\mu}.
\end{equation*}
 for $\epsilon>0$. Thus $\Folge{f}$ is a Cauchy sequence for convergence in measure, so we can find $f\in\MeasbFnct{X, \mathcal{A}}$ such that $f_{n}\nmC f$ by Proposition~\ref{f-is-banach-space}. Proposition~\ref{conv-in-measure-subseq} tells us that we can find a subsequence $(f_{n_{k}})_{k\in\Nat}$ such that $f_{n_{k}}\aeC f$. But we do not yet know that $f\in\cLp{\mu}$. We infer $\lim_{k\to \infty} |f_{n_{k}}-f|^{p}\to 0$ outside a set of measure zero. Thus we obtain from Fatou's Lemma Proposition~\ref{fatous-lemma} for every $n\in\Nat$
 \begin{equation*}
   \dInt{|f-f_{n}|^{p}}{\mu}\leq \liminf_{k\to \infty}\dInt{|f_{n_{k}} - f_{n}|^{p}}{\mu}.
 \end{equation*}
Thus $f - f_{n}\in\cLp{\mu}$ for all $n\in\Nat$, and from $f = (f - f_{n}) + f_{n}$  we infer $f\in\cLp{\mu}$, since $\cLp{\mu}$ is closed under addition. We see also that $\|f-f_{n}\|_{p}\to 0$, as $n\to \infty$. 

2.
If the measure space is $\sigma$-finite, we may write $\dInt{f}{\mu}$ as $\lim_{n\to \infty}\dInt[A_{n}]{f}{\mu}$, where $\mu(A_{n})<\infty$ for an increasing sequence $\Folge{A}$ of measurable sets with $\bigcup_{n\in\Nat}A_{n} = X$. Since the restriction to each $A_{n}$ yields a finite measure space, where the result holds, it is not difficult to see that completeness holds for the whole space as well. Specifically, given $\epsilon>0$, there exists $n_{0}\in\Nat$ so that for all $n, m\geq n_{0}$
\begin{equation*}
  \|f_{n}-f_{m}\|_{p} \leq \|f_{n}-f_{m}\|_{p}^{(n)} + \epsilon
\end{equation*}
holds, with $\|g\|_{p}^{(n)} := \bigl(\dInt{|g|^{p}}{\mu_{n}}\bigr)^{1/p}$, and $\mu_{n}: B \mapsto \mu(B\cap A_{n})$ as the measure $\mu$ localized to $A_{n}$. Then $\|f_{n}-f\|_{p}^{(n)}\to 0$, from which we obtain $\|f_{n}-f\|_{p}\to 0$. Hence completeness is also valid for the $\sigma$-finite case. 
\EndProof

\BeginExample{lp-spaces-discrete}
Let $|\cdot |$ be the counting measure on $(\Nat, \PowerSet{\Nat})$, then this is a $\sigma$-finite measure space. Define 
\begin{align*}
  \ell_{p} & := \rLp{|\cdot |}, 1 \leq p < \infty,\\
\ell_{\infty} & := L_{\infty}(|\cdot |).
\end{align*}
Then $\ell_{p}$ is the set of all real sequences $\Folge{x}$ with $\sum_{n\in\Nat}|x_{n}|^{p}<\infty$, and $\Folge{x}\in\ell_{\infty}$ iff $\sup_{n\in\Nat}|x_{n}|<\infty$. Note that we do not need to pass to equivalence classes, since $|A| = 0$ iff $A= \emptyset$. These spaces are well known and well studied; they will not be considered further. 
\EndExample

The case $p = 2$ deserves particular attention, since the norm is in this case obtained from the inner product
\begin{equation*}
  (f, g) := \dInt{f\cdot g}{\mu}.
\end{equation*}
In fact, linearity of the integral shows that 
\begin{equation*}
  (\alpha\cdot f + \beta\cdot g, h) = \alpha\cdot (f, h) + \beta\cdot (g, h)
\end{equation*}
holds, commutativity of multiplications yields $(f, g) = (g, f)$, finally it is clear that $(f, f)\geq 0$ always holds. If we have $f\in\cLp[2]{\mu}$ with $f=_{\mu}0$, then we know that also $(f, f)=0$, thus $(f, f)=0$ iff $f = 0$ in $\rLp[2]{\mu}$. 

Thus we obtain from Proposition~\ref{all-lp-are-banach}

\BeginCorollary{l2-is-hilbert}
$\rLp{\mu}$ is a Hilbert space with the inner product $(f, g) := \dInt{f\cdot g}{\mu}.$ \QED
\EndCorollary

This will have some interesting consequences, which we will explore in Section~\ref{sec:radon-nikodym}.

Before doing so, we show that the step functions belonging to $\rLpS$ are dense.

\BeginCorollary{step-is-dense}
Given $1 \leq p < \infty$, the set 
\begin{equation*}
D := \{f\in\StepFnct{X, \mathcal{A}}\mid \mu(\{x\in X \mid f(x)\not=0\})<\infty\}
\end{equation*}
is dense in $\rLp{\mu}$ with respect to $\|\cdot \|_{p}$. 
\EndCorollary

\BeginProof
The proof makes use of the fact that the step functions are dense with
respect to pointwise convergence: we'll just have to filter out those
functions which are in $\rLp{\mu}$.  Assume that $f\in\cLp{\mu}$ with
$f\geq 0$, then there exists by Proposition~\ref{ApproxStepFncts} an
increasing sequence $\Folge{g}$ of step functions with $f(x) =
\lim_{n\to \infty}f_{n}(x)$. Because $0\leq g_{n}\leq f$, we conclude
$g_{n}\in D$, and we know from Lebesgue's Dominated Convergence
Theorem~\ref{lebesgue-dominated-convergence} that $\|f-g_n\|_{p}\to
0$. Thus every non-negative element of $\cLp{\mu}$ can be approximated
through elements of $D$ in the $\|\cdot \|_{p}$-norm. In the general case,
decompose $f = f^{+} - f^{-}$ and apply the argument to both summands
separately.
\EndProof

Because the rationals form a countable and dense subset of the reals,
we take all step functions with rational coefficients, and obtain

\BeginCorollary{step-is-dense-separable}
$\rLp{\mu}$ is a separable Banach space for $1 \leq p < \infty$. \QED
\EndCorollary

Note that we did exclude the case $p=\infty$; in fact, $\rLp[\infty]{\mu}$ is usually not a separable Banach space, as this example shows.

\BeginExample{linf-not-separable}
Let $\lambda$ be Lebesgue measure on the Borel sets of the unit interval $[0, 1]$. Put $f_{t}:=\chi_{[0, t]}$ for $0\leq t \leq 1$, then $f_{t}\in\rLp[\infty]{\lambda}$ for all $t$, and we have 
$
\infNorm{f_{s} - f_{t}}{\lambda} = 1
$
for $0 < s < t < 1.$ Let 
$
K_{t} := \{f\in \rLp[\infty]{\lambda} \mid \infNorm{f - f_{t}}{\lambda} < 1/2\},
$
thus $K_{s}\cap K_{t} = \emptyset$ for $s\not= t$ (if $g \in K_{s}\cap K_{t}$, then 
$
\infNorm{f_{s} - f_{t}}{\lambda} \leq \infNorm{g - f_{t}}{\lambda} + \infNorm{f_{s} - g}{\lambda} < 1
$). On the other hand, each $K_{t}$ is open, so if we have a countable subset $D\subseteq \rLp[\infty]{\lambda}$, then $K_{t}\cap D = \emptyset$ for uncountably many $t$. Thus $D$ cannot be dense. But this means that $\rLp[\infty]{\lambda}$ is not separable.  
\EndExample

This is the first installment on the properties of $\rLpS$-spaces. We
will be back with a general discussion in Section~\ref{sec:cont-lin-fnct-lp} after having
explored the Lebesgue-Radon-Nikodym Theorem as a valuable tool in
general, and for our discussion.

%%% Local Variables: 
%%% mode: latex
%%% TeX-master: "../Mskr3"
%%% End: 

%\Input{\Folder/RadonNikodym}
%spell checked - 24Aug14
\Subsubsection{The Lebesgue-Radon-Nikodym Theorem}
\label{sec:radon-nikodym}

The Hilbert space structure of the $\rLpS[2]$ spaces will now be used for decomposing a measure into an absolutely and a singular part with respect to another measure, and for constructing a density. This construction requires a more general study of the relationship between two measures. 

We even go a bit beyond that and define absolute continuity and
singularity as a relationship of two arbitrary additive set
functions. This will be specialized fairly quickly to a relationship
between finite measures, but this added generality will turn out to be
beneficial nevertheless, as we will see.
{
% beware of local macros
\def\alpha{\rho}
\def\beta{\zeta}
\BeginDefinition{abs-cont-sing}
Let $(X, {\cal A})$ be a measurable space with two additive set functions $\alpha, \beta: {\cal A}\to \Real$. 
\begin{enumerate}
\item $\alpha$ is said to be \emph{\index{absolute continuity}absolutely continuous} with respect to $\beta$  ($\isEquiv{\alpha}{\beta}{\absCont}$) iff $\alpha(E) = 0$ for every $E\in{\cal A}$ for which $\beta(A) = 0$.
\item $\alpha$ is said to be \emph{concentrated} on $A\in{\cal A}$ iff $\alpha(E) = \alpha(E\cap A)$ for all $E\in{\cal A}$.
\item $\alpha$ and $\beta$ are called \emph{\index{mutual singular}mutually singular} ($\isEquiv{\alpha}{\beta}{\bot}$) iff there exists a pair of disjoint sets $A$ and $B$ such that $\alpha$ is concentrated on $A$ and $\beta$ is concentrated on $B$. 
\end{enumerate}
\EndDefinition

If two additive set functions are mutually singular, they live on
disjoint measurable sets in the same measurable space. These are
elementary properties.

\BeginLemma{elem-abs-cont-sing}
Let $\alpha_{1}, \alpha_{2}, \beta: {\cal A}\to \Real$ additive set functions, then we have for $a_{1}, a_{2}\in\Real$
\begin{enumerate}
\item \label{e-a-c-s-1} If $\isEquiv{\alpha_{1}}{\beta}{\bot}$, and  $\isEquiv{\alpha_{2}}{\beta}{\bot}$, then $\isEquiv{a_{1}\cdot \alpha_{1} + a_{2}\cdot \alpha_{2}}{\beta}{\bot}$.
\item \label{e-a-c-s-2} If $\isEquiv{\alpha_{1}}{\beta}{\absCont}$, and  $\isEquiv{\alpha_{2}}{\beta}{\absCont}$, then $\isEquiv{a_{1}\cdot \alpha_{1} + a_{2}\cdot \alpha_{2}}{\beta}{\absCont}$.
\item \label{e-a-c-s-3} If $\isEquiv{\alpha_{1}}{\beta}{\absCont}$ and $\isEquiv{\alpha_{2}}{\beta}{\bot}$, then $\isEquiv{\alpha_{1}}{\alpha_{2}}{\bot}$.
\item \label{e-a-c-s-4} If $\isEquiv{\alpha}{\beta}{\absCont}$ and $\isEquiv{\alpha}{\beta}{\bot}$, then $\alpha = 0$.
\end{enumerate}
\EndLemma

\BeginProof
1.  For proving~\ref{e-a-c-s-1}, note that we can find a measurable
set $B$ and sets $A_{1}, A_{2}\in {\cal A}$ with $B\cap(A_{1}\cup
A_{2})=\emptyset$ with $\beta(E) = \beta(E\cap B)$ and $\alpha_{i}(E)
= \alpha_{i}(E\cap A_{i})$ for $i = 1, 2$. By additivity, we obtain
$(a_{1}\cdot \alpha_{1}+a_{2}\cdot \alpha_{2})(E) = (a_{1}\cdot
\alpha_{1}+a_{2}\cdot \alpha_{2})(E\cap(A_{1}\cup A_{2}))$. Property~\ref{e-a-c-s-2} is obvious.

2.
$\alpha_{2}$ is concentrated on $A_{2}$, $\beta$ is concentrated on $B$ with $A \cap B =\emptyset$, hence $\beta(E\cap A_{2}) = 0$, thus $\alpha_{1}(E \cap A_{2}) = 0$ for all $E\in{\cal A}$. Additivity implies $\alpha_{1}(E) = \alpha_{1}\bigl(E\cap (X\setminus A_{2})\bigr)$, so $\alpha_{1}$ is concentrated on $X\setminus A_{2}$. This proves~\ref{e-a-c-s-3}. For proving~\ref{e-a-c-s-4}, note that $\isEquiv{\alpha}{\beta}{\absCont}$ and $\isEquiv{\alpha}{\beta}{\bot}$ imply $\isEquiv{\alpha}{\alpha}{\bot}$ by property~\ref{e-a-c-s-3}, which implies $\alpha=0$.  
\EndProof
}

We specialize these relations now to finite measures on ${\cal A}$. Absolute continuity can be expressed in a different way, which makes the concept more transparent.

\BeginLemma{equiv-abs-cont}
Given measures $\mu$ and $\nu$ on a measurable space $(X, {\cal A})$, these conditions are equivalent:
\begin{enumerate}
\item\label{equiv-abs-cont:1} $\isEquiv{\mu}{\nu}{\absCont}$.
\item\label{equiv-abs-cont:2} For every $\epsilon>0$ there exists $\delta>0$ such that $\nu(A)<\delta$ implies $\mu(A)<\epsilon$ for all measurable sets $A\in{\cal A}$.
\end{enumerate}
\EndLemma

\BeginProof
\labelImpl{equiv-abs-cont:1}{equiv-abs-cont:2}:
Assume that we can find $\epsilon>0$ so that there exist sets $A_{n}\in{\cal A}$ with $\nu(A_{n})<2^{-n}$ but $\mu(A_{n})\geq \epsilon$. Then we have $\mu(\bigcup_{k\geq n}A_{k})\geq \epsilon$ for all $n \in\Nat$, consequently, by monotone convergence, also $\mu\bigl(\bigcap_{n\in\Nat}\bigcup_{k\geq n}A_{k}\bigr) \geq \epsilon$. On the other hand, $\nu(\bigcup_{k\geq n}A_{k}) \leq \sum_{k\geq n}2^{-k} = 2^{-n+1}$ for all $n \in\Nat$, so by monotone convergence again, $\nu\bigl(\bigcap_{n\in\Nat}\bigcup_{k\geq n}A_{k}\bigr) =0.$ Thus $\isEquiv{\mu}{\nu}{\absCont}$ does not hold. 

\labelImpl{equiv-abs-cont:2}{equiv-abs-cont:1}:
Let $\nu(A) = 0$, then $\mu(A) \leq \epsilon$ for every $\epsilon>0$, hence $\isEquiv{\mu}{\nu}{\absCont}$ is true.
\EndProof

Given this equivalence, absolute continuity could have been be defined
akin to the well-known $\epsilon$-$\delta$ definition of
continuity for real functions. Then the name becomes a bit more
descriptive.

Given two measures $\mu$ and $\nu$, one, say $\mu$, can be decomposed
uniquely as a sum $\mu_{a}+\mu_{s}$ such that
$\isEquiv{\mu_{a}}{\nu}{\absCont}$ and $\isEquiv{\mu_{s}}{\nu}{\bot}$,
additionally $\isEquiv{\mu_{s}}{\mu_{a}}{\bot}$ holds. This is stated
and proved in the following theorem, which actually shows much more,
viz., that there exists a \emph{density} $h$ of $\mu_{a}$ with respect
to $\nu$. This means that $\mu_{a}(A) = \dInt[A]{h}{\nu}$ holds for all
$A\in{\cal A}$. What this is will be described now also in greater
detail. Before entering into formalities, it is noted that the
decomposition is usually called the \emph{\index{Lebesgue
    decomposition}Lebesgue decomposition} of $\mu$ with respect to
$\nu$, and that the density $h$ is usually called the
\emph{\index{Radon-Nikodym derivative}Radon-Nikodym derivative} of
$\mu_{a}$ with respect to $\nu$ and denoted by \index{$d\mu/d\nu$}$d\mu/d\nu$.

The proof both for the existence of Lebesgue decomposition and of the
Radon-Nikodym derivative is done in one step. The beautiful proof
given below was proposed by von Neumann, see~\cite[6.9]{Rudin}. Here we go:

\BeginTheorem{lebesgue-radon-nikodym}
Let $\mu$ and $\nu$ be finite measures on $(X, {\cal A})$. 
\begin{enumerate}
\item \label{l-r-n-1} There exists a unique pair $\mu_{a}$ and $\mu_{s}$ of finite measures on $(X,{\cal A})$ such that $\mu = \mu_{a} + \mu_{s}$ with $\isEquiv{\mu_{a}}{\nu}{\absCont}$, $\isEquiv{\mu_{a}}{\nu}{\bot}$. In addition, $\isEquiv{\mu_{a}}{\mu_{s}}{\bot}$ holds.
\item \label{l-r-n-2} There exists a unique $h\in\rLp[1]{\nu}$ such that 
  \begin{equation*}
    \mu_{a}(A) = \dInt[A]{h}{\nu}
  \end{equation*}
for all $A\in{\cal A}$. 
\end{enumerate}
\EndTheorem

The line of attack will be as follows: we show that $f\mapsto
\dInt{f}{\mu}$ is a continuous linear functional on the Hilbert space
$\rLp[2]{\mu+\nu}$. By the representation for these functionals on
Hilbert spaces, we can express this functional through some function
$g\in\cLp[2]{\mu+\nu}$, hence $\dInt{f}{\mu} = \dInt{f\cdot
  g}{(\mu+\nu)}$ (note the way the measures $\mu$ and $\mu+\nu$
interact by exploiting the integral with respect to $\mu$ as a linear
functional on $\rLp[2]{\mu}$). A closer investigation of $g$ will then
yield the sets we need for the decomposition, and permit constructing
the density $h$.

\BeginProof
1.
Define the finite measure $\phi := \mu+\nu$ on ${\cal A}$; note that $\dInt{f}{\phi} = \dInt{f}{\mu} +\dInt{f}{\nu}$ holds for all measurable $f$ for which the sum on the right hand side is defined; this follows from Levi's Theorem~\ref{beppo-levi} (for $f\geq 0$) and from additivity (for general $f$). We show first that $L: f\mapsto \dInt{f}{\mu}$ is a continuous linear operator on $\rLp[2]{\phi}$. In fact, 
\begin{equation*}
 \bigl|\dInt{f}{\mu}\bigr| \leq \dInt{|f|}{\phi} =\dInt{|f|\cdot 1}{\phi} \leq \bigl(\dInt{}{|f|^{2}}\bigr)^{1/2}\cdot \sqrt{\phi(X)}
\end{equation*}
by Schwarz's inequality (Lemma~\ref{schwarz-inequality}). Thus 
\begin{equation*}
  \sup_{\|f\|_{2}\leq 1} |L(f)|  \leq \sqrt{\phi(X)} < \infty.
\end{equation*}
Hence $L$ is continuous (Exercise~\ref{ex-linear-cont-is-bounded}), thus by Theorem~\ref{repr-for-lin-op} there exists $g\in\cLp[2]{\mu}$ such that 
\begin{equation}
\label{eq-l-r-n}
  L(f) = \dInt{f\cdot g}{\phi}
\end{equation}
for all $f\in\rLp[2]{\mu}$.

2.
Let $f=\chi_{A}$ for $A\in{\cal A}$, then we obtain $\dInt[A]{g}{\phi} = \mu(A) \leq \phi(A)$ from (\ref{eq-l-r-n}). This yields $0\leq g\leq 1$ $\phi$-a.e.; we can change $g$ on a set of $\phi$-measure $0$ to the effect that $0\leq g(x) \leq 1$ holds for all $x\in X$. This will not affect the representation in~(\ref{eq-l-r-n}). 

We know that  
\begin{equation}
\label{eq-l-r-n-1}
  \dInt{(1-g)\cdot f}{\mu} = \dInt{f\cdot g}{\nu}
\end{equation}
holds for all $f\in\rLp[2]{\phi}$. Put
\begin{align*}
  A & := \{x\in X \mid 0\leq g(x) < 1\},\\
 B & := \{x\in X\mid g(x) = 1\},
\end{align*}
then $A, B\in{\cal A}$, and we define for $E\in{\cal A}$
\begin{align*}
  \mu_{a}(E) & := \mu(E\cap A),\\ 
\mu_{s}(E) & := \mu(E\cap B).
\end{align*}
If $f=\chi_{B}$, then we obtain from~(\ref{eq-l-r-n-1}) $\nu(B) = \dInt[B]{g}{\nu} = \dInt[B]{0}{\mu} = 0$, thus $\nu(B) = 0$ so that $\isEquiv{\mu_{s}}{\nu}{\bot}$. 

3.
Replace for a fixed $E\in{\cal A}$ in~(\ref{eq-l-r-n-1}) the function $f$ by $(1 + g + \dots + g^{n})\cdot \chi_{E}$, then we have
\begin{equation*}
  \dInt[E]{(1-g^{n+1})}{\mu} = \dInt[E]{g\cdot (1 + g +\dots + g^{n})}{\nu}.
\end{equation*}
Look at the integrand on the right hand side: it equals zero on $B$, and increases monotonically to $1$ on $A$, hence $\lim_{n\to \infty}\dInt[E]{(1-g^{n+1})}{\mu} = \mu(E\cap A) = \mu_{a}(E)$. This provides a bound for the left hand side for all $n\in\Nat$. The integrand on the left hand side converges monotonically to some function $0\leq h\in\cLp[1]{\nu}$ with 
$\lim_{n\to \infty}\dInt[E]{g\cdot (1 + g+ \dots + g^{n})}{\nu} = \dInt[E]{h}{\nu}$ by Levi's Theorem~\ref{beppo-levi}. Hence we have
\begin{equation*}
  \dInt[E]{h}{\nu} = \mu_{a}(E)
\end{equation*}
 for all $E\in{\cal A}$, in particular $\isEquiv{\mu_{a}}{\nu}{\absCont}$.

4.
Assume that we can find another pair $\mu_{a}'$ and $\mu_{s}'$ with $\isEquiv{\mu_{a}'}{\nu}{\absCont}$ and $\isEquiv{\mu_{s}'}{\nu}{\bot}$ and $\mu = \mu_{a}'+\mu_{s}'$. Then we have $\mu_{a}-\mu_{a}' = \mu_{s}'-\mu_{s}$ with $\isEquiv{\mu_{a}-\mu_{a}'}{\nu}{\absCont}$ and $\isEquiv{\mu_{s}'-\mu_{s}}{\nu}{\bot}$ by Lemma~\ref{elem-abs-cont-sing}, hence $\mu_{s}-\mu_{s}'=0$, again by Lemma~\ref{elem-abs-cont-sing}, which implies $\mu_{a}-\mu_{a}' = 0$. So the decomposition is unique. From this, uniqueness of the density $h$ in inferred. 
\EndProof

We obtain as a consequence the well-know Radon-Nikodym Theorem:

\BeginTheorem{radon-nikodym}
Let $\mu$ and $\nu$ be finite measures on $(X, {\cal A})$ with
$\isEquiv{\mu}{\nu}{\absCont}$. Then there exists a unique
$h\in\rLp[1]{\mu}$ with $\mu(A) = \dInt[A]{h}{\nu}$ for all $A\in {\cal
  A}$. Moreover, $f\in\rLp[1]{\mu}$ iff $f\cdot h\in\rLp[1]{\nu}$, in this case
\begin{equation*}
  \dInt{f}{\mu} = \dInt{f\cdot h}{\nu}.
\end{equation*}
$h$ is called the \emph{\index{Radon-Nikodym derivative}Radon-Nikodym derivative} of $\mu$ with
respect to $\nu$ and sometimes denoted by $d\mu/d\nu$. 
\EndTheorem

\BeginProof
Write $m = \mu_{a}+\mu_{s}$, where $\mu_{a}$ and $\mu_{s}$ are the
Lebesgue decomposition of $\mu$ with respect to $\nu$ by Theorem~\ref{lebesgue-radon-nikodym}. Since
$\isEquiv{\mu_{s}}{\nu}{\bot}$, we find $\mu_{s}=0$, so that
$\mu_{a}=\mu$. Then apply the second part of
Theorem~\ref{lebesgue-radon-nikodym} to $\mu$. This accounts for the first part. The second part follows from this by an approximation through step functions according to Corollary~\ref{step-is-dense}. 
\EndProof

Note that the Radon-Nikodym Theorem gives a one-to-one correspondence
between finite measures $\mu$ such that $\isEquiv{\mu}{\nu}{\absCont}$
and the Banach space $\rLp[1]{\nu}$. 

Theorem~\ref{lebesgue-radon-nikodym} can be extended to complex
measures; we will comment on this after the Jordan Decomposition will
be established in Proposition~\ref{jordan-decomposition-exists}.

Both constructions have, as one might suspect, a plethora of
applications. We will not discuss the Lebesgue decomposition further
but rather focus on the Radon-Nikodym Theorem and discuss two
applications, viz., identifying the dual space of the $\rLpS$-spaces
for $p<\infty$, and disintegrating a measure on a product space.

Before we do this, we have a look at integration by substitution, a
technique well-known from Calculus. The multi-dimensional case has
been hinted at on page~\pageref{change of variables_calculus}, we deal
here with the one-dimensional case. The approach displays a pretty
interplay of integrating with respect to an image measure, and the
Radon-Nikodym Theorem, which should not be missed.

We prepare the stage with an auxiliary statement, which is of interest of its own.

\BeginLemma{pre-for-int-by-subst}
Let $(X, {\cal A}, \mu)$ and $(Y, {\cal B}, \rho)$ be finite measure spaces, $\phi: X\to Y$ be measurable and onto such that $\rho_{*}(\Bild{\phi}{A}) = 0$, whenever $\mu(A) = 0$. Put $\nu := \FinM{\phi}(\mu)$. Then there exists a measurable function $w: X\to \pReal$ such that 
\begin{enumerate}
\item $f\in \rLp[1]{\rho}$ iff $(f\circ g)\cdot w\in\rLp[1]{\mu}$.
\item $\dInt[Y]{f(y)}{\rho(y)} = \dInt{(f\circ \phi)(x)\cdot g(x)}{\mu(x)}$ for all $f\in \rLp[1]{\rho}$. 
\end{enumerate}
\EndLemma

\BeginProof
We show first that $\rho\absCont\nu$, from which we obtain a
derivative. This is used then through the change of variable formula
for obtaining the desired result. 

In fact, assume that $\nu(B) = 0$ for some $B\in{\cal B}$,
equivalently, $\mu(\InvBild{\phi}{B}) = 0$. Thus by assumption
$0 = \rho_{*}(\Bild{\phi}{\InvBild{\phi}{B}}) = \rho(B)$, since $B =
\Bild{\phi}{\InvBild{\phi}{B}}$ due to $\phi$ being onto. Thus we find $g_{1}: Y\to \pReal$ such that 
$f\in \rLp[1]{\rho}$ iff $f\cdot g_{1}\in\rLp[1]{\nu}$ 
and
$
\dInt[Y]{f}{\rho} = \dInt[Y]{f\cdot g_{1}}{\nu}.
$
Since $\nu=\FinM{\phi}(\mu)$, we obtain from Corollary~\ref{int-image-measure} that
$
\dInt[Y]{f}{\rho} = \dInt{(f\circ \phi)\cdot (g_{1}\circ \phi)}{\mu}.
$
Thus putting $g := g_{1}\circ \phi$, the assertion follows. 
\EndProof

The r\^ole of $\nu$ as the image measure is interesting here. It just
serves as a kind of facilitator, but it remains in the
background. Only the measures $\rho$ and $\mu$ are acting, the image
measure is used only for obtaining the Radon-Nikodym derivative, and
for converting its integral to an integral with respect to its
preimage through change of variables.

We specialize things now to intervals on the real line and make
restrictive assumptions on $\phi$. Then ---~voil\`{a}!~--- the well known
formula on integration by substitution will result.

But first a more general consequence of
Lemma~\ref{pre-for-int-by-subst} is to be presented.  We will be working with Lebesgue
measure on intervals of the reals. Here we assume that $\phi: [\alpha,
\beta]\to [a, b]$ is continuous with the additional property that
$\lambda(A) = 0$ implies $\lambda_{*}(\Bild{\phi}{A}) = 0$ for all
$A\subseteq \Borel{[\alpha, \beta]}$. This class of functions is
generally known as \emph{absolutely continuous} and discussed in great
detail in~\cite[Section 18, Theorem (18.25)]{Hewitt-Stromberg}. We
obtain from Lemma~\ref{pre-for-int-by-subst}

\BeginCorollary{int-by-subst-prelim}
Let $[\alpha, \beta]\subseteq \Real$ be a closed interval, $\phi:
[\alpha, \beta]\to [a, b]$ be a surjective and absolutely continuous
function. Then there exists a Borel measurable function $w: [\alpha,
\beta]\to \Real$ such that
\begin{enumerate}
\item $f\in\rLp[1]{[a, b], \lambda}$ iff $(f\circ \phi)\cdot w\in\rLp[1]{[\alpha, \beta],\lambda}$ 
\item $\int_{a}^{b}f(x)\ dx = \int_{\alpha}^{\beta}(f(\phi(t))\cdot w(t)\ dt$.
\end{enumerate}
\EndCorollary

\BeginProof
The assertion follows from Lemma~\ref{pre-for-int-by-subst} by specializing $\mu$ and $\rho$ to $\lambda$. 
\EndProof

If we restrict $\phi$ further, we obtain even more specific
informations about the function $w$. The following proof shows how we
exploit the properties of $\phi$, viz., being monotone and having a
continuous first derivative, through the definition of the integral as
a limit of approximations on a system on subintervals which get
smaller and smaller. The subdivisions in the domain are then related to the one in the range of $\phi$, the relationship is done through Lagrange's Theorem which brings in the derivative. But see for yourself: 
  
\BeginProposition{integration-by-substitution}
Assume that $\phi: [\alpha, \beta]\to [a, b]$ is continuous and monotone
with a continuous first derivative such that $\phi(\alpha) = a$ and
$\phi(\beta) = b$. Then $f$ is Lebesgue integrable over $[a, b]$ iff
$(f\circ \phi)\cdot \phi'$ is Lebesgue integrable over $[\alpha,
\beta]$, and
\begin{equation*}
  \int_{a}^{b}f(x)\ dx = \int_{\alpha}^{\beta}f(\phi(z))\cdot \phi'(z)\ dz
\end{equation*}
holds. 
\EndProposition

We follow~\cite[Nr. 316]{Fichtenholz} in his proof. The basic idea is
to approximate the integral through step functions, which are obtained
by subdividing the interval $[\alpha, \beta]$ into sub intervals, and
to refine the subdivisions, using uniform continuity both of $\phi$
and $\phi'$ on its compact domain. So this is a fairly classical
proof.

\BeginProof
0.  We may assume that $f \geq 0$, otherwise we decompose $f = f^{+} -
f^{-}$ with $f^{+}, f^{-} \geq 0$. Also we assume that $f$ is bounded
by some constant $L$, otherwise we establish the property for $f\wedge
n$ with $n\in\Nat$, and, letting $n\to \infty$, appeal to Levi's
Theorem~\ref{beppo-levi}. Moreover we assume that $\phi$ is
increasing.

1.  The interval $[\alpha, \beta]$ is subdivided through $\alpha =
z_{0} < z_{1} < \dots < z_{n} = \beta$; put $x_{i} := \phi(z_{i})$,
then $a = x_{0} \leq x_{1} \leq \dots \leq x_{n} = b$, and $\Delta
z_{i} := z_{i+1} - z_{i}$, and $\Delta x_{i} := x_{i+1} -
x_{i}$. Let $\ell := \max_{i=1, \dots, n-1}\Delta z_{i}$, then if
$\ell\to 0$, the maximal difference $\max_{i=1, \dots,
  n-1}\Delta x_{i}$ tends to $0$ as well, because $\phi$ is
uniformly continuous. This is so since the interval $[\alpha, \beta]$
is compact.

For approximating the integral
$\int_{\alpha}^{\beta}f(\phi(z))\cdot \phi'(z)\ dz$ we select
$\zeta_{i}$ from each interval $[z_{i}, z_{i+1}]$ and write
\begin{equation*}
  S := \sum_{i}f(\phi(\zeta_{i}))\cdot \phi'(\zeta_{i})\cdot \Delta z_{i}.
\end{equation*}
Put $\xi_{i} := \phi(\zeta_{i})$, hence $x_{i} \leq \xi_{i}\leq
x_{i+1}$. By Lagrange's Formula\footnote{Recall that \emph{\index{Lagrange's Formula}Lagrange's Formula} says the following: Assume that $g$ is continuous on the interval $[c, d]$ with a continuous derivative $g'$ on the open interval $]c, d[$. Then there exists $t\in]c, d[$ such that $g(d)-g(c) = g'(t)\cdot (d-c)$.} there exists $\tau_{i}\in[z_{i},
z_{i+1}]$ such that $\Delta x_{i} = \phi'(\tau_{i})\cdot \Delta z_{i}$,
so that we can write as an approximation to the integral
$\int_{a}^{b}f(x)\ dx$ the sum
\begin{align*}
  s  := & \sum_{i}f(\xi_{i})\cdot \Delta z_{i}\\
 = &\sum_{i}f(\xi_{i}) \cdot \phi(\tau_{i})\cdot \Delta z_{i}\\
 = &\sum_{i}f(\phi(\zeta_{i})) \cdot \phi'(\tau_{i})\cdot \Delta z_{i}.
\end{align*}
If $\ell\to 0$, we know that $s\to \int_{a}^{b}f(x)\ dx$ and $S \to
\int_{\alpha}^{\beta}f(\phi(z))\cdot \phi'(z)\ dz$, so that we have to
get a handle at the difference $|S-s|$. We claim that this difference
tends to zero, as $\ell\to 0$. Given $\epsilon>0$, we find $\delta>0$
such that $|\phi'(\zeta_{i}) - \phi'(\tau_{i})| < \epsilon$, provided
$\ell < \delta$. This is so because $\phi'$ is continuous, hence
uniformly continuous. But then we obtain by telescoping
\begin{equation*}
  |S-s| \leq \sum_{i}|f(\phi(\zeta_{i}))|cdot |\phi'(\zeta_{i}) - \phi'(\tau_{i})|\cdot \Delta z_{i} < L\cdot (\beta-\alpha)\cdot \epsilon. 
\end{equation*}
Thus the difference vanishes, and we obtain indeed the equality claimed above. 
\EndProof

%%% Local Variables: 
%%% mode: latex
%%% TeX-master: "../Mskr3"
%%% End: 

\Subsubsection{Continuous Linear Functionals on $\rLpS$}
\label{sec:cont-lin-fnct-lp}

After all these preparations, we will investigate now continuous
linear functionals on the $\rLpS$-spaces and show that the map
$f\mapsto \dInt{f}{\mu}$ plays an important r\^ole in identifying
them. For full generality with respect to the functional concerned we
introduce signed measures here and show that they may be obtained in a
fairly specific way from the (unsigned) measures considered so far.

But before entering into this discussion, some general remarks. If $V$
is a real vector space with a norm $\|\cdot \|$, then a map $\Lambda:
V\to \Real$ is a \emph{\index{linear functional}linear functional} on
$V$ iff it is compatible with the vector space structure, i.e., iff
$\Lambda(\alpha\cdot x+\beta\cdot y) = \alpha\cdot \Lambda(x) +
\beta\cdot \Lambda(y)$ holds for all $x, y\in V$ and all $\alpha,
\beta\in\Real$. If $\Lambda\not= 0$, the range of $\Lambda$ is
unbounded, so $\sup_{x\in V}|\Lambda(x)|=\infty$. Consequently it is
difficult to assign to $\Lambda$ something like the $\sup$-norm for
characterizing continuity. It turns out, however, that we may
investigate continuity through the behavior of $\Lambda$ on the unit
ball of $V$, so we define
\begin{equation*}
  \|\Lambda\| := \sup_{\|x\|\leq1}|\Lambda(x)|
\end{equation*}
Call $\Lambda$ \emph{\index{bounded}bounded} iff $\|\Lambda\|\leq\infty$. Then $\Lambda$ is continuous iff $\Lambda$ is bounded, see Exercise~\ref{ex-linear-cont-is-bounded}. 

Now let $\mu$ be a finite measure with $p$ and $q$ conjugate to each other. Define for $g\in\rLp[q]{\mu}$ the linear functional
\begin{equation*}
  \Lambda_{g}(f) := \dInt{f\cdot g}{\mu}
\end{equation*}
on $\rLp{\mu}$, then we know from Hölder's inequality in Proposition~\ref{hoelder-inequality} that
\begin{equation*}
  \|\Lambda_{g}\| \leq \sup_{\|f\|_{p}\leq 1} \dInt{|f\cdot g|}{\mu}  \leq \|g\|_{q}
\end{equation*}
That was easy. But what about the converse? Given a bounded linear functional $\Lambda$ on $\rLp{\mu}$, does there exists $g\in\rLp[q]{\mu}$ with $\Lambda = \Lambda_{g}$? It is immediate that this will not work in general, since $\Lambda_{g}(f)\geq 0$, provided $f\geq 0$, so we have to assume that $\Lambda$ maps positive functions to a non-negative value. Call $\Lambda$ \emph{\index{linear functional!positive}positive} iff this is the case.

Summarizing, we consider maps $\Lambda: \cLp{\mu}\to \Real$ with these properties:
\begin{description}
\item[Linearity:] $\Lambda(\alpha\cdot x+\beta\cdot y) = \alpha\cdot \Lambda(x) + \beta\cdot \Lambda(y)$ holds for all $x, y\in V$ and all $\alpha, \beta\in\Real$.
\item[Boundedness:] $\|\Lambda\| := \sup_{\|f\|_{p}=1}|\Lambda(f)|
  \leq \infty$ (hence $|\Lambda(f)| \leq \|\Lambda\|\cdot \|f\|_{p}$
  for all $f$).
\item[Positiveness:] $f\geq 0 \Rightarrow \Lambda(f)\geq 0$ (note that
  $f\geq 0$ means $f'\geq0$-almost everywhere with respect to $\mu$ for each representative $f'$ of $f$ by our convention). 
\end{description}

We will first work on this restricted problem, and then we will expand the answer. This will require a slight generalization: we will talk about signed measures rather than about measures. 

Let's jump right in:

\BeginTheorem{lq-dual-to-lp}
Assume that $\mu$ is a finite measure on $(X, {\cal A})$, $1\leq p <\infty$, and that $\Lambda$ is a bounded positive linear functional on $\rLp{\mu}$. Then there exists a unique $g\in\rLp[q]{\mu}$ such that 
\begin{equation*}
  \Lambda(f) = \dInt{f\cdot g}{\mu}
\end{equation*}
holds for each $f\in\rLp{\mu}$. In addition, $\|\Lambda\| = \|g\|_{q}$.
\EndTheorem

This is our line of attack: We will first see that we obtain from
$\Lambda$ a finite measure $\nu$ on ${\cal A}$ such that
$\isEquiv{\nu}{\mu}{\absCont}$. The Radon-Nikodym Theorem will then
give us a density $g := d\nu/d\mu$ which will turn out to be the
function we are looking for. This is shown by separating the cases
$p=1$ and $p>1$.

\BeginProof
1.
Define for $A\in{\cal A}$ 
\begin{equation*}
  \nu(A) := \Lambda(\chi_A).
\end{equation*}
Then $A\subseteq B$ implies $\chi_{A}\leq \chi_{B}$, hence
$\Lambda(\chi_{A})\leq \Lambda(\chi_{B})$. Because $\Lambda$ is
monotone, hence $\nu$ is monotone. Since $\Lambda$ is linear, we have
$\nu(\emptyset) = 0$, and $\nu$ is additive. Let $\Folge{A}$ be an
increasing sequence of measurable sets with $A :=
\bigcup_{n\in\Nat}A_{n}$, then $\chi_{A\setminus A_{n}}\to 0$, and thus
\begin{equation*}
  \nu(A)-\nu(A_{n}) = \|\chi_{A\setminus A_{n}}\|_{p}^{p} =
  \Lambda(\chi_{A\setminus A_{n}})^{p}\to 0,
\end{equation*}
since $\Lambda$ is continuous. Thus $\Lambda$ is a finite measure on
${\cal A}$ (note $\nu(X) = \Lambda(1)<\infty$). If $\mu(A) = 0$, we see that
$\chi_{A}=_{\mu}0$, thus
$\Lambda(\chi_{A}) = 0$ (we are dealing with the $=_{\mu}$-class of
$\chi_{A}$), so that $\nu(A) = 0$. Thus
$\isEquiv{\nu}{\mu}{\absCont}$, and the Radon-Nikodym
Theorem~\ref{radon-nikodym} tells
us that there exists $g\in\rLp[1]{\mu}$ with 
\begin{equation*}
  \Lambda(\chi_{A}) = \nu(A) = \dInt[A]{g}{\mu}
\end{equation*}
for all $A\in {\cal A}$. Since the integral as well as $\Lambda$ are
linear, we obtain from this
\begin{equation*}
  \Lambda(f) = \dInt{f\cdot g}{\mu}
\end{equation*}
for all step functions $f$. 

2. 
We have to show that $g\in\rLp[q]{\mu}$. Consider these cases.
\begin{description}
\item[Case $p=1$:] We have for each $A\in{\cal A}$ 
\begin{equation*}
  \bigl|\dInt[A]{g}{\mu}\bigr| \leq |\Lambda(\chi_{A})| \leq
  \|\Lambda\|\cdot \|\chi_{A}\|_{1} = \|\Lambda\|\cdot \mu(A)
\end{equation*}
But this implies $|g(x)| \leq_{\mu} \|\Lambda\|$, thus $\|g\|_{\infty}\leq
\|\Lambda\|$.
\item[Case $1 < p < \infty$:] Let $t = \chi_{\{x\in X\mid g(x) \geq
    0\}}-\chi_{\{x\in X\mid g(x) < 0\}}$, then $|g| = t\cdot g$, and
  $t$ is measurable, since $g$ is. Define $A_{n} := \{x\in X \mid
  |g(x)| \leq n\}$, and put $f := \chi_{A_{n}}\cdot |g|^{q-1}\cdot
  t$. Then
  \begin{align*}
    |f|^{p}\cdot \chi_{A_{n}} & = |g|^{(q-1)\cdot p}\cdot \chi_{A_{n}}\\ & =
    |g|^{q}\cdot \chi_{A_{n}},\\
\chi_{A_{n}}\cdot (f\cdot g) & = \chi_{A_{n}}\cdot |g|^{q-1}\cdot t\cdot
g\\ & = \chi_{A_{n}}\cdot |g|^{q}\cdot t,
  \end{align*}
thus
\begin{equation*}
  \dInt[A_{n}]{|g|^{q}}{\mu} = \dInt[A_{n}]{f\cdot g}{\mu} = \Lambda(f)
  \leq \|\Lambda\|\cdot \bigl(\dInt[A_{n}]{|g|^{q}}{\mu}\bigr)^{1/p}
\end{equation*}
Since $1 - 1/p = 1/q$, dividing by the factor $\|\Lambda\|$ and
raising the result by $q$ yields
\begin{equation*}
  \dInt[E_{n}]{|g|^{q}}{\mu}\leq \|\Lambda\|^{q}.
\end{equation*}
By Lebesgue's Dominated Convergence
Theorem~\ref{lebesgue-dominated-convergence} we obtain that
$\|g\|_{q}\leq \|\Lambda\|$ holds, hence $g\in\rLp[q]{\mu}$, and
$\|g\|_{q}= \|\Lambda\|$. 
\end{description}
The proof is completed now by the observation that $\Lambda(f) =
\dInt{f\cdot g}{\mu}$ holds for all step functions $f$. Since both
sides of this equation represent continuous functions, and since the
step functions are dense in $\rLp{\mu}$ by
Corollary~\ref{step-is-dense}, the equality holds on all of $\rLp{\mu}$. 
\EndProof

This representation holds only for positive linear functions; what
about the rest? It turns out that we need to extend our notion of
measures to signed measures, and that a very similar statement
holds for signed measures (of course we would have to explain what the integral of a
signed measure is, but this will work out very smoothly). So what we will
do next is to define signed measures, and to relate them to the
measures with which we have worked until now. We follow essentially Halmos'
exposition~\cite[§ 29]{Halmos}.

\BeginDefinition{signed-measure}
A map $\mu: {\cal A}\to \Real$ is said to be a \emph{\index{signed measure}signed measure} iff
$\mu$ is $\sigma$-additive, i.e., iff $\mu(\bigcup_{n\in\Nat}A_{n}) =
\sum_{n\in\Nat}\mu(A_{n})$, whenever $\Folge{A}$ is a sequence of
mutually disjoint sets in ${\cal A}$. 
\EndDefinition

Clearly, $\mu(\emptyset) = 0$, since a signed measure $\mu$ is finite, so the
distinguishing feature is the absence of monotonicity. 
It turns out, however, that we can partition the whole space $X$ into
a positive and a negative part, that restricting $\mu$ to these parts
will yield a measure each, and that $\mu$ can be written in this way
as the difference of two measures. 

Fix a signed measure $\mu$. Call $N\in{\cal A}$ a \emph{\index{signed measure!positive set}negative set} iff
$\mu(A\cap N) \leq 0$ for all $A\in{\cal A}$; a positive set is
defined accordingly.  It is immediate that the difference of two
negative sets is a negative set again, and that the union of a
disjoint sequence of negative sets is a negative set as well. Thus the
union of a sequence of negative sets is negative again. 

\BeginProposition{jordan-decomposition-exists}
Let $\mu$ be a signed measure on ${\cal A}$. Then there exists a pair
$X^{+}$ and $X^{-}$ of disjoint measurable sets such that $X^{+}$ is a
positive set, $X^{-}$ is a negative set. Then $\mu^{+}(B) := \mu(B\cap
X^{+})$ and $\mu^{-}(B) := -\mu(B\cap X^{-})$ are finite measures on
${\cal A}$ such that $\mu = \mu^{+}-\mu^{-}$. The pair $\mu^{+}$ and
$\mu^{-}$  is called the
\emph{\index{Jordan decomposition}Jordan Decomposition} of the signed
measure $\mu$.  
\EndProposition

\BeginProof
1.
Define 
\begin{equation*}
\alpha :=\inf\{\mu(A) \mid  A\in{\cal A} \text{ is
  negative}\}>-\infty.
\end{equation*}
Assume that $\Folge{A}$ is a sequence of
measurable sets with $\mu(A_{n})\to \alpha$, then we know that $A :=
\bigcup_{n\in\Nat}A_{n}$ is negative again with $\alpha = \mu(A)$. In
fact, put $B_{1} := A_{1}$, $B_{n+1} := A_{n+1}\setminus B_{n}$, then
each $B_{n}$ is negative, we have 
\begin{equation*}
\mu(A) =
\mu(\bigcup_{n\in\Nat}B_{n}) = \sum_{n\in\Nat}\mu(B_{n}) = \lim_{n\to
  \infty}\mu(A_{n})
\end{equation*} 
by telescoping. 

2.
We claim that  
\begin{equation*}
X^{+} := X\setminus A
\end{equation*}
is a positive set. In fact,
assume that this is not true ---~now this becomes the tricky part~--- then there exists $E_{0}\subseteq X^{+}$
with $\mu(E_{0})<0$. $E_{0}$ cannot be a negative set, because
otherwise $A\cup E_{0}$ would be a negative set with $\mu(A\cup E_{0})
= \mu(A)+\mu(E_{0}) < \alpha$, which is contradicts the construction
of $A$. Let $k_{1}$ be the smallest positive integer  such that
$E_{0}$ contains a measurable set $E_{1}$ with $\mu(E_{1})\geq
1/k_{1}$. Now look at $E_{0}\setminus E_{1}$. We have
\begin{equation*}
\mu(E_{0}\setminus E_{1}) = \mu(E_{0})-\mu(E_{1}) \leq
\mu(E_{0})-\mu(E_{1}) \leq \mu(E_{0}) - 1/k_{1}< 0.
\end{equation*}
We may repeat
the same consideration now for $E_{0}\setminus E_{1}$; let $k_{2}$ be
the smallest positive integer such that $E_{0}\setminus E_{1}$
contains a measurable set $E_{2}$ with $\mu(E_{2})\geq 1/k_{2}$. This
produces a sequence of disjoint measurable sets $\Folge{E}$ with
\begin{equation*}
E_{n+1}\subseteq E_{0}\setminus(E_{1}\cup\dots\cup E_{n}),
\end{equation*} 
and since
$\sum_{n\in\Nat}\mu(E_{n})$ is finite (because
$\bigcup_{n\in\Nat}E_{n}\in{\cal A}$, and $\mu$ takes only finite
values), we infer that $\lim_{n\to \infty} 1/k_{n} = 0$. 

3.
Let $F\subseteq F_{0} := E_{0}\setminus\bigcup_{n\in\Nat}E_{n}$, and assume
that $\mu(F) \geq 0$. Let $\ell$ be the largest positive integer with
$\mu(F)\geq 1/\ell$. Since $k_{n}\to 0$, as $n\to \infty$, we 
find $m\in\Nat$ with $1/\ell \geq 1/k_{m}$. Since $F\subseteq
E_{0}\setminus(E_{1}\cup\dots\cup E_{m})$, this yields a
contradiction. But $F_{0}$ is disjoint from $A$, and since 
\begin{equation*}
\mu(F_{0}) = \mu(E_{0}) - \sum_{n\in\Nat}\mu(E_{n}) \leq \mu(E_{0}) <
0,
\end{equation*}
we have arrived at a contradiction. Thus $\mu(E_{0})\geq 0$. 

4. Now define $\mu^{+}$ and $\mu^{-}$ as the traces of $\mu$ on
$X^{+}$ and $X^{-} := A$, resp., then the assertion follows. 
\EndProof

It should be noted that the decomposition of $X$ into $X^{+}$ and
$X^{-}$ is not unique, but the decomposition of $\mu$ into $\mu^{+}$
and $\mu^{-}$ is. Assume that $X_{1}^{+}$ with $X_{1}^{-}$ and
$X_{2}^{+}$ with $X_{2}^{-}$ are two such decompositions. Let
$A\in{\cal A}$, then we have $A\cap (X_{1}^{+}\setminus X_{2}^{+})
\subseteq A\cap X_{1}^{+}$, hence $\mu(A\cap(X_{1}^{+}\setminus
X_{2}^{+}) \geq 0$; on the other hand, $A\cap (X_{1}^{+}\setminus X_{2}^{+})
\subseteq A\cap X_{2}^{-}$, thus $\mu(A\cap(X_{1}^{+}\setminus
X_{2}^{+}) \leq 0$, so that we have $\mu(A\cap(X_{1}^{+}\setminus
X_{2}^{+}) = 0$, which implies $\mu(A\cap X_{1}^{+}) = \mu(A\cap
X_{2}^{+})$. Thus uniqueness of $\mu^{+}$ and $\mu^{-}$ follows. 
 
Given a signed measure $\mu$ with a Jordan decomposition $\mu^{+}$ and
$\mu^{-}$, we define a (positive) measure $|\mu| := \mu^{+}+\mu^{-}$;
$|\mu|$ is called the \emph{\index{signed measure!total
    variation}total variation} of $\mu$.  It is clear that $|\mu|$ is
a finite measure on ${\cal A}$. A set $A\in{\cal A}$ is called a
\emph{$\mu$-\index{signed measure!nullset}nullset} iff $\mu(B) = 0$
for every $B\in{\cal A}$ with $B\subseteq A$; hence $A$ is a
$\mu$-nullset iff $A$ is a $|\mu|$-nullset iff $|\mu|(A) = 0$. In this
way, we can define that a property holds $\mu$-everywhere also for
signed measures, viz., by saying that it holds $|\mu|$-everywhere (in
the traditional sense). Also the relation
$\isEquiv{\mu}{\nu}{\absCont}$ of \index{absolute continuity!signed measure}absolute continuity \index{signed
  measure!absolute continuity} between the signed measure $\mu$ and
the positive measure $\nu$ can be redefined as saying that each
$\nu$-nullset is a $\mu$-nullset. Thus $\isEquiv{\mu}{\nu}{\absCont}$
is equivalent to $\isEquiv{|\mu|}{\nu}{\absCont}$ and to both
$\isEquiv{\mu^{+}}{\nu}{\absCont}$ and
$\isEquiv{\mu^{-}}{\nu}{\absCont}$. For the derivatives, it is easy to
see that
\begin{equation*}
  \frac{d\mu}{d\nu}   = \frac{d\mu^{+}}{d\nu}-\frac{d\mu^{-}}{d\nu}
  \text{ and }
  \frac{d|\mu|}{d\nu}   = \frac{d\mu^{+}}{d\nu}+\frac{d\mu^{-}}{d\nu}
\end{equation*}
hold. 

We define integrability of a measurable function through $|\mu|$ by putting
\begin{equation*}
  \cLp{|\mu|} := \cLp{\mu^{+}}\cap\cLp{\mu^{-}},
\end{equation*}
and define $\rLp{\mu}$ again as the set of equivalence classes. 

These observations provide a convenient entry point into discussing
complex measures. Call $\mu: {\cal A}\to \mathbb{C}$ a (complex)
measure iff $\mu$ is $\sigma$-additive, i.e., iff
$\mu(\bigcup_{n\in\Nat}A_{n}) = \sum_{n\in\Nat}\mu(A_{n})$ for each
sequence $\Folge{A}$ of mutually disjoint sets in ${\cal A}$. Then it
can be easily shown that $\mu$ can be written as $\mu = \mu_{r}+i\cdot
\mu_{c}$ with (real) signed measures $\mu_{r}$ and $\mu_{c}$, which in
turn have a Jordan decomposition and consequently a total variation
each. In this way the $\rLpS$-spaces can be defined also for complex
measures and complex measurable functions; the reader is referred
to~\cite{Rudin} or~\cite{Hewitt-Stromberg} for further information.~---
 
Returning to the main current of the discussion, we are able to state the general representation of continuous
linear functionals on an $\rLp{\mu}$-space. We need only to sketch the
proof, mutatis mutandis, since the main work has already been done in
the proof of Theorem~\ref{lq-dual-to-lp}. 

\BeginTheorem{lq-dual-to-lp-general}
Assume that $\mu$ is a finite measure on $(X, {\cal A})$, $1\leq p
<\infty$, and that $\Lambda$ is a bounded linear functional
on $\rLp{\mu}$. Then there exists a unique $g\in\rLp[q]{\mu}$ such
that
\begin{equation*}
  \Lambda(f) = \dInt{f\cdot g}{\mu}
\end{equation*}
holds for each $f\in\rLp{\mu}$. In addition, $\|\Lambda\| = \|g\|_{q}$.
\EndTheorem

\BeginProof
$\nu(A) := \Lambda(\chi_{A})$ defines a signed measure on ${\cal A}$
with $\isEquiv{\nu}{\mu}{\absCont}$. Let $h$ be the Radon-Nikodym
derivative of $\nu$ with respect to $\mu$, then $h\in\rLp[q]{\mu}$ and 
\begin{equation*}
  \Lambda(f) = \dInt{f\cdot h}{\mu}
\end{equation*}
are shown as above.
\EndProof

It should be noted that Theorem~\ref{lq-dual-to-lp-general} holds also
for $\sigma$-finite measures, and that it is true for $1 < p < \infty$
in the case of general (positive) measures, see,
e.g.,\cite[§ VII.3]{Elstrodt} for a discussion. 

The case of continuous linear functionals for the space
$\rLp[\infty]{\mu}$ is considerably more
involved. Example~\ref{linf-not-separable} indicates already that
these spaces play a special r\^ole. Looking back at the discussion
above, we found that for $p<\infty$ the map $ A\mapsto
\dInt[A]{|f|^{p}}{\mu} $ yields a measure, and this measure was
instrumental through the Radon-Nikodym Theorem for providing the
factor which could be chosen to represent the linear functional. This
argument, however, is not available for the case $p=\infty$, since we
are not working there with a norm which is derived from an integral. It
can be shown, however, that continuous linear functional have an
integral representation with respect to finitely additive set
functions; in fact,\cite[Theorem 20.35]{Hewitt-Stromberg} or \cite[Theorem IV.8.16 ]{Dunford-Schwartz} show  that
the continuous linear functionals on $\rLp[\infty]{\mu}$ are in a
one-to-one correspondence with all finitely additive set functions
$\xi$ such that $\isEquiv{\xi}{\mu}{\absCont}$. Note that this
requires an extension of integration to not necessarily
$\sigma$-additive set functions.

\Subsubsection{Disintegration}
\label{sec:disintegration}
One encounters occasionally the situation the need to decompose a
measure on a product of two spaces. Consider this scenario. Given a
measurable space $(X, {\cal A})$ as an input , $(Y, {\cal B})$ as an
output space, let $(\mu\otimes K)(B) = \dInt{K(x)(D_{x})}{\mu(x)}$ be
the probability for $\langle x_{1}, x_{2}\in B\in{\cal A}\otimes{\cal
  B}$ with $\mu$ as the initial distribution and $K: (X, {\cal
  A})\Trans (Y, {\cal B})$ as the transition law (think of an epidemic
which is set off according to $\mu$ and propagates according to
$K$). Assume that you want to reverse the process: Given $F\in{\cal
  B}$, you put $\nu(F) := \SubProb{\pi_{Y}(\mu\otimes K)}(F) =
(\mu\otimes K)(X\times F)$, so this is the probability that your
process hits an element of $F$. Can you find a stochastic relation $L:
(Y, {\cal B})\Trans (X, {\cal A})$ such that $(\mu\otimes K)(B) =
\dInt{L(x)(B^{y})}{\nu(y)}$? The relation $L$ is the
\emph{\index{converse}converse} of $K$ given $\mu$. It is probably not
particularly important that the measure on the product has the shape
$\mu\otimes K$, so we state the problem in such a way that we are
given a measure on a product of two measurable spaces, and the
question is whether we can decompose it into the product of a
projection onto one space, and a stochastic relation between the spaces.

This problem is of course easiest dealt with
when one can deduce that the measure is the product of
measures on the coordinate spaces; probabilistically, this
would correspond to the distribution of two independent random
variables. But sometimes one is not so lucky, and there is
some hidden dependence, or one simply cannot assess the degree of
independence. Then one has to live with a somewhat weaker
result: in this case one can decompose the measure
into a measure on one component and a transition probability.
This will be made specific in the discussion to follow.

Because it will not cost substantially more attention, we will treat the question
a bit more generally. Let $(X, \mathcal{A})$, $(Y, \mathcal{B}),$ and $(Z, \mathcal{C})$
be measurable spaces, assume that $\mu \in \SubProb{X, \mathcal{A}}$, and
let $f: X \rightarrow Y$ and $g: X \rightarrow Z$ be measurable maps. Then
$\mu_{f} := \SubProb{f}(\mu)$ and $\mu_{g}:= \SubProb{g}(\mu)$
define subprobabilities on $(Y, \mathcal{B})$ resp.
$(Z, \mathcal{C})$. $\mu_f$ and $\mu_g$ can be interpreted as the probability distribution of $f$ resp.
$g$ under $\mu$.

We will show that we can represent the joint distribution as
\begin{equation*}
\mu(\{x \in X \mid f(x) \in B, g(x) \in C\})
= \dInt[B]{K(y)(C)}{\mu_f(y)},
\end{equation*}
where $K: (Y, {\cal B})\Trans (Z, {\cal C})$ is a stochastic
relation. This will require $Z$ to be a Polish space with $\mathcal{C}
= \Borel{Z}$.

Let us see how this corresponds to the initially stated problem.
Suppose $X := Y \times Z$ with $\mathcal{A} = \mathcal{B}\otimes \mathcal{C}$,
and let $f := \pi_Y$, $g := \pi_Z$, then
\begin{eqnarray*}
\mu_f(B) & =& \mu(B \times Z),\\
\mu_g(C) & =& \mu(Y \times Z),\\
\mu(B \times C) & = & \mu(\{x \in X \mid f(x) \in B, g(x) \in C\}).
\end{eqnarray*}
% and
% \begin{equation*}
% \mu(\{x \in X \mid f(x) \in B, g(x) \in C\}) = \mu(B \times C).
% \end{equation*}
Granted that we have established the decomposition, we can then
write
\begin{equation*}
\mu(B \times C) = \dInt[B]{K(y)(C)}{\mu_f(y)};
\end{equation*}
thus we have decomposed the probability on the product into a probability
on the first component, and, conditioned on the value the first component
may take, a probability on the second factor.

\BeginDefinition{reg-cond-prob}
Using the notation from above, $K$ is called a \emph{regular conditional
distribution of $g$ given $f$}\index{conditional distribution!regular} iff
\begin{equation*}
\mu(\{x \in X \mid f(x) \in B, g(x) \in C\})
= \int_B K(y)(C)\ \mu_f(dy)
\end{equation*}
holds for each $B \in \mathcal{B}, C \in \mathcal{C}$, where $K: (Y,
{\cal B})\Trans (C, {\cal C})$ is a stochastic relation on $(X,
\mathcal{A})$ and $(Z, \mathcal{C})$.  If only $y\mapsto
K(y)(C)$ is $\mathcal{B}$-measurable for all $C \in \mathcal{C}$, then
it will be called a \emph{conditional distribution of $g$ given
  $f$}.\index{conditional distribution}
\EndDefinition

The existence of regular conditional distribution will be established,
provided $Z$ is Polish with $\mathcal{C} = \Borel{Z}$. This will be
accomplished in several steps: first the existence of a conditional
distribution will be shown using the well known Radon-Nikodym
Theorem. The latter construction will then be scrutinized. It will
turn out that there exists a set of measure zero outside of which the
conditional distribution behaves like a regular one, but at first
sight only on an algebra of sets, not on the entire
$\sigma$-algebra. But don't worry, the second step will apply a
classical extension argument and yield a regular conditional
distribution on the Borel sets, just as we want it.  The proofs are
actually a kind of a round trip through the first principles of
measure theory, where the Radon-Nikodym Theorem together with the
classical Hahn Extension Theorem are the main vehicles. It displays
also some nice and helpful proof techniques.

We fix  $(X, \mathcal{A})$, $(Y, \mathcal{B}),$ and $(Z, \mathcal{C})$
as measurable spaces, assume that $\mu \in \SubProb{X, \mathcal{A}}$, and
take $f: X \rightarrow Y$ and $g: X \rightarrow Z$ to be measurable maps.
The measures $\mu_f := \SubProb{f}(\mu)$ and $\mu_g := \SubProb{g}(\mu)$ are defined as above as the distribution
of $f$ resp. $g$ under $\mu$.

The existence of a conditional distribution of $g$ given $f$ is established first,
and it is shown that it is essentially unique.

\BeginLemma{simple-reg-cond-distr}
Using the notation from above, then
\begin{enumerate}
  \item there exists a conditional distribution $K_0$ of $g$ given $f$,
  \item if there
   is another conditional distribution $K_0'$ of $g$ given $f$, then
   there exists for any $C \in \mathcal{C}$
   a set $N_C \in \mathcal{B}$ with $\mu_f(N_C) = 0$ such that
   $
   K_0(y)(C) = K_0'(C)
   $
   for all $y \notin C$.
\end{enumerate}
\EndLemma

\BeginProof
1.
Fix $C \in \mathcal{C}$, then
\begin{equation*}
\varpi_C(B) := \mu(\InvBild{f}{B} \cap \InvBild{g}{C})
\end{equation*}
defines a subprobability measure $\varpi_C$ on $\mathcal{B}$ which
is absolutely continuous with respect to $\mu_g$, because
$\mu_g(B) = 0$ implies $\varpi_C(B) = 0$. The
Radon-Nikodym Theorem~\ref{radon-nikodym}
now gives a density $h_C \in \MeasbFnct{Y, \mathcal{B}}$ with
\begin{equation*}
\varpi_C(B) = \dInt[B]{h_C}{\mu_f}
\end{equation*}
for all $B \in \mathcal{B}$. Setting
$
K_0(y)(C) := h_C(y)
$
yields the desired conditional distribution.

2.
Suppose $K_0'$ is another conditional distribution of $g$
given $f$, then we have
\begin{equation*}
\forall B \in \mathcal{B}: \dInt[B]{K_0(y)(C)}{\mu_f(y)} = \dInt[B]{K_0(y)(C)}{\mu_f(y)},
\end{equation*}
for all $C \in \mathcal{C}$, which implies that the set on which $K_0(\cdot)(C)$
disagrees with $K_0'(\cdot)(C)$ is $\mu_f$-null.
\EndProof

Essential uniqueness may strengthened if the $\sigma$-algebra $\mathcal{C}$
is countably generated, and if the conditional distribution is regular.

\BeginLemma{uniq-count-gen}
Assume that $K$ and $K'$ are regular conditional distributions of $g$
given $f$, and that $\mathcal{C}$ has a countable generator. Then there exists
a set $N \in \mathcal{B}$ with $\mu_f(N) = 0$ such that
$
K(y)(C) = K'(y)(C)
$
for all $C \in \mathcal{C}$ and all $y \notin N$.
\EndLemma

\BeginProof
If $\mathcal{C}_0$ is a countable generator of $\mathcal{C}$, then
\begin{equation*}
\mathcal{C}_f := \{\bigcap \mathcal{\mathcal{E}}\mid
\mathcal{E} \subseteq \mathcal{C}_0\text{ is finite}\}
\end{equation*}
is a countable generator of $\mathcal{C}$ well, and $\mathcal{C}_f$ is closed under
finite intersections; note that $Z \in \mathcal{C}_f$.
Construct for $D \in \mathcal{C}_f$ the set $N_D \in \mathcal{B}$ outside of which
$K(\cdot)(D)$ and $K'(\cdot)(D)$ coincide, and define
\begin{equation*}
N := \bigcup_{D \in \mathcal{C}_f} N_D \in \mathcal{B}.
\end{equation*}
Evidently, $\mu_f(N) = 0.$
We claim that $K(y)(C) = K'(y)(C)$ holds for all $C \in \mathcal{C}$, whenever $y \notin N$.
In fact, fix $y \notin N$, and let
\begin{equation*}
\mathcal{C}_1 := \{C \in \mathcal{C} \mid K(y)(C) = K'(y)(C)\},
\end{equation*}
then $\mathcal{C}_1$ contains $\mathcal{C}_f$ by construction,
and is a $\pi$-$\lambda$-system.  This is so since it is
closed under complements and countable disjoint unions. Thus
$
\mathcal{C} = \sigma(\mathcal{C}_f) \subseteq \mathcal{C}_1,
$
by the $\pi$-$\lambda$-Theorem~\ref{Pi-Lambda},
and we are done.
\EndProof

We will show now that a regular conditional distribution of $g$ given
$f$ exists. This will be done through several steps, given the
construction of a conditional distribution $K_0$:
\begin{dingautolist}{192}
  \item \label{step-disint-A} 
A set $N_a \in \mathcal{B}$ is constructed with $\mu_f(N_a) = 0$
  such that $K_0(y)$ is additive on a countable generator $\mathcal{C}_z$ for $\mathcal{C}$.
  \item \label{step-disint-B} 
We construct a set $N_z \in B$ with $\mu_f(N_z) = 0$
  such that $K_0(y)(Z) \leq 1$ for $y \notin N_z$.
  \item \label{step-disint-C} 
For each element $G$ of $\mathcal{C}_z$ we will find a
  set $N_G \in \mathcal{B}$ with $\mu_f(N_G) = 0$ such that
  $K_0(y)(G)$ can be approximated from inside through compact
  sets, whenever $y \notin N_G$.
  \item \label{step-disint-D} 
Then we will combine all these sets of $\mu_f$-measure zero
  to produce a set $N \in \mathcal{B}$ with $\mu_f(N) = 0$ outside of which
  $K_0(y)$ is a premeasure on the generator $\mathcal{C}_z$, hence can be extended to a measure
  on all of $C$.
\end{dingautolist}
Well, this looks like a full program, so let us get on with it.

\BeginTheorem{reg-cond-distr-exists}
Given measurable spaces $(X, \mathcal{A})$ and $(Y, \mathcal{B})$, a Polish space $Z$,
a subprobability $\mu \in \SubProb{X, \mathcal{A}}$, and measurable maps
$f: X \rightarrow Y$, $g: X \rightarrow Z$, there exists a regular
conditional distribution $K$ of $g$ given $f$. $K$ is uniquely determined
up to a set of $\mu_f$-measure zero.
\EndTheorem

\BeginProof
0.
Since $Z$ is a Polish space, its topology has a countable base. We infer from
Lemma~\ref{DenseWithBorel} that $\Borel{Z}$ has a countable generator
$\mathcal{C}$. Then the Boolean algebra $\mathcal{C}_1$ generated by $\mathcal{C}$ is
also a countable generator of $\Borel{Z}$.

1.
Given $C_n \in \mathcal{C}_1$, we find by
Proposition~\ref{tight-on-polish} a sequence
$(E_{n, k})_{k \in \Nat}$ of compact sets in $Z$ with
\begin{equation*}
E_{n, 1} \subseteq E_{n, 2}\subseteq E_{n, 3} \ldots \subseteq C_n
\end{equation*}
such that
\begin{equation*}
\mu_g(C_n) = \sup_{k \in \Nat} \mu_g(E_{n, k}).
\end{equation*}
Then the Boolean algebra $\mathcal{C}_z$ generated by
$
\mathcal{C} \cup \{E_{n, k} \mid n, k \in \Nat\}
$
is also a countable generator of $\Borel{Z}$.

2.
From the construction of the conditional distribution of $g$ given $f$ we
infer that for disjoint $C_1, C_2 \in \mathcal{C}_z$
\begin{eqnarray*}
\dInt[Y]{K_0(y)(C_1 \cup C_2)}{\mu_f(y)}
& = &
\mu(\{x \in X \mid f(x) \in B, g(x) \in C_1 \cup C_2\})\\
& = & \mu(\{x \in X \mid f(x) \in B, g(x) \in C_1\}) +\\
&&\phantom{\mu(\{x \in X \mid}\mu(\{x \in X \mid f(x) \in B, g(x) \in C_2\})\\
& = &
\dInt[Y]{K_0(y)(C_1)}{\mu_f(y)} + \dInt[Y]{K_0(y)(C_2)}{\mu_f(y)}.
\end{eqnarray*}
Thus there exists $N_{C_1, C_2} \in \mathcal{B}$ with
$
\mu_f(N_{C_1, C_2}) = 0
$
such that
\begin{equation*}
K_0(y)(C_1 \cup C_2) = K_0(y)(C_1) + K_0(y)(C_2)
\end{equation*}
for $y \notin N_{C_1, C_2}$. Because $\mathcal{C}_z$ is countable,
we may deduce (by taking the union of $N_{C_1, C_2}$ over all pairs $C_1, C_2$)
that there exists a set $N_a \in \mathcal{B}$ such that $K_0$ is additive outside
$N_a$, and $\mu_f(N_a) = 0$. This accounts for part~\ref{step-disint-A}
in the plan above.

3.
By the previous arguments it is easy to construct
a set $N_z \in \mathcal{B}$ with $\mu_f(N_z) = 0$
such that $K_0(y)(Z) \leq 1$ for $y \notin N_z$ (part~\ref{step-disint-B}).

4.
Because
\begin{align*}
  \dInt[Y]{K_0(y)(C_n)}{\mu_f(y)} &= \mu(\InvBild{f}{Y} \cap \InvBild{g}{C_n}) \\
    &=  \mu_g(C_n) \\
    &=  \sup_{k \in \Nat} \mu_g(E_{n, k})\\
    & =   \sup_{k \in \Nat} \int_Y K_0(y)(E_{n, k})\ \mu_f(dy)&&
    \text{ Levi's Theorem~\ref{beppo-levi}}\\
    & =   \dInt[Y]{\sup_{k \in \Nat}  K_0(y)(E_{n, k})}{\mu_f(y)}
\end{align*}
we find for each $n \in \Nat$ a set $N_n \in \mathcal{B}$ with
\begin{equation*}
\forall y \notin N_n: K_0(y)(C_n) = \sup_{k \in \Nat} K_0(y)(E_{n, k})
\end{equation*}
and $\mu_f(N_n) = 0$. This accounts for part~\ref{step-disint-C}.

5.
Now we may begin to work on part~\ref{step-disint-D}. Put
\begin{equation*}
N := N_a \cup N_z \cup \bigcup_{n \in \Nat} N_n,
\end{equation*}
then $N \in B$ with $\mu_f(N) = 0$. We claim that $K_0(y)$ is a
premeasure on $\mathcal{C}_z$ for each $y \notin N$. It is clear that
$K_0(y)$ is additive on $\mathcal{C}_z$, hence monotone, so
merely $\sigma$-additivity has
to be demonstrated: let $(D_\ell)_{\ell \in \Nat}$ be a sequence
in $\mathcal{C}_z$ that is monotonically decreasing with
\begin{equation*}
\eta := \inf_{\ell \in \Nat} K_0(y)(D_\ell) > 0,
\end{equation*}
then we have to show that
\begin{equation*}
\bigcap_{\ell \in \Nat} D_\ell \not= \emptyset.
\end{equation*}
We approximate the sets $D_\ell$ now by compact sets,
so we assume that $D_\ell = C_{n_\ell}$ for some $n_\ell$ (otherwise
the sets are compact themselves). By construction we find for
each $\ell \in \Nat$ a compact set
$
E_{n_\ell, k_\ell}\subseteq C_\ell
$
with
\begin{equation*}
K_0(y)(C_{n_\ell}\setminus E_{n_\ell, k_\ell}) < \eta\cdot2^{\ell +1},
\end{equation*}
then
\begin{equation*}
E_r := \bigcap_{i=\ell}^r E_{n_\ell, k_\ell} \subseteq C_{n_r} = D_r
\end{equation*}
defines a decreasing sequence of compact sets with
\begin{equation*}
K_0(y)(E_r) \geq K_0(y)(C_{n_r}) - \sum_{i=\ell}^r K_0(y)(E_{n_\ell, k_\ell}) > \eta/2,
\end{equation*}
thus $E_r \not= \emptyset$.
Since $E_r$ is compact and decreasing, we know that the sequence has a nonempty
intersection (otherwise one of the $E_r$ would already be empty). We may infer
\begin{equation*}
\bigcap_{\ell \in \Nat} D_\ell \supseteq \bigcap_{r \in \Nat} E_r \not= \emptyset.
\end{equation*}

6.
The classic Hahn Extension
Theorem~\SetCite{Theorem 1.113}
now tells us that there exists a unique
extension of $K_0(y)$ from $\mathcal{C}_z$ to a measure $K(y)$ on
$\sigma(\mathcal{C}_z) = \Borel{Z}$, whenever $y \notin N$. If, however,
$y \in N$, then we define $K(y) := \nu$, where $\nu \in \SubProb{Z}$ is
arbitrary. Because
\begin{equation*}
\dInt[B]{K(y)(C)}{\mu_f(y)} = \dInt[B]{K_0(y)(C)}{\mu_f(y)}
=
\mu(\{x \in X \mid f(x) \in B, g(x) \in C\})
\end{equation*}
holds for $C \in \mathcal{C}_z$, the $\pi$-$\lambda$-Theorem~\ref{Pi-Lambda}
asserts that this equality is valid for all $C \in \Borel{Z}$ as well.

Measurability of $y \mapsto K(y)(C)$ needs to be shown, and then we are
done. We do this by the principle of good sets: put
\begin{equation*}
\mathcal{E} := \{C \in \Borel{Z} \mid y \mapsto K(y)(C)\text{ is } \mathcal{B}-\text{measurable}\}.
\end{equation*}
Then $\mathcal{E}$ is a $\sigma$-algebra, and $\mathcal{E}$ contains the generator $\mathcal{C}_z$ by
construction, thus $\mathcal{E} = \Borel{Z}$.
\EndProof

The scenario in which the space $X = Y \times Z$ with a measurable
space $(Y, \mathcal{B})$ and a Polish space $Z$ with
$\mathcal{A} = \mathcal{B} \otimes \Borel{Z}$ with $f$ and $g$
as projections deserves particular attention. In this case
we decompose a measure on $A$ into its projection onto $Z$
and a conditional distribution for the projection onto
$Z$ given the projection onto $Y$. This is sometimes called
the \emph{disintegration}\index{disintegration} of a
measure $\mu \in \SubProb{Y \times Z}$.

We state the corresponding proposition explicitly, since one needs it
usually in this specialized form.

\BeginProposition{Disintegration-exists}
Given a measurable space $(Y, \mathcal{B})$ and a Polish space $Z$, there
exists for every subprobability
$
\mu \in \SubProb{Y \times Z, \mathcal{B} \otimes \Borel{Z}}
$
a regular conditional distribution of $\pi_Z$ given $\pi_Y$, that is,
a stochastic relation $K: (Y, {\cal B})\Trans (Z, \Borel{Z})$ such
that 
\begin{equation*}
  \mu(E) = \dInt[Y]{K(y)(E_{y})}{\SubProb{\pi_{Y}}(\mu)(y)}
\end{equation*}
for all $E\in {\cal B}\otimes\Borel{Z}$.
\QED
\EndProposition

The construction is done with a Polish as one of the factors. The
proof shows that it is indeed tightness which saves the days, since
otherwise it would be difficult to make sure that the condition
distribution constructed above is $\sigma$-additive. In fact, examples
show that this assumption is in fact necessary:~\cite{Kellerer}
constructs a product measure on spaces which fail to be Polish, for
which no disintegration exists.
%%% Local Variables: 
%%% mode: latex
%%% TeX-master: "../Mskr3"
%%% End: 

%\Input{\Folder/BibNotes}
\Subsection{Bibliographic Notes}
\label{sec:probabs-bib-notes}

Most topics of this chapter are fairly standard, hence there are
plenty of sources to mention. One of my favourite texts is the rich
compendium written by~\cite{Bogachev}. The discussion on Souslin's
operation $\sA({\cal A})$ on a $\sigma$-algebra ${\cal A}$ is heavily
influenced by Srivastava's representation~\cite{Srivastava} of this
topic, but see also~\cite{Parthasarathy,Arveson,Kellerer}. The measure
extension is taken from~\cite{Lubin}, following a suggestion by
S. M. Srivastava; the extension of a stochastic relation is
from~\cite{EED-PDL-TR}. The approach to integration centering around
B. Levi's Theorem is taken mostly from the elegant representation
by Doob~\cite[Chapter~VI]{Doob}, see
also~\cite[Kapitel~IV]{Elstrodt}. The introduction of the Daniell
integral follows essentially~\cite[Sec. 7.8]{Bogachev}, see
also~\cite{Kellerer}. The logic \texttt{CSL} is defined and
investigated in terms of model checking
in~\cite{Baier+Haverkort+Hermanns+Katoen}, the stochastic
interpretation is taken from~\cite{EED-HennessyMilner}, see
also~\cite{Desharnais-Panangaden-ContTime}. The Hutchinson metric is
discussed in detail in Edgar's monograph~\cite{Edgar}, from which the
present proof of Proposition~\ref{weak-is-hutchinson} is taken. There
are many fine books on Banach spaces, Hilbert spaces and the
application to $L^{p}$ spaces; my sources are~\cite{Doob,Halmos,Rudin,
  Dunford-Schwartz,Loomis, Schaefer}. The exposition of projective
limits, and of disintegration follows basically~\cite[Chapter
V]{Parthasarathy} with an occasional glimpse at~\cite{Bogachev}.

%%% Local Variables: 
%%% mode: latex
%%% TeX-master: "../Mskr3"
%%% End: 

%\Input{\Folder/Exercises}
%spell checked - 24Aug14
\Subsection{Exercises}
\label{sec:probs-exercises}
\Exercise{ex-weak-generator}{
Assume that ${\cal A} = \sigma({\cal A}_{0})$. Show that the weak-*-$\sigma$-algebra $\schwach{{\cal A}}$ on $\FinM{X, {\cal A}}$ is the initial $\sigma$-algebra with respect to $\{ev_{A}\mid A\in{\cal A}_{0}\}$. 

Show also that both $\SubProb{X, {\cal A}}$ and $\Prob{X, {\cal A}}$
are measurable subsets of $\FinM{X, {\cal A}}$.  }
\Exercise{ex-cont-baire-meas}{ Let $(X, \tau)$ be a topological, and
  $(Y, d)$ a metric space. Each continuous function $X\to Y$ is also
  Baire measurable.  } \Exercise{ex-support-nghbhood}{ Let $(X, d)$ be
  a separable metric space, $\mu\in\FinM{X, \Borel{X}}$. Show that $x
  \in\supp(\mu)$ iff $\mu(U)>0$ for each open neighborhood $U$ of $x$.
} \Exercise{ex-unif-impl-ae}{ Let $(X, {\cal A}, \mu)$ be a finite
  measure space. Show that norm convergence in $L_{\infty}(X, {\cal
    A},\mu)$ implies convergence almost everywhere ($f_{n}\aeC f$,
  provided $\infNorm{f_{n}-f}{\mu}\to 0$). Give an example showing
  that the converse is false.  } 
\Exercise{ex-extend-top-and-sigma}{
If ${\cal A}$ is a $\sigma$-algebra on $X$ and $B\subseteq X$ with $A\not\in{\cal A}$, then 
\begin{equation*}
\{(A_{1}\cap B)\cup(A_{2}\cap(X\setminus B)) \mid  A_{1}, A_{2}\in {\cal A}\}
\end{equation*}
is the smallest $\sigma$-algebra $\sigma({\cal A}\cup\{B\})$ on $X$ containing ${\cal A}$ and $B$. If $\tau$ is a topology on $X$ with $H\not\in\tau$, then 
\begin{equation*}
\{G_{1}\cup (G_{2}\cap H) \mid  G_{1}, G_{2}\in \tau\}
\end{equation*}
is the smallest topology $\tau_{H}$ on $X$ containing $\tau$ and $H$. Show that $\Borel{\tau_{H}} = \sigma({\cal A}\cup\{H\})$
}
\Exercise{ex-extend-measure}{
Let $(X, {\cal A}, \mu)$ be a finite measure space, $B\not\in{\cal A}$, and 
$\beta := \alpha\cdot \mu_{*}(B) + (1 - \alpha)\cdot \mu^{*}(B)$ with $0\leq \alpha \leq 1$. Then there exists a measure  $\nu$ on $\sigma({\cal A}\cup\{B\})$ which extends $\mu$ such that $\nu(B) = \beta.$ (Hint: Exercise~\ref{ex-extend-top-and-sigma}).
}
\Exercise{ex-meas-approx-from-below}{
  Given the measurable space $ (X, \mathcal{A}) $ and $f \in
  \MeasbFnct{X, \mathcal{A}}$ with $f \geq 0$. Show that there exists
  a decreasing sequence $\Folge{f}$ of step functions $f_n \in
  \MeasbFnct{X, \mathcal{A}}$ with
  \begin{equation*}
  f(x) = \inf_{n \in \Nat} f_n(x)
\end{equation*}
  for all $x \in X.$
}
\Exercise{ex-prod-meas}{
Let $f_{i}: X_{i}\to Y_{i}$ be ${\cal A}_{i}$-${\cal B}_{i} $-measurable maps for $i\in I$. Show that 
\begin{equation*}
  f:
  \begin{cases}
    \prod_{i\in I}X_{i}&\to \prod_{i\in I}Y_{i}\\
(x_{i})_{i\in I}& \mapsto (f_{i}(x_{i}))_{i\in I}
  \end{cases}
\end{equation*}
is $\bigotimes_{i\in I}{\cal A}_{i}$-$\bigotimes_{i\in I}{\cal B}_{i}$-measurable. Conclude that the kernel of $f$
\begin{equation*}
  \Kern{f} := \{\langle x, x'\rangle \mid f(x) = f(x')\}
\end{equation*}
is a measurable subset of $Y\times Y$, whenever $f: (X, {\cal A})\to (Y, {\cal B})$ is measurable, and ${\cal B}$ is separable. 
}
\Exercise{ex-graph-meas}{
Let $f: X\to Y$ be ${\cal A}$-${\cal B}$ measurable, and assume that ${\cal B}$ is separable. Show that the graph of $f$
\begin{equation*}
  \Graph{f} := \{\langle x, f(x)\rangle \mid x\in X\}
\end{equation*}
is a measurable subset of ${\cal A}\otimes{\cal B}$. 
}
\Exercise{char-fnct}{
  Let $\chi_{A}$ be the indicator function of set $A$. Show that
  \begin{enumerate}
  \item $A\subseteq B$ iff $\chi_{A}\leq\chi_{B}$,
  \item $\chi_{\bigcup_{n\in\Nat}A_{n}} = \sup_{n\in\Nat}\chi_{A_{n}}$
    and $\chi_{\bigcap_{n\in\Nat}A_{n}} = \inf_{n\in\Nat}\chi_{A_{n}}$ 
  \item $\chi_{A\Delta B} = |\chi_{A} - \chi_{B}| = \chi_{A}+ \chi_{B}
    \ (\mathrm{mod}\ 2)$. Conclude that the
    power set $(\PowerSet{X}, \Delta)$ is a commutative group with
    $A\Delta A = \emptyset$. 
  \item $\bigl(\bigcup_{n\in\Nat}A_{n}\bigr)\Delta
    \bigl(\bigcup_{n\in\Nat}B_{n}\bigr)\subseteq \bigcup_{n\in\Nat}(A_{n}\Delta B_{n})$
  \end{enumerate}
  \Exercise{ex-complete-psudo-metric}{ Let $(X, {\cal A}, \mu)$ be a
    finite measure space, put $d(A, B) := \mu(A\Delta B)$ for $A,
    B\in{\cal A}$. Show that $({\cal A}, d)$ is a complete pseudo
    metric space.  } }
\Exercise{ex-equiv-compl-metric}{
Let $(X, d)$ be a metric space. Show that 
\begin{equation*}
  D(x, y) := \frac{d(x, y)}{1 + d(x, y)}
\end{equation*}
defines a metric on $X$ such that the metric spaces are homeomorphic as topological spaces. The $(X, d)$ is complete iff $(X, D)$ is. 
}
\Exercise{ex-meas-non-regular}{
(This Exercise draws heavily on Exercises~\ref{ex-extend-top-and-sigma} and~\ref{ex-extend-measure}). 
Let $X := [0, 1]$ with $\lambda$ as the Lebesgue measure on the Borel set of $X$. There exists a set $B\subseteq X$ with $\lambda_{*}(B) = 0$ and $\lambda^{*}(B) = 1$~\SetCite{Lemma 1.141}, so that $B\not\in\Borel{X}$. 
\begin{enumerate}
\item Show that $(X, \tau_{B})$ is a Hausdorff space with a countable base, where $\tau_{B}$ is the smallest topology containing the interval topology on $[0, 1]$ and $B$ (see Exercise~\ref{ex-extend-top-and-sigma}).
\item Extend $\lambda$ to a measure $\mu$ with $\alpha = 1/2$ in Exercise~\ref{ex-extend-measure}.
\item Show that $\inf\{\mu(G) \mid G\supseteq X\setminus B\text{ and } G\text{ is $\tau_{B}$-open}\} = 1$, but $\mu(X\setminus B) = 1/2$. Thus $\mu$ is not regular (since $(X, \tau_{B})$ is not a metric space).
\end{enumerate}
}
\Exercise{ex-limit-is-measurable}{ Prove
  Proposition~\ref{limit-is-measurable}.  }
\Exercise{ex-kernel-yields}{ Let $K: (X, {\cal A})\Trans (Y, {\cal
    B})$ be a transition kernel.
  \begin{enumerate}
  \item Assume that $f\in\MeasbFnctP{Y, \mathcal{B}}$ is
    integrable with respect to $K(x)$ for all $x\in X$. Show that
    \begin{equation*}
      K(f)(x) := \dInt{f}{K(x)}
    \end{equation*}
    defines a measurable function $K(f): X \to \Real_{+}$.
  \item Assume that $x\mapsto
    K(x)(Y)$ is bounded. Define for $B\in{\cal B}$ 
    \begin{equation*}
      \comp{K}(\mu)(B) := \dInt{K(x)(B)}{\mu(x)}.
    \end{equation*}
Show that $\comp{K}: \SubProb{X, {\cal A}}\to\SubProb{Y, {\cal B}}$ is
$\schwach{X, {\cal A}}$-$\schwach{Y, {\cal B}}$-measurable
(see~\CategCite{Example 1.99}). 
  \end{enumerate}
} 
\Exercise{ex-int-meas-times-kernel}{
Let $\mu\in\SubProb{X,{\cal A}}$ be s subprobability measure on $(X, {\cal A})$, and let $K: (X, {\cal A})\Trans (Y, {\cal B})$ be a stochastic relation. Assume that $f: X\times Y\to \Real$ is bounded and measurable. Show that
\begin{equation*}
  \dInt[X\times Y]{f}{\mu\otimes K} = \dInt{\bigl(\dInt[Y]{f_{x}}{K(x)}\bigr)}{\mu(x)}
\end{equation*}
($\mu\otimes K$ is defined in Example~\ref{measure-times-kernel} on page~\pageref{measure-times-kernel}).
}
\Exercise{ex-joint-measurability}{
Let $(X, {\cal A})$ and $(Y, {\cal B})$ be measurable spaces and $D\in{\cal A}\otimes{\cal B}$. Show that the map
\begin{equation*}
  \begin{cases}
    \FinM{Y, {\cal B}}\times X & \to \Real\\
\langle \nu, x\rangle & \mapsto \nu(D_{x})
  \end{cases}
\end{equation*}
is $\schwach{Y, {\cal B}}\otimes{\cal A}$-$\Borel{\Real}$-measurable (the weak-*-$\sigma$-algebra $\schwach{Y, {\cal B}}$ has been defined in Section~\ref{sec:sigma-algebra-on-measures}). 
}
\Exercise{non-pushout}{
  Show that the category of analytic spaces with measurable maps is
  not closed under taking pushouts. \textbf{Hint}: Show that the
  pushout of $\Faktor{X}{\alpha_{1}}$ and $\Faktor{X}{\alpha_{2}}$ is
  $\Faktor{X}{(\alpha_{1}\cup\alpha_{2})}$ for equivalence relations
  $\alpha_{1}$ and $\alpha_{2}$ on a Polish space $X$. Then use
  Proposition~\ref {RemainsAnalytic} and
  Example~\ref{inters-non-smooth}.  }
\Exercise{ex-weaktop-discr-space}{
Let $S := \{1,\dots, n\}$ for some $n\in \Nat$. Show that the weak topology on $\FinM{S, \PowerSet{S}}$ can be identified with the Euclidean topology on $(\pReal)^{n}$.
}
\Exercise{ex-trans-syst-eff-fnct}{
Let $S := \{1, \dots, n\}$ for some $n\in \Nat$ be the finite state space of some transition system $\to_{S}$. Given $s\in S$, let $R(s) := \{s'\in S \mid s\to_{S}s'\}$ be the set of all successors to $s$, put
\begin{align*}
\kappa(s) & := \{\sum_{s'\in R(s)}\alpha_{s'}\cdot \delta_{s'}\mid 0\geq\alpha_{s'}\text{ rational}, \sum_{s'\in R(s)}\alpha_{s'}\leq 1\},\\
P(s) & := \{A\in\schwach{S, \PowerSet{S}} \mid \kappa(s)\subseteq A\}.
\end{align*}
Show that the set 
$
\{\langle s, q\rangle \mid H^{q}\in P(s)\}
$
is a member of $\PowerSet{S}\otimes\Borel{[0, 1]}$ 
for any $H\in \schwach{S, \PowerSet{S}}\otimes\Borel{[0, 1]}$.

This construction is of interest in the analysis of stochastic non-determinism. 
}
%%% Local Variables: 
%%% mode: latex
%%% TeX-master: "../Mskr3"
%%% End: 
\Exercise{ex-weak-bisim}{Let $X$ and $Y$ be Polish spaces with a
  transition kernel $K: X\Trans Y$. The equivalence relations $\alpha$
  on $X$ and $\beta$ on $Y$ are assumed to be smooth with determining
  sequences $\Folge{A}$ resp. $\Folge{B}$ of Borel sets. Put ${\cal I}_{\alpha} := \sigma(\{A_{n}\mid n\in\Nat\}$ and ${\cal J}_{\beta} := \sigma(\{B_{n}\mid n\in\Nat\}$. Show that the following statements are equivalent
  \begin{enumerate}
  \item $K: (X,{\cal I}_{\alpha})\Trans (Y,{\cal J}_{\beta})$
is a transition kernel.
\item $(\alpha, \beta)$ is a congruence for $K$. 
\item $\alpha\subseteq \Kern{\SubProb{\fMap{\beta}\circ K}}$.
\item There exists a transition kernel $K': (X, {\cal
    I}_{\alpha})\Trans (Y, {\cal J}_{\beta})$ such that $(i_{\alpha},
  j_{\beta}): K\to K'$ is a morphism, where the measurable maps
  $i_{\alpha}: (X, \Borel{X})\to (X, {\cal I}_{\alpha})$ and
  $j_{\beta}: (Y, \Borel{Y})\to (Y, {\cal I}_{\beta})$ are given by
  the respective identities.
\end{enumerate}
}  
% \Exercise{ex-weak-bisim}{Let $X$ and $Y$ be Polish spaces with a
%   transition kernel $K: X\Trans Y$. The equivalence relations $\alpha$
%   on $X$ and $\beta$ on $Y$ are assumed to be smooth with determining
%   sequences $\Folge{A}$ resp. $\Folge{B}$. Show that the following statemenets are equivalent
%   \begin{enumerate}
%   \item $K: \bigl(X, \sigma(\{A_{n}\mid n\in\Nat\})\bigr)\Trans \bigl(Y, \sigma(\{B_{n}\mid
% n\in\Nat\})\bigr)$
% is a transition kernel.
% \item $(\alpha, \beta)$ is a congruence for $K$. 
% \item $\alpha\subseteq \Kern{\SubProb{\fMap{\beta}\circ K}}$.
% \end{enumerate}
% }  
\Exercise{ex-atoms-in-sigma-algebra}{
\def\sM{\ensuremath{\mathtt{S}_{X}}}
Let $\sM$ be the set of all smooth equivalence relations on the Polish space $X$, which is ordered by inclusion. Then $\sM$ is closed under countable infima, and $\Delta_{X}\subseteq \rho \subseteq \nabla_{X}$, where $\nabla_{X} := X\times X$ is the universal relation. 
\begin{enumerate}
\item $\rho\mapsto\{A\in \Borel{X}\mid A\text{ is } \rho-invariant\}$
is an order reversing bijection between $\sM$ and the countably generated sub-$\sigma$-algebras of $\Borel{X}$ such that $\Delta_{X}\mapsto \Borel{X}$ and $\nabla_{X}\mapsto \{\emptyset, X\}$.
\item Define for $x, x'\in X$ with $x\not=x'$ the equivalence relation
  $\theta_{x, x'} := \Delta_{X}\cup\{\langle x, x'\rangle, \langle x', x\rangle\}$. Then $\theta_{x, x'}$ is an atom of $\sM$. Describe the $\sigma$-algebra of $\theta_{x, x'}$-invariant Borel sets.
\item Define for the Borel set $B$ with $\emptyset\not=B\not=X$ the equivalence relation $\tau_{B}$ through $\isEquiv{x}{x'}{\tau_{B}}$ iff $\{x, x'\}\subseteq B$ or $\{x, x'\}\cap B=\emptyset$ for all $x, x'\in X$. Then $\tau_{B}$  is an anti-atom in $\sM$ (i.e., an atom in the reverse order). Describe the $\sigma$-algebra of $\tau_{B}$-invariant Borel sets.
\item Show that for each $\rho\in\sM$ there exists a countable family $\Folge{\beta}$ of anti-atoms with $\rho=\bigwedge_{n\in \Nat} \beta_{n}$.
\item Show that $\tau_{B}\wedge\theta_{x, x'}= \Delta_{X}$ and $\tau_{B}\vee\theta_{x, x'}= \nabla_{X}$, whenever $B$ is a Borel set with $\emptyset\not=B\not=X$ and $x\in B$, $x'\not\in B$.
\end{enumerate}
} 
\Exercise{ex-selection-algebra}{
  Let $Y$ be a Polish space, $F: X\to \Closed(Y)$ and ${\cal L}$ and
  algebra of sets on $X$. We assume that $\weak{F}(G)\in{\cal
    L}_{\sigma}$ for each open $G\subseteq Y$ ($\weak{F}$ is defined
  on page~\pageref{set-valued-measurability}). Show that there exists
  a map $s: X\to Y$ such that $s(x)\in F(x)$ for all $x\in X$ such
  that $\InvBild{s}{B}\in{\cal L}_{\sigma}$ for each
  $B\in\Borel{Y}$. \textbf{Hint}: Modify the proof for
  Theorem~\ref{meas-selections-exist} suitably.  }
\Exercise{step-fnct-val}{ Given a finite measure space $(X, {\cal A},
  \mu)$, let $f = \sum_{i=1}^{n}\alpha_{i}\cdot\chi_{A_{i}}$ be a step
  function with $A_{1}, \dots, A_{n}\in{\cal A}$ and coefficients
  $\alpha_{1}, \dots, \alpha_{n}$. Show that
\begin{equation*}
  \sum_{i=1}^{n}\alpha_{i}\cdot\mu(A_{i}) = \sum_{\gamma>0}\gamma\cdot\mu(\{x\in X\mid f(x) = \gamma\}).
\end{equation*}
}
\Exercise{ex-vietoris-top}{
\def\KK{\mathfrak{C}}
Let $(X, d)$ be a metric space, and define 
\begin{equation*}
\KK(X) := \{C\subseteq X\mid \emptyset\not=C\text{ is compact}\}.
\end{equation*}
Given $C_{1}, C_{2}\in \KK(X)$, define the \emph{\index{Hausdorff distance}Hausdorff distance} of $C_{1}$ and $C_{2}$ through 
\begin{equation*}
  D_{H}(C_{1}, C_{2}) := \max\{\sup_{x\in C_{2}}d(x, C_{1}), \sup_{x\in C_{1}}d(x, C_{2})\}. 
\end{equation*}
  \begin{enumerate}
  \item Show that $D_{H}(C_{1}, C_{2}) < \epsilon$ iff $C_{1}\subseteq C_{2}^{\epsilon} \text{ and } C_{2}\subseteq C_{1}^{\epsilon}$ (the $\epsilon$-neighborhood $B^{\epsilon}$ if a set is defined on page~\pageref{a-top-basis-without-bound}).
  \item Show that $(\KK(x), D_{H})$ is a metric space.
  \item If $(X, d)$ has a countable dense subset, so has $(\KK(X), D_{H})$. 
  \item Let $(Y, {\cal B})$ be a measurable space and assume that $X$ is compact. Show that $F: Y\to \KK(X)$ is ${\cal B}$-$\Borel{\KK(X)}$ measurable iff $F$ is measurable (as a relation, in the sense of Definition~\ref{set-valued-measurability} on page~\pageref{set-valued-measurability}). 
  \end{enumerate}
}
\Exercise{ex-for-jenny}{
Given the plane $E := \{\langle x_{1}, x_{2}, x_{3}\rangle\in \Real^{3}\mid  2\cdot x_{1} + 4\cdot x_{2}-7\cdot x_{3} = 12\}$, determine the point in $E$ which is closest to $\langle 4, 2, 0\rangle$ in the Euclidean distance.
}
\Exercise{ex-linear-cont-is-bounded}{
Let $(V, \|\cdot \|)$ be a real normed space, $L: V\to \Real$ be linear. Show that $L$ is continuous iff $L$ is bounded, i.e., iff 
$
\sup_{\|v\|\leq 1}|L(v)|<\infty.
$ 
}
\Exercise{ex-linear-cont-is-bounded-cont}{
Let $(V, \|\cdot \|)$ be a real normed space, and define 
\begin{equation*}
  V^{*} := \{L: V\to \Real\mid L\text{ is linear and continuous}\},
\end{equation*}
the \emph{\index{space!dual}dual space} of $V$. Then $V^{*}$ is a vector space. Show that 
\begin{equation*}
  \|L\|:=\sup_{\|v\|\leq 1}|L(v)|
\end{equation*}
defines a norm on $V^{*}$ with which $(V^{*}, \|\cdot \|)$ is a Banach space. 
}
\Exercise{ex-dual-hilbert}{
Let $H$ be a Hilbert space, then $H^{*}$ is isometrically isomorphic to $H$.
}
\Exercise{ex-double-dual}{
Let $(V, \|\cdot \|)$ be a real normed space, and define 
\begin{equation*}
  \pi(x)(L) := L(x)
\end{equation*}
for $x\in V$ and $L\in V^{*}$. 
\begin{enumerate}
\item Show that $\pi(x)\in V^{**}$, and that $x\mapsto\pi(x)$ defines
  a continuous map $V\to V^{**}$.
\item Given $x\in V$ with $x\not=0$, there exists $L\in V^{*}$ with $\|L\|=1$ and $L(x)=\|x\|$ (use the Hahn-Banach Theorem~\SetCite{Theorem 1.55}).
\item Show that $\pi$ is an isometry (thus a normed space can be embedded isometrically into its bidual).
\end{enumerate}

}
\Exercise{ex-parallelogram-law}{
Given a real vector space $V$.
\begin{enumerate}
\item Let $(\cdot , \cdot )$ be an inner product on $V$. Show that
  \begin{equation*}
    \|x + y\|^{2} + \|x - y\|^{2} = 2\cdot \|x\|^{2} + 2\cdot \|y\|^{2}
  \end{equation*}
  always holds. This equation is known as the \emph{parallelogram
    law}: The sum of the squares of the diagonals is the sum of the
  squares of the sides in a parallelogram. 
\item Assume, conversely, that $\|\cdot \|$ is a norm for which the parallelogram law holds. Show that 
  \begin{equation*}
    (x, y) := \frac{\|x+y\|^{2}-\|x-y\|^{2}}{4}
  \end{equation*}
defines an inner product on $V$. 
\end{enumerate}
\Exercise{ex-factor-linear}{
Let $H$ be a Hilbert space, $L: H\to \Real$ be a continuous linear map with $L\not=0$. Relating $Kern(L)$ and $\Kern{L}$, show that $\Faktor{H}{Kern(L)}$ and $\Real$ are isomorphic as vector spaces. 
}
}

%%% Local Variables: 
%%% mode: latex
%%% TeX-master: "../Mskr3"
%%% End: 

%%% Local Variables: 
%%% mode: latex
%%% TeX-master: "../Mskr3"
%%% End: 

%\bibliographystyle{alpha}
\newpage
\addcontentsline{toc}{section}{References}
%\bibliography{EED,LitDB,Haskell}

\newpage
\addcontentsline{toc}{section}{Index}
%\printindex
\begin{theindex}

  \item $H_{\gamma}$, 103
  \item $\FinM{X, {\cal A}}$, 10
  \item $\SigmaM{X, {\cal       A}}$, 10
  \item $\aeC$, 22
  \item $\mu $-a.e., 21
  \item $\mu $-essentially bounded, 21
  \item $\mu $-null set, 56
  \item $\nmC$, 24
  \item $\pi $-$\lambda $-Theorem., 5
  \item $\schwach{X, {\cal A}}$, 10
  \item $\sigma $-algebra
    \subitem countable-cocountable, 28
    \subitem countably generated, 27
    \subitem separable, 29
  \item $\sigma$-algebra
    \subitem final, 8
    \subitem initial, 8
    \subitem product, 9
    \subitem sum, 9
    \subitem trace, 8
  \item $\supp(\mu)$, 16
  \item $\tau $-regular, 15
  \item $d\mu/d\nu$, 128
  \item ${\cal  A}$-cover, 54
  \item ${\cal L}_{\infty}(\mu)$, 21

  \indexspace

  \item  $\betaSenza{{\cal A}}$, 10

  \indexspace

  \item absolute continuity, 126
    \subitem signed measure, 136
  \item analytic set, 41
  \item atom, 31

  \indexspace

  \item Baire   sets, 4
  \item bisimilar, 107
  \item Borel   sets, 4
  \item Borel measurability, 11
  \item Borel sets, 4
  \item bounded, 132

  \indexspace

  \item change of variables
    \subitem calculus, 72
    \subitem image measure, 71
  \item closure of a set, 35
  \item complete
    \subitem measure space, 56
  \item completion
    \subitem $M$, 57
    \subitem $\mu $, 56
    \subitem universal, 57
  \item conditional distribution, 138
    \subitem regular, 138
  \item congruence, 110
  \item conjugate numbers, 122
  \item continuity, 32
  \item convergence
    \subitem almost everywhere, 22
    \subitem in measure, 24
    \subitem pointwise, 22
    \subitem uniform, 22
  \item converse, 137
  \item convex, 119
  \item CSL, 90
    \subitem next operator, 90
    \subitem path quantifier, 90
    \subitem steady-state, 90
    \subitem until operator, 90
  \item cut
    \subitem horizontal, 35
    \subitem vertical, 35
  \item cylinder sets, 85

  \indexspace

  \item diameter, 16
  \item disintegration, 142
  \item distance
    \subitem L\IeC {\'e}vy-Prohorov, 14

  \indexspace

  \item equivalence relation
    \subitem invariant set, 48
    \subitem smooth, 46
      \subsubitem determining sequence, 46

  \indexspace

  \item function
    \subitem indicator, 19
    \subitem step, 19

  \indexspace

  \item graph
    \subitem of a map, 42

  \indexspace

  \item H\"older's inequality, 123
  \item Hausdorff distance, 146
  \item Hilbert cube, 37
  \item Hilbert space, 119
  \item hit
    \subitem $\sigma$-algebra, 5
    \subitem measurable, 66
  \item homeomorphism, 33
  \item Hutchinson metric, 103

  \indexspace

  \item indicator function, 19
  \item inequality
    \subitem H\"older, 122, 123
    \subitem Minkowski, 123
    \subitem Schwarz, 118
  \item inner product, 118
  \item invariant set, 48

  \indexspace

  \item Jordan decomposition, 134

  \indexspace

  \item kernel
    \subitem linear functional, 119
    \subitem Markov, 11
    \subitem transition, 10
  \item Kuratowski's trap, 36

  \indexspace

  \item Lagrange's Formula, 131
  \item Lebesgue decomposition, 128
  \item linear functional, 132
    \subitem positive, 132
  \item logic
    \subitem continuous time stochastic, 89
    \subitem CSL, 90
    \subitem modal, 11

  \indexspace

  \item Mackey, 30
  \item map
    \subitem continuous, 32
    \subitem graph, 42
    \subitem kernel, 31
    \subitem measurable, 7
    \subitem semicontinuous, 7
  \item Markov property, 92
  \item measurable
    \subitem rectangle, 9
    \subitem relation, 65
    \subitem selector, 65
    \subitem set-valued map, 65
  \item measure
    \subitem projective limit, 87, 89
    \subitem projective system, 87
    \subitem tight, 102
  \item Minkowski's inequality, 123
  \item morphism
    \subitem measure spaces, 111
    \subitem stochastic relations, 106
  \item mutual singular, 127

  \indexspace

  \item next operator
    \subitem CSL, 90

  \indexspace

  \item orthogonal
    \subitem complement, 119
    \subitem vector, 119
  \item oscillation, 35

  \indexspace

  \item path quantifier
    \subitem CSL, 90
  \item principle of good sets, 7
  \item pseudo-norm, 21

  \indexspace

  \item quotient object, 112

  \indexspace

  \item Radon-Nikodym derivative, 128, 129

  \indexspace

  \item Schwarz inequality, 118
  \item separable, 29
  \item set
    \subitem analytic, 41
    \subitem clopen, 38
    \subitem co-analytic, 41
    \subitem cylinder, 9
  \item signed   measure
    \subitem absolute continuity, 136
  \item signed measure, 134
    \subitem nullset, 136
    \subitem positive set, 134
    \subitem total variation, 135
  \item Souslin scheme, 51
    \subitem regular, 51
  \item space
    \subitem analytic, 44
    \subitem dual, 147
    \subitem measurable, 4
    \subitem metric
      \subsubitem complete, 33
    \subitem Polish, 33
      \subsubitem Borel sets, 4
    \subitem standard Borel, 51
    \subitem topological
      \subsubitem Baire sets, 4
      \subsubitem base, 32
      \subsubitem Borel sets, 4
      \subsubitem compact, 38
      \subsubitem continuity, 32
      \subsubitem homeomophism, 33
      \subsubitem subbase, 32
  \item steady-state
    \subitem CSL, 90
  \item stochastic relations, 11
  \item support, 16

  \indexspace

  \item theorem
    \subitem Alexandrov, 37, 89
    \subitem Blackwell-Mackey, 49
    \subitem Egorov, 24
    \subitem Kuratowski and Ryll-Nardzewski, 65
    \subitem Kuratowski Isomorphism, 44
    \subitem Lebesgue Dominated Convergence, 71
    \subitem Lubin, 60
    \subitem Lusin, 43
    \subitem Riesz Representation, 76
    \subitem Souslin, 44
    \subitem unique structure, 45
    \subitem von Neumann Selection, 59
  \item topology
    \subitem Alexandrov, 12
    \subitem Baire   sets, 4
    \subitem base, 32
    \subitem Borel   sets, 4
    \subitem compact, 38
    \subitem initial, 33
    \subitem product, 33
    \subitem subbase, 32
    \subitem subspace, 33
    \subitem sum, 33
    \subitem weak, 98

  \indexspace

  \item universal set, 35
  \item until operator
    \subitem CSL, 90

  \indexspace

  \item vector lattice, 73

  \indexspace

  \item weak topology, 98

\end{theindex}

\end{document}